\documentclass[11pt]{article}

\usepackage{apacite}
\usepackage{natbib}
\usepackage{epsfig,lscape}
\usepackage{amssymb,amsfonts,amsmath,amsthm}
\usepackage{rotating}
\usepackage{bbm}

\usepackage{threeparttable}
\usepackage{booktabs}
\usepackage{subcaption}
\usepackage{multirow}
\usepackage{graphicx}
\usepackage{comment}
\usepackage{multirow,graphicx}

\usepackage{hyperref}

\def\Var{\mathrm{Var}}

\def\exp{\mathrm{exp}}

\def\T{ {\mathrm{\scriptscriptstyle T}} }

\def\bbR{\mathbb R}

\def\argmin{\mathrm{argmin}}

\def\MSE{\mathrm{MSE}}
\def\V{\mathrm{V}}
\def\B{\mathrm{B}}
\def\P{\mathrm{P}}
\def\A{\mathrm{A}}
\def\out{\mathrm{out}}

\def\op{\mathrm{op}}

\def\In{\mathrm{in}}
\def\poly{\mathrm{Poly}}
\def\T{\mathrm{T}}
\def\div{\mathrm{div}}
\def\Im{\mathrm{Im}}
\def\Re{\mathrm{Re}}

\def\Diag{\mathrm{Diag}}
\def\Tr{\mathrm{Tr}}

\newenvironment{prf}
{\noindent \textbf{Proof.}}{\hfill $\Box$ \vspace{.1in}}

\newtheorem{thm}{Theorem}
\newtheorem{lem}{Lemma}
\newtheorem{pro}{Proposition}
\newtheorem{cor}{Corollary}
\newtheorem{ass}{Assumption}

\theoremstyle{definition}

\theoremstyle{definition}
\newtheorem{rem}{Remark}

\allowdisplaybreaks

\begin{document}

\begin{titlepage}

\begin{center}
{\Large On Ridge Estimation in High-dimensional Rotationally Sparse Linear Regression}

\vspace{.1in} {\large Libin Liang\footnotemark[1] and Zhiqiang Tan\footnotemark[1]}

\vspace{.1in}
\today
\end{center}

\footnotetext[1]{Department of Statistics, Rutgers University. Address: 110 Frelinghuysen Road,
Piscataway, NJ 08854. E-mails: ll866@stat.rutgers.edu, ztan@stat.rutgers.edu.
The authors thank Qiyang Han for helpful comments related to error approximation formulas in Section \ref{sec:4_1}.}

\paragraph{Abstract.}
Recently, deep neural networks have been found to nearly interpolate training data but still generalize well in various applications.
To help understand such a phenomenon,
it has been of interest to analyze the ridge estimator and its interpolation limit in high-dimensional regression models.
For this motivation, we study the ridge estimator in a rotationally sparse setting of high-dimensional linear regression,
where the signal of a response is aligned with a small number, $d$, of covariates with large or spiked variances,
compared with the remaining covariates with small or tail variances,
\textit{after} an orthogonal transformation of the covariate vector.
We establish high-probability upper and lower bounds on the out-sample and in-sample prediction errors in two distinct regimes
depending on the ratio of the effective rank of tail variances over the sample size $n$.
The separation of the two regimes enables us to exploit relevant concentration inequalities
and derive concrete error bounds without making any oracle assumption or independent components assumption on covariate vectors.
Moreover, we derive sufficient and necessary conditions which indicate that the prediction errors of ridge estimation can be of the order $O(\frac{d}{n})$
if and only if the gap between the spiked and tail variances are sufficiently large.
We also compare the orders of optimal out-sample and in-sample prediction errors and find that, remarkably,
the optimal out-sample prediction error may be significantly smaller than the optimal in-sample one.
Finally, we present numerical experiments which empirically confirm our theoretical findings.


\end{titlepage}

\section{Introduction}

Over-parameterized models, in particular, deep neural networks, have been successful when trained without penalty or with a mild penalty in various applications.
This phenomenon appears to be at odds with conventional statistical thinking that complex models tend to overfit the training data and generalize poorly
without proper regularization \citep{belkin2019reconciling,zhang_2016}.
Considerable research has been devoted to studying such phenomena;
 see \cite{Bartlett_2020}, \cite{BuneaFlorentinaStrimas2020}, \cite{HasteTrevorAndrea},  and \cite{Tsigler_ridge_2023} among others.

We consider high-dimensional linear regression as a basic example of over-parameterized models
and study ridge estimation under a rotationally sparse setting, which will be introduced shortly.
Suppose that
the training data, $(y_1,x_1), \ldots, (y_n,x_n)$, are i.i.d.~from the following model:
\begin{align}
    y_i = x_i^\T\theta^* + \epsilon_i, \quad i=1,\ldots,n   , \label{eq:1}
\end{align}
where $\theta^* \in \mathbb{R}^p$,
$y_i \in \bbR$ is a response variable, $x_i\in\bbR^p$ is a covariate vector satisfying $E(x_i)=0$ and $\Var(x_i)=\Sigma$ (assumed to be non-singular),
and $\epsilon_i$ is a noise variable, independent of $x_i$ and satisfying $E(\epsilon_i)=0$ and $\Var(\epsilon_i)=\sigma^2$
for $i=1,\ldots,n$.
In addition, assume that $\Sigma^{-1/2}x_i$ is a sub-gaussian random vector with sub-gaussian norm $\sigma_x$ for $i=1,\ldots,n$.
See Supplement Section \ref{sec_s0} for the definition of the sub-gaussian norm.

We are interested in the ridge estimator defined as
\begin{align}
    \hat{\theta}(\tau) &= \argmin_{\theta\in\mathbb{R}^p}
    \left\{ \frac{1}{n} \Vert Y-X\theta\Vert^2 + \tau \Vert\theta\Vert^2 \right\} , \label{eq:2}
\end{align}
where $Y = (y_1,\ldots,y_n)^\T \in \bbR^p$, $X = (x_1,\ldots,x_n)^\T \in \bbR^{n\times p}$,
$\tau \ge 0$ is a tuning parameter, and $\|\cdot\|$ denotes the $L_2$ norm.
In the case of $\tau=0$, the estimator $\hat{\theta}(0)$ is defined as the limit of $\hat{\theta}(\tau)$ as $\tau\to 0+$,
and called a min-norm interpolator.

As a measure of prediction performance, we investigate both the out-sample and in-sample mean squared errors (MSE) of $\hat\theta(\tau)$,
defined conditionally on $X$ as follows.
\begin{itemize}
    \item  Out-sample error
    \begin{align}
        \MSE_{\out} = E[(x^\T_0(\hat{\theta}(\tau)-\theta^*))^2|X] = E[\Vert\hat{\theta}(\tau)-\theta^*\Vert_\Sigma^2|X], \notag
    \end{align}
    where $x_0\in \mathbb{R}^p$ is a new covariate vector independent of $X$, and
    $\| b\|_M = (b^\T M b)^{1/2}$ for any positive semi-definite matrix $M$ and any vector $b$ of the suitable dimension.

    \item In-sample error
    \begin{align}
        \MSE_{\In}=E[\frac{1}{n}\Vert X(\hat{\theta}(\tau) - \theta^*)\Vert^2|X]=E[\Vert\hat{\theta}(\tau)-\theta^*\Vert_{\hat{\Sigma}}^2|X] ,\notag
    \end{align}
    where $\hat{\Sigma} =X^\T X /n$ is the (uncentered) sample covariance matrix of $X$.
\end{itemize}

An important property of the ridge estimator is that its out-sample and in-sample prediction errors are invariant to
an orthonormal transformation of the covariate vectors $x_i$.
By this property, we assume without loss of generality that the covariance matrix $\Sigma$ is diagonal, i.e.,
 $\Sigma=\Diag(\lambda_1,\ldots,\lambda_p)$, where $\lambda_1\geq \ldots \geq \lambda_p>0$ are the eigenvalues of $\Sigma$.

As motivated by the recent literature \citep{Bartlett_2020, BuneaFlorentinaStrimas2020, HasteTrevorAndrea},
we assume that only a few covariates (or a few directions of the covariate vector before orthogonalization) with large variances contain most of the information
about the response. Formally, we assume that the mean response is aligned with a small number of covariates with large variances,
for instance, the first $1\le d <p $ covariates:
\begin{align}
&\frac{\Vert\theta^*_{(d+1):p}\Vert^2_{\Sigma_{(d+1):p}}}{\Vert\theta^*_{1:d}\Vert^2_{\Sigma^{-1}_{1:d}}} \approx 0 .  \notag
\end{align}
We denote as $\theta_{1:d}^*$ the first $d$ entries of $\theta^*$ and as $\theta_{(d+1):p}^*$ the remaining $p-d$ entries of $\theta^*$,
which are called the spiked part and tail part respectively.
In addition, we denote $\Sigma_{1:d}=\Diag(\lambda_1,\ldots, \lambda_d)$ and $\Sigma_{(d+1):p}=\Diag(\lambda_{d+1},\ldots, \lambda_p)$.
We refer to such a setting as a rotationally sparse setting, because a small number of coefficients $\theta_{1:d}^*$ are nonzero
whereas the remaining ones $\theta_{(d+1):p}^*$ are zero or close to zero, \textit{after} an orthogonal transformation of the covariates.
From both the recent literature and our results later, a sufficiently large gap between $\lambda_{d}$ and $\lambda_{d+1}$
may be necessary and sufficient for ridge estimation to achieve meaningful prediction performance.

The rotationally sparse setting fundamentally differs from the directly sparse setting of regression models,
where a small number of coefficients, for example, $s$, are nonzero whereas the remaining $p-s$ ones are zero or close to zero,
with the original covariate vectors. Such a sparse structure depends on the particular coordinate system for the covariate vectors
and would be lost after an orthogonal transformation.
For the directly sparse setting, Lasso estimation is known to be effective in achieving (in-sample) prediction errors in the order $O( s \frac{\log p}{n})$,
that is, $O(\frac{s}{n})$ up to a logarithmic factor of $p$ under suitable conditions (including compatibility or restricted eigenvalue conditions)
(e.g., \cite{Bickel_lasso}; \cite{book_lasso}).
Hence one of the interesting questions which motivate our work is
to investigate plausible conditions for ridge estimation to achieve prediction errors in the order $O(\frac{d}{n})$
in the rotationally sparse linear regression.

\vspace{.05in}
\textbf{Our work.}\;
We study the out-sample and in-sample prediction errors of the ridge estimator in high-dimensional rotationally sparse linear regression.
Classical analysis of ridge estimation deals with the in-sample MSE in a fixed design (with fixed $X$), whereas
recent work on the ridge estimator and the min-norm interpolator, as reviewed later, has been focused on the out-sample MSE in a random design (with random $X$).
Nevertheless, it is also of interest to study the in-sample MSE in a random design.
In fact, the standard analysis of Lasso estimation concerns the in-sample MSE either in a fixed design or in a random design
(e.g., \cite{Bickel_lasso}; \cite{book_lasso}).
The in-sample MSE in a random design also plays an important role in analyzing de-biased Lasso estimation
\citep{Zhang_Zhang_2014, van_de_geer_lasso} and high-dimensional estimation of average treatment effects
\citep{Chernozhukov_2018, tan2020model}.

\begin{table} 
  \centering
\caption{Summary of our main results in the rotationally sparse setting}
  \begin{minipage}{1\linewidth}
    \centering
    \caption*{Regime I (small or moderate TER): $r_d(\Sigma)\lesssim n$}
\tiny
    \begin{tabular}{|l|l|l|}
\hline
Range of $\tau$      & Out-sample error [Theorem \ref{thm1}]     & In-sampe error [Theorem \ref{thm2}]       \\ \hline
$\tau < \lambda_{d+1}$     &    $\MSE_{\out} \gtrsim \sigma^2(\frac{d}{n}+\frac{r_d(\Sigma^2)}{n}) $   &   $\MSE_{\In}\gtrsim \sigma^2(\frac{d}{n} +\frac{r^2_d(\Sigma)}{n^2})$
\\ \hline
\multirow{5}{*}{$\lambda_{d+1}\leq \tau \leq \lambda_{d}$}          &       & Assume $\tau\gg \lambda_{d+1}$:
\\
&$\MSE_{\out} \asymp \Vert\theta^*_{1:d}\Vert^2_{\Sigma^{-1}_{1:d}}\tau^2$&$\Vert\theta^*_{1:d}\Vert^2_{\Sigma^{-1}_{1:d}}\tau^2+\sigma^2(\frac{d}{n}+\frac{\lambda_{d+1}^2}{\tau^2}\frac{r_d(\Sigma)}{n})$\\
&$~~~~~~~~~~~~~~+ \sigma^2(\frac{d}{n}+\frac{\lambda_{d+1}^2}{\tau^2}\frac{r_d(\Sigma^2)}{n}$& $~~~~~~~~~~~~~\gtrsim \MSE_{\In} \gtrsim$  \\
&&$~~~~~~~~~~~~~~~~~~~~~\Vert\theta^*_{1:d}\Vert^2_{\Sigma^{-1}_{1:d}}\tau^2+\sigma^2(\frac{d}{n} +\frac{\lambda_{d+1}^2}{\tau^2}\frac{r_d^2(\Sigma)}{n^2})$\\
\cline{1-3}
\multirow{2}{*}{$\tau>\lambda_d$}         &  \multirow{2}{*}{$\MSE_{\out} \gtrsim \Vert\theta^*_{1:d}\Vert^2_{\Sigma^{-1}_{1:d}}\lambda_d^2$}     &  Assume $\tau\gg \lambda_{d+1}$:
\\
&&$\MSE_{\In} \gtrsim\Vert\theta^*_{1:d}\Vert^2_{\Sigma^{-1}_{1:d}}\lambda_d^2$\\
\hline
\end{tabular}
  \end{minipage}
  \vspace{1em}

  \begin{minipage}{1\linewidth}
    \centering
    \caption*{Regime II (large TER): $r_d(\Sigma)\geq c_x n$}
   \tiny
    \begin{tabular}{|l|l|l|}
\hline
Range of $\tau$ & Out-sample error [Theorem \ref{thm3}] & In-sample error [Theorem \ref{thm4}] \\ \hline
\multirow{4}{*}{$\tau+\lambda_{d+1}\frac{r_d(\Sigma)}{n}\lesssim \lambda_d$}      &   $\MSE_{\out} \asymp \Vert\theta^*_{1:d}\Vert^2_{\Sigma^{-1}_{1:d}}(\tau+\lambda_{d+1}\frac{r_d(\Sigma)}{n})^2$    &Assume $\tau\gg \lambda_{d+1}\frac{r_d(\Sigma)}{n}$: \\
&$~~~~~~~~~~~~~~~~~~~~~~~+\sigma^2(\frac{d}{n}+ \frac{\lambda_{d+1}^2}{(\tau+\lambda_{d+1}\frac{r_d(\Sigma)}{n})^2} \frac{r_d(\Sigma^2)}{n} )$&$\MSE_{\In} \asymp \Vert\theta_{1:d}^*\Vert^2_{\Sigma_{1:d}^{-1}}(\tau+\lambda_{d+1}\frac{r_d(\Sigma)}{n})^2 $\\
&&$~~~~~~~~~~~~~~~~~~~~~~~+\sigma^2(\frac{d}{n}+\frac{\lambda_{d+1}^2}{(\tau+\lambda_{d+1}\frac{r_d(\Sigma)}{n})^2}\frac{r^2_d(\Sigma)}{n^2})$\\
\hline
\multirow{3}{*}{$\tau+\lambda_{d+1}\frac{r_d(\Sigma)}{n}\gtrsim \lambda_d$}      &  \multirow{3}{*}{$\MSE_{\out} \gtrsim\Vert\theta^*_{1:d}\Vert^2_{\Sigma^{-1}_{1:d}}\lambda_d^2$}    &  Assume $\tau\gg \lambda_{d+1}\frac{r_d(\Sigma)}{n}$:   \\
          &      &  \multirow{2}{*}{$\MSE_{\In} \gtrsim\Vert\theta^*_{1:d}\Vert^2_{\Sigma^{-1}_{1:d}}\lambda_d^2$}    \\
       &&   \\
\hline
\end{tabular}
  \end{minipage}
 \begin{tablenotes}
      \small
      \item Note: $c_x$ is a constant depending only on the sub-gaussian norm $\sigma_x$ of covariate vectors.
    \end{tablenotes}
    \label{table:1a}
\end{table}

The main findings of our work can be summarized as follows.
First, we establish high-probability bounds on the out-sample and in-sample MSEs of the ridge estimator (see Table \ref{table:1a})
in two distinct, albeit possibly overlapping, regimes, called small or moderate TER and large TER.\
The two regimes are defined by whether the ratio $\frac{r_d(\Sigma)}{n}$ is small or large, where $r_d(\Sigma)$ is the tail effective rank (TER)
\begin{align*}
r_d(\Sigma) =\frac{\sum_{j>d}\lambda_j}{\lambda_{d+1}} .
\end{align*}
Such a quantity is also central in the related analyses of ridge estimation and min-norm interpolation \citep{Bartlett_2020,Tsigler_ridge_2023}.
A main difference between the two regimes is that in the large TER regime,
the prediction errors of the ridge estimator can sometimes be controlled for a small ridge parameter $\tau$,
including $\tau=0$ corresponding to the min-norm interpolator.
From a technical perspective, the separation of the two regimes enables us to exploit rerlevant concentration inequalities
and derive concrete error bounds without making any oracle assumption or independent components assumption on covariate vectors
as used in \cite{Tsigler_ridge_2023}.

\begin{table} 
\caption{Summary of conditions on $\frac{\lambda_{d+1}}{\lambda_d}$ to achieve $O(\frac{d}{n})$ prediction errors in the rotationally sparse setting}
\tiny
\centering
\begin{tabular}{|l|l|l|l|}
 \hline
& & $\MSE_\out=O(\frac{d}{n})$  & $\MSE_\In=O(\frac{d}{n})$ \\ [0.5ex]
\hline
\multirow{4}{*}{$r_d(\Sigma)\lesssim n$} &\multirow{2}{*}{Sufficient Condition} &\multirow{2}{*}{$\frac{\lambda_{d+1}}{\lambda_{d}}\lesssim  \sqrt{\frac{d}{n}}\min\{1,\sqrt{\frac{d}{r_d(\Sigma^2)}}\}$}& \multirow{2}{*}{$\frac{\lambda_{d+1}}{\lambda_{d}} \lesssim  \sqrt{\frac{d}{n}}\min\{1, \sqrt{\frac{d}{r_d(\Sigma)}}\}$}\\
&&&\\ \cline{2-4}
&\multirow{2}{*}{Necessary Condition}&Assume $n\gg d$ and $r_d( \Sigma^2) \gg d$:&Assume $n\gg d$ and $ \frac{r_d(\Sigma)}{n} \sqrt{\frac{n}{d}}\gg 1$:\\
&&$\frac{\lambda_{d+1}}{\lambda_d} \lesssim \sqrt{\frac{d}{n}}\sqrt{\frac{d}{r_d(\Sigma^2)}} $&$\frac{\lambda_{d+1}}{\lambda_{d}}\lesssim  \frac{d}{r_d(\Sigma)}$\\
\hline
 \multirow{4}{*}{$r_d(\Sigma)\geq c_x n$}&\multirow{2}{*}{Sufficient Condition} &\multirow{2}{*}{$\frac{\lambda_{d+1}}{\lambda_{d}} \lesssim  \sqrt{\frac{d}{n}}\min\{\sqrt{\frac{d}{r_d(\Sigma^2)}},\frac{n}{r_d(\Sigma)}\}
 $}&\multirow{2}{*}{$\frac{\lambda_{d+1}}{\lambda_{d}}\lesssim \frac{d}{r_d(\Sigma)}$}\\
 &&&\\ \cline{2-4}
 &\multirow{2}{*}{Necessary Condition}&Assume $n\gg d$: &Assume $n\gg d$:\\
 &&$\frac{\lambda_{d+1}}{\lambda_{d}} \lesssim \sqrt{\frac{d}{n}}\min\{\sqrt{\frac{d}{r_d(\Sigma^2)}}, \frac{n}{r_d(\Sigma)}\} $&$\frac{\lambda_{d+1}}{\lambda_{d}}\lesssim \frac{d}{r_d(\Sigma)}$\\
 \hline
\end{tabular}
\begin{tablenotes}
      \small
      \item Note: $c_x$ is a constant depending only on the sub-gaussian norm $\sigma_x$ of covariate vectors.
    \end{tablenotes}
\label{table:1b}
\end{table}

\begin{table} 
\caption{Summary of optimal $\MSE$ in the rotationally sparse setting}
\tiny
\centering
\begin{tabular}{|l|l|l|}
 \hline
& Out-sample error  & In-sample error \\ [0.5ex]
\hline
\multirow{2}{*}{$r_d(\Sigma)\lesssim n$} &Assume $\lambda_d\gtrsim \lambda_{d+1}\sqrt{\frac{n}{r_d(\Sigma^2)}}:$& Assume $r_d(\Sigma)\asymp n$ and $\lambda_{d} \gg\lambda_{d+1}$:\\
&$\MSE^*_\out\asymp \max\{\frac{\lambda_{d+1}}{\lambda_d}\sqrt{\frac{r_d(\Sigma^2)}{n}},\frac{d}{n}\}$&$\MSE_\In^*\asymp \max\{\frac{\lambda_{d+1}}{\lambda_d},\frac{d}{n}\}$\\
\hline
 \multirow{2}{*}{$r_d(\Sigma)\geq c_x n$}&&Assume $\lambda_d\gg \lambda_{d+1}\frac{r_d(\Sigma)}{n}$:\\
 &$\MSE_\out^*\asymp \max\{\frac{\lambda_{d+1}}{\lambda_d}\sqrt{\frac{r_d(\Sigma^2)}{n}},\frac{\lambda_{d+1}^2}{\lambda_d^2}\frac{ r_d(\Sigma)^2}{n^2},\frac{d}{n}\}$&$\MSE_\In^*\asymp \max\{\frac{\lambda_{d+1}}{\lambda_d}\frac{r_d(\Sigma)}{n},\frac{d}{n}\}$\\
 \hline
\end{tabular}
\begin{tablenotes}
      \small
      \item Note: $c_x$ is a constant depending only on the sub-gaussian norm $\sigma_x$ of covariate vectors. $\MSE_\out^*$ denotes the $\MSE_\out$ with optimal $\tau$ and $\MSE_\In^*$ denotes the $\MSE_\In$ with optimal $\tau$.
    \end{tablenotes}
\label{table:1c}
\end{table}

\normalsize

Second, from our error bounds, we derive sufficient and necessary conditions on the ratio $\frac{\lambda_{d+1}}{\lambda_d}$
together with the choice of ridge parameter $\tau$ such that
the out-sample and in-sample MSEs is of the order $O(\frac{d}{n})$ respectively (see Table \ref{table:1b}).
All of these conditions are determined in the simple form that the ratio $\frac{\lambda_{d+1}}{\lambda_d}$ is sufficiently small,
i.e., the gap between the spiked and tail variances is sufficiently large.
In other words, our results indicate that ridge estimation can achieve prediction errors in the order $O(\frac{d}{n})$ for a suitable choice of $\tau$
if and only if the gap between the spiked and tail variances is sufficiently large in the rotationally sparse linear regression.
These results can be seen to serve as a counterpart to existing theory for Lasso estimation to achieve prediction errors in the order $O( s \frac{\log p}{n})$
under suitable conditions including compatibility conditions on $\Sigma$.

Third, from our error bounds depending on the ridge parameter $\tau$, we also derive the optimal orders of out-sample and in-sample MSEs
obtained respectively with the optimal choices of $\tau$ (see Table \ref{table:1c}). The optimal orders of prediction errors may be greater than  $O(\frac{d}{n})$.
Remarkably, we find that if $\frac{\lambda_{d+1}}{\lambda_d}$ is sufficiently small,
then the optimal out-sample MSE is, up to a constant factor, smaller than the optimal in-sample MSE in both the regime of small or moderate TER
(under some technical conditions)
and the regime of large TER.\ We also identify specific conditions under which the optimal out-sample MSE
is significantly smaller than the optimal in-sample MSE (see Remarks \ref{in-out-com-small_TER} and \ref{in-out-com-large_TER} for details, and Figure \ref{figure2} for numerical results).
This phenomenon seems to be surprising: out-sample MSEs may be usually considered to be no smaller than in-sample MSEs.

\vspace{.05in}
\textbf{Related works.}\;
There is a large and growing literature on prediction properties of ridge estimators and min-norm interpolators.
See \cite{Tsigler_ridge_2023}, Section 9, for a recent review.
We discuss directly related works to ours, in addition to the earlier discussion.

\cite{HsuKakadeZhang2011} allowed sub-gaussian covariate vectors and
studied the out-sample MSE of the ridge estimator when the ridge parameter $\tau$ is large enough such that the effective dimension,
$\sum_{j=1}^p \frac{\lambda_j}{\lambda_j + \tau}$, is small compared with the sample size.
For this reason, their error bounds are not applicable to a small ridge parameter or a min-norm interpolator.

\cite{HasteTrevorAndrea} derived out-sample error approximation formulas for the ridge estimator and the min-norm interpolator
using random matrix theory. They also showed that the deviation between the out-sample error and the approximation formula is upper bounded by the order of $n^{-\frac{1}{2}}$ for the ridge estimator (with a ridge tuning parameter bounded away from 0)
and is upper bounded by the order of $n^{-\frac{1}{7}}$ for the min-norm interpolator.
Compared to our results, the independent components assumption and boundedness of $\frac{p}{n}$ are assumed in \cite{HasteTrevorAndrea}.
Moreover, the orders of their deviation bounds may be much larger than $\frac{d}{n}$, so that combining the approximation formulas and the deviation bounds
may lead to less sharp out-sample error bounds than ours in the rotationally sparse setting.
Despite these differences, it can be shown that the orders of the approximation formulas in \cite{HasteTrevorAndrea} match the orders of our error bounds,
which are obtained without the independent components assumption or boundedness of $\frac{p}{n}$
in the rotationally sparse setting. See Section \ref{sec:4_1} for details.

\cite{Bartlett_2020} studied the min-norm interpolator and gave upper bounds of the out-sample error variance and bias
and a lower bound of the out-sample error variance (but not bias).
However, their results rely on the independent components assumption. Moreover, although the tail effective rank is involved,
their out-sample error bounds are obtained in terms of the overall $\Vert\theta^*\Vert^2$,
regardless of how the mean response is aligned differently with the spiked and tail parts of covariate vectors,
which are essential to the rotationally sparse setting.

\cite{Tsigler_ridge_2023} provided upper and lower bounds of both the out-sample error variance and bias,
while exploiting the decomposition of the spiked and tail parts of covariate vectors.
However, although the variance and bias upper bounds in \cite{Tsigler_ridge_2023} are obtained with sub-Gaussian covariate vectors
instead of the independent components assumption,
an oracle assumption is required on some random matrix from covariate vectors.
Moreover, their variance lower bound is obtained under the independent components assumption,
and their bias lower bound is provided in terms of the expectation with respect to a prior distribution on $\theta^*$
under an extra oracle assumption on covariate vectors.
By comparison, our error upper bounds match those in \cite{Tsigler_ridge_2023} for the ridge tuning parameter in suitable ranges,
and all our upper and lower bounds are obtained with sub-Gaussian covariate vectors without
making any oracle assumption or independent components assumption. See Section \ref{sec_4_2_main} for details.

\cite{BuneaFlorentinaStrimas2020} studied the min-norm interpolator in a latent factor model as follows:
\begin{align}
y_i=\beta^\T z_i +\xi_i, \quad x_i=Az_i + e_i, \quad i=1,\ldots,n, \label{eq:bunea-model}
\end{align}
where 
$\beta\in\mathbb{R}^d$, $A\in\mathbb{R}^{p\times d}$,
$z_i\in\mathbb{R}^d$ is a latent feature vector,
$\xi_i \in \mathbb{R}$ and $e_i\in \mathbb{R}^p$ are mean-zero noises,
and $(z_i, \xi_i, e_i)$ are mutually independent for each $i$.
The matrix $\Sigma= \Var(x_i)$ can be expressed as $\Sigma = A \Sigma_Z A^\T + \Sigma_E$,
where $\Sigma_Z = \Var(z_i)$ and $\Sigma_E = \Var(e_i)$.
The latent factor model can be seen to share a similar structure as our rotationally sparse linear regression model.
From our comparison in Section \ref{sec:4_3_main},
the upper bound of the out-sample MSE in \cite{BuneaFlorentinaStrimas2020} is obtained for the min-norm interpolator in the large TER regime,
and is less sharp than our result which gives the order of out-sample MSE
(i.e., matching upper and lower bounds up to a constant factor)
except in the trivial situation where the out-sample MSE is bounded away from zero.
In addition, the analysis of \cite{BuneaFlorentinaStrimas2020} assumes that
the whiten noises, $\Sigma_E^{-1/2} e_i$, has independent components.
Such an assumption of independent components is avoided in our analysis.

\section{Assumptions and notation}

We formulate the following assumptions to facilitate our theoretical analysis. Let $d$ be the dimension of the spike part satisfying $0<d<p$.

\begin{ass}[Low dimension of spiked part]
\label{ass:3} Suppose that $d\leq n$ and  $\frac{d}{n}$ is small enough such that
\begin{align*}
      \eta_1=  C_0 \sigma_x^2 \sqrt{\frac{d}{n}}   < \frac{1}{2} ,
   \end{align*}
\end{ass}
\noindent where $C_0$ is an absolute constant from Lemma \ref{lem_02}--\ref{lem_06}.

As shown in \cite{Bartlett_2020}, the tail effective rank (TER) is important for analyzing benign linear regression.
For the covariance matrix $\Sigma$, define
\begin{align*}
r_d(\Sigma)=\frac{\sum_{j>d}\lambda_j}{\lambda_{d+1}},
\end{align*}
where $\lambda_1 \geq \ldots \geq \lambda_p >0$ are the eigenvalues of $\Sigma$.
We refer to $r_d(\Sigma)$ as TER, because it pertains to the tail eigenvalues of $\Sigma$.
We also use the following related quantity:
\begin{align*}
r_d(\Sigma^2)=\frac{\sum_{j>d}\lambda_j^2}{\lambda_{d+1}^2}.
\end{align*}
It can be easily verified that $r_d(\Sigma^2)\leq r_d(\Sigma)$.

The following two assumptions describe two regimes of TER, in terms of the ratio $\frac{r_d(\Sigma)}{n}$.
The magnitude of $\frac{r_d(\Sigma)}{n}$ affects the behavior of the out-sample error and in-sample error.

\begin{ass}[Small or moderate TER]
\label{ass:3a}
$\frac{r_d(\Sigma)}{n} \leq C_1$, where $C_1>0$ is a constant.
\end{ass}
\begin{ass}[Large TER]
\label{ass:4}
$\frac{r_d(\Sigma)}{n}$ is large enough such that
\begin{align*}
   \eta_2= C_0\sigma_x^2\sqrt{\frac{4n^2}{r_d(\Sigma)^2}+\frac{2n}{r_d(\Sigma)}}\leq \frac{1}{2}.
\end{align*}
Alternatively, it is sufficient to assume that
$  \frac{r_d(\Sigma)}{n} \geq c_x $
for some $c_x$ depending only on $\sigma_x$. For instance, $c_x$ can be $\max\{4\sqrt{2}C_0\sigma_x^2,16C^2_0\sigma_x^4\}$.
\end{ass}

Separating the two regimes above is desirable for theoretical analysis, because
it enables us to establish concrete results and avoid making
any oracle assumption or independent components assumption on covariate vectors as used in \cite{Tsigler_ridge_2023}.
See Section \ref{sec_4_2_main} for further information.
The two regimes above are not contained by each other. For instance, $\frac{r_d(\Sigma)}{n}$ satisfies Assumption \ref{ass:3a} but not Assumption \ref{ass:4} if $\frac{r_d(\Sigma)}{n}\ll 1$,
whereas $\frac{r_d(\Sigma)}{n}$ satisfies Assumption \ref{ass:4} but not Assumption \ref{ass:3a} if $\frac{r_d(\Sigma)}{n}\gg 1$. Overlapping of the two regimes is possible, for instance, $\frac{r_d(\Sigma)}{n} \asymp 1$ and $\frac{r_d(\Sigma)}{n}$ satisfies both Assumptions \ref{ass:3a} and \ref{ass:4}.

The following assumption describes the rotationally sparse setting in terms of the relative magnitudes of
$\Vert\theta^*_{(d+1):p}\Vert^2_{\Sigma_{(d+1):p}}$ and $\Vert\theta^*_{1:d}\Vert^2_{\Sigma_{1:d}^{-1}}$.

\begin{ass}[Rotational Sparsity]
\label{ass1}

$\newline$
(i) [Applied with small or moderate TER]. For some $0<\delta_{1}< 1$,
\begin{align*}
\frac{\Vert\theta^*_{(d+1):p}\Vert^2_{\Sigma_{(d+1):p}}}{\Vert\theta^*_{1:d}\Vert^2_{\Sigma_{1:d}^{-1}}}  \leq  \frac{\delta_{1}}{4(1+\sigma_x^2)} \lambda_{d+1}^2 .
\end{align*}

\noindent (ii) [Applied with large TER]. For some $0<\delta_{2}< 1$,
\begin{align*}
\frac{\Vert\theta^*_{(d+1):p}\Vert^2_{\Sigma_{(d+1):p}}}{\Vert\theta^*_{1:d}\Vert^2_{\Sigma_{1:d}^{-1}}}  \leq    \frac{\delta_{2}}{4(1+\sigma_x^2)} (\frac{1}{\lambda_d}+ \frac{4n}{\sum_{j>d}\lambda_j})^{-2} ,
\end{align*}
which can be equivalently stated as
\begin{align*}
\frac{\Vert\theta^*_{(d+1):p}\Vert^2_{\Sigma_{(d+1):p}}}{\Vert\theta^*_{1:d}\Vert^2_{\Sigma_{1:d}^{-1}}}  \leq
 \frac{\delta_{2}}{4(1+\sigma_x^2)} (\frac{\lambda_{d+1}}{\lambda_d}+ \frac{4n}{r_d(\Sigma)} )^{-2} \lambda_{d+1}^2 .
\end{align*}
\end{ass}
Assumption \ref{ass1}(i) is specified for the small or moderate regime, whereas Assumption \ref{ass1}(ii) is specified for the regime large TER regime.
Assumption \ref{ass1}(ii) provides a much weaker condition than Assumption \ref{ass1}(i) on
$\frac{\Vert\theta^*_{(d+1):p}\Vert^2_{\Sigma_{(d+1):p}}}{\Vert\theta^*_{1:d}\Vert^2_{\Sigma_{1:d}^{-1}}}$
if $ \lambda_d \gg \lambda_{d+1}$ and $ r_d(\Sigma) \gg n$.

\vspace{.05in}
\textbf{Notation.} Given two positive sequence $\{a_k\}$ and $\{b_k\}$,  $a_k\lesssim b_k$ ($a_k\gtrsim b_k$) indicates that there exist constants $c>0$ and $K\geq 1$ such that $a_k \leq c b_k$ ($a_k\geq cb_k$) for all $k\geq K$. We also denote $a_k = O(b_k)$ if $a_k \lesssim b_k$.
Moreover, $a_k\asymp b_k$ indicates that both $a_k\lesssim b_k$ and $a_k \gtrsim b_k$;
and $a_k \ll b_k$ (or $a_k \gg b_k $) indicates $\lim_{k\rightarrow \infty} \frac{a_k}{b_k} = 0$ (or $\lim_{k\rightarrow \infty} \frac{a_k}{b_k} = \infty$).
Finally, $\poly_{deg}(x)$ denotes a polynomial of $x$ with positive bounded coefficients and the highest order equal to $deg$.

\section{Main results}
\label{main_results}

The standard formula of the ridge estimator for $n>p$ is
\begin{align*}
    \hat{\theta}(\tau) = (X^TX + n\tau)^{-1} X^T Y.
\end{align*}
However, this formula does not hold for $\tau=0$ when $X^TX$ is not invertible with $n<p$. In the high-dimensional setting,
the ridge estimator $\hat{\theta}(\tau)$ can be expressed as
\begin{align}
    \hat{\theta}(\tau) = X^T (
    XX^\T+n\tau)^{-1}Y  .  \label{eq:3}
\end{align}
See, for example, Appendix B in \cite{Tsigler_ridge_2023}.

With the expression (\ref{eq:3}) and the assumption that $X$ and $\epsilon$ are independent, the out-sample and in-sample errors can be decomposed into bias and variance as follows:
\begin{itemize}
    \item Out-sample error
    \small
    \begin{align}
        &\textrm{MSE}_{\out} = \underbrace{\Vert(I_p - X^\T( XX^\T+n\tau I_n)^{-1}X)\theta^*\Vert_{\Sigma}^2}_{\B_{\out}} \notag \\
        &~~~~~~~~~~~~~~~~~+   \underbrace{\sigma^2\Tr(( XX^\T+n\tau I_n )^{-1}X\Sigma X^\T(XX^\T+n\tau I_n)^{-1})}_{\V_{\out}} , \label{eq:3.3}
    \end{align}
    \normalsize
    \item In-sample error
    \small
    \begin{align}
        &\textrm{MSE}_{\In} = \underbrace{\Vert(I_p - X^\T( XX^\T+n\tau I_n)^{-1}X)\theta^*\Vert_{\hat{\Sigma}}^2}_{\B_{\In}} \notag \\
        &~~~~~~~~~~~~~~~~~+   \underbrace{\sigma^2\Tr(( XX^\T+n\tau I_n)^{-1}X\hat{\Sigma} X^\T(XX^\T+n\tau I_n)^{-1})}_{\V_{\In}} . \label{eq:3.4}
    \end{align}
    \normalsize
\end{itemize}
We present our main results about the out-sample and in-sample errors in the small or moderate TER regime in Section \ref{main_results_1}
and the large TER regime in Section \ref{main_results_2}.

\subsection{Regime I: Small or moderate TER}
\label{main_results_1}

Consider the regime of small or moderate TER such that $r_d(\Sigma)\leq C_1 n$, as stated in Assumption \ref{ass:3a}.
The following is our main result about $\MSE_{\out}$, the out-sample MSE.
Let $\A_0$ be a constant satisfying $A_0\geq 1$ and $(1+\A_0)^2<4\delta_1^{-1}$,
with $0 < \delta_1 <1$ from Assumption \ref{ass1}(i).
For $\A_0^{-1}\lambda_{d+1}\leq \tau\leq \A_0 \lambda_{d}$, we determine the order of $\MSE_{\out}$ (including upper and lower bounds).
For $\tau \leq \A_0^{-1} \lambda_{d+1}$ and $\tau \geq \A_0 \lambda_d$, we give lower bounds of $\MSE_{\out}$ through, respectively, the variance and bias terms.

\begin{thm}[Out-sample error with small or moderate TER]
\label{thm1}
Under Assumption \ref{ass:3}, \ref{ass:3a} and \ref{ass1}(i), for any $\nu$ satisfying $0<\nu < \frac{1}{2}\min\{1,\sigma_x^2\}$ and any $\A_0$ satisfying $A_0\geq 1$ and $(1+\A_0)^2<4\delta_1^{-1}$,
the following inequalities hold uniformly in the range of $\tau$ stated
with probability at least $1-2\exp\{-\frac{\nu^2 n}{C_0^2\sigma_x^4}\}-2\exp\{-\frac{\nu^2 n}{ C_0\sigma_x^4}\}-18\exp\{-\frac{n}{C_0}\}$:
\begin{align*}
 & (i)~ \MSE_{\out} \geq M_1 \underbrace{\sigma^2(\frac{d}{n}+ \frac{r_d( \Sigma^2) }{n })}_{\underline{\V}_\out} \quad \text{for } \tau \leq \A_0^{-1} \lambda_{d+1} ,  \\
 & (ii)~ M_2(\underbrace{\Vert\theta_{1:d}^*\Vert^2_{\Sigma_{1:d}^{-1}} \tau^2}_{\overline{\B}_\out } + \underbrace{\sigma^2(\frac{d}{n}+\frac{\lambda_{d+1}^2}{\tau^2}\frac{r_d( \Sigma^2)}{n}  )}_{\overline{\V}_\out } ) \geq \MSE_{\out} \geq \notag\\  &~~~~~~~~~~~~~~~~~~~~~~M_1(\underbrace{\Vert\theta_{1:d}^*\Vert^2_{\Sigma_{1:d}^{-1}} \tau^2}_{\underline{\B}_\out } + \underbrace{\sigma^2(\frac{d}{n}+\frac{\lambda_{d+1}^2}{\tau^2}\frac{r_d( \Sigma^2)}{n }  )}_{\underline{\V}_\out })
  \quad \text{for } \A_0^{-1}\lambda_{d+1} \leq \tau \leq \A_0\lambda_d,  \\
 & (iii)~\MSE_{\out} \geq M_1 \underbrace{\Vert\theta_{1:d}^*\Vert^2_{\Sigma_{1:d}^{-1}}\lambda_d^2}_{\underline{\B}_\out } \quad \text{for }  \tau \geq \A_0 \lambda_d ,
\end{align*}
where
$M_1,M_2>0$ are constants, depending only on $(\sigma_x, \eta_1, C_1, \delta_1, \nu, \A_0)$,
and $\underline{B}_{\out}$ and  $\underline{V}_{\out}$ represent lower bounds of out-sample bias and variance
and $\overline{B}_{\out}$ and $\overline{V}_{\out}$ represent upper bounds of out-sample bias and variance,
up to the constant $M_1$ or $M_2$.
\end{thm}

The following corollary provides simple conditions for achieving $\MSE_{\out}=O(\frac{d}{n})$ in the regime of small or moderate TER.

\begin{cor}[Conditions for $\MSE_{\out} = O(\frac{d}{n})$ with small or moderate TER]
\label{cor1}
In the setting of Theorem 1, assume further
that $\sigma^2\asymp 1$ and $\Vert\theta_{1:d}^*\Vert^2_{\Sigma_{1:d}^{-1}}\lambda_d^2\asymp 1$.

(i) A sufficient condition for $\MSE_{\out} = O(\frac{d}{n})$ with a probability approaching 1 as $n\to\infty$ is that $\frac{\lambda_{d+1}}{\lambda_{d}}\lesssim  \sqrt{\frac{d}{n}}\min\{1,\sqrt{\frac{d}{r_d(\Sigma^2)}}\}$ and
the ridge parameter $\tau$ is chosen in the range $\A_0^{-1}\lambda_{d+1} \leq \tau \leq \A_0 \lambda_{d+1}$ if
 $r_d(\Sigma^2) \leq d$ or $\A_0^{-1}\lambda_{d+1}\max\{\frac{1}{c}\sqrt{\frac{r_d( \Sigma^2)}{d}},1\} \leq \tau \leq \A_0\lambda_{d}\min\{c\sqrt{\frac{d}{n}}, 1\}$ if $r_d(\Sigma^2) > d$, where $c$ is a constant satisfying $c\geq 1$ and $\frac{\lambda_{d+1}}{\lambda_{d}}\leq c  \sqrt{\frac{d}{n}}\min\{1,\sqrt{\frac{d}{r_d(\Sigma^2)}}\}$.

(ii) Suppose that $n\gg d$ and $r_d( \Sigma^2) \gg d$. Then a necessary condition for $\MSE_{\out} = O(\frac{d}{n})$ with a probability bounded away from 0
is that $\frac{\lambda_{d+1}}{\lambda_d} \lesssim \sqrt{\frac{d}{n}}\sqrt{\frac{d}{r_d(\Sigma^2)}}$ and
the ridge parameter $\tau$ is chosen in the range $\sqrt{\frac{r_d( \Sigma^2)}{d} } \lambda_{d+1} \lesssim \tau \lesssim \sqrt{\frac{d}{n}}\lambda_d$.

\noindent The sufficient and necessary conditions become matched,  $\frac{\lambda_{d+1}}{\lambda_d} \lesssim \sqrt{\frac{d}{n}}\sqrt{\frac{d}{r_d(\Sigma^2)}}$,
in the case where $n\gg d$ and $r_d(\Sigma^2)\gg d$ in addition to the assumptions stated.
\end{cor}

Next, we give our main result about $\MSE_{\In}$, the in-sample MSE, in the regime of small or moderate TER stated in Assumption \ref{ass:3a}.
As before, let $\A_0$ be a constant satisfying $A_0\geq 1$ and $(1+\A_0)^2<4\delta_1^{-1}$, with $0 < \delta_1 <1$ from Assumption \ref{ass1}(i). For $\A_0^{-1}\lambda_{d+1}\leq \tau\leq \A_0 \lambda_{d}$, we derive
upper and lower bound of $\MSE_{\In}$.
For $\tau \leq \A_0^{-1} \lambda_{d+1}$, we give a lower bound through variance. For $\tau \geq \A_0 \lambda_d$, we give a lower bound of $\MSE_{\In}$ through
the sum of bias and variance terms.

\begin{thm}[In-sample error with small or moderate TER]
\label{thm2}
Under Assumption \ref{ass:3}, \ref{ass:3a} and \ref{ass1}(i), for any $\nu$ satisfying $0<\nu <\frac{1}{4}\min\{1,\sigma_x^2\}$ and any $A_0$  satisfying $A_0\geq 1$ and $(1+\A_0)^2<4\delta_1^{-1}$, the following inequalities hold  uniformly in the range of $\tau$ stated
with probability at least $1-2\exp\{-\frac{\nu^2 n}{C_0^2\sigma_x^4}\}-2 \exp\{-\frac{\nu^2 n}{C_0\sigma_x^4} \}-8\exp\{-\frac{n}{C_0}\}$:
\begin{align*}
   & (i)~\MSE_{\In} \geq M_1 \underbrace{\sigma^2(\frac{d}{n}+\frac{r_d^2(\Sigma)}{n^2})}_{\underline{V}_\In} \quad \text{for }   \tau \leq \A_0^{-1} \lambda_{d+1},   \\
   &(ii)~M_2(\underbrace{\Vert\theta^*_{1:d}\Vert^2_{\Sigma_{1:d}^{-1}} \tau^2}_{\overline{\B}_\In } +  \underbrace{\sigma^2(\frac{d}{n}+\frac{\lambda^2_{d+1}}{\tau^2}\frac{r_d(\Sigma)}{n})}_{\overline{\V}_\In } )  \geq \MSE_{\In} \geq \notag \\ &~~~~~~~~~~M_1(\underbrace{\kappa_1(\tau)\Vert\theta^*_{1:d}\Vert^2_{\Sigma_{1:d}^{-1}} \tau^2}_{\underline{\B}_\In } +  \underbrace{\sigma^2(\frac{d}{n} +\frac{\lambda_{d+1}^2}{\tau^2}\frac{r_d^2(\Sigma)}{n^2} )}_{\underline{\V}_\In } )\quad \text{for }  \A_0^{-1} \lambda_{d+1} \leq \tau \leq \A_0 \lambda_d,  \\
   &(iii)~\MSE_{\In} \geq M_1(\underbrace{ \kappa_1(\tau) \Vert\theta^*_{1:d}\Vert^2_{\Sigma_{1:d}^{-1}} \lambda_d^2}_{\underline{\B}_\In } +  \underbrace{\sigma^2\frac{\lambda_{d+1}^2}{\tau^2}\frac{r_d^2(\Sigma)}{n^2} }_{\underline{\V}_\In } )  \quad \text{for }   \tau \geq \A_0 \lambda_d ,
\end{align*}
where $\kappa_1(\tau)=\max\{1-(\frac{2 C_0\sigma_x^2(2+C_1)\lambda_{d+1}}{\tau}(1+16(2C_0\sigma_x^2+1)(1+C_1)\frac{\sqrt{\delta_1}}{1-\sqrt{\delta_1}}) + 64\frac{\sqrt{\delta_1}}{1-\sqrt{\delta_1}}),0\}$, and $M_1,M_2>0$ are constants depending only on $(\sigma_x, \eta_1, C_1, \delta_1, \nu, \A_0)$. The terms $\underline{B}_{\In}$ and  $\underline{V}_{\In}$ represent lower bounds of in-sample bias and variance
and $\overline{B}_{\In}$ and $\overline{V}_{\In}$ represent upper bounds of in-sample bias and variance,
up to the constant $M_1$ or $M_2$.
\end{thm}

By the definition of $\kappa_1(\tau)$, the bias term $\underline{\B}_\In$ is activated in the lower bound of $\MSE_{\In}$ only when $\frac{\lambda_{d+1}}{\tau} $ and $\frac{\sqrt{\delta_1}}{1-\sqrt{\delta_1}}$ are small enough.
This can be explained from our proof strategy as follows (see Section \ref{sec_6_2} for details).
The in-sample bias $\Vert \hat{\theta}(\tau)-\theta^*\Vert^2_{\hat{\Sigma}}$ can be expressed as the sum of $\Vert \hat{\theta}(\tau)_{1:d}-\theta_{1:d}^*\Vert^2_{\hat{\Sigma}_{1:d}}$, $2 (\hat{\theta}(\tau)_{1:d}-\theta_{1:d}^*)\hat{\Sigma}_{(1:d),(d+1):p}(\hat{\theta}(\tau)_{(d+1):p}-\theta_{(d+1):p}^*)$ and $\Vert \hat{\theta}(\tau)_{(d+1):p}-\theta_{(d+1):p}^*\Vert^2_{\hat{\Sigma}_{(d+1):p}}$ and only the interaction term $2 (\hat{\theta}(\tau)_{1:d}-\theta_{1:d}^*)\hat{\Sigma}_{(1:d),(d+1):p}\cdot(\hat{\theta}(\tau)_{(d+1):p}-\theta_{(d+1):p}^*)$ can be negative.
When $\frac{\lambda_{d+1}}{\tau}$ and $\frac{\sqrt{\delta_1}}{1-\sqrt{\delta_1}}$ are small enough,
the bias from the spiked part, $\Vert \hat{\theta}(\tau)_{1:d}-\theta_{1:d}^*\Vert^2_{\hat{\Sigma}_{1:d}}$,
can be shown to dominate the interaction term.
Then a lower bound on the bias from the spiked part, which can be deduced in a convenient manner,
also provides a lower bound on the overall bias.

From Theorem \ref{thm2}, we deduce the following simple conditions for achieving $\MSE_{\In}=O(\frac{d}{n})$ in the regime of small or moderate TER.

\begin{cor}[Conditions for $\MSE_{\In} = O(\frac{d}{n})$ with small or moderate TER]
\label{cor2}
In the setting of Theorem \ref{thm2}, assume further
that $\sigma^2\asymp 1$ and $\Vert\theta_{1:d}^*\Vert^2_{\Sigma_{1:d}^{-1}}\lambda_d^2\asymp 1$.

(i) A sufficient condition for $\MSE_{\In}=O(\frac{d}{n})$ with a probability approaching 1 as $n\to\infty$ is that $\frac{\lambda_{d+1}}{\lambda_{d}} \lesssim  \sqrt{\frac{d}{n}}\min\{1, \sqrt{\frac{d}{r_d(\Sigma)}}\}$ and
the ridge parameter $\tau$ is chosen in the range $\A_0^{-1}\lambda_{d+1} \leq \tau \leq \A_0 \lambda_{d+1}$ if
 $r_d(\Sigma) \leq d$ or $\A_0^{-1}\lambda_{d+1}\max\{\frac{1}{c}\sqrt{\frac{r_d( \Sigma)}{d}},1\} \leq \tau \leq \A_0\lambda_{d}\min\{c\sqrt{\frac{d}{n}}, 1\}$ if $r_d(\Sigma) > d$, where $c$ is a constant satisfying $c\geq 1$ and $\frac{\lambda_{d+1}}{\lambda_{d}}\leq c   \sqrt{\frac{d}{n}}\min\{1,\sqrt{\frac{d}{r_d(\Sigma)}}\}$.

(ii) Suppose that $n\gg d$, $ \frac{r_d(\Sigma)}{n} \sqrt{\frac{n}{d}}\gg 1$  and
$64\frac{\sqrt{\delta_1}}{1-\sqrt{\delta_1}}<1$. Then a necessary condition for $\MSE_{\In} = O(\frac{d}{n})$ with a probability bounded away from 0 is that $\frac{\lambda_{d+1}}{\lambda_{d}}\lesssim  \frac{d}{r_d(\Sigma)}$ and the ridge parameter $\tau$ is chosen in the range $\lambda_{d+1}\frac{r_d(\Sigma)}{n}\sqrt{\frac{n}{d}} \lesssim \tau \lesssim \lambda_d\sqrt{\frac{d}{n}}$.

\noindent The sufficient and necessary conditions become matched, $\frac{\lambda_{d+1}}{\lambda_{d}}\lesssim  \frac{d}{r_d(\Sigma)}$,
in the case where $n\gg d$ and $r_d(\Sigma)\asymp n$ in addition to the assumptions stated.
\end{cor}

From Theorems \ref{thm1} and \ref{thm2}, we derive the order of $\MSE_\out$ with an optimal choice $\tau$, denoted as $\MSE_\out^*$, and the order of $\MSE_\In$ with an optimal choice $\tau$, denoted as $\MSE_\In^*$. The following corollary gives the orders of $\MSE_{\out}^*$ and $\MSE_{\In}^*$ in the small or moderate TER regime.

\begin{cor}[Optimal error orders with small or moderate TER]
\label{cor5a}
Suppose that Assumption \ref{ass:3}, \ref{ass:3a} and \ref{ass1}(i) are satisfied and further $\sigma^2\asymp 1$, $\Vert\theta_{1:d}^*\Vert^2_{\Sigma_{1:d}^{-1}}\lambda_d^2\asymp 1$, $r_d(\Sigma)\asymp n$,
$\lambda_d \gtrsim  \lambda_{d+1}\sqrt{\frac{n}{r_d(\Sigma^2)}}$, $\lambda_d \gg \lambda_{d+1}$, and $64\frac{\sqrt{\delta_1}}{1-\sqrt{\delta_1}}<1$. Then
\newline
(i) $\MSE_{\out}^*\asymp\max\{\frac{\lambda_{d+1}}{\lambda_d}\sqrt{\frac{r_d(\Sigma^2)}{n}}, \frac{d}{n}\}$ with a probability approaching to 1 and the optimal $\tau$ is chosen as $\tau =  \sqrt{\lambda_d\lambda_{d+1}\sqrt{\frac{r_d(\Sigma^2)}{n}}}\min\{\sqrt{cA_0^{-2}}, \frac{A_0\lambda_d}{\sqrt{\lambda_d\lambda_{d+1}\sqrt{\frac{r_d(\Sigma^2)}{n}}}}\}$ where $c$ is a constant satisfying $\lambda_{d+1}\sqrt{\frac{n}{r_d(\Sigma^2)}}\leq c \lambda_{d}$.
\newline
(ii) $\MSE^*_{\In}\asymp \max\{\frac{\lambda_{d+1}}{\lambda_d},\frac{d}{n}\}$ with a probability approaching to 1 and the optimal $\tau$ is chosen as $\tau \asymp  \sqrt{\lambda_{d+1}\lambda_d}$.

\noindent Therefore $\MSE_{\out}^* \lesssim \MSE_{\In}^*$ with a probability approaching to 1, by noting $r_d(\Sigma^2) \le r_d(\Sigma) \asymp n$.
\end{cor}

The additional conditions $r_d(\Sigma)\asymp n$, $\lambda_d\gtrsim \lambda_{d+1}\sqrt{\frac{n}{r_d(\Sigma^2)}}$ and $\lambda_d\gg \lambda_{d+1}$ can be explained as follows.
First, $r_d(\Sigma)\asymp n$ and $\lambda_d\gg \lambda_{d+1}$ are important for determining the order of $\MSE_\In^*$. In fact,
the order of $\V_\In$ for $A_0^{-1}\lambda_{d+1}\leq \tau\leq A_0\lambda_d$ can be determined from Theorem \ref{thm2}(ii)
only under $r_d(\Sigma)\asymp n$, due to the difference between $\underline{\V}_\In$ and $\overline{\V}_\In$.
Note that
the sum of the in-sample bias and the tail part of in-sample variance,
$\Vert\theta^*_{1:d}\Vert^2_{\Sigma_{1:d}^{-1}} \tau^2+\sigma^2\frac{\lambda_{d+1}^2}{\tau^2}$,
reaches the minimum order of $\frac{\lambda_{d+1}}{\lambda_d}$ by the choice $\tau \asymp \sqrt{\lambda_d\lambda_{d+1}}$.
The condition $\lambda_d\gg \lambda_{d+1}$ ensures that this choice of $\tau$
is large enough so that $\kappa_1(\tau)$ is activated.
Second, $\lambda_d\gtrsim \lambda_{d+1}\sqrt{\frac{n}{r_d(\Sigma^2)}}$ is important for determining the order of $\MSE_\out^*$ because under $\lambda_d\gtrsim \lambda_{d+1}\sqrt{\frac{n}{r_d(\Sigma^2)}}$, the order of $\V_\out$ for
$\tau\leq A_0^{-1}\lambda_{d+1}$ can be shown to be larger than $\max\{\frac{\lambda_{d+1}}{\lambda_d}\sqrt{\frac{r_d(\Sigma^2)}{n}}, \frac{d}{n}\}$
and then the range $\tau\leq A_0^{-1}\lambda_{d+1}$ can be ruled out
when optimizing $\MSE_{\out}$. See the proof of Corollary \ref{cor5a} in Supplement Section \ref{sec_II_2} for details.

\begin{rem}
\label{in-out-com-small_TER}

In the setting of Corollary \ref{cor5a}, we observe that the gap between $\MSE_\In^*$ and $\MSE_\out^*$ can be significantly large, for example,
if further $\frac{\lambda_{d+1}}{\lambda_d}\gtrsim \frac{d}{n}\sqrt{\frac{n}{r_d(\Sigma^2)}}$ and $\frac{n}{r_d(\Sigma^2)}\gg 1$. In this case,
both $\MSE_\out$ and $\MSE_\In$ do not achieve $O(\frac{d}{n})$, and
 $\MSE_{\out}^*\asymp \frac{\lambda_{d+1}}{\lambda_d}\sqrt{\frac{r_d(\Sigma^2)}{n}}$ and
 $\MSE^*_{\In}\asymp \frac{\lambda_{d+1}}{\lambda_d}$ by Corollary \ref{cor5a}.
With $\frac{n}{r_d(\Sigma^2)}\gg 1$, there can be a substantial gap between $\MSE_\In^*$ and $\MSE_\out^*$.
\end{rem}

In the setting of Corollary \ref{cor5a}, we point out that
the advantage of $\MSE_\out^*$ over $\MSE_\In^*$
can be attributed to $\V_\out \lesssim  \V_\In$ for $A_0^{-1}\lambda_{d+1}\leq \tau\leq A_0\lambda_d$ under $r_d(\Sigma)\asymp n$.
See the proof of Corollary \ref{cor5a} in Supplement Section \ref{sec_II_2} for details.
In fact, in the setting of Corollary \ref{cor5a}, the optimal choices of $\tau$ for both $\MSE_\out$ and $\MSE_\In$ are chosen from the range $A_0^{-1}\lambda_{d+1}\leq \tau \leq A_0\lambda_d$, because both $\MSE_\out$ and $\MSE_\In$ are lower bounded by large variance for $\tau\leq A_0^{-1}\lambda_{d+1}$ and by large bias for $\tau\geq A_0\lambda_d$. Moreover, as seen from the proof,
the optimal order of $\MSE_\In$, $\max\{\frac{\lambda_{d+1}}{\lambda_d},\frac{d}{n}\}$, in the setting of Corollary \ref{cor5a}
can be achieved only when $\underline{\B}_\In$ is activated. With $\underline{\B}_\In$ activated, the orders of $\B_\out$ and $\B_\In$ are the same for $A_0^{-1}\lambda_{d+1}\leq \tau\leq A_0\lambda_d$.
By comparison, $\V_\out$ can be shown to be smaller than $\V_\In$ up to a constant factor
for $A_0^{-1}\lambda_{d+1}\leq \tau \leq A_0 \lambda_d$ under $r_d(\Sigma)\asymp n$:
 \begin{align}
     \underbrace{\sigma^2(\frac{d}{n}+ \frac{\lambda_{d+1}^2}{\tau^2}\frac{r_d(\Sigma^2)}{n})}_{\text{order~of~$\V_\out$~in~Theorem~\ref{thm1}(ii)}} \lesssim \underbrace{\sigma^2(\frac{d}{n}+\frac{\lambda_{d+1}^2}{\tau^2})}_{\text{order~of~$\V_\In$~in~Theorem~\ref{thm2}(ii)}} .  \notag
 \end{align}
Hence the advantage of $\MSE_\out^*$ stems from the smaller order of $\V_\out$ in the setting of Corollary \ref{cor5a}.

\subsection{Regime II: Large TER}
\label{main_results_2}

The second regime we investigate is when $\frac{r_d(\Sigma)}{n}$ is large enough such that Assumption \ref{ass:4} is satisfied. In this regime,
as shown below,
the smallest eigenvalue of
$\tau I_n + n^{-1} X_{(d+1):p}X_{(d+1):p}^\T $ is lower bounded away from 0 for any $\tau \ge 0$, so that $\MSE_\out$ can
sometimes be controlled even for $\tau=0$. Let $A_0$ be any positive constant. For the ridge parameter $\tau\ge 0$ satisfying
$\tau +\lambda_{d+1}\frac{r_d(\Sigma)}{n}\leq A_0 \lambda_d$, we
determine the order of $\MSE_{\out}$ (including upper and lower bounds). For $\tau +\lambda_{d+1}\frac{r_d(\Sigma)}{n}\geq A_0 \lambda_d$, we give a lower bound of $\MSE_{\out}$ through the bias term.

 \begin{thm}[Out-sample error with large TER]
\label{thm3}
Under Assumption \ref{ass:3}, \ref{ass:4} and \ref{ass1}(ii), for any $\nu$ satisfying $0< \nu < \frac{1}{4}\min\{1,\sigma_x^2\}$ and $\frac{r_d(\Sigma)\nu^2}{C_0^2\sigma_x^4}>1$ and any $A_0>0$,
the following inequalities hold  uniformly in the range of $\tau$ stated
with probability at least $1-2n \exp\{-\frac{\nu \sqrt{r_d(\Sigma)}}{C_0\sigma_x^2}\}-2\exp\{-\frac{\nu^2 n}{C_0^2\sigma_x^4}\}-2\exp\{-\frac{\nu^2 n}{C_0\sigma_x^4}\}-16 \exp\{-\frac{n}{C_0}\}$:
\begin{align}
& (i)~M_2 (\underbrace{\Vert\theta_{1:d}^*\Vert^2_{\Sigma_{1:d}^{-1}}(\tau+ \lambda_{d+1}\frac{r_d(\Sigma)}{n})^2}_{\overline{\B}_\out }
+ \underbrace{\sigma^2(\frac{d}{n}+ \frac{\lambda_{d+1}^2}{(\tau+\lambda_{d+1}\frac{r_d(\Sigma)}{n})^2}\frac{r_d(\Sigma^2)}{n}) }_{\overline{\V}_\out }) \geq \MSE_{\out} \geq \notag \\
&~~~M_1 ( \underbrace{\Vert\theta_{1:d}^*\Vert^2_{\Sigma_{1:d}^{-1}}(\tau+ \lambda_{d+1}\frac{r_d(\Sigma)}{n})^2}_{\underline{\B}_\out }
+ \underbrace{\sigma^2(\frac{d}{n}+ \frac{\lambda_{d+1}^2}{(\tau+\lambda_{d+1}\frac{r_d(\Sigma)}{n})^2}\frac{r_d(\Sigma^2)}{n} )  }_{\underline{\V}_\out }) \quad \text{for } \tau +\lambda_{d+1}\frac{r_d(\Sigma)}{n} \leq A_0\lambda_d, \notag \\
&(ii)~ \MSE_{\out} \geq M_1 \underbrace{\Vert\theta_{1:d}^*\Vert^2_{\Sigma_{1:d}^{-1}}\lambda_d^2}_{\underline{\B}_\out }, \quad \text{for }  \tau+\lambda_{d+1}\frac{r_d(\Sigma)}{n} \geq A_0 \lambda_d ,   \notag
\end{align}
where $M_1,M_2>0$ are constants depending only on $(\sigma_x, \eta_2, \delta_2, \nu, \A_0)$.
\end{thm}

\begin{rem}
\label{rem_small_large_TER}
The error bounds in Theorems \ref{thm1} and \ref{thm3} are derived from the same set of algebraic bounds
but then by applying relevant high-probability inequalities to control random quantities for different ranges of $\tau$
under the two regimes of TER (see Section \ref{sec_6_1}).
In fact, after ignoring the range choice of $\tau$, the error bounds (i)--(ii) in Theorem \ref{thm3} appear similar to (ii)--(iii) in Theorem \ref{thm1}
except with $\tau+ \lambda_{d+1}\frac{r_d(\Sigma)}{n}$ in place of $\tau$.
In the overlapping case of small or moderate TER regime and large TER regime, i.e., $\frac {r_d(\Sigma)}{n}\asymp 1$
and $\frac {r_d(\Sigma)}{n}$ satisfies Assumption \ref{ass:3a} and \ref{ass:4},
the result from Theorem \ref{thm1} (ii), $\MSE_\out \asymp \Vert\theta_{1:d}^*\Vert^2_{\Sigma_{1:d}^{-1}} \tau^2+\sigma^2(\frac{d}{n}+\frac{\lambda_{d+1}^2}{\tau^2}\frac{r_d( \Sigma^2)}{n } )$, and the result from Theorem \ref{thm3} (i), $\MSE_\out \asymp \Vert\theta_{1:d}^*\Vert^2_{\Sigma_{1:d}^{-1}}(\tau+ \lambda_{d+1}\frac{r_d(\Sigma)}{n})^2+\sigma^2(\frac{d}{n}+ \frac{\lambda_{d+1}^2}{(\tau+\lambda_{d+1}\frac{r_d(\Sigma)}{n})^2}\frac{r_d(\Sigma^2)}{n} )$,
are equivalent to each other for $\tau$ in the range $A_0^{-1}\lambda_{d+1} \leq \tau \leq A_0\lambda_d$,
where $\tau \asymp \tau+ \lambda_{d+1}\frac{r_d(\Sigma)}{n}$ with $r_d(\Sigma)\asymp n$.
However, the error bound in Theorem \ref{thm3} (i) remains applicable, but that in Theorem \ref{thm1} (ii) does not apply, to
small $\tau$ satisfying $0\le \tau \le A_0^{-1} \lambda_{d+1}$ including $\tau=0$.
\end{rem}

The following corollary provides simple conditions for achieving $\MSE_{\out} = O(\frac{d}{n})$ in the regime of large TER.

\begin{cor}[Conditions for $\MSE_{\out} = O(\frac{d}{n})$ with large TER]
\label{cor3}
In the setting of Theorem \ref{thm3}, assume further that $\sigma^2\asymp 1$ and
$\Vert\theta_{1:d}^*\Vert^2_{\Sigma_{1:d}^{-1}}\lambda_d^2\asymp 1$.

(i) A sufficient condition for $\MSE_{\out}=O(\frac{d}{n})$ with a probability approaching 1 as $n\to\infty$ is
that $\frac{\lambda_{d+1}}{\lambda_{d}} \lesssim  \sqrt{\frac{d}{n}} \min \{ \sqrt{\frac{d}{r_d(\Sigma^2)}}, \frac{n}{r_d(\Sigma)}\}$
and the ridge parameter $\tau$ is chosen satisfying
$\tau + \lambda_{d+1}\frac{r_d(\Sigma)}{n} \lesssim \sqrt{\frac{d}{n}}\lambda_d$ if $\frac{n\sqrt{r_d(\Sigma^2)}}{\sqrt{d} r_d(\Sigma)}\leq 1$ or $\sqrt{\frac{r_d( \Sigma^2)}{d} } \lambda_{d+1} \lesssim \tau +\lambda_{d+1}\frac{r_d(\Sigma)}{n}\lesssim \sqrt{\frac{d}{n}}\lambda_{d}$ if $\frac{n\sqrt{r_d(\Sigma^2)}}{\sqrt{d} r_d(\Sigma)}> 1$.

(ii) Suppose that $n\gg d$. Then a necessary condition for $\MSE_{\out}=O(\frac{d}{n})$ with a probability bounded away from 0 is that
$\frac{\lambda_{d+1}}{\lambda_{d}} \lesssim  \sqrt{\frac{d}{n}} \min \{ \sqrt{\frac{d}{r_d(\Sigma^2)}}, \frac{n}{r_d(\Sigma)}\}$ and $\tau$ is chosen satisfying $\tau+\lambda_{d+1}\frac{r_d(\Sigma)}{n}\lesssim \sqrt{\frac{d}{n}}\lambda_d$ if $\frac{n\sqrt{r_d(\Sigma^2)}}{\sqrt{d}r_d(\Sigma)}\leq 1$ or $\sqrt{\frac{r_d( \Sigma^2)}{d} } \lambda_{d+1} \lesssim \tau +\lambda_{d+1}\frac{r_d(\Sigma)}{n}\lesssim
           \sqrt{\frac{d}{n}}\lambda_{d}$ if $\frac{n\sqrt{r_d(\Sigma^2)}}{\sqrt{d}r_d(\Sigma)}> 1$.

\noindent The sufficient and necessary conditions become matched, $\frac{\lambda_{d+1}}{\lambda_{d}} \lesssim  \sqrt{\frac{d}{n}} \min \{ \sqrt{\frac{d}{r_d(\Sigma^2)}}, \frac{n}{r_d(\Sigma)}\}$, in the case where $n\gg d$ in addition to the assumptions stated.
\end{cor}

Next, we give our main result about $\MSE_{\In}$ in the regime of large TER stated in Assumption \ref{ass:4}.
Let $A_0$ be any positive constant. We derive upper and lower bounds of $\MSE_{\In}$ in the case of $\tau+\lambda_{d+1}\frac{r_d(\Sigma)}{n}\leq A_0 \lambda_d$,
and a lower bound of $\MSE_{\In}$ in the case of $\tau+\lambda_{d+1}\frac{r_d(\Sigma)}{n}\geq A_0 \lambda_d$, through the sum of bias and variance terms.

\begin{thm}[In-sample error with large TER]
\label{thm4}
    Under Assumption \ref{ass:3}, \ref{ass:4} and \ref{ass1}(ii), for any $\nu$ satisfying $0<\nu < \frac{1}{4}$ and $\frac{r_d(\Sigma)\nu^2}{C_0^2\sigma_x^4}>1$ and any $A_0>0$, the following inequalities hold  uniformly in the range of $\tau$ stated
    with probability at least $1-2n \exp\{-\frac{\nu\sqrt{r_d(\Sigma)}}{C_0\sigma_x^2}\}-2 \exp\{-\frac{\nu^2 n}{C_0^2\sigma_x^4}\}-2 \exp\{-\frac{\nu^2 n}{C_0\sigma_x^4}\}-12\exp\{-\frac{n}{C_0}\}$:
    \begin{align*}
&(i)~M_2(\underbrace{\Vert\theta_{1:d}^*\Vert^2_{\Sigma_{1:d}^{-1}}(\tau+\lambda_{d+1}\frac{r_d(\Sigma)}{n})^2}_{\overline{\B}_\In }+\underbrace{\sigma^2(\frac{d}{n}+\frac{\lambda_{d+1}^2}{(\tau+\lambda_{d+1}\frac{r_d(\Sigma)}{n})^2}\frac{r^2_d(\Sigma)}{n^2}) }_{\overline{\V}_\In} )\geq  \MSE_{\In} \geq \notag\\
      & M_1(\underbrace{\kappa_{2}(\tau)\Vert\theta_{1:d}^*\Vert^2_{\Sigma_{1:d}^{-1}}(\tau+\lambda_{d+1}\frac{r_d(\Sigma)}{n})^2}_{ \underline{\B}_\In } +\underbrace{ \sigma^2(\frac{d}{n}+\frac{\lambda_{d+1}^2}{(\tau+\lambda_{d+1}\frac{r_d(\Sigma)}{n})^2}\frac{r^2_d(\Sigma)}{n^2})  }_{\underline{\V}_\In }) \quad \text{for }  \tau+\lambda_{d+1}\frac{r_d(\Sigma)}{n}
      \leq A_0 \lambda_d, \\
      &(ii)~\MSE_{\In} \geq M_1(\underbrace{\kappa_{2}(\tau)\Vert\theta_{1:d}^*\Vert^2_{\Sigma_{1:d}^{-1}}\lambda_d^2}_{\underline{\B}_\In } + \underbrace{ \sigma^2\frac{\lambda_{d+1}^2}{(\tau+\lambda_{d+1}\frac{r_d(\Sigma)}{n})^2}\frac{r^2_d(\Sigma)}{n^2}  }_{\underline{\V}_\In }) \quad \text{for }  \tau+\lambda_{d+1}\frac{r_d(\Sigma)}{n}
      \geq A_0 \lambda_d,
    \end{align*}
where $\kappa_2(\tau)=\max\{1-(16\frac{\lambda_{d+1}\frac{r_d(\Sigma)}{n}}{\tau+\lambda_{d+1}\frac{r_d(\Sigma)}{n}}(1+112\frac{\sqrt{\delta_2}}{1-\sqrt{\delta_2}})+64\frac{\sqrt{\delta_2}}{1-\sqrt{\delta_2}}),0\}$ and $M_1,M_2>0$ are constants depending only on $(\sigma_x, \eta_2, \delta_2, \nu, \A_0)$.
\end{thm}

Similarly to $\kappa_1(\tau)$ in Theorem \ref{thm3}, the definition of $\kappa_2(\tau)$ indicates that the bias term $\underline{\B}_\In$ is activated in the lower bound of $\MSE_\In$
only when $\frac{\lambda_{d+1}}{(\tau+\lambda_{d+1}\frac{r_d(\Sigma)}{n})}\frac{r_d(\Sigma)}{n}$ and $\frac{\sqrt{\delta_2}}{1-\sqrt{\delta_2}}$ are small enough.
In this case, the bias from the spiked part, $\Vert \hat{\theta}(\tau)_{1:d}-\theta_{1:d}^*\Vert^2_{\hat{\Sigma}_{1:d}}$,
 can be shown to dominate the interaction term between the spiked part and the tail part.

\begin{rem}
Similarly to out-sample error bounds discussed in Remark \ref{rem_small_large_TER},
the error bounds in Theorems \ref{thm2} and \ref{thm4} are also derived from the same set of algebraic bounds
but then by applying relevant high-probability inequalities to control random quantities for different ranges of $\tau$
under the two regimes of TER (see Section \ref{sec_6_2}).
After ignoring the range choice of $\tau$, the error bounds (i)--(ii) in Theorem \ref{thm4} appear similar to (ii)--(iii) in Theorem \ref{thm2}
except with $\tau+ \lambda_{d+1}\frac{r_d(\Sigma)}{n}$ in place of $\tau$
and with the additional difference that $ \frac{r^2_d(\Sigma)}{n^2}$ are involved both $\overline{\V}_\In$ and $\underline{\V}_\In$ in
Theorem \ref{thm4} (i), but not in Theorem \ref{thm2} (ii).
In the overlapping case of small or moderate TER regime and large TER regime, i.e., $\frac {r_d(\Sigma)}{n}\asymp 1$
and $\frac {r_d(\Sigma)}{n}$ satisfies Assumption \ref{ass:3a} and \ref{ass:4},
the result from Theorem \ref{thm2} (ii) with $\kappa_1(\tau)$ is activated reduces to $\MSE_\out \asymp \Vert\theta_{1:d}^*\Vert^2_{\Sigma_{1:d}^{-1}} \tau^2+\sigma^2(\frac{d}{n}+\frac{\lambda_{d+1}^2}{\tau^2} )$, and the result from Theorem \ref{thm4} (i) with $\kappa_2(\tau)$ is activated
reduces to $\MSE_\out \asymp \Vert\theta_{1:d}^*\Vert^2_{\Sigma_{1:d}^{-1}}(\tau+ \lambda_{d+1}\frac{r_d(\Sigma)}{n})^2+\sigma^2(\frac{d}{n}+ \frac{\lambda_{d+1}^2}{(\tau+\lambda_{d+1}\frac{r_d(\Sigma)}{n})^2} )$,
and the two results are equivalent to each other, for $\tau$ in the range $A_0^{-1}\lambda_{d+1} \leq \tau \leq A_0\lambda_d$,
where $\tau \asymp \tau+ \lambda_{d+1}\frac{r_d(\Sigma)}{n}$ with $r_d(\Sigma)\asymp n$.
However, the error bound in Theorem \ref{thm4} (i) remains applicable, but that in Theorem \ref{thm2} (ii) does not apply, to
small $\tau$ satisfying $0\le \tau \le A_0^{-1} \lambda_{d+1}$ including $\tau=0$.
\end{rem}

From Theorem \ref{thm4}, we deduce the following simple conditions for achieving $\MSE_{\In}=O(\frac{d}{n})$ in the regime of large TER.

\begin{cor}[Conditions for $\MSE_{\In} = O(\frac{d}{n})$ with large TER]
\label{cor4}
In the setting of Theorem \ref{thm4}, assume further
that $\sigma^2\asymp 1$ and $\Vert\theta_{1:d}^*\Vert^2_{\Sigma_{1:d}^{-1}}\lambda_d^2\asymp 1$.

(i) A sufficient condition for $\MSE_{\In}=O(\frac{d}{n})$ with a probability approaching to 1 as $n\to\infty$ is $\frac{\lambda_{d+1}}{\lambda_{d}}\lesssim \frac{d}{r_d(\Sigma)}$ and the ridge parameter $\tau$ is chosen such that $\lambda_{d+1}\frac{r_d(\Sigma)}{n}\sqrt{\frac{n}{d}} \lesssim \tau +\lambda_{d+1}\frac{r_d(\Sigma)}{n} \lesssim  \lambda_d\sqrt{\frac{d}{n}}$.

(ii) Suppose that $n\gg d$ and $64\frac{\sqrt{\delta_2}}{1-\sqrt{\delta_2}}<1$. Then a necessary condition for $\MSE_{\In} = O(\frac{d}{n})$ with a probability bounded away from 0 is
$\frac{\lambda_{d+1}}{\lambda_{d}}\lesssim \frac{d}{r_d(\Sigma)}$ and the ridge parameter $\tau$ is chosen in the range $\lambda_{d+1}\frac{r_d(\Sigma)}{n}\sqrt{\frac{n}{d}} \lesssim \tau+\lambda_{d+1}\frac{r_d(\Sigma)}{n}\lesssim \lambda_d\sqrt{\frac{d}{n}}$.

\noindent The sufficient and necessary conditions become matched, $\frac{\lambda_{d+1}}{\lambda_{d}}\lesssim \frac{d}{r_d(\Sigma)}$,
in the case where $n\gg d$ in addition to the assumptions stated.
\end{cor}

From Theorem \ref{thm3} and \ref{thm4}, we derive the orders of $\MSE^*_\out$ and $\MSE^*_\In$, i.e., $\MSE_\out$ and $\MSE_\In$ with optimal
choices of $\tau$ respectively, in large TER regime.

\begin{cor}[Optimal error orders with large TER]
\label{cor5}
Suppose that Assumption \ref{ass:3}, \ref{ass:4} and \ref{ass1}(ii) are satisfied and further
$\sigma^2\asymp 1$, $\Vert\theta_{1:d}^*\Vert^2_{\Sigma_{1:d}^{-1}}\lambda_d^2\asymp 1$, $\lambda_d \gg \lambda_{d+1}\frac{r_d(\Sigma)}{n}$ and $64\frac{\sqrt{\delta_2}}{1-\sqrt{\delta_2}}<1$. Then
\newline
(i) The order of $\MSE_{\out}^*$ is $\max\{\frac{\lambda_{d+1}}{\lambda_d}\sqrt{\frac{r_d(\Sigma^2)}{n}},\frac{\lambda_{d+1}^2}{\lambda_d^2}\frac{ r_d(\Sigma)^2}{n^2},\frac{d}{n}\}$ with a probability approaching to 1 and the optimal $\tau$ is chosen as $\tau=0$ if $\frac{n\sqrt{r_d(\Sigma^2)}}{\sqrt{d}r_d(\Sigma)}\leq \frac{\lambda_{d+1}}{\lambda_d}\frac{r_d(\Sigma)}{n}$ or satisfying $\tau+\lambda_{d+1}\frac{r_d(\Sigma)}{n}\asymp\sqrt{\lambda_d\lambda_{d+1}\sqrt{\frac{r_d(\Sigma^2)}{n}}}$ if $\frac{n\sqrt{r_d(\Sigma^2)}}{\sqrt{d}r_d(\Sigma)}> \frac{\lambda_{d+1}}{\lambda_d}\frac{r_d(\Sigma)}{n}$.
\newline
(ii) The order of $\MSE^*_{\In}$ is $\max\{\frac{\lambda_{d+1}}{\lambda_d}\frac{r_d(\Sigma)}{n},\frac{d}{n}\}$ with a probability approaching to 1 and the optimal $\tau$  is chosen satisfying $\tau+\lambda_{d+1}\frac{r_d(\Sigma)}{n} \asymp \sqrt{\lambda_d\lambda_{d+1}\frac{r_d(\Sigma)}{n}}$.

\noindent Therefore $\MSE_{\out}^* \lesssim \MSE_{\In}^*$ with a probability approaching to 1 because $\frac{\lambda_{d+1}}{\lambda_d}\sqrt{\frac{r_d(\Sigma^2)}{n}}\lesssim \frac{\lambda_{d+1}}{\lambda_d}\frac{r_d(\Sigma)}{n}$ by noting $r_d(\Sigma^2)\leq r_d(\Sigma)$ and $r_d(\Sigma)\gtrsim n$ (Assumption \ref{ass:4}) and because $\frac{\lambda^2_{d+1}}{\lambda_d^2}\frac{r_d(\Sigma)^2}{n^2}\lesssim \frac{\lambda_{d+1}}{\lambda_d}\frac{r_d(\Sigma)}{n}$ by noting $\lambda_d\gg \lambda_{d+1}\frac{r_d(\Sigma)}{n}$.
\end{cor}

The additional condition $\lambda_d\gg \lambda_{d+1}\frac{r_d(\Sigma)}{n}$ is important for determining the order of $\MSE_\In^*$.
In fact, the sum of the in-sample bias and the tail part of in-sample variance, $\Vert\theta^*_{1:d}\Vert^2_{\Sigma_{1:d}^{-1}} (\tau+\lambda_{d+1}\frac{r_d(\Sigma)}{n})^2+\sigma^2\frac{\lambda_{d+1}^2}{(\tau+\lambda_{d+1}\frac{r_d(\Sigma)}{n})^2}\frac{r_d^2(\Sigma)}{n^2}$,
reach the minimum order of $\frac{\lambda_{d+1}}{\lambda_d}\frac{r_d(\Sigma)}{n}$ by the choice $\tau$ satisfying $\tau+\lambda_{d+1}\frac{r_d(\Sigma)}{n}\asymp \sqrt{\lambda_d\lambda_{d+1}\frac{r_d(\Sigma)}{n}}$.
The condition $\lambda_d\gg \lambda_{d+1}\frac{r_d(\Sigma)}{n}$
ensures that this choice of $\tau$ is large enough so that $\kappa_2(\tau)$ is activated.
See the proof of Corollary \ref{cor5} in Supplement Section \ref{sec_II_2} for details.

\begin{rem}
\label{in-out-com-large_TER}
In the setting of Corollary \ref{cor5}, we observe that
the gap between $\MSE_\In^*$ and $\MSE_\out^*$ can be significantly large, for example, in two cases,
$\frac{\lambda_{d+1}}{\lambda_d}\gtrsim \frac{d}{\sqrt{nr_d(\Sigma^2)}} \min \{ 1, \frac{n\sqrt{r_d(\Sigma^2)}}{\sqrt{d}r_d(\Sigma)}\}$ and $\frac{\lambda_{d+1}}{\lambda_d}\gtrsim \frac{n\sqrt{n r_d(\Sigma^2)}}{r_d(\Sigma)^2}$ hold or $\frac{\lambda_{d+1}}{\lambda_d}\gtrsim \frac{d}{\sqrt{nr_d(\Sigma^2)}} \min \{ 1, \frac{n\sqrt{r_d(\Sigma^2)}}{\sqrt{d}r_d(\Sigma)}\}$,  $\frac{\lambda_{d+1}}{\lambda_d}\lesssim \frac{n\sqrt{n r_d(\Sigma^2)}}{r_d(\Sigma)^2}$ and $\frac{r_d(\Sigma)}{\sqrt{n r_d(\Sigma^2)}}\gg 1$ hold. In the first case, $\MSE_\In^*\asymp \frac{\lambda_{d+1}}{\lambda_d}\frac{r_d(\Sigma)}{n}$ and $\MSE_\out^*\asymp \frac{\lambda_{d+1}^2}{\lambda_d^2}\frac{r_d(\Sigma)^2}{n^2}$ by Corollary \ref{cor5}, and hence $\frac{\MSE_\In^*}{\MSE_\out^*}=\frac{\lambda_d}{\lambda_{d+1}}\frac{n}{r_d(\Sigma)}\gg 1$ with $\lambda_d \gg \lambda_{d+1}\frac{r_d(\Sigma)}{n}$. In the second case, $\MSE_\In^*\asymp \frac{\lambda_{d+1}}{\lambda_d}\frac{r_d(\Sigma)}{n}$ and $\MSE_\out^*\asymp \frac{\lambda_{d+1}}{\lambda_d}\sqrt{\frac{r_d(\Sigma^2)}{n}}$ by Corollary \ref{cor5}, and hence $\frac{\MSE_\In^*}{\MSE_\out^*}=\frac{r_d(\Sigma)}{\sqrt{n r_d(\Sigma^2)}}\gg 1$.
\end{rem}

In the setting of Corollary \ref{cor5}, we point out that
the advantage of $\MSE_\out^*$ over $\MSE_\In^*$ can be attributed to $\V_\out\lesssim \V_\In$ for $\tau+\lambda_{d+1}\frac{r_d(\Sigma)}{n}\leq A_0\lambda_d$.
See the proof of Corollary \ref{cor5} in Supplement Section \ref{sec_II_2} for details. In fact,
in the setting of Corollary \ref{cor5}, the optimal $\tau$ for both $\MSE_\out$ and $\MSE_\In$ are chosen from the range $\tau +\lambda_{d+1}\frac{r_d(\Sigma)}{n}\leq A_0\lambda_{d}$, because both $\MSE_\out$ and $\MSE_\In$ are lower bounded by large bias for $\tau +\lambda_{d+1}\frac{r_d(\Sigma)}{n}\geq  A_0\lambda_{d}$.
Moreover, as seen from the proof, the optimal order of $\MSE_\In$, $\max\{\frac{\lambda_{d+1}}{\lambda_d}\frac{r_d(\Sigma)}{n},\frac{d}{n}\}$, in the setting of Corollary \ref{cor5} can be achieved only when $\underline{\B}_\In$ is activated.
With $\underline{\B}_\In$ activated, the orders of $\B_\out$ and $\B_\In$ are the same for $\tau+\lambda_{d+1}\frac{r_d(\Sigma)}{n}\leq A_0\lambda_d$.
By comparison,
$\V_\out$ can be shown to be smaller than $\V_\In$ up to a constant factor for $\tau+\lambda_{d+1}\frac{r_d(\Sigma)}{n}\leq A_0\lambda_d$ in the large TER regime:
\begin{align}
    \underbrace{\sigma^2(\frac{d}{n}+\frac{\lambda_{d+1}^2}{(\tau+\lambda_{d+1}\frac{r_d(\Sigma)}{n})^2}\frac{r_d(\Sigma^2)}{n}) }_{\text{order~of~$\V_\out$~in~Theorem~\ref{thm3}(i)} }\lesssim \underbrace{\sigma^2(\frac{d}{n}+ \frac{\lambda_{d+1}^2}{(\tau+\lambda_{d+1}\frac{r_d(\Sigma)}{n})^2}\frac{r^2_d(\Sigma)}{n^2})}_{\text{order~of~$\V_\In$~in~Theorem~\ref{thm4}(i)} } . \notag
\end{align}
Hence the advantage of $\MSE_\out^*$ stems from the smaller order of $\V_\out$ in the setting of Corollary \ref{cor5}.

\section{Connection and comparison with existing results}

\subsection{Comparison with error approximation formulas}

\label{sec:4_1}

We compare our results to error approximation formulas, obtained under an independence assumption on whiten covariates in ridge linear regression
(see Assumption \ref{ass:5} below).
We first review the results of \cite{HasteTrevorAndrea} about approximation formulas for out-sample error. To facilitate the comparison, we also derive and justify the approximation formulas for in-sample error. Then we show that in the rotationally sparse setting,
the orders of out-sample and in-sample error approximation formulas match those from our results for the ridge tuning parameter in suitable ranges.

For ridge linear regression,
\cite{HasteTrevorAndrea} gave the following approximation formulas for the out-sample bias and variance:
\begin{align}
        \mathcal{B}_\out(\tau, \hat{H}_n, \hat{G}_n, \gamma)&= \tau^2\Vert\theta^*\Vert^2(1+\gamma m_{n,1}(-\tau))\int \frac{s}{[\tau + (1-\gamma + \gamma \lambda m_n(-\tau))s]^2}d \hat{G}_n(s) , \label{eq:s2_2a}  \\
            \mathcal{V}_\out(\tau, \hat{H}_n, \gamma) &= \sigma^2\gamma \int \frac{s^2(1-\gamma + \gamma \lambda^2 m_n^\prime(-\tau))}{[\tau + s(1-\gamma+\gamma \tau m_n(-\tau))]^2} d \hat{H}_n(s) ,\label{eq:s2_2b}
\end{align}
where $\gamma=\frac{p}{n}$, $\hat{H}_n(s)=p^{-1}\sum_{j=1}^p1_{\{s\geq \lambda_j\}}$, $\hat{G}_n(s)=\sum_{i=j}^p (\langle\theta^*, v_j\rangle^2/\Vert\theta^*\Vert^2) 1_{\{s\geq \lambda_j\}}$, $v_1,...,v_p$ are the eigenvectors of $\Sigma$, $m_n(z)$ is determined by solving the following equation
\begin{align}
     m_n(z)&=\int \frac{1}{s[1-\gamma-\gamma z m_n(z)]-z}d \hat{H}_n(s) , \label{eq:s2_1}
\end{align}
and $m_{n,1}(z)$ is calculated by
\begin{align}
    m_{n,1}(z)&=\frac{\int \frac{ s^2[1-\gamma-\gamma z m_n(z)]  }{[s[1-\gamma-\gamma z m_n(z)]-z]^2} d \hat{H}_n(s) }{1-\gamma \int \frac{zs}{[s[1-\gamma -\gamma z m_n(z)]-z]^2}d \hat{H}_n(s)} . \label{eq:s2_2}
\end{align}
Consider the following assumption on the whiten covariates, defined as $z_i=\Sigma^{-1/2}x_i$, the variance matrix $\Sigma$, and the ratio $\frac{p}{n}$.
\begin{ass}
\label{ass:5}
$\newline$
(i) Each vector $z_i=(z_{i1},\cdots, z_{ip})^\T$ has independent entries with $E(z_{ij})=0$, $E(z_{ij}^2)=1$ and $E( |z_{ij}|^k )\leq C_k < \infty$ for all $k\geq 2$.
\newline
(ii) $ \lambda_1 \leq M$ and $\int s^{-1}\hat{H}_n(s)ds < M$.
\newline
(iii) $|1-\frac{p}{n}|\geq \frac{1}{M}$, $\frac{1}{M}\leq \frac{p}{n}\leq M$ .
\end{ass}
\cite{HasteTrevorAndrea} showed that under Assumption \ref{ass:5} and assuming $\max\{\tau,\lambda_p\} >\frac{1}{M} $ and $n^{-2/3+1/M}<\tau < M$, for any constant $D>0$ and $\delta>0$, with probability at least $1-C(M,D,\delta)n^{-D}$,
\begin{align*}
    |\mathcal{B}_\out(\tau, \hat{H}_n, \hat{G}_n, \gamma) - \B_\out| < \frac{C(M)\Vert \theta^*\Vert^2}{\tau n^{(1-\delta)/2}} ,\quad
    |\mathcal{V}_\out(\tau, \hat{H}_n, \hat{G}_n, \gamma) - \V_\out| < \frac{C(M)}{\tau^2 n^{(1-\delta)/2}},
\end{align*}
where $C(M,D,\delta)$ is a constant depending only on $(M, D,\delta)$, and $C(M)$ is a constant depending only on $M$.

To facilitate the comparison between our results and error approximation formulas, we also give approximation formulas for in-sample bias and variance as follows:
\begin{align}
    \mathcal{B}_{\In}(\tau;\hat{H}_n,\hat{G}_n,\gamma)&=\tau^2||\theta^*||^2(\gamma\tau^2 m_n^\prime(-\tau) + 1-\gamma)\int \frac{s}{[\tau+(1-\gamma+\gamma\lambda m_n(-\tau))s]^2}d\hat{G}_n(s) , \label{eq:s2_2c}\\
    \mathcal{V}_{\In}(\tau;\hat{H}_n,\gamma)&=\sigma^2\gamma(1 - 2\tau m_n(-\tau)+\tau^2 m_n^\prime(-\tau)) . \label{eq:s2_2d}
\end{align}
We establish the convergence of $\mathcal{B}_{\In}(\tau;\hat{H}_n,\hat{G}_n,\gamma)$ and $\mathcal{V}_{\In}(\tau;\hat{H}_n,\gamma)$ in the following theorem.

\begin{thm}[Convergence of in-sample error approximation formulas]
\label{thm5}
Under Assumption \ref{ass:5}, further assume that $\tau >\frac{1}{M}$ and $n^{-2/3+1/M}<\tau < \frac{M}{2}$. Then for any $D>0$, $\delta>0$, with probability at least $1-C(M,D, \delta)n^{-D}$,
\begin{align*}
    |\mathcal{B}_{\In}(\tau;\hat{H}_n,\hat{G}_n,\gamma) - \B_{\In}| &\leq  C(M)\max\{\frac{1}{\tau^{2/3}n^{(1-\delta)/3}},\frac{8M}{\tau n^{(1-\delta)/2}}\} , \\
    |\mathcal{V}_{\In}(\tau;\hat{H}_n,\gamma) - \V_{\In}| &\leq \sigma^2 C(M)(\max\{\frac{1}{\tau^{2/3}n^{(1-\delta)/3}},\frac{8M}{\tau n^{(1-\delta)/2}}\}  + \frac{1}{n^{(1-\delta)/2}} ),
\end{align*}
where $C(M,D,\delta)$ is a constant depending only on $(M, D, \delta)$, and $C(M)$ is a constant depending only on $M$.
\end{thm}

Given any $\tau >0, \gamma = \frac{p}{n}>0$ and $\tilde{\lambda}=(\lambda_1,\ldots,\lambda_p)$ with $\lambda_j>0$ for $1\leq j\leq p$
, we define $\alpha > 1$ as a solution to the equation
\begin{align}
    \frac{1}{\alpha}&= 1-\gamma\frac{1}{p}\sum_{j=1}^p\frac{1}{1+\frac{\alpha \tau}{\lambda_j}}  .\label{eq:4.1.0}
\end{align}
The approximation formulas above can be equivalently expressed as follows.
These formulas can also be calculated using a distributional approximation method in \cite{HanShen2022} under the independent components assumption, for which the discussion is deferred to Supplement Section~\ref{sec_supp_II4}.

\begin{cor}[Equivalent expressions of error approximation formulas]
\label{cor5b}
With $\alpha$ defined in (\ref{eq:4.1.0}), we have
\begin{align}
        \mathcal{B}_{\out}(\tau;\hat{H}_n,\hat{G}_n,\gamma) &= (1-\frac{1}{n}\sum_{j=1}^p \frac{\lambda_j^2}{(\lambda_j+\alpha \tau)^2})^{-1}\frac{1}{n}\sum_{j=1}^p \frac{\alpha^2\tau^2 \lambda_j \theta_j^{*2} }{(\lambda_j +\alpha\tau)^2} ,  \label{eq:4.1.1}  \\
        \mathcal{V}_{\out}(\tau;\hat{H}_n,\gamma) &=(1-\frac{1}{n}\sum_{j=1}^p \frac{\lambda_j^2}{(\lambda_j+\alpha \tau)^2})^{-1}(\frac{1}{n}\sum_{j=1}^p\frac{\lambda_j^2}{(\lambda_j + \alpha \tau)^2})\sigma^2,  \label{eq:4.1.2} \\
        \mathcal{B}_{\In}(\tau;\hat{H}_n,\hat{G}_n,\gamma) &=  \frac{1}{\alpha^2}(1-\frac{1}{n}\sum_{j=1}^p \frac{\lambda_j^2}{(\lambda_j+\alpha \tau)^2})^{-1}\frac{1}{n}\sum_{j=1}^p \frac{\alpha^2\tau^2 \lambda_j \theta_j^{*2} }{(\lambda_j +\alpha\tau)^2}  ,    \label{eq:4.1.3} \\
       \mathcal{V}_{\In}(\tau;\hat{H}_n,\gamma)&=  (1-\frac{2}{\alpha} + \frac{(1-\frac{1}{n}\sum_{j=1}^p\frac{\lambda_j^2}{(\lambda_j + \alpha \tau)^2})^{-1} }{\alpha^2} )\sigma^2 .\label{eq:4.1.4}
\end{align}
\end{cor}

Next, we study the orders of error approximation formulas (\ref{eq:4.1.1})-(\ref{eq:4.1.4})
in the high-dimensional rotationally sparse setting,
and compare them with our results, which are obtained without requiring independence of the whiten covariates.
The first result is about the small or moderate TER regime.

\begin{cor}[Matching error approximation formulas with small or moderate TER]
\label{cor7_c}
~\newline
\indent (i) Suppose that $\frac{d}{n} < 1$, $r_d(\Sigma)\lesssim n$, and $\Vert\theta^*_{(d+1):p}\Vert^2_{\Sigma_{(d+1):p}}\lesssim  \Vert\theta^*_{1:d}\Vert^2_{\Sigma_{1:d}^{-1}}\lambda_{d+1}^2$.
For $\lambda_{d+1} \lesssim   \tau \lesssim \lambda_{d}$, we have
\begin{align}
    \mathcal{B}_{\out}(\tau;\hat{H}_n,\hat{G}_n,\gamma) +
    \mathcal{V}_{\out}(\tau;\hat{H}_n,\gamma) &\asymp  \Vert\theta^*_{1:d}\Vert_{\Sigma^{-1}_{1:d}}^2 \tau^2+\sigma^2(\frac{d}{n}+\frac{\lambda_{d+1}^2}{\tau^2}\frac{r_d(\Sigma^2)}{n}).  \label{eq:4.1.11}
\end{align}

(ii) Suppose further that $r_d(\Sigma)\asymp n$. For $\lambda_{d+1} \lesssim   \tau \lesssim \lambda_{d}$, we have
\begin{align}
    \mathcal{B}_{\In}(\tau;\hat{H}_n,\hat{G}_n,\gamma) +
    \mathcal{V}_{\In}(\tau;\hat{H}_n,\gamma) &\asymp  \Vert\theta^*_{1:d}\Vert_{\Sigma^{-1}_{1:d}}^2 \tau^2+\sigma^2(\frac{d}{n}+\frac{\lambda_{d+1}^2}{\tau^2}) . \label{eq:4.1.11a}
\end{align}
\end{cor}

For comparison, we notice that the conditions, $\frac{d}{n} < 1$, $r_d(\Sigma)\lesssim n$ and $\Vert\theta^*_{(d+1):p}\Vert^2\lesssim  \Vert\theta^*_{1:d}\Vert^2_{\Sigma_{1:d}^{-1}}\lambda_{d+1}^2$
correspond to, respectively, Assumption \ref{ass:3}, \ref{ass:3a} and \ref{ass1}(i) used in Theorems \ref{thm1} and \ref{thm2}.
The order of the approximation formula (\ref{eq:4.1.11}) matches Theorem \ref{thm1}(ii).
The order of the approximation formula (\ref{eq:4.1.11a}) matches Theorem \ref{thm2}(ii) when $r_d(\Sigma)\asymp n$ and $\kappa_1(\tau)$ is activated.

Then we study the large TER regime and the results are summarized as follows.

\begin{cor}[Matching error approximation formulas with large TER]
\label{cor7_d}
~$\newline$
\indent (i) Suppose that  $\frac{d}{n}<\frac{1}{5}$, $r_d(\Sigma) > c n$ for some $c>10$, $\Vert\theta^*_{(d+1):p}\Vert^2_{\Sigma_{(d+1):p}}\lesssim  \Vert\theta^*_{1:d}\Vert^2_{\Sigma_{1:d}^{-1}}(\frac{1}{\lambda_d}+\frac{n}{\sum_{j>d}\lambda_j})^{-2}$. For
    $\lambda_{d} \gtrsim  \tau + \lambda_{d+1}\frac{r_d(\Sigma)}{n}$ and $\tau>0$, we have
\begin{align}
       \mathcal{B}_{\out}(\tau;\hat{H}_n,\hat{G}_n,\gamma) +
    \mathcal{V}_{\out}(\tau;\hat{H}_n,\gamma) &\asymp \Vert\theta_{1:d}^*\Vert^2_{\Sigma_{1:d}^{-1}}(\tau+\lambda_{d+1}\frac{r_d(\Sigma)}{n})^2 + \sigma^2(\frac{d}{n}+\frac{\lambda_{d+1}^2}{(\tau+\lambda_{d+1}\frac{r_d(\Sigma)}{n})^2}\frac{r_d(\Sigma^2)}{n}) .\label{eq:4.1.6}
\end{align}

(ii) Suppose further that $\tau > \lambda_{d+1}\frac{r_d(\Sigma)}{n}$. For
    $\lambda_{d} \gtrsim  \tau + \lambda_{d+1}\frac{r_d(\Sigma)}{n}$ and $\tau>0$, we have
\begin{align}
    \mathcal{B}_{\In}(\tau;\hat{H}_n,\hat{G}_n,\gamma) +
    \mathcal{V}_{\In}(\tau;\hat{H}_n,\gamma) &\asymp \Vert\theta_{1:d}^*\Vert^2_{\Sigma_{1:d}^{-1}}(\tau+\lambda_{d+1}\frac{r_d(\Sigma)}{n})^2 + \sigma^2(\frac{d}{n}+\frac{\lambda_{d+1}^2}{(\tau+\lambda_{d+1}\frac{r_d(\Sigma)}{n})^2}\frac{r^2_d(\Sigma)}{n^2}). \label{eq:4.1.7}
\end{align}
\end{cor}

The conditions, $\frac{d}{n}<\frac{1}{5}$, $r_d(\Sigma) > c n$ for some $c>10$, $\Vert\theta^*_{(d+1):p}\Vert^2\lesssim  \Vert\theta^*_{1:d}\Vert^2_{\Sigma_{1:d}^{-1}}(\frac{1}{\lambda_d}+\frac{n}{\sum_{j>d}\lambda_j})^{-2}$
correspond to, respectively, Assumption \ref{ass:3}, \ref{ass:4} and \ref{ass1}(ii) used in Theorems \ref{thm3} and \ref{thm4}.
The order of the approximation formula (\ref{eq:4.1.6}) matches Theorem \ref{thm3}(i).
The order of the approximation formula (\ref{eq:4.1.7}) matches Theorem \ref{thm4}(i) when $\kappa_2(\tau)$ is activated.

\subsection{Comparison with \cite{Tsigler_ridge_2023}}

\label{sec_4_2_main}

In this section, we compare our results to \cite{Tsigler_ridge_2023},
where the non-asymptotic out-sample error bounds are studied for high-dimensional ridge regression.
We make a comparison in both upper bounds and lower bounds. In the following, denote $A_d = X_{(d+1):p}X_{(d+1):p}^\T+n\tau I_n$.
The conditional number of $A_d$ is defined $\mu_1(A_d)/\mu_n(A_d)$, where $\mu_1 (A_d)$ and $\mu_n(A_d)$ are the maximum eigenvalue and minimum eigenvalue of $A_d$ repectively.

\vspace{.05in}
\textbf{Upper bound of $\MSE_\out$.} Our upper bounds of $\MSE_\out$ match the results in \cite{Tsigler_ridge_2023} for most ridge tuning parameters, but our result avoids making any oracle condition on covariate vectors as used in \cite{Tsigler_ridge_2023}. In fact, for any $d$ small enough compared to $n$, given that the conditional number of $A_d$ is controlled by $L$, it is shown in \cite{Tsigler_ridge_2023} that with a high probability:
\begin{align}
    \B_{\out}/c &\leq \Vert \theta^*_{1:d} \Vert_{\Sigma_{1:d}^{-1}}(\tau + \lambda_{d+1}\frac{r_d(\Sigma)}{n})^2
    +  \Vert \theta^*_{(d+1):p}\Vert^2_{\Sigma_{(d+1):p}} ,\label{eq:4_1_0_c}\\
    \V_\out /c &\leq \frac{d}{n} + \frac{\lambda_{d+1}^2}{(\tau+\lambda_{d+1}\frac{r_d(\Sigma)}{n})^2}\frac{r_d(\Sigma^2)}{n} ,\label{eq:4_1_1_c}
\end{align}
where $c$ is a constant depending on the $\sigma_x$ and $L$.
When $r_d(\Sigma) \lesssim n$ (i.e., in the small or moderate TER regime), our result (\ref{p1_2b})--(\ref{p1_2}) in Proposition \ref{pro1} matches (\ref{eq:4_1_0_c})--(\ref{eq:4_1_1_c}) for $\lambda_{d+1}\leq \tau\leq \lambda_{1}$.
When $r_d(\Sigma)>c_x n$ for some $c_x$ depending only on $\sigma_x$ (i.e., in the large TER regime), our result (\ref{p7_2b})--(\ref{p7_2}) in Proposition \ref{pro7} matches (\ref{eq:4_1_0_c})--(\ref{eq:4_1_1_c}) for $\tau+\lambda_{d+1}\frac{r_d(\Sigma)}{n}\leq \lambda_{1}$.

To control the the conditional number of $A_d$, \cite{Tsigler_ridge_2023} requires an oracle small-ball assumption on covariate vectors:
$\sum_{j>d}x_{ij}^2>c(\sum_{j>d}\lambda_j+n\tau)$ for all $i=1,\ldots,n$ and some $c$ satisfying $0<c<1$ with a high probability.

Instead of requiring an extra oracle assumption, our analysis achieves the control of the conditional number of $A_d$ using concentration inequalities
specifically in the two TER regimes for certain ranges of $\tau$,
which are summarized below. See Section \ref{sec_6_1_3} and Supplement Section \ref{bound_A_d} for details.
We derive an upper bound of $\mu_1(A_d)$, using concentration properties of matrix operator norms based on sub-gaussian covariate vectors,
for $\tau\ge \lambda_{d+1}$ in the small or moderate TER regime and
for $\tau\ge 0$ in the large TER regime.
To obtain a lower bound of $\mu_n(A_d)$, we handle the small or moderate TER regime and the large TER regime separately.
In small or moderate TER regime, we use a trivial lower bound: $\mu_n(A_d)\geq n\tau$. In the large TER regime,
we use  concentration properties of quadratic forms of sub-gaussian random vectors (\cite{Zajkowski_2020}, Corollary 2.8)
to show that the oracle small-ball assumption in \cite{Tsigler_ridge_2023} is satisfied for $\tau \ge 0$.
Then we derive a lower bound of $\mu_n(A_d)$ following similar reasoning as in \cite{Tsigler_ridge_2023}.

\vspace{.05in}
\textbf{Lower bound of $\V_\out$.} Our lower bound of $\V_\out$ also matches the result in \cite{Tsigler_ridge_2023} for a certain range of ridge tuning parameters,
 but our result does not require the independence of the components of the whiten covariate vector as assumed in \cite{Tsigler_ridge_2023}.
In fact, given that the components of whiten $x_i$ are independent, for any $d$ small enough compared to $n$,
it is shown in \cite{Tsigler_ridge_2023} that with a high probability,
\begin{align}
    \V_{\out}/c \geq \frac{1}{n }\sum_{j=1}^p\min\{1,\frac{\lambda_j^2}{\lambda_{d+1}^2(\frac{\sum_{l>d}\lambda_l+n\tau}{n\lambda_{d+1}}+1)^2}\} . \label{eq:sec4_2_1}
\end{align}
where $c$ is a constant depending on $\sigma_x$.
When $r_d(\Sigma) \lesssim n$ (i.e., in the small or moderate TER regime), for $\lambda_{d+1}\leq \tau \leq \lambda_d$,
(\ref{eq:sec4_2_1}) can be shown to reduce to
\begin{align}
  \V_{\out} \gtrsim \frac{d}{n} + \frac{\lambda_{d+1}^2}{\tau^2}\frac{r_d(\Sigma^2)}{n} , \notag
\end{align}
which matches our result in (\ref{p3_1}) of Proposition \ref{pro3} for small or moderate TER regime.
When $r_d(\Sigma)>c_x n$ for some $c_x$ depending on $\sigma_x$ (i.e., in the large TER regime),
for $ \lambda_{d+1}\frac{r_d(\Sigma)}{n} +\tau\leq \lambda_d$, (\ref{eq:sec4_2_1}) can be shown to reduce to
\begin{align}
\V_{\out} \gtrsim \frac{d}{n} + \frac{\lambda_{d+1}^2}{(\tau+\lambda_{d+1}\frac{r_d(\Sigma)}{n})^2}\frac{r_d(\Sigma^2)}{n}  , \notag 
\end{align}
which matches our result in (\ref{p9_1}) of Proposition \ref{pro9} for large TER regime.

Instead of requiring the assumption of independent components in whiten covariate vectors as in \cite{Tsigler_ridge_2023},
we derive a lower bound of $\V_\out$ using concentration properties of sub-gaussian random vectors (see Section \ref{sec_6_1} for further information).

\vspace{.05in}
\textbf{Lower bound of $\B_\out$.} The lower bound of $\B_\out$ in \cite{Tsigler_ridge_2023} is provided
as a probability lower bound on the expectation of $\B_\out$ with respect to a prior distribution on $\theta^*$
under an extra oracle assumption on certain modifications of matrix $A = XX^\T+ n\tau I_n$.
Our lower bound on $\B_\out$ is a direct probability bound without assuming a prior distribution on $\theta^*$
and any extra oracle assumption on covariate vectors, but focuses on the rotationally sparse setting.

\subsection{Comparison with \cite{BuneaFlorentinaStrimas2020}}

\label{sec:4_3_main}

For the min-norm interpolator ($\tau = 0$), we compare Theorem \ref{thm3} with Theorem 16 in \cite{BuneaFlorentinaStrimas2020},
which are both obtained in the large TER regime. See Supplement Section \ref{sec:approximation_com} for details.
Note that models (\ref{eq:1}) and (\ref{eq:bunea-model}) are related via
$\Sigma =A \Sigma_Z A^\T + \Sigma_E$ and $\theta^*=(\Sigma_E+A\Sigma_Z A^T)^{-1}A\Sigma_Z\beta$.
To facilitate the comparison, we let $\Sigma_Z = I_d$,
$\Sigma_E =\Diag(\underbrace{\lambda_{d+1},\ldots,\lambda_{d+1}}_{\text{$d$ entries}},\lambda_{d+1},\ldots,\lambda_p)\in \mathbb{R}^{p\times p}$, and
\begin{align*}
A & =\begin{pmatrix}\Diag(\sqrt{\lambda_1-\lambda_{d+1}},\ldots,\sqrt{\lambda_d-\lambda_{d+1}}) \\
0_{(p-d)\times d}
\end{pmatrix}\in \mathbb{R}^{p\times d},
\end{align*}
such that $\Sigma = \Diag(\lambda_1,\ldots,\lambda_d)\in\mathbb{R}^{p\times p}$ and
\begin{align*}
\theta^*=(\Diag(\frac{\sqrt{\lambda_1-\lambda_d}}{\lambda_1},\ldots, \frac{\sqrt{\lambda_d-\lambda_{d+1}}}{\lambda_{d}})\beta,\underbrace{0,\ldots,0}_{\text{$p-d$ entries}})^T\in\mathbb{R}^p.
\end{align*}
We further assume that $\lambda_1\asymp \lambda_d$, $\lambda_d  > c_1 \lambda_{d+1}$ and $r_d(\Sigma)>c_2 d$ for some $c_1>1$  and $c_2>1$.
In this setting, Theorem \ref{thm3} with $\tau=0$ shows that with a high probability,
\begin{align*}
  &  \MSE_\out \asymp \underbrace{\Vert \theta_{1:d}^*\Vert^2_{\Sigma_{1:d}} \frac{\lambda_{d+1}^2}{\lambda_d^2} \frac{r^2_d(\Sigma)}{n^2}}_{\B_\out} + \underbrace{\sigma^2(\frac{d}{n}+ \frac{n r_d(\Sigma^2)}{r^2_d(\Sigma)})}_{\V_\out}, \quad \text{for } \lambda_{d+1}\frac{r_d(\Sigma)}{n}\leq \lambda_{d}, \\
  &  \MSE_\out \gtrsim  \underbrace{\Vert \theta_{1:d}^*\Vert^2_{\Sigma_{1:d}} }_{\B_\out} , \quad
 \text{for } \lambda_{d+1}\frac{r_d(\Sigma)}{n}>\lambda_d.
\end{align*}
After ignoring the $\log(n)$ factor, Theorem 16 in \cite{BuneaFlorentinaStrimas2020} gives that with a high probability,
\begin{align*}
    \MSE_\out \lesssim \underbrace{\Vert \theta_{1:d}^*\Vert^2_{\Sigma_{1:d}} \frac{\lambda_{d+1}}{\lambda_d}\frac{r_d(\Sigma)}{n}}_{\B_\out} +\underbrace{\sigma^2 (\frac{d}{n}+\frac{n}{r_d(\Sigma)})}_{\V_\out}.
\end{align*}
Hence for $\lambda_d\geq \lambda_{d+1}\frac{r_d(\Sigma)}{n}$, our result gives the order of out-sample MSE which is sharper than the upper bound in \cite{BuneaFlorentinaStrimas2020}. For $\lambda_d < \lambda_{d+1}\frac{r_d(\Sigma)}{n}$, our lower bound indicates
that the out-sample MSE is larger than $\Vert \theta_{1:d}^*\Vert^2_{\Sigma_{1:d}}$ up to a constant, and accordingly
the upper bound in \cite{BuneaFlorentinaStrimas2020} is larger than $\Vert \theta_{1:d}^*\Vert^2_{\Sigma_{1:d}}$ up to a constant.

\section{Proofs of main results (Theorems \ref{thm1}--\ref{thm4})}
\label{sec:5}

We provide proofs of the main results (Theorems \ref{thm1}--\ref{thm4}), depending on auxiliary bounds on $\B_{\out}$, $\B_{\In}$, $\V_{\out}$, $\V_{\In}$,
for which the proofs are outlined in Section \ref{sec:6}.
Without loss of generality, we only consider the case of $\A_0=1$ involved in the ranges of the ridge parameter $\tau$.

\subsection{Proof of Theorem \ref{thm1}}
\label{sec:5_1}

We provide auxiliary bounds for the out-sample squared bias and variance $\B_{\out}$ and $\V_{\out}$ under the small or moderate TER regime (Assumption \ref{ass:3a}).

\begin{pro}[Upper bound of out-sample error with small or moderate TER]
\label{pro1}
Under Assumption \ref{ass:3} and \ref{ass:3a}, for any $\nu$ satisfying
$0<\nu <\frac{1}{2}$, the following inequalities hold uniformly in the range of $\tau$ stated with probability at least $1-2\exp\{-\frac{\nu^2 n}{C_0^2\sigma_x^4}\}-18 \exp\{-\frac{n}{C_0}\}$: for $\tau \geq \lambda_{d+1}$,
\begin{align}
    \B_{\out}&\leq   \frac{(1+C_1)^3(1+\nu+\eta_1)^2\poly_6(\sigma_x)}{(1-\nu-\eta_1)^4}(\Vert\theta^*_{1:d}\Vert^2_{\Sigma_{1:d}^{-1}}(\frac{1}{\lambda_1}+\frac{1}{\tau})^{-2} + \Vert \theta_{(d+1):p}^*\Vert^2_{\Sigma_{(d+1):p}}) ,\label{p1_2b} \\
    \V_{\out} &\leq  \frac{(1+C_1)^2\poly_6(\sigma_x)}{(1-\nu-\eta_1)^4}\sigma^2(\frac{d}{n}+\frac{\lambda_{d+1}^2}{\tau^2}\frac{r_d( \Sigma^2)}{n} ). \label{p1_2}
\end{align}
Further with Assumption \ref{ass1}(i),
\begin{align}
    \B_{\out}&\leq   \frac{(1+C_1)^3(1+\nu+\eta_1)^2\poly_6(\sigma_x)}{(1-\nu-\eta_1)^4}\Vert\theta^*_{1:d}\Vert^2_{\Sigma_{1:d}^{-1}}(\frac{1}{\lambda_1}+\frac{1}{\tau})^{-2}.  \label{p1_1}
\end{align}
\end{pro}

\begin{pro}[Lower bound of $\B_\out$ with small or moderate TER]
\label{pro2}
Under Assumption \ref{ass:3} and \ref{ass1}(i), for any $\nu$ satisfying $0<\nu <\frac{1}{2}$, the following inequality holds uniformly in the range of $\tau$ stated with probability at least $1-2 \exp\{-\frac{\nu^2 n}{C_0^2\sigma_x^4}\}-2 \exp\{-\frac{n}{C_0}\}$: for $\tau \geq \lambda_{d+1}$,
\begin{align}
    &\B_{\out} \geq \frac{(1-\sqrt{\delta_1})^2}{(1+\nu +\eta_1)^2}\Vert\theta^*_{1:d}\Vert^2_{\Sigma_{1:d}^{-1}}(\frac{1}{\lambda_d}+\frac{1}{\tau})^{-2} .\label{p2_1}
\end{align}
\end{pro}

\begin{pro}[Lower bound of $\V_\out$ with small or moderate TER]
\label{pro3}
Under Assumption \ref{ass:3} and \ref{ass:3a}, for any $\nu$ satisfying $0<\nu < \frac{1}{2}\min\{1,\sigma_x^2\}$, the following inequalities hold uniformly in the range of $\tau$ stated with probability at least $1-2 \exp\{-\frac{\nu^2 n}{C_0\sigma_x^4}\}-2 \exp\{-\frac{\nu^2 n}{C_0^2\sigma_x^4}\}-12\exp\{-\frac{n}{C_0}\}$: for $\tau \leq \lambda_{d+1}$,
\begin{align}
   \V_{\out}\geq \frac{(1-\nu-\eta_1)^2(\frac{1}{2}-\nu)}{(1+C_1)^2(1+\nu+\eta_1)^4\poly_4(\sigma_x)}\frac{1}{1+\frac{2C_0\sigma_x^2}{\frac{1}{2}-\eta_1}}\sigma^2(\frac{d}{n}+\frac{r_d( \Sigma^2)}{n}), \label{p3_2}
\end{align}
and for $\lambda_{d} \geq \tau \geq \lambda_{d+1}$,
\begin{align}
   \V_{\out}\geq \frac{(1-\nu -\eta_1)^2(\frac{1}{2}-\nu)}{(1+C_1)^2(1+\nu +\eta_1)^4\poly_4(\sigma_x)}\frac{1}{1+\frac{2C_0\sigma_x^2}{\frac{1}{2}-\eta_1}}\sigma^2(\frac{d}{n}+\frac{\lambda_{d+1}^2}{\tau^2}\frac{r_d( \Sigma^2)}{n}).  \label{p3_1}
\end{align}
\end{pro}

Theorem \ref{thm1} can be deduced by combining the bounds for $\B_{\out}$ and $\V_{\out}$ above.
The probability control is determined from the intersection of the relevant events included in the propositions.
\begin{itemize}
    \item  If $\tau\leq \lambda_{d+1}$, then the lower bound for $\V_\out$ in Theorem \ref{thm1}(i) is obtained from (\ref{p3_2}) in Proposition~\ref{pro3}.
    \item If $\lambda_{d+1}\leq \tau\leq \lambda_d$, then the upper bounds for $\B_\out$ and $\V_\out$ in Theorem \ref{thm1}(ii) are obtained from (\ref{p1_1}) and (\ref{p1_2}) in Proposition \ref{pro1}, and the lower bounds for $\B_\out$ and $\V_\out$ in Theorem \ref{thm1}(ii) are obtained from (\ref{p2_1}) in Proposition \ref{pro2} and (\ref{p3_1}) in Proposition~\ref{pro3}.
    \item If $\tau \geq \lambda_d$, then the lower bound for $\B_\out$ in Theorem \ref{thm1}(iii) is obtained from (\ref{p2_1}) in Proposition~\ref{pro2}.
\end{itemize}

\subsection{Proof of Theorem \ref{thm2}}
\label{sec:5_2}

We provide auxiliary bounds for the in-sample squared bias and variance $\B_{\In}$ and $\V_{\In}$ under the small or moderate TER regime (Assumption \ref{ass:3a}).

\begin{pro}[Upper bound of in-sample error with small or moderate TER]
\label{pro4}
Under Assumption \ref{ass:3} and \ref{ass:3a}, for any $\nu$ satisfying $0<\nu < \frac{1}{2}$, the following inequalities hold uniformly in the range of $\tau$ stated with probability at least $1-2 \exp\{-\frac{\nu^2 n}{C_0^2\sigma_x^4}\}-2 \exp\{-\frac{\nu^2 n}{C_0\sigma_x^4} \}-8\exp\{-\frac{n}{C_0}\}$: for $\tau \geq \lambda_{d+1}$,
\begin{align}
    \B_{\In} &\leq  \frac{(1+C_1)^4\poly_8(\sigma_x)}{(1-\nu - \eta_1)^2}    (\Vert\theta^*_{1:d}\Vert_{\Sigma_{1:d}^{-1}}^2 (\frac{1}{\lambda_1}+\frac{1}{\tau})^{-2} + \Vert \theta^{*}_{(d+1):p}\Vert^2_{\Sigma_{(d+1):p}}) ,\label{p4_2b}  \\
    \V_{\In} &\leq (1+C_1)^2\poly_4(\sigma_x)\sigma^2(\frac{d}{n} +\frac{\lambda^2_{d+1}}{\tau^2}\frac{r_d(\Sigma)}{n}  ) . \label{p4_2}
\end{align}
Further with Assumption \ref{ass1}(ii),
\begin{align}
  \B_{\In} &\leq  \frac{(1+C_1)^4\poly_8(\sigma_x)}{(1-\nu - \eta_1)^2}    \Vert\theta^*_{1:d}\Vert_{\Sigma_{1:d}^{-1}}^2 (\frac{1}{\lambda_1}+\frac{1}{\tau})^{-2}  . \label{p4_1}
\end{align}
\end{pro}

\begin{pro}[Lower bound of $\B_\In$ with small or moderate TER]
\label{pro5}
Under Assumption \ref{ass:3}, \ref{ass:3a} and \ref{ass1}(i), for any $\nu$ satisfying $0<\nu < \frac{1}{4}$, the following inequality holds uniformly in the range of $\tau$ stated with probability at least $1-2 \exp\{-\frac{\nu^2 n}{C_0^2\sigma_x^4}\}-8\exp\{-\frac{n}{C_0}\}$: for $\tau\geq \lambda_{d+1}$,
\begin{align}
    \B_{\In} \geq \kappa_1(\tau)\frac{(1-\sqrt{\delta_1})^2}{(1+\nu+\eta_1)^2}\Vert\theta^*_{1:d}\Vert^2_{\Sigma_{1:d}^{-1}}(\frac{1}{\lambda_d}+\frac{1}{\tau})^{-2} , \label{p5_1}
\end{align}
where $\kappa_1(\tau)=\max\{1-(\frac{2 C_0\sigma_x^2(2+C_1)\lambda_{d+1}}{\tau}(1+16(2C_0\sigma_x^2+1)(1+C_1)\frac{\sqrt{\delta_1}}{1-\sqrt{\delta_1}}) + 64\frac{\sqrt{\delta_1}}{1-\sqrt{\delta_1}}),0\}$.
\end{pro}

\begin{pro}[Lower bound of $\V_\In$ with small or moderate TER]
\label{pro6}
Under Assumption \ref{ass:3}, \ref{ass:3a}, \ref{ass1}(i), then for any $\nu$ satisfying $0<\nu < \frac{1}{2}\min\{1,\sigma_x^2\}$, the following inequalities hold uniformly in the range of $\tau$ stated with probability at least $1-2 \exp\{-\frac{\nu^2 n}{C_0^2\sigma_x^4}\}-2 \exp\{-\frac{\nu^2 n}{C_0\sigma_x^4} \}-6\exp\{-\frac{n}{C_0}\}$: \\
(i) for $\tau \leq \lambda_{d+1}$,
\begin{align}
\V_{\In} \geq  \frac{(1-\nu)^2}{\poly_4(\sigma_x)(1+C_1)^2(1+\frac{1}{(1-\nu-\eta_2)^2})^2}\sigma^2(\frac{d}{n} +\frac{r_d^2(\Sigma)}{n^2})  , \label{p6_1}
\end{align}
(ii) for $\lambda_{d+1}\leq \tau \leq \lambda_d$,
\begin{align}
    \V_{\In} \geq \frac{(1-\nu)^2}{\poly_4(\sigma_x)(1+C_1)^2(1+\frac{1}{(1-\nu-\eta_2)^2})^2}\sigma^2(\frac{d}{n} +\frac{\lambda_{d+1}^2}{\tau^2}\frac{r_d^2(\Sigma)}{n^2}) , \label{p6_2}
\end{align}
(iii) for $\tau \geq \lambda_d$,
\begin{align}
     V_{\In} &\geq \frac{(1-\nu)^2}{\poly_4(\sigma_x)(1+C_1)^2} \sigma^2\frac{\lambda_{d+1}^2}{\tau^2}\frac{r_d^2(\Sigma)}{n^2} . \label{p6_3}
\end{align}
\end{pro}

Theorem \ref{thm2} can be deduced by combining the bounds for $\B_{\In}$ and $\V_{\In}$ above.
The probability control is determined from the intersection of the relevant events included in the propositions.
\begin{itemize}
    \item If $\tau \leq \lambda_{d+1}$, then the lower bound for $\V_{\In}$ in Theorem \ref{thm2}(i) is obtained from (\ref{p6_1}) in Proposition~\ref{pro6}.
    \item If $\lambda_{d+1}\leq \tau\leq \lambda_d$, then the upper bounds for $\B_\In$ and $\V_\In$ in Theorem \ref{thm2}(ii) are obtained from (\ref{p4_1}) and (\ref{p4_2}) in Proposition \ref{pro4}, and the lower bounds for $\B_\In$ and $\V_\In$ in Theorem \ref{thm2}(ii) are obtained from (\ref{p5_1}) in Proposition \ref{pro5} and (\ref{p6_2}) in Proposition \ref{pro6}.
    \item If $\tau\geq\lambda_d$, then the lower bounds for $\B_{\In}$ and $\V_{\In}$ in Theorem \ref{thm2}(iii) are obtained from (\ref{p5_1}) in Proposition \ref{pro5} and (\ref{p6_3}) in Proposition \ref{pro6}.
\end{itemize}

\subsection{Proof of Theorem 3}
\label{sec:5_3}

We provide auxiliary bounds for the out-sample squared bias and variance $\B_{\out}$ and $\V_{\out}$ under the large TER regime (Assumption \ref{ass:4}).

\begin{pro}[Upper bound of out-sample error with large TER]
\label{pro7}
Under Assumption \ref{ass:3}, \ref{ass:4}, for any $\nu$ satisfying
$0<\nu < \frac{1}{2}$ and $\frac{r_d(\Sigma)\nu^2 }{C_0^2\sigma_x^4}>1$, the following inequalities hold uniformly in the range of $\tau$ stated with probability at least $1-2 \exp\{-\frac{\nu^2 n}{C_0^2\sigma_x^4}\}-2n\exp\{-\frac{\nu \sqrt{r_d(\Sigma)}}{C_0\sigma_x^2} \}-16 \exp\{-\frac{n}{C_0}\}$: for $\tau \geq 0$,
\begin{align}
    \B_{\out}&\leq   \frac{(1+\nu +\eta_1)^2(1+\nu + \eta_2)^2\poly_4(\sigma_x)}{(1-\nu-\eta_1)^4(1-\nu-\eta_2)^2}(\Vert\theta^*_{1:d}\Vert^2_{\Sigma_{1:d}^{-1}}(\frac{1}{\lambda_1}+\frac{1}{\tau+\lambda_{d+1}\frac{r_d(\Sigma)}{n}})^{-2} +\Vert\theta^*_{(d+1):p}\Vert_{\Sigma_{(d+1):p}}^2), \label{p7_2b} \\
    \V_{\out} &\leq  \frac{(1+\nu+\eta_2)^2\poly_2(\sigma_x)}{(1-\nu-\eta_1)^4(1-\nu-\eta_2)^2}\sigma^2(\frac{d}{n}+\frac{\lambda_{d+1}^2}{(\tau+\lambda_{d+1}\frac{r_d(\Sigma)}{n})^2}\frac{r_d(\Sigma^2)}{n} ) . \label{p7_2}
\end{align}
Further with Assumption \ref{ass1}(ii),
\begin{align}
    \B_{\out}&\leq   \frac{(1+\nu +\eta_1)^2(1+\nu + \eta_2)^2\poly_4(\sigma_x)}{(1-\nu-\eta_1)^4(1-\nu-\eta_2)^2}\Vert\theta^*_{1:d}\Vert^2_{\Sigma_{1:d}^{-1}}(\frac{1}{\lambda_1}+\frac{1}{\tau+\lambda_{d+1}\frac{r_d(\Sigma)}{n}})^{-2}  .\label{p7_1}
\end{align}
\end{pro}

\begin{pro}[Lower bound of $\B_\out$ with large TER]
\label{pro8}
Under Assumption \ref{ass:3}, \ref{ass:4} and \ref{ass1}(ii), for any $\nu$ satisfying $0<\nu<\frac{1}{4}$ and
 $\frac{r_d(\Sigma)\nu^2 }{C_0^2\sigma_x^4}>1$, the following inequality holds uniformly in the range of $\tau$ stated with probability at least $1-2 \exp\{-\frac{\nu n}{C_0^2\sigma_x^4}\}-2n \exp\{-\frac{\nu \sqrt{r_d(\Sigma)}}{C_0\sigma_x^2} \}-6 \exp\{-\frac{n}{C_0}\}$: for $\tau \geq 0$,
\begin{align}
    &\B_{\out} \geq \frac{(1-\sqrt{\delta_2})^2(1-\nu-\eta_2)^2}{(1+\nu+\eta_1)^2}\Vert\theta^*_{1:d}\Vert_{\Sigma_{1:d}^{-1}}^2(\frac{1}{\lambda_d}+\frac{1}{\tau +\lambda_{d+1}\frac{r_d(\Sigma)}{n}})^{-2}. \label{p8_1}
\end{align}
\end{pro}

\begin{pro}[Lower bound of $\V_\out$ with large TER]
\label{pro9}
Under Assumption \ref{ass:3} and \ref{ass:4}, for any $\nu$ satisfying $0<\nu < \frac{1}{2}\min\{1,\sigma_x^2\}$ and
 $\frac{r_d(\Sigma)\nu^2 }{C_0^2\sigma_x^4}>1$, the following inequality holds uniformly in the range of $\tau$ stated with probability $1-2 \exp\{-\frac{\nu^2 n}{C_0\sigma_x^4}\}-2 \exp\{-\frac{\nu^2 n}{C_0^2\sigma_x^4}\}-2n \exp\{-\frac{\nu  \sqrt{r_d(\Sigma)}}{C_0\sigma_x^2} \}-10\exp\{-\frac{n}{C_0}\}$: for $\tau+\frac{\sum_{j>d}\lambda_j}{n}\leq \lambda_d$,
\begin{align}
 \V_{\out} \geq \frac{(1-\nu-\eta_1)^2(1-\nu-\eta_2)^2(\frac{1}{2}-\nu)}{4(1+\nu+\eta_1)^4(1+\nu+\eta_2)^2}\frac{8C_0\sigma_x^2(\frac{1}{2}-\eta_1)}{1+8C_0\sigma_x^2(\frac{1}{2}-\eta_1)} \sigma^2(\frac{d}{n}+\frac{\lambda_{d+1}^2}{(\tau+\lambda_{d+1}\frac{r_d(\Sigma)}{n})^2}\frac{r_d(\Sigma^2)}{n}) . \label{p9_1}
\end{align}
\end{pro}

Theorem \ref{thm3} can be deduced by combining the bounds for $\B_{\out}$ and $\V_{\out}$ above.
The probability control is determined from the intersection of the relevant events included in the propositions.
\begin{itemize}
    \item If $\tau+\lambda_{d+1}\frac{r_d(\Sigma)}{n}\leq \lambda_d$, then the upper bound for $\B_\out$ and $\V_\out$ in Theorem \ref{thm3}(i) are obtained from (\ref{p7_1}) and (\ref{p7_2}) in Proposition \ref{pro7}, and the lower bounds for $\B_\out$ and $\V_\out$ in Theorem \ref{thm3}(i) are obtained from (\ref{p8_1}) in Proposition \ref{pro8} and (\ref{p9_1}) in Proposition \ref{pro9}.
    \item If $\tau+\lambda_{d+1}\frac{r_d(\Sigma)}{n}\geq \lambda_d$, then the lower bound for $\B_{\out}$ in Theorem \ref{thm3}(ii) is obtained from (\ref{p8_1}) in Proposition \ref{pro8}.
\end{itemize}

\subsection{Proof of Theorem 4}
\label{sec:5_4}

We provide auxiliary bounds for the in-sample squared bias and variance $\B_{\In}$ and $\V_{\In}$ under the large TER regime (Assumption \ref{ass:4}).

\begin{pro}[Upper bound of in-sample error with large TER]
\label{pro10}
Given Assumption \ref{ass:3} and \ref{ass:4}, for any $\nu$ satisfying $0<\nu < \frac{1}{2}$ and $\frac{r_d(\Sigma)\nu^2}{C_0^2\sigma_x^4}>1$, the following inequalities hold uniformly in the range of $\tau$ stated with probability at least $1-2 \exp\{-\frac{\nu^2 n}{C_0\sigma_x^4}\}-2 \exp\{-\frac{\nu^2 n}{C_0^2\sigma_x^4}\}-2n \exp\{-\frac{\nu \sqrt{r_d(\Sigma)}}{C_0\sigma_x^2} \}-12 \exp\{-\frac{n}{C_0}\}$: for $\tau \geq 0$,
\begin{align}
    \B_{\In} &\leq  \frac{(1+\nu+\eta_2)^2\poly_6(\sigma_x)}{(1-\nu-\eta_1)^2(1-\nu-\eta_2)^2}(\Vert\theta^*_{1:d}\Vert_{\Sigma_{1:d}^{-1}}^2 (\frac{1}{\lambda_1}+\frac{1}{\tau+\lambda_{d+1}\frac{r_d(\Sigma)}{n}})^{-2} + \Vert \theta^*_{(d+1):p}\Vert^2_{\Sigma_{(d+1):p}}) ,\label{p10_2b} \\
    \V_{\In} &\leq  \frac{(1+\nu+\eta_2)^2\poly_4(\sigma_x)}{(1-\nu-\eta_2)^2}\sigma^2(\frac{d}{n} + \frac{\lambda_{d+1}^2}{(\tau+\lambda_{d+1}\frac{r_d(\Sigma)}{n})^2}\frac{r^2_d(\Sigma)}{n^2}). \label{p10_2}
\end{align}
Further with Assumption \ref{ass1}(ii),
\begin{align}
  \B_{\In} &\leq  \frac{(1+\nu+\eta_2)^2\poly_6(\sigma_x)}{(1-\nu-\eta_1)^2(1-\nu-\eta_2)^2}\Vert\theta^*_{1:d}\Vert_{\Sigma_{1:d}^{-1}}^2 (\frac{1}{\lambda_1}+\frac{1}{\tau+\lambda_{d+1}\frac{r_d(\Sigma)}{n}})^{-2}.   \label{p10_1}
\end{align}
\end{pro}

\begin{pro}[Lower bound of $\B_\In$ with large TER]
\label{pro11}
Under Assumption \ref{ass:3}, \ref{ass:4} and \ref{ass1}(ii), for any $\nu$ satisfying $0<\nu < \frac{1}{4}$ and $\frac{r_d(\Sigma)\nu^2}{C_0^2\sigma_x^4}>1$, the following inequality holds uniformly in the range of $\tau$ stated with probability at least $1-2 \exp\{-\frac{\nu^2 n}{C_0^2\sigma_x^4}\}-2n \exp\{-\frac{\nu \sqrt{r_d(\Sigma)}}{C_0\sigma_x^2} \}-6 \exp\{-\frac{n}{C_0}\}$: for $\tau\geq 0$,
\begin{align}
    \B_{\In} \geq   \frac{(1-\sqrt{\delta_2})^2(1-\nu-\eta_2)^2}{(1+\nu+\eta_1)^2}\kappa_2(\tau)\Vert\theta^*_{1:d}\Vert^2_{\Sigma_{1:d}^{-1}}(\frac{1}{\lambda_d}+\frac{1}{\tau+\lambda_{d+1}\frac{r_d(\Sigma)}{n}})^{-2},  \label{p11_1}
\end{align}
where $\kappa_2(\tau)=\max\{1-(16\frac{\lambda_{d+1}\frac{r_d(\Sigma)}{n}}{\tau+\lambda_{d+1}\frac{r_d(\Sigma)}{n}}(1+112\frac{\sqrt{\delta_2}}{1-\sqrt{\delta_2}})+64\frac{\sqrt{\delta_2}}{1-\sqrt{\delta_2}}),0\}$.
\end{pro}

\begin{pro}[Lower bound of $\V_\In$ with large TER]
\label{pro12}
Under Assumption \ref{ass:3} and \ref{ass:4}, for any $\nu$ satisfying $0<\nu < \frac{1}{2}$ and $\frac{r_d(\Sigma)\nu^2}{C_0^2\sigma_x^4}>1$, the following inequalities hold uniformly in the range of $\tau$ stated with probability at least $1-2 \exp\{-\frac{\nu^2 n}{C_0^2\sigma_x^4}\}-2n \exp\{-\frac{\nu \sqrt{r_d(\Sigma)}}{C_0\sigma_x^2} \}-4 \exp\{-\frac{n}{C_0}\}$: for $\tau + \lambda_{d+1}\frac{r_d(\Sigma)}{n}\leq \lambda_d$,
\begin{align}
    \V_{\In} \geq  \frac{(1-\nu-\eta_2)^2}{2(1+\nu+\eta_2)^2(1+\frac{1}{(1-\nu-\eta_1)^2})}\sigma^2(\frac{d}{n} + \frac{\lambda_{d+1}^2}{(\tau+\lambda_{d+1}\frac{r_d(\Sigma)}{n})^2}\frac{r^2_d(\Sigma)}{n^2}) ,  \label{p12_2}
\end{align}
and for $\tau +\lambda_{d+1}\frac{r_d(\Sigma)}{n} \geq \lambda_{d}$,
\begin{align}
    V_{\In} &\geq  \frac{(1-\nu-\eta_2)^2}{2(1+\nu+\eta_2)^2}\sigma^2\frac{\lambda_{d+1}^2}{(\tau+\lambda_{d+1}\frac{r_d(\Sigma)}{n})^2}\frac{r^2_d(\Sigma)}{n^2}. \label{p12_1}
\end{align}
\end{pro}

Theorem \ref{thm4} can be deduced by combining the bounds for $\B_{\In}$ and $\V_{\In}$ above.
The probability control is determined from the intersection of the relevant events included in the propositions.
\begin{itemize}
    \item If $\tau+\lambda_{d+1}\frac{r_d(\Sigma)}{n} \leq \lambda_d$, then the upper bounds for $\B_\In$ and $\V_\In$ in Theorem \ref{thm4}(i) are obtained from (\ref{p10_1}) and (\ref{p10_2}) of Proposition \ref{pro10}, and the lower bounds for $\B_\In$ and $\V_\In$ in Theorem \ref{thm4}(i) are obtained from (\ref{p11_1}) and (\ref{p12_2}) of Proposition \ref{pro12}.
    \item If $\tau+\lambda_{d+1}\frac{r_d(\Sigma)}{n} \geq \lambda_d$, then the lower bounds for $\B_\In$ and $\V_\In$ are obtained from (\ref{p11_1}) in Proposition \ref{pro11} and (\ref{p12_1}) in Proposition \ref{pro12}.
\end{itemize}

\section{Proof outlines of auxiliary results (Propositions \ref{pro1}--\ref{pro12})}
\label{sec:6}

We provide proof outlines of the auxiliary bounds (Proposition \ref{pro1}--\ref{pro12}) used in the proofs in Section \ref{sec:5}.
In Section \ref{sec_6_1}, we discuss the results for out-sample error including Propositions \ref{pro1}--\ref{pro3} and Propositions \ref{pro7}--\ref{pro9}.
In Section \ref{sec_6_2}, we discuss the results for in-sample error including Propositions \ref{pro4}--\ref{pro6} and Propositions \ref{pro10}--\ref{pro12}.
We introduce the following notation.
\begin{itemize}
 \item $A = XX^\T+ n\tau I_n$,  $A_d = X_{(d+1):p}X_{(d+1):p}^\T+n\tau I_n$.
 \item $X_{1:d}$ denotes the the matrices comprised of the
first $d$ columns of $X$ and $X_{(d+1):p}$ denotes the the matrices comprised of the
last $p-d$ columns of $X$.
\item $\hat{\Sigma}_{1:d}=\frac{X_{1:d}^TX_{1:d}}{n}$, $\hat{\Sigma}_{(d+1):p}=\frac{X_{(d+1):p}^TX_{(d+1):p}}{n}$ and $\hat{\Sigma}_{1:d,(d+1):p}=\frac{X_{1:d}^TX_{(d+1):p}}{n}$.
    \item $H_d=\Sigma^{-1/2}_{1:d}X_{1:d}^\T$, $\hat{H}_d = \hat{\Sigma}^{-1/2}_{1:d}X_{1:d}^\T$.
    \item $M_d = X_{(d+1):p}\Sigma_{(d+1):p}X_{(d+1):p}^\T$, $\hat{M}_d = X_{(d+1):p}\hat{\Sigma}_{(d+1):p}X_{(d+1):p}^\T$.
    \item $\mu_j(M)$ is the $j$-th largest eigenvalue of symmetric semi-positive definite matrix $M$.
    \item $c$, $c_1$, $c_2$, $c_3$ are absolute constants that may differ from line to line.
    \item $c(\sigma_x)$, $c_1(\sigma_x)$, $c_2(\sigma_x)$, $c_3(\sigma_x)$, $c_4(\sigma_x)$ are constants
related to $\sigma_x$ which may differ from line to line.
\item w.h.p. indicates that an event holds with probability at least $1-\frac{n}{c}\exp(-c\sqrt{n})$ for a constant $c$ and sample size $n$.
\end{itemize}

The out-sample squared bias and variance can be decomposed as follows:
\begin{align*}
    \B_{\out} &= \underbrace{ \Vert\theta^*_{1:d} - X_{1:d}^\T A^{-1}X\theta^*\Vert_{\Sigma_{1:d}}^2 }_{\B_{\out,1}}+ \underbrace{\Vert\theta^*_{(d+1):p} - X_{d+1:p}^\T A^{-1}X\theta^*\Vert_{\Sigma_{d+1:p}}^2}_{\B_{\out,2}} , \\
    \V_{\out} &=  \underbrace{\sigma^2\Tr(A^{-1}X_{1:d}\Sigma_{1:d} X_{1:d}^\T A^{-1})}_{\V_{\out,1}} + \underbrace{\sigma^2\Tr( A^{-1}X_{(d+1):p}\Sigma_{(d+1):p}X_{(d+1):p}^\T A^{-1})}_{\V_{\out,2}}.
\end{align*}
Similarly, the in-sample squared bias and variance can be shown to satisfy
\begin{align}
    \B_{\In} &= \underbrace{\Vert\theta^*_{1:d} - X_{1:d}^\T A^{-1}X\theta^*\Vert_{\hat{\Sigma}_{1:d}}^2 }_{\B_{\In,1}} + \underbrace{\Vert\theta^*_{(d+1):p} - X_{d+1:p}^\T A^{-1}X\theta^*\Vert_{\hat{\Sigma}_{d+1:p}}^2 }_{\B_{\In,2}} \notag \\
    &+\underbrace{2(\theta^{*T}_{1:d}-\theta^{*T}X^\T A^{-1}X_{1:d})\hat{\Sigma}_{1:d, (d+1):p}(\theta_{(d+1):p}^* - X_{d+1:p}^\T A^{-1}X\theta^*)}_{\B_{\In,12}} ,\notag  \\
    \V_{\In} &\leq \underbrace{2\sigma^2\Tr(A^{-1}X_{1:d}\hat{\Sigma}_{1:d} X_{1:d}^\T A^{-1})}_{\V_{\In,1}} +  \underbrace{2\sigma^2\Tr( A^{-1}X_{(d+1):p}\hat{\Sigma}_{(d+1):p}X_{(d+1):p}^\T A^{-1})}_{\V_{\In,2}} . \notag
\end{align}
where the Cauchy--Schwartz inequality is used to bound $\V_{\In}$.

\subsection{Derivation of bounds for out-sample error}

\label{sec_6_1}

To derive the bounds for out-sample error, we first build algebraic bounds of the out-sample error. Then we provide intermediate bounds of the out-sample error by controlling some random quantities in the algebraic bounds. We deduce the final bounds mainly by further controlling the extreme eigenvalues $\mu_1(A_d)$ and $\mu_n(A_d)$ (from Lemma \ref{lem3a}) in the intermediate bounds and incorporating Assumption \ref{ass1} (rotational sparsity) to control quantities related to $\Vert \theta_{(d+1):p}^*\Vert_{\Sigma_{(d+1):p}}$.
The details of our derivation are presented in Supplement Section \ref{bounds_pro}.

\subsubsection{Algebraic bounds of out-sample error}

\label{sec:alge_bound_out}

The first step of deriving the bounds is to build algebraic bounds for the out-sample bias and variance.
The lemmas below are inspired by Lemma 27 and 28 in \cite{Tsigler_ridge_2023}.
\setcounter{lem}{0}
\begin{lem}[Algebraic upper bounds of out-sample error] Given invertible $\Sigma_{1:d}$, we have
\label{lem1aa}
\begin{align}
    B_{\out}&\leq  2\Vert\theta^*_{1:d}\Vert^2_{\Sigma_{1:d}^{-1}} (\frac{1}{\lambda_1}+\frac{\mu_n(H_dH_d^T)}{\mu_1(A_d)})^{-2}+ \frac{2\mu_1(H_dH_d^T)}{\mu_n^2(H_dH_d^T)} \frac{\mu_1^2(A_d)}{\mu_n^2(A_d)} \Vert X_{(d+1):p}\theta^*_{(d+1):p}\Vert^2\notag \\
    &+3(\Vert M_d\Vert \frac{\mu_1(H_dH_d^T)}{\mu^2_n(A_d)} \Vert\theta^*_{1:d}\Vert_{\Sigma^{-1}}^2 (\frac{1}{\lambda_1}+\frac{\mu_n(H_dH_d^T)}{\mu_1(A_d)})^{-2}  +\Vert\theta^*_{(d+1):p}\Vert^2_{\Sigma_{(d+1):p}} \notag\\
    &+ \frac{\Vert M_d\Vert}{\mu_n(A)^2}\Vert X_{(d+1):p}\theta^*_{(d+1):p}\Vert^2 ) ,\label{lem_prf_03a}   \\
    \V_{\out} &\leq   \sigma^2 \frac{\mu_1^2(A_d)}{\mu_n^2(A_d)}\frac{\Tr(H_d^\T H_d)}{\mu_d(H_dH_d^\T )^2} +  \sigma^2 \frac{\Tr(M_d)}{\mu_n(A_d)^2} . \notag
\end{align}
\end{lem}

\begin{lem}[Algebraic lower bounds of out-sample error]
\label{lem2a}
Given invertible $\Sigma_{1:d}$, we have
\begin{align*}
    \V_{\out} &\geq   \sigma^2\frac{d\mu_d(H_dH_d^\T )}{\mu_1(A_d)^2}(\frac{1}{\lambda_d}+\frac{\mu_1(H_dH_d^\T )}{\mu_n(A_d)})^{-2} +  \max\{0,\sigma^2\frac{\Tr(M_d)-d\mu_1(M_d)}{\mu_1(A_d)^2} \}.
\end{align*}
Further if $\Vert\theta^*_{1:d}\Vert_{\Sigma_{1:d}^{-1}}\geq \frac{\mu_1(H_dH_d)^{1/2}\Vert X_{(d+1):p}\theta^*_{(d+1):p}\Vert}{\mu_n(A_d)}$, then
\begin{align*}
     \B_{\out}&\geq  (1 - \frac{\mu_1(H_dH_d^T)^{1/2}\Vert X_{(d+1):p}\theta^*_{(d+1):p}\Vert}{\mu_n(A_d)\Vert\theta^*_{1:d}\Vert_{\Sigma_{1:d}^{-1}}})^2\Vert\theta^*_{1:d}\Vert^2_{\Sigma_{1:d}^{-1}}(\frac{1}{\lambda_d}+\frac{\mu_1(H_dH_d^\T )}{\mu_n(A_d)})^{-2}.
\end{align*}
\end{lem}

\subsubsection{Intermediate bounds of out-sample error}

\label{sec_6_1_2}

Based on the sub-gaussianity of the covariate vectors, the following random quantities in Lemmas \ref{lem1aa} and \ref{lem2a}
can be controlled with high probability. See Supplement Lemma \ref{lem3} (i), (ii), (iii) and (iv).

\begin{itemize}
    \item Bounds of eigenvalues of $H_dH_d^\T$: w.h.p.
    \begin{align*}
\mu_1(H_dH_d^\T )\leq c_1  n ,~~~
\mu_d( H_d H_d^\T ) \geq c_2  n.
    \end{align*}
    \item Bounds of traces of random matrices: w.h.p.
    \begin{align*}
& \Tr(H_d^\T H_d)= \Tr(X_{1:d}\Sigma_{1:d}^{-1}X_{1:d}) \leq c_1(\sigma_x) nd ,\\
& \Tr(M_d)=\Tr(X_{(d+1:p)}\Sigma_{(d+1):p}X^\T _{(d+1):p}) \leq c_2(\sigma_x) n\sum_{j>d}\lambda_j^2, \\
&\Tr(M_d)=\Tr(X_{(d+1:p)}\Sigma_{(d+1):p}X^\T _{(d+1):p}) \geq c_1 n\sum_{j>d}\lambda_j^2.
\end{align*}
\item Bounds of norms of random matrix and vector: w.h.p.
\begin{align*}
   \mu_1(M_d)=\Vert M_d \Vert \leq c_1(\sigma_x) (n\lambda_{d+1}^2+\sum_{j>d}\lambda_j^2) ,~~~
   \Vert X_{(d+1):p}\theta^*_{(d+1):p}\Vert^2&\leq c_2(\sigma_x) n \Vert \theta_{(d+1):p}^*\Vert_{\Sigma_{(d+1):p}}^2.
    \end{align*}
\end{itemize}

Substituting the probability bounds above into the algebraic bounds in Lemmas \ref{lem1aa} and \ref{lem2a} yields
the intermediate bounds below with high probability.
$\newline$
\textbf{Upper bounds:}
\begin{align}
     \B_{\out} &\leq c_1
\Vert\theta^*_{1:d}\Vert_{\Sigma_{1:d}^{-1}}^2
(\frac{1}{\lambda_1} + \frac{n}{\mu_1(A_d)})^{-2}
    +c_1(\sigma_x) \frac{\mu_1^2(A_d)}{\mu_n^2(A_d)} \Vert\theta^*_{(d+1):p}\Vert_{\Sigma_{(d+1):p}}^2 \notag \\
    &+c_2(\sigma_x) \frac{(n^2\lambda_{d+1}^2+n\sum_{j>d}\lambda_j^2)}{\mu^2_n(A_d)} \Vert\theta^*_{1:d}\Vert_{\Sigma^{-1}}^2 (\frac{1}{\lambda_1} + \frac{n}{\mu_1(A_d)})^{-2} \notag \\
    &+(3+c_3(\sigma_x) \frac{(n^2\lambda_{d+1}^2+n\sum_{j>d}\lambda_j^2)}{\mu_n(A_d)^2})\Vert\theta^*_{(d+1):p}\Vert_{\Sigma_{(d+1):p}}^2 ,\notag \\
    \V_{\out} &\leq c_1(\sigma_x)\frac{\mu^2_1(A_d)}{\mu_n^2(A_d)}\sigma^2\frac{d}{n} + c_2(\sigma_x)\sigma^2 \frac{n\sum_{j>d}\lambda_j^2}{\mu_n(A_d)^2}. \notag
\end{align}
\textbf{Lower bounds:}
\begin{align}
    \V_{\out} &\geq  c_1\sigma^2\frac{nd}{\mu_1(A_d)^2}(\frac{1}{\lambda_d}+\frac{n}{\mu_n(A_d)})^{-2} + \sigma^2 \max\{0,c_2\frac{n\sum_{j>d}\lambda_j^2}{\mu_1(A_d)^2}(1-c_1(\sigma_x)(\frac{d}{r_d(\Sigma^2)}+\frac{d}{n}))\} . \notag
\end{align}
Further if $\Vert\theta^*_{1:d}\Vert_{\Sigma_{1:d}^{-1}}\geq \frac{\mu_1(H_dH_d)^{1/2}\Vert X_{(d+1):p}\theta^*_{(d+1):p}\Vert}{\mu_n(A_d)}$, then
\begin{align}
    \B_{\out} &\geq c(1 - \frac{c_1(\sigma_x)n\Vert\theta^*_{(d+1):p}\Vert_{\Sigma_{(d+1):p}}}{\mu_n(A_d)\Vert\theta^*_{1:d}\Vert_{\Sigma_{1:d}^{-1}}})^2\Vert\theta^*_{1:d}\Vert_{\Sigma_{1:d}^{-1}}^2(\frac{1}{\lambda_d}+\frac{n}{\mu_n(A_d)})^{-2}  . \notag
\end{align}

\subsubsection{Final bounds for out-sample error}
\label{sec_6_1_3}

The eigenvalues $\mu_1(A_d)$ and $\mu_n(A_d)$ in the intermediate bounds above can be controlled with high probability
respectively in the small or moderate TER regime and the large TER regime. See Supplement Lemma \ref{lem3a}.

\begin{itemize}
    \item Under the small or moderate TER regime (Assumption \ref{ass:3a}), for $\tau\geq \lambda_{d+1}$,
    \begin{align*}
     \mu_1(A_d) \leq c(\sigma_x) n\tau ,
    ~~~~~~~~ \mu_n(A_d) \geq n\tau.
    \end{align*}
    \item Under the large TER regime (Assumption \ref{ass:4}), for $\tau\geq 0$,
    \begin{align*}
        \mu_1(A_d) \leq c_1 n(\tau+\lambda_{d+1}\frac{r_d(\Sigma)}{n}),~~~~~~~~\mu_n(A_d) \geq c_2 n(\tau+\lambda_{d+1}\frac{r_d(\Sigma)}{n}).
    \end{align*}
\end{itemize}

We deduce the final bounds for the out-sample error from the intermediate bounds as follows, mainly by applying the probability bounds on $\mu_1(A_d)$ and $\mu_n(A_d)$ and incorporating Assumption \ref{ass1} (rotational sparsity).

\textbf{Upper bounds of $\B_\out$, $\V_\out$ (Proposition \ref{pro1},\ref{pro7}).} We first substitute the bounds of $\mu_1(A_d)$ and $\mu_n(A_d)$ into the intermediate bounds. Then we incorporate Assumption \ref{ass1} to control quantities related to $\Vert\theta^*_{(d+1):p}\Vert_{\Sigma_{(d+1):p}}$.
This leads to the final upper bounds for $\B_\out$ and $\V_\out$.

\textbf{Lower bound of $\B_{\out}$ (Proposition \ref{pro2},\ref{pro8}).} We first show that, w.h.p.
\begin{align*}
\Vert\theta^*_{1:d}\Vert_{\Sigma_{1:d}^{-1}}\geq \frac{\mu_1(H_dH_d)^{1/2}\Vert X_{(d+1):p}\theta^*_{(d+1):p}\Vert}{\mu_n(A_d)},
\end{align*}
if $\tau\geq \lambda_{d+1}$ under Assumption \ref{ass:3a} (small or modereate TER) or if $\tau\geq 0$ under Assumption \ref{ass:4} (large TER).
Then we substitute the bounds of $\mu_n(A_d)$ into the intermediate lower bound of $\B_\out$ and incorporate
Assumption \ref{ass1} to control quantities related to $\Vert\theta^*_{(d+1):p}\Vert_{\Sigma_{(d+1):p}}^2$. This leads to the final lower bound for $\B_\out$.

\textbf{Lower bound of $\V_\out$ (Proposition \ref{pro3},\ref{pro9}).} The intermediate lower bound of the out-sample variance in
Section~\ref{sec_6_1_2} is
\begin{align}
    \V_{\out} &\geq  c\sigma^2\frac{nd}{\mu_1(A_d)^2}(\frac{1}{\lambda_d}+\frac{n}{\mu_n(A_d)})^{-2} + \sigma^2 \max\{0,c_2\frac{n\sum_{j>d}\lambda_j^2}{\mu_1(A_d)^2}(1-c_1(\sigma_x)(\frac{d}{r_d(\Sigma^2)}+\frac{d}{n}))\}.  \label{eq6_1_3}
\end{align}
To deduce the final lower bound for $\V_{\out}$, our strategy is as follows. We first derive a lower bound for the first term $\sigma^2\frac{nd}{\mu_1(A_d)^2}(\frac{1}{\lambda_d}+\frac{n}{\mu_n(A_d)})^{-2}$. Then we discuss two complementary cases.
See Supplement Section \ref{sec_I_1_4} for details.
The first case is that $\frac{r_d(\Sigma^2)}{d}$ is upper bounded by a constant (possibly depending on $\sigma_x$).
The second case is that $\frac{r_d(\Sigma^2)}{d}$ is large enough such that $(1-c_1(\sigma_x)(\frac{d}{r_d(\Sigma^2)}+\frac{d}{n}))>0$.
Lastly, we show that in these two cases, $\V_\out$ satisfies lower bounds of the same order, which gives the final lower bound for $\V_\out$.

For small or moderate TER, after substituting the bounds of $\mu_1(A_d)$ and $\mu_n(A_d)$ into the first term
of the intermediate lower bound in (\ref{eq6_1_3}), we have for $\tau\geq \lambda_{d+1}$, w.h.p.
\begin{align*}
    \sigma^2\frac{nd}{\mu_1(A_d)^2}(\frac{1}{\lambda_d}+\frac{n}{\mu_n(A_d)})^{-2} \geq c(\sigma_x) \sigma^2\frac{d}{n\tau^2}(\frac{1}{\lambda_d}+\frac{1}{\tau})^{-2},
\end{align*}
which implies that for $\lambda_{d+1}\leq \tau\leq \lambda_d$, w.h.p.
\begin{align}
    \sigma^2\frac{nd}{\mu_1(A_d)^2}(\frac{1}{\lambda_d}+\frac{n}{\mu_n(A_d)})^{-2} \geq c(\sigma_x) \sigma^2\frac{d}{n} .\label{eq:6_1}
\end{align}
Then we discuss two cases which are complementary to each other.
The first case is that $\frac{r_d(\Sigma^2)}{d}$ is upper bounded by a constant $c(\sigma_x)$.
In this case, we have for $\tau\geq \lambda_{d+1}$,
\begin{align*}
    d &\geq \frac{1}{c(\sigma_x)}\frac{\sum_{j>d}\lambda_j^2}{\lambda_{d+1}^2}
      \geq \frac{1}{c(\sigma_x)} \frac{\sum_{j>d}\lambda_j^2}{\tau^2}  ,
\end{align*}
and hence (allowing that $c(\sigma_x)$ below may vary from the previous line)
\begin{align}
 \frac{d}{n} \geq c(\sigma_x)(\frac{d}{n}+\frac{\sum_{j>d}\lambda_j^2}{n\tau^2}).  \label{eq:6_1a}
\end{align}
Combining (\ref{eq6_1_3}), (\ref{eq:6_1}) and (\ref{eq:6_1a}) shows that for $\lambda_{d+1}\leq \tau\leq \lambda_d$, we have, w.h.p.
\begin{align*}
    \V_{\out}\geq c \sigma^2\frac{nd}{\mu_1(A_d)^2}(\frac{1}{\lambda_d}+\frac{n}{\mu_n(A_d)})^{-2} \geq c_1(\sigma_x)\sigma^2\frac{d}{n} \geq  c_2(\sigma_x)\sigma^2(\frac{d}{n}+\frac{\sum_{j>d}\lambda_j^2}{n\tau^2}).
\end{align*}
The second case is that $\frac{r_d(\Sigma^2)}{d}$ is large enough such that $1-c_1(\sigma_x)(\frac{d}{r_d(\Sigma^2)}+\frac{d}{n})>0$.
In this case, after substituting the upper bound of $\mu_1(A_d)$ into the second term
of the intermediate lower bound in (\ref{eq6_1_3}), we have for $\tau \geq \lambda_{d+1}$, w.h.p.
\begin{align}
    &\sigma^2 \max\{0,\frac{n\sum_{j>d}\lambda_j^2}{\mu_1(A_d)^2}(1-c_1(\sigma_x)(\frac{d}{r_d(\Sigma^2)}+\frac{d}{n}))\} \geq c(\sigma_x)\frac{\sum_{j>d}\lambda_j^2}{n\tau^2}. \label{eq:6_1b}
\end{align}
Then combining (\ref{eq6_1_3}), (\ref{eq:6_1}) and (\ref{eq:6_1b}) yields that for $ \lambda_{d+1} \leq \tau \leq \lambda_d$, w.h.p.
\begin{align*}
    \V_{\out}\geq  c(\sigma_x)\sigma^2(\frac{d}{n}+\frac{\sum_{j>d}\lambda_j^2}{n\tau^2}) .
\end{align*}
In conclusion of the two cases, it holds that for $\lambda_{d+1}\leq \tau \leq \lambda_d$, w.h.p.
\begin{align*}
    \V_{\out}\geq  c(\sigma_x)\sigma^2(\frac{d}{n}+\frac{\lambda_{d+1}^2}{\tau^2}\frac{r_d(\Sigma^2)}{n}).
\end{align*}
For $\tau\leq \lambda_{d+1}$, the lower bound of $\V_{\out}$ follows by the monotonicity of variance: $\V_\out$ for $\tau \leq \lambda_{d+1}$
is no smaller than $\V_\out$ for $\tau =\lambda_{d+1}$. Hence for $\tau\leq \lambda_{d+1}$,
\begin{align*}
    \V_{\out}\geq  c(\sigma_x)\sigma^2(\frac{d}{n}+\frac{r_d(\Sigma^2)}{n}).
\end{align*}

The lower bound of $\V_\out$ in the large TER regime can be derived similarly to the small or moderate TER regime.
For succinctness, we omit the associated details.

\subsection{Derivation of bounds for in-sample error}

\label{sec_6_2}

Our strategy for deriving the bounds for in-sample error is similar to that for out-sample error.
The details of our derivation are presented in Supplement Section \ref{sec_s2}.

\subsubsection{Algebraic bounds of in-sample error}

\label{sec_6_2_1}

Similarly to the out-sample error, we first
give the algebraic bounds for in-sample bias and variance.

\begin{lem}[Algebraic upper bounds of in-sample error]
\label{lem3aa}
Given invertible $\hat{\Sigma}_{1:d}$, we have
\begin{align}
     \B_{\In}&\leq  2\Vert\theta^*_{1:d}\Vert^2_{\hat{\Sigma}_{1:d}^{-1}} (\frac{1}{\lambda_1}+\frac{n^2}{\mu_1^2(A_d)})^{-2} + \frac{2}{n} \frac{\mu_1^2(A_d)}{\mu_n^2(A_d)} \Vert X_{(d+1):p}\theta^*_{(d+1):p}\Vert^2  \notag \\
    &+3(\Vert \hat{M}_d\Vert \frac{n}{\mu^2_n(A_d)} \Vert\theta^*_{1:d}\Vert_{\hat{\Sigma}^{-1}}^2 (\frac{1}{\lambda_1}+\frac{n}{\mu_1(A_d)})^{-2}  +\Vert\theta^*_{(d+1):p}\Vert^2_{\hat{\Sigma}_{(d+1):p}} + \frac{\Vert \hat{M}_d\Vert}{\mu_n(A)^2}\Vert X_{(d+1):p}\theta^*_{(d+1):p}\Vert^2 ), \notag  \\
    \V_{\In} &\leq   2\sigma^2 \frac{\mu_1^2(A_d)}{\mu_n^2(A_d)}\frac{d}{n} +  2\sigma^2 \frac{\frac{1}{n}\Tr(X_{(d+1):p}X_{(d+1):p}^\T)\mu_1(X_{(d+1):p}X_{(d+1):p}^\T )}{\mu_n(A_d)^2} .\notag
\end{align}
\end{lem}
\begin{lem}[Algebraic lower bounds of in-sample error]
\label{lem4aa}
Given invertible $\hat{\Sigma}_{1:d}$, we have
\begin{align*}
    \V_\In
    &\geq \sigma^2 \frac{1}{2}(\frac{1}{n}\sum_{i=1}^d \frac{\mu^2_i(X_{1:d}X_{1:d}^\T )}{(\mu_i(X_{1:d}X_{1:d}^\T )+n\tau)^2} +\frac{1}{n}\sum_{i=1}^n \frac{\mu_i^2(X_{(d+1):p}X_{(d+1):p}^\T )}{(\mu_i(X_{(d+1):p}X_{(d+1):p}^\T )+n\tau)^2}).
\end{align*}
Further if $\Vert\theta^*_{1:d}\Vert_{\hat{\Sigma}_{1:d}^{-1}}\geq \frac{n^{1/2}\Vert X_{(d+1):p}\theta^*_{(d+1):p}\Vert}{\mu_n(A_d)}$, then
\begin{align*}
     \B_{\In}&\geq  \max\{0,1-\frac{|\B_{\In,12}|}{\B_{\In}}\}(1 - \frac{n^{1/2}\Vert X_{(d+1):p}\theta^*_{(d+1):p}\Vert}{\mu_n(A_d)\Vert\theta^*_{1:d}\Vert_{\hat{\Sigma}_{1:d}^{-1}}})^2\Vert\theta^*_{1:d}\Vert^2_{\hat{\Sigma}_{1:d}^{-1}}(\frac{1}{\lambda_d}+\frac{n}{\mu_n(A_d)})^{-2}.
\end{align*}
\end{lem}

\subsubsection{Intermediate bounds of in-sample error}

In addition to the probability bounds in Section~\ref{sec_6_1_2}, the following probability bounds can be obtained about random quantities in Lemmas \ref{lem3aa} and \ref{lem4aa}. See Supplement Lemma \ref{lem3} (i), (v) and (vii).

\begin{itemize}
    \item Bounds of $\Vert\theta^*_{1:d}\Vert_{\hat{\Sigma}^{-1}}^2$, w.h.p.
    \begin{align*}
c_1 \Vert\theta^*_{1:d}\Vert^2_{\Sigma_{1:d}^{-1}}\leq \Vert\theta^*_{1:d}\Vert_{\hat{\Sigma}^{-1}_{1:d}}^2 \leq c_2 \Vert\theta^*_{1:d}\Vert^2_{\Sigma_{1:d}^{-1}}.
    \end{align*}
    \item Bounds of $\mu_1(X_{(d+1):p}X_{(d+1):p}^T)$ and $\Vert\hat{M}_d\Vert$, w.h.p.
    \begin{align*}
\mu_1(X_{(d+1):p}X_{(d+1):p}^\T )&\leq  c_1(\sigma_x)(n\lambda_{d+1}+\sum_{j>d}\lambda_j),~~~\Vert \hat{M}_d \Vert\leq  c_2(\sigma_x)\frac{(n\lambda_{d+1}+\sum_{j>d}\lambda_j)^2 }{n}.
\end{align*}
\item Bounds of $\Tr(X_{(d+1):p}X_{(d+1):p}^T)$, w.h.p.
\begin{align*}
    \Tr(X_{(d+1):p}X_{(d+1):p}^T) \leq c n\sum_{j>d}\lambda_j.
\end{align*}
\end{itemize}

Substituting the probability bounds in Section~\ref{sec_6_1_2} and above into the algebraic bounds in Lemmas \ref{lem3aa} and \ref{lem4aa} yields
the intermediate bounds below, w.h.p.
$\newline$
\textbf{Upper bounds:}
\begin{align}
    \B_{\In}&\leq c_1\Vert\theta^*_{1:d}\Vert_{\Sigma_{1:d}^{-1}}^2(\frac{1}{\lambda_1}+\frac{n}{\mu_1(A_d)})^{-2}
    + c_1(\sigma_x) \frac{\mu_1^2(A_d)}{\mu_n^2(A_d)} \Vert \theta^*_{(d+1):p}\Vert_{\Sigma_{(d+1):p}}^2 \notag \\
    &+ c_2(\sigma_x)\frac{(n\lambda_{d+1}+\sum_{j>d}\lambda_j)^2}{\mu^2_n(A_d)}  \Vert\theta^*_{1:d}\Vert_{\Sigma^{-1}}^2  (\frac{1}{\lambda_1}+\frac{n}{\mu_1(A_d)})^{-2}\notag \\
    &+(c_3(\sigma_x)+c_4(\sigma_x)\frac{(n\lambda_{d+1}+\sum_{j>d}\lambda_j)^2}{\mu_n(A_d)^2})\Vert\theta^*_{(d+1):p}\Vert_{\Sigma_{(d+1):p}}^2 ,\notag \\
    \V_{\In}&\leq c_3 \sigma^2 \frac{\mu_1^2(A_d)}{\mu_n^2(A_d)} \frac{d}{n} + c_3(\sigma_x) \sigma^2 \frac{(\sum_{j>d}\lambda_j)(n\lambda_{d+1}+\sum_{j>d}\lambda_j)}{\mu_n(A_d)^2} .  \notag
\end{align}
\textbf{Lower bounds:}
$\newline$
Further if $\Vert\theta^*_{1:d}\Vert_{\hat{\Sigma}_{1:d}^{-1}}\geq \frac{n^{1/2}\Vert X_{(d+1):p}\theta^*_{(d+1):p}\Vert}{\mu_n(A_d)}$, then
\begin{align}
    \B_{\In} &\geq c_1\max\{0,1-\frac{|\B_{\In,12}|}{\B_{\In}}\} (1 - \frac{c_1(\sigma_x)n\Vert\theta^*_{(d+1):p}\Vert_{\Sigma_{(d+1):p}}}{\mu_n(A_d)\Vert\theta^*_{1:d}\Vert_{\Sigma_{1:d}^{-1}}})^2\Vert\theta^*_{1:d}\Vert_{\Sigma_{1:d}^{-1}}^2(\frac{1}{\lambda_d}+\frac{n}{\mu_n(A_d)})^{-2} .  \notag
\end{align}

\subsubsection{Final bounds for in-sample error}

We deduce the final bounds for the in-sample error from the intermediate bounds as follows, mainly by applying the probability bounds on $\mu_1(A_d)$ and $\mu_n(A_d)$ (from Section \ref{sec_6_1_3}) and incorporating Assumption \ref{ass1} (rotational sparsity).

\textbf{Upper bounds of $\B_\In$, $\V_\In$ (Proposition \ref{pro4},\ref{pro10}).} We first substitute the bounds of $\mu_1(A_d)$ and $\mu_n(A_d)$ into the intermediate bounds. Then we incorporate Assumption \ref{ass1} to control quantities related to $\Vert\theta^*_{(d+1):p}\Vert_{\Sigma_{(d+1):p}}^2$.
This leads to the final upper bounds for $\B_\In$ and $\V_\In$.

\textbf{Lower bound of $\B_{\In}$ (Proposition \ref{pro5},\ref{pro11}).} We first show that w.h.p.
\begin{align*}
\Vert\theta^*_{1:d}\Vert_{\hat{\Sigma}_{1:d}^{-1}}\geq \frac{n^{1/2}\Vert X_{(d+1):p}\theta^*_{(d+1):p}\Vert}{\mu_n(A_d)},
\end{align*}
if $\tau\geq \lambda_{d+1}$ under Assumption \ref{ass:3a} (small or modereate TER) or if $\tau\geq 0$ under Assumption \ref{ass:4} (large TER).
Then we show that w.h.p.
\begin{align*}
\max\{0,1-\frac{|\B_{\In,2}|}{\B_{\In,1}}\}\geq \kappa_1(\tau),
\end{align*}
if $\tau\geq \lambda_{d+1}$ under Assumption \ref{ass:3a} (small or moderate regime TER regime) and Assumption \ref{ass1}(i) , or
\begin{align*}
\max\{0,1-\frac{|\B_{\In,2}|}{\B_{\In,1}}\}\geq \kappa_2(\tau),
\end{align*}
if $\tau\geq 0$ under Assumption \ref{ass:4} (large TER regime) and Assumption \ref{ass1}(ii).
See Theorems \ref{thm2} and \ref{thm4} for the definition of $\kappa_1(\tau)$ and $\kappa_2(\tau)$, and Supplement Lemma \ref{lem6a} for details.
We substitute the bounds of $\mu_n(A_d)$ into the intermediate lower bound of $\B_\In$ and incorporate Assumption \ref{ass1} to control quantities related to $\Vert\theta^*_{(d+1):p}\Vert_{\Sigma_{(d+1):p}}$. This leads to the final lower bound of $\B_{\In}$.

\textbf{Lower bound of $\V_{\In}$ (Proposition \ref{pro6},\ref{pro12}).} The algebraic lower bound of in-sample variance
from Section~\ref{sec_6_2_1} is
\begin{align}
    \V_{\In}
    &\geq \sigma^2 \frac{1}{2}(\frac{1}{n}\sum_{i=1}^d \frac{\mu^2_i(X_{1:d}X_{1:d}^\T )}{(\mu_i(X_{1:d}X_{1:d}^\T )+n\tau)^2} +\frac{1}{n}\sum_{i=1}^n \frac{\mu_i^2(X_{(d+1):p}X_{(d+1):p}^\T )}{(\mu_i(X_{(d+1):p}X_{(d+1):p}^\T )+n\tau)^2}) .  \label{sec_6_2_3_1}
\end{align}
By the concentration of $X_{1:d}^\T X_{1:d}$ (see Supplement Lemma \ref{lem3} (i)), we have w.h.p.
\begin{align*}
   \mu_1(X_{1:d}^\T X_{1:d}) \geq \cdots \geq \mu_d(X_{1:d}^\T X_{1:d}) \geq  c\lambda_d.
\end{align*}
If $\tau \leq \lambda_d$, then the first term in the algebraic bound satisfies w.h.p.
\begin{align*}
   \frac{1}{n}\sum_{i=1}^d \frac{\mu^2_i(X_{1:d}X_{1:d}^\T )}{(\mu_i(X_{1:d}X_{1:d}^\T )+n\tau)^2} \geq c\frac{d}{n}.
\end{align*}
For the second term of the algebraic lower bound, we first give the algebraic bound
\begin{align*}
    \sigma^2\frac{1}{n}\sum_{i=1}^n \frac{\mu_i^2(X_{(d+1):p}X_{(d+1):p}^\T )}{(\mu_i(X_{(d+1):p}X_{(d+1):p}^\T )+n\tau)^2}
    &\geq \sigma^2\frac{1}{n^2} \frac{\Tr(X_{(d+1):p}X_{(d+1):p}^\T )^2}{\mu_1(A_d)^2} .
\end{align*}
By the control of $\mu_1(A_d)$ in Section \ref{sec_6_1_3} and the fact that w.h.p.~(see Supplement Lemma \ref{lem3} (vii))
\begin{align*}
    \Tr(X_{(d+1):p}X_{(d+1):p}^\T ) \geq c n\sum_{j>d}\lambda_j,
\end{align*}
the second term in the algebraic bound (\ref{sec_6_2_3_1}) can be lower bounded as follows.
\begin{itemize}
    \item In small or moderate TER regime, under Assumption \ref{ass:3a}, for $\tau\geq \lambda_{d+1}$, w.h.p.
    \begin{align*}
        \sigma^2\frac{1}{n}\sum_{i=1}^n \frac{\mu_i^2(X_{(d+1):p}X_{(d+1):p}^\T )}{(\mu_i(X_{(d+1):p}X_{(d+1):p}^\T )+n\tau)^2}
    &\geq  c(\sigma_x) \sigma^2\frac{\lambda_{d+1}^2}{\tau^2}\frac{r_d^2(\Sigma)}{n^2}.
    \end{align*}
    \item  In large TER regime, under Assumption \ref{ass:4}, for $\tau\geq 0$, w.h.p.
    \begin{align*}
        \sigma^2\frac{1}{n}\sum_{i=1}^n \frac{\mu_i^2(X_{(d+1):p}X_{(d+1):p}^\T )}{(\mu_i(X_{(d+1):p}X_{(d+1):p}^\T )+n\tau)^2}
    &\geq  c \sigma^2\frac{\lambda_{d+1}^2}{(\tau+\lambda_{d+1}\frac{r_d(\Sigma)}{n})^2}\frac{r^2_d(\Sigma)}{n^2}.
    \end{align*}
\end{itemize}

\noindent Combining the preceding bounds on the two terms of (\ref{sec_6_2_3_1}) gives the lower bound of $\V_{\In}$.

\section{Numerical studies}

We present numerical results in support of our theoretical results, including
 the sufficient conditions and necessary conditions for $\MSE=O(\frac{d}{n})$ and
 the conditions for when $\MSE_\out^*$ can be much smaller than $\MSE_\In^*$ as described in Remark \ref{in-out-com-small_TER} and \ref{in-out-com-large_TER}.

\subsection{Data generation and MSE calculation}

We first generate the covariance matrix $\Sigma$ and coefficient vector $\theta^*$.

\textbf{Generating $\Sigma$.} Given $0<d<p$ and $\rho < 1$, we generate a diagonal covariance matrix $\Sigma$ as follows. We let $\lambda_i = 1$ for $i=1,\ldots, d$ and let $\lambda_i = \rho$ for $i=d+1,\ldots, p$ unless otherwise stated. Here $\rho$ represents the gap between the spiked and tail eigenvalues of $\Sigma$.

\textbf{Generating $\theta^*$.} Given the covariance matrix $\Sigma$, we generate $\theta^*\in \mathbb{R}^p$ as follows. We let $\theta^*_{1:d}=\frac{1}{\sqrt{d}}$. To generate $\theta_{(d+1):p}^*$, we first generate $\beta_{(d+1):p}\sim N_{p-d}(0,I)$ and then let
\begin{align*}
\theta_{(d+1):p}^*=\left\{ \begin{array}{ll}
\frac{\beta_{(d+1):p}}{\Vert \beta_{(d+1):p}\Vert_{\Sigma_{(d+1):p}}}\sqrt{0.01||\theta^*_{1:d}||^2_{\Sigma_{1:d}^{-1}}\lambda_{d+1}^2},
 & \text{if } r_d(\Sigma)< 10 n,\\
\frac{\beta_{(d+1):p}}{\Vert \beta_{(d+1):p}\Vert_{\Sigma_{(d+1):p}}}\sqrt{0.01||\theta^*_{1:d}||^2_{\Sigma_{1:d}^{-1}} (\frac{1}{\lambda_{d}}+ \frac{n}{\sum_{i>d}\lambda_i})^{-2}}, & \text{if } r_d(\Sigma)\geq 10 n.
\end{array} \right.
\end{align*}
In the numerical study, we consider $r_d(\Sigma)<10n$ as the small or moderate TER regime and consider $r_d(\Sigma)\geq 10n$ as the large TER regime. From the generating process above, $\Vert \theta_{1:d}^*\Vert_{\Sigma_{1:d}^{-1}}^2\lambda^2_d=1$ and $||\theta_{(d+1):p}||^2_{\Sigma_{(d+1):p}}$ satisfies rotational sparsity Assumption \ref{ass1}(i) in small or moderate TER regime or satisfies rotational sparsity Assumption \ref{ass1}(ii) in large TER regime.

\textbf{Generating $x_i$ and $y_i$.} Given $\Sigma$ and $\theta^*$ from above, we generate data $x_i$ and $y_i$ for $i=1,\ldots,n$ as follows. We sample $z_{1i}\sim unif(\sqrt{p}S^{p-1})$ and $z_{2i}\sim N_p(0, I)$, where $S^{p-1}$ is the spherical surface with radius 1 in $\mathbb{R}^p$. Then we let $z_i=\frac{\sqrt{2}}{2}z_{1i}+\frac{\sqrt{2}}{2}z_{2i}$ and $x_i=\Sigma^{1/2}z_i$. By the generating process,
$z_i\in\mathbb{R}^p$ is an isotropic random vector with dependent components. Then we sample $\epsilon_i\sim N(0,1)$ and generate $y_i = x_i^\T \theta^*+\epsilon_i$ by model (\ref{eq:1}).

With $x_i$ and $y_i$ generated from above, for different ridge parameters $\tau$, we calculate $\MSE_{\out}$ or $\MSE_{\In}$ according to (\ref{eq:3.3}) or (\ref{eq:3.4}), where $y_i$'s are averaged out. We report the $\MSE_{\out}$ and $\MSE_{\In}$ based on the average of 10 repeated runs of data generation.

\subsection{Experiment settings}

\subsubsection{Study of conditions for $\MSE=O(\frac{d}{n})$}

\label{sec_7_2_1}

We use the following settings to study the sufficient conditions and the necessary conditions for $\MSE=O(\frac{d}{n})$ in Corollary \ref{cor1}, \ref{cor2}, \ref{cor3} and \ref{cor4}. We focus on the scenarios where the sufficient condition matches the necessary condition up to a constant, that is,
$\MSE=O(\frac{d}{n})$ if and only if the ratio $\frac{\lambda_{d+1}}{\lambda_{d}}$ is smaller than or equal to a certain threshold.

\textbf{Study of Corollary \ref{cor1}.} Given the small or moderate TER regime, $n \gg d$ and $r_d(\Sigma^2)\gg d$, from Corollary \ref{cor1},
the sufficient condition for $\MSE_\out=O(\frac{d}{n})$ is $\frac{\lambda_{d+1}}{\lambda_{d}}\lesssim  \sqrt{\frac{d}{n}}\sqrt{\frac{d}{r_d(\Sigma^2)}}$ and the necessary condition for $\MSE_\out=O(\frac{d}{n})$ is $\frac{\lambda_{d+1}}{\lambda_{d}}\lesssim \sqrt{\frac{d}{n}}\sqrt{\frac{d}{r_d(\Sigma^2)}}$. To embody this condition, we set $d=5, n=1500, p=1500$ and $\rho=[0.1,1,10]\times\sqrt{\frac{d}{n}}\sqrt{\frac{d}{r_d(\Sigma^2)}}$ such that $r_d(\Sigma^2)=1495$.

\textbf{Study of Corollary \ref{cor2}.} Given the small or moderate TER regime, $n\gg d$ and $r_d(\Sigma)\asymp n$, from Corollary \ref{cor2}, the sufficient condition for $\MSE_\In=O(\frac{d}{n})$ is $\frac{\lambda_{d+1}}{\lambda_{d}}\lesssim \frac{d}{r_d(\Sigma)}$ and the necessary condition for $\MSE_\In=O(\frac{d}{n})$ is $\frac{\lambda_{d+1}}{\lambda_{d}}\lesssim \frac{d}{r_d(\Sigma)}$.
To embody this condition, we set $d=5, n=1500$, $p=1500$ and $\rho=[0.1,1,10]\times\frac{d}{r_d(\Sigma)}$ such that $r_d(\Sigma)=1495$.

\textbf{Study of Corollary \ref{cor3}.}
Given the large TER regime and $n\gg d$, from Corollary \ref{cor3}, the sufficient condition for $\MSE_\out=O(\frac{d}{n})$ and the necessary condition for $\MSE_\out=O(\frac{d}{n})$ are the same. The condition is $\frac{\lambda_{d+1}}{\lambda_{d}}\lesssim \frac{\sqrt{nd}}{r_d(\Sigma)}$ if $\frac{n\sqrt{r_d(\Sigma^2)}}{\sqrt{d}r_d(\Sigma)}\leq 1$ and the condition is
$\frac{\lambda_{d+1}}{\lambda_{d}}\lesssim \frac{d}{\sqrt{n r_d(\Sigma^2)}}$ if $\frac{n\sqrt{r_d(\Sigma^2)}}{\sqrt{d}r_d(\Sigma)}> 1$. To embody the first condition, we set $d=5$, $ n=50$, $p=1500$ such that $\frac{n\sqrt{r_d(\Sigma^2)}}{\sqrt{d}r_d(\Sigma)}=0.577\leq 1$ and $\rho=[0.1,1,10]\times\frac{\sqrt{nd}}{r_d(\Sigma)}$. To embody the second condition, we set $d=5, n=150$, $p=1500$ such that $\frac{n\sqrt{r_d(\Sigma^2)}}{\sqrt{d}r_d(\Sigma)}=1.732> 1$ and $\rho=[0.1,1,10]\times\frac{d}{\sqrt{n r_d(\Sigma^2)}}$.

\textbf{Study of Corollary \ref{cor4}.}
Given the large TER regime and $n\gg d$, from Corollary \ref{cor4}, the sufficient condition for $\MSE_\In=O(\frac{d}{n})$ and the necessary condition for $\MSE_\In=O(\frac{d}{n})$ are the same, and the condition is $\frac{\lambda_{d+1}}{\lambda_{d}}\lesssim \frac{d}{r_d(\Sigma)}$. To embody this condition, we set $d=5$, $ n=150,p=1500$ and $\rho=[0.1,1,10]\times\frac{d}{r_d(\Sigma)}$.

\subsubsection{Study of conditions for $\MSE_\out^*$ much smaller than $\MSE_\In^*$}
\label{sec_7_2_2}

We use the following settings to study the conditions for when $\MSE_\out^*$ can be much smaller than $\MSE_\In^*$, which is discussed in Remark \ref{in-out-com-small_TER} and \ref{in-out-com-large_TER}.

\noindent (i) In the small or moderate TER regime, $\MSE_\out^*$ can be much smaller than $\MSE_\In^*$ if
\begin{align*}
r_d(\Sigma)\asymp n, \; \frac{\lambda_{d+1}}{\lambda_d}\gtrsim \frac{d}{n}\sqrt{\frac{n}{r_d(\Sigma^2)}}, \;
\frac{n}{r_d(\Sigma^2)}\gg 1.
\end{align*}
To embody this condition, we set $d=2$, $n=300$, $p=15000$, $\lambda_{1}=\cdots=\lambda_d=1$, $\lambda_{d+1}=\cdots = \lambda_{11d}=\rho$, $\lambda_{11d+1}=\cdots =\lambda_p=0.02\rho$ and $\rho = \frac{d}{n}\sqrt{\frac{n}{r_d(\Sigma^2)}}$ such that $r_d(\Sigma)=319.56$, $\frac{\lambda_{d+1}}{\lambda_d}=\frac{d}{n}\sqrt{\frac{n}{r_d(\Sigma^2)}}$ and $\frac{n}{r_d(\Sigma^2)}=11.54$.

\noindent (ii) In the large TER regime, $\MSE_\out^*$ can be much smaller than $\MSE_\In^*$ if
\begin{align*}
    \frac{\lambda_{d+1}}{\lambda_d}\gtrsim \frac{d}{\sqrt{nr_d(\Sigma^2)}} \min \{ 1, \frac{n\sqrt{r_d(\Sigma^2)}}{\sqrt{d}r_d(\Sigma)}\},\; \frac{\lambda_{d+1}}{\lambda_d}\gtrsim \frac{n\sqrt{nr_d(\Sigma^2)}}{r_d(\Sigma)^2}.
\end{align*}
To embody this condition, we set $d=2$, $n=150$, $p=15000$ and $\rho = \frac{d}{\sqrt{nr_d(\Sigma^2)}} \min \{ 1, \frac{n\sqrt{r_d(\Sigma^2)}}{\sqrt{d}r_d(\Sigma)}\}$ such that $\frac{\lambda_{d+1}}{\lambda_d}=\frac{d}{\sqrt{nr_d(\Sigma^2)}} \min \{ 1, \frac{n\sqrt{r_d(\Sigma^2)}}{\sqrt{d}r_d(\Sigma)}\}$, $\frac{\lambda_{d+1}}{\lambda_d}=0.0011549$ and $\frac{n\sqrt{nr_d(\Sigma^2)}}{r_d(\Sigma)^2}=0.0010002$.

\noindent (iii) In the large TER regime, $\MSE_\out^*$ can be much smaller than $\MSE_\In^*$ if
\begin{align*}
    \frac{\lambda_{d+1}}{\lambda_d}\gtrsim \frac{d}{\sqrt{nr_d(\Sigma^2)}} \min \{ 1, \frac{n\sqrt{r_d(\Sigma^2)}}{\sqrt{d}r_d(\Sigma)}\},\; \frac{\lambda_{d+1}}{\lambda_d}\lesssim \frac{n\sqrt{nr_d(\Sigma^2)}}{r_d(\Sigma)^2},\;\frac{r_d(\Sigma)^2}{n r_d(\Sigma^2)}\gg 1.
\end{align*}
To embody this condition, we set $d=2$, $n=300$, $p=15000$ and $\rho = \frac{d}{\sqrt{nr_d(\Sigma^2)}} \min \{ 1, \frac{n\sqrt{r_d(\Sigma^2)}}{\sqrt{d}r_d(\Sigma)}\}$ such that $\frac{\lambda_{d+1}}{\lambda_d}=\frac{d}{\sqrt{nr_d(\Sigma^2)}} \min \{ 1, \frac{n\sqrt{r_d(\Sigma^2)}}{\sqrt{d}r_d(\Sigma)}\}$, $\frac{\lambda_{d+1}}{\lambda_d}=0.0009429$, $\frac{n\sqrt{nr_d(\Sigma^2)}}{r_d(\Sigma)^2}=0.0028290$ and $\frac{r_d(\Sigma)^2}{n r_d(\Sigma^2)}=50$.

\vspace{.2in}
\begin{figure}[!h]  
  \begin{subfigure}[t]{.32\textwidth}
    \centering
    \includegraphics[width=\linewidth]{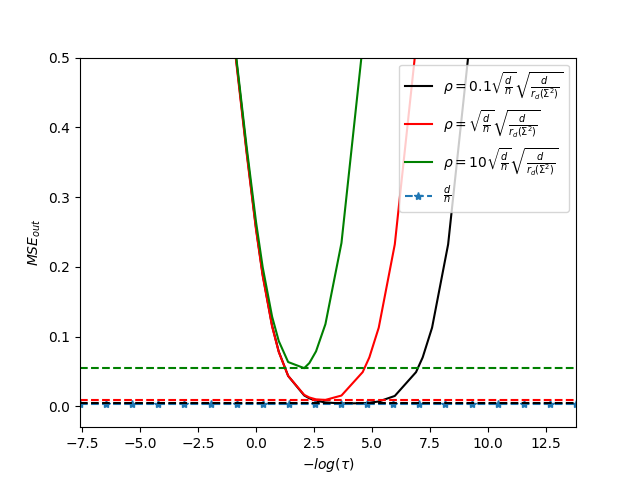}
    \caption{\centering \tiny Corollary \ref{cor1}: Small or moderate TER, \newline $n\gg d$, $r_d(\Sigma^2)\gg d$. }
  \end{subfigure}
  \begin{subfigure}[t]{.32\textwidth}
    \centering
    \includegraphics[width=\linewidth]{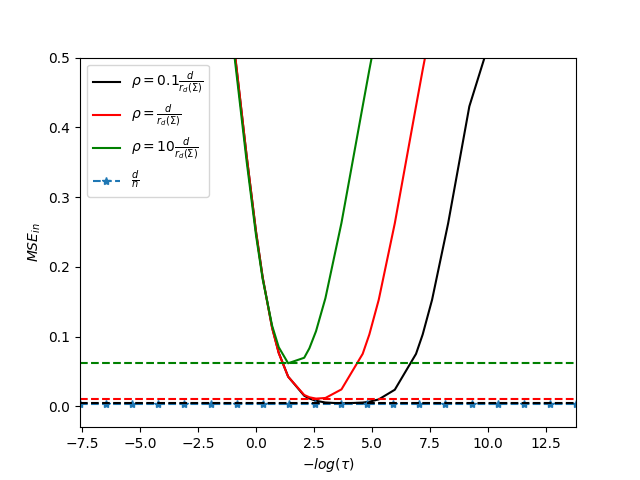}
    \caption{\centering \tiny
    Corollary \ref{cor2}: Small or moderate TER,\newline $n\gg d$ and $r_d(\Sigma)\asymp n$.}
  \end{subfigure}\\
  \begin{subfigure}[t]{.32\textwidth}
    \centering
    \includegraphics[width=\linewidth]{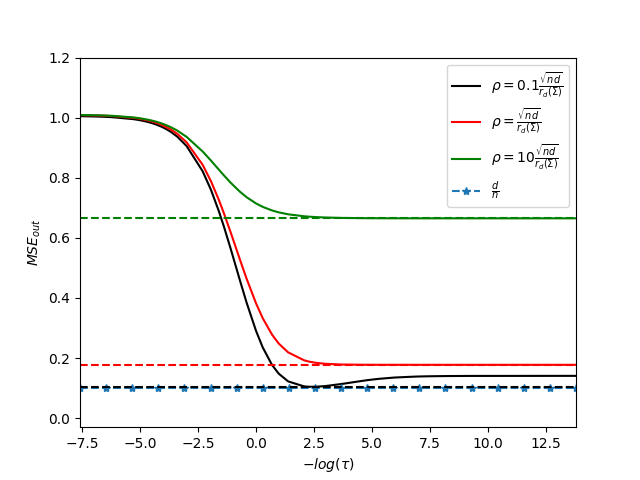}
    \caption{\centering \tiny Corollary \ref{cor3}: Large TER,\newline $n\gg d$ and $\frac{n\sqrt{r_d(\Sigma^2)}}{\sqrt{d}r_d(\Sigma)}\leq 1$.}
  \end{subfigure}
  \hfill
  \centering
  \begin{subfigure}[t]{.32\textwidth}
    \centering
    \includegraphics[width=\linewidth]{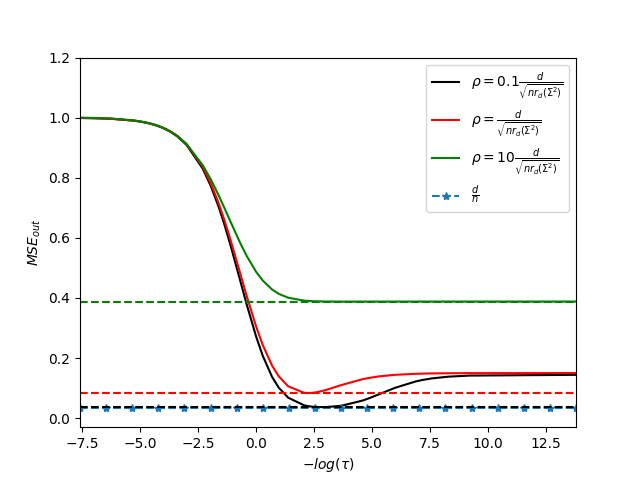}
    \caption{\centering \tiny Corollary \ref{cor3}: Large TER, \newline $n\gg d$ and $\frac{n\sqrt{r_d(\Sigma^2)}}{\sqrt{d}r_d(\Sigma)}> 1$.}
  \end{subfigure}
  \begin{subfigure}[t]{.32\textwidth}
    \centering
    \includegraphics[width=\linewidth]{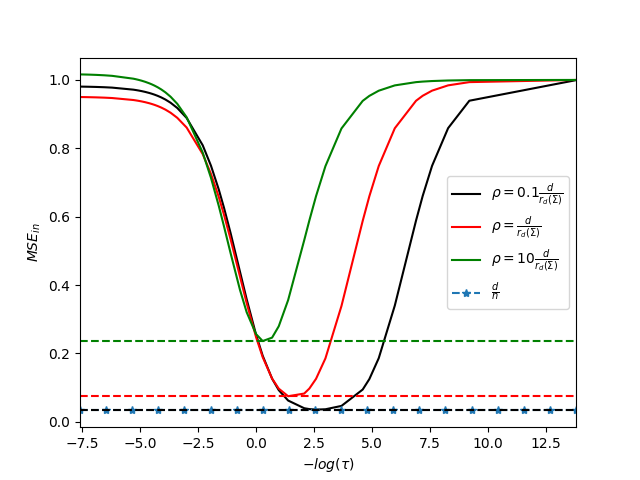}
    \caption{\centering \tiny Corollary \ref{cor4}: Large TER,\newline $n\gg d$.}
  \end{subfigure}
  \caption{Study of conditions for $\MSE=O(\frac{d}{n})$.}
  \label{figure1}
\end{figure}

\begin{figure}[!h]  
  \begin{subfigure}[t]{.32\textwidth}
    \centering
    \includegraphics[width=\linewidth]{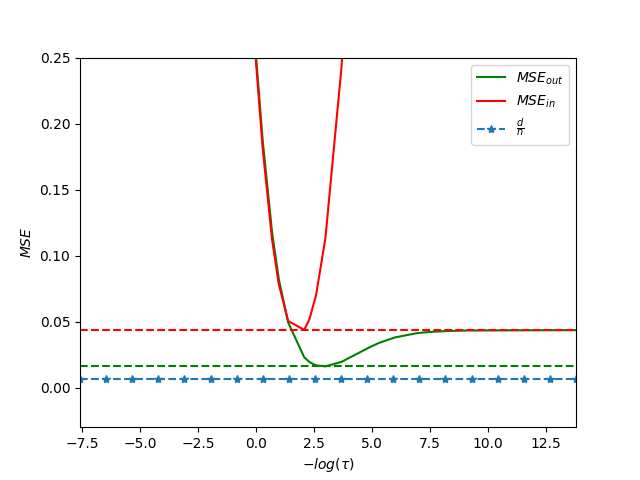}
    \caption{\centering \tiny Small or moderate TER: \newline
 $r_d(\Sigma)\asymp n$, $ \frac{\lambda_{d+1}}{\lambda_d}\gtrsim \frac{d}{n}\sqrt{\frac{n}{r_d(\Sigma^2)}}$,
$\frac{n}{r_d(\Sigma^2)}\gg 1$. }
  \end{subfigure}
  \begin{subfigure}[t]{.32\textwidth}
    \centering
    \includegraphics[width=\linewidth]{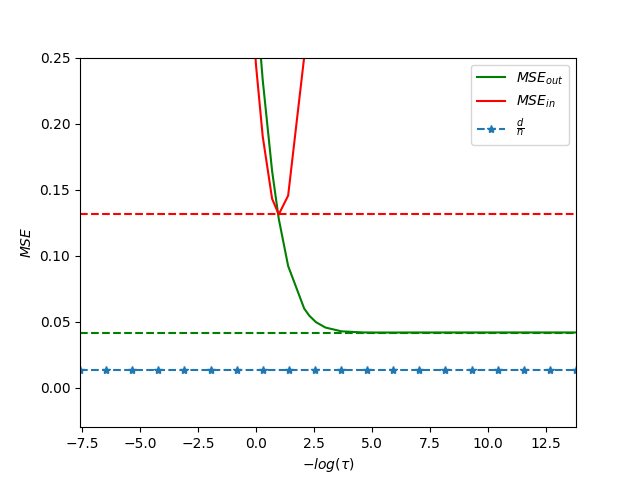}
    \caption{\centering \tiny
    Large TER: $\frac{\lambda_{d+1}}{\lambda_d}\gtrsim \frac{d}{\sqrt{nr_d(\Sigma^2)}} \min \{ 1, \frac{n\sqrt{r_d(\Sigma^2)}}{\sqrt{d}r_d(\Sigma)}\}$, $\frac{\lambda_{d+1}}{\lambda_d}\gtrsim \frac{n\sqrt{nr_d(\Sigma^2)}}{r_d(\Sigma)^2}$.}
  \end{subfigure}
  \begin{subfigure}[t]{.32\textwidth}
    \centering
    \includegraphics[width=\linewidth]{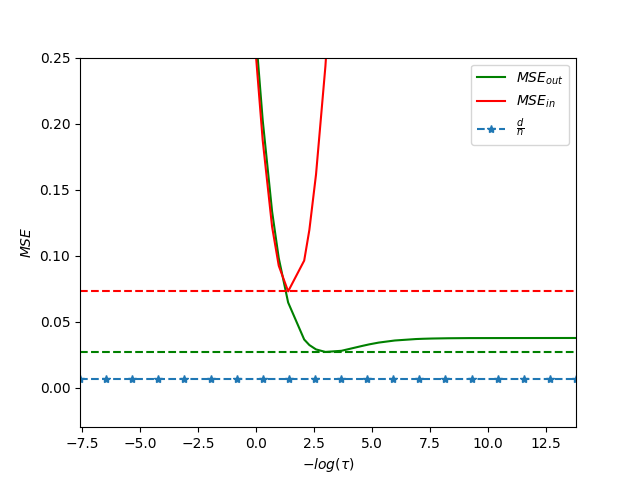}
    \caption{\centering \tiny Large TER: $\frac{\lambda_{d+1}}{\lambda_d}\gtrsim \frac{d}{\sqrt{nr_d(\Sigma^2)}} \min \{ 1, \frac{n\sqrt{r_d(\Sigma^2)}}{\sqrt{d}r_d(\Sigma)}\}$, $\frac{\lambda_{d+1}}{\lambda_d}\lesssim \frac{n\sqrt{nr_d(\Sigma^2)}}{r_d(\Sigma)^2}$ and $\frac{r_d(\Sigma)^2}{n r_d(\Sigma^2)}\gg 1$.}
  \end{subfigure}
  \hfill
  \centering
  \caption{Study of conditions for $\MSE_\out^*$ much smaller than $\MSE_\In^*$.}
  \label{figure2}
\end{figure}

\subsection{Results}

The numerical results are summarized in Figure \ref{figure1} and \ref{figure2}. From Figure \ref{figure1}, we see that when the
ratio $\frac{\lambda_{d+1}}{\lambda_{d}}$ is related to the thresholds described in Section \ref{sec_7_2_1} by a pre-factor equal to $0.1$ or 1, but not 10,
the MSEs with near optimal choices of $\tau$ are close to $\frac{d}{n}$, which gives numerical support to our conditions for $\MSE=O(\frac{d}{n})$ discussed in Corollary \ref{cor1}, \ref{cor2}, \ref{cor3} and \ref{cor4}. From Figure \ref{figure2}, we see that $\MSE^*_{\out}$ is much smaller than $\MSE^*_{\In}$,
each associated with the optimal choices of $\tau$, in the settings described in Section \ref{sec_7_2_2} which embody the conditions from our theory for when $\MSE_\out^*$ can be much smaller than $\MSE_\In^*$ in Section \ref{main_results_1} and \ref{main_results_2}.

\bibliographystyle{apacite}
\bibliography{reference}



\clearpage

\setcounter{page}{1}

\setcounter{section}{0}
\setcounter{equation}{0}

\setcounter{figure}{0}
\setcounter{table}{0}

\renewcommand{\theequation}{S\arabic{equation}}
\renewcommand{\thesection}{\Roman{section}}

\renewcommand\thefigure{S\arabic{figure}}
\renewcommand\thetable{S\arabic{table}}

\setcounter{lem}{0}
\renewcommand{\thelem}{S\arabic{lem}}

\begin{center}
{\Large Supplementary Material for}

{\Large ``On Ridge Estimation in High-dimensional Rotationally Sparse Linear Regression''}

\end{center}

\section{Definition of sub-gaussianity}

\label{sec_s0}

A random variable $z\in \mathbb{R}$ is sub-gaussian if it has a finite sub-gaussian norm
\begin{align*}
    \Vert z\Vert_{\psi_2} = \inf\{t>0:\mathbb{E}\exp(z^2/t^2)\leq 2\}.
\end{align*}
The sub-gaussian norm of a random vector $Z\in\mathbb{R}^p$ is
\begin{align*}
    \Vert Z\Vert_{\psi_2} = \sup_{s\neq 0}\Vert \frac{ \langle s, Z\rangle}{\Vert s\Vert} \Vert_{\psi_2}.
\end{align*}

\section{Proofs of main results}

\label{sec_s1}

We provide proofs of Propositions \ref{pro1}--\ref{pro12} in Section \ref{sec:5}.
For convenience, we re-state the following notation from Section \ref{sec:6}.
\begin{itemize}
 \item $A = XX^\T+ n\tau I_n$,  $A_d = X_{(d+1):p}X_{(d+1):p}^\T+n\tau I_n$.
 \item $X_{1:d}$ denotes the the matrices comprised of the
first $d$ columns of $X$ and $X_{(d+1):p}$ denotes the the matrices comprised of the
last $p-d$ columns of $X$.
\item $\hat{\Sigma}_{1:d}=\frac{X_{1:d}^TX_{1:d}}{n}$, $\hat{\Sigma}_{(d+1):p}=\frac{X_{(d+1):p}^TX_{(d+1):p}}{n}$ and $\hat{\Sigma}_{1:d,(d+1):p}=\frac{X_{1:d}^TX_{(d+1):p}}{n}$.
    \item $H_d=\Sigma^{-1/2}_{1:d}X_{1:d}^\T$, $\hat{H}_d = \hat{\Sigma}^{-1/2}_{1:d}X_{1:d}^\T$.
    \item $M_d = X_{(d+1):p}\Sigma_{(d+1):p}X_{(d+1):p}^\T$, $\hat{M}_d = X_{(d+1):p}\hat{\Sigma}_{(d+1):p}X_{(d+1):p}^\T$.
    \item $\mu_j(M)$ is the $j$-th largest eigenvalue of symmetric semi-positive definite matrix $M$.
\end{itemize}

\subsection{Proof of the bounds for out-sample error}
\label{bounds_pro}

\subsubsection{Algebraic bounds of the out-sample error}

The bias and variance of the out-sample error can be decomposed as follows:
\begin{align*}
    \B_{\out} &= \underbrace{ \Vert\theta^*_{1:d} - X_{1:d}^\T A^{-1}X\theta^*\Vert_{\Sigma_{1:d}}^2 }_{\B_{\out,1}}+ \underbrace{\Vert\theta^*_{(d+1):p} - X_{d+1:p}^\T A^{-1}X\theta^*\Vert_{\Sigma_{d+1:p}}^2}_{\B_{\out,2}} , \\
    \V_{\out} &=  \underbrace{\sigma^2\Tr(A^{-1}X_{1:d}\Sigma_{1:d} X_{1:d}^\T A^{-1})}_{\V_{\out,1}} + \underbrace{\sigma^2\Tr( A^{-1}X_{(d+1):p}\Sigma_{(d+1):p}X_{(d+1):p}^\T A^{-1})}_{\V_{\out,2}}.
\end{align*}

The following algebraic bounds are the foundation of the upper bounds and lower bounds for the out-sample error. The algebraic upper bounds below are mainly inspired by Lemma 27 and 28 in \cite{Tsigler_ridge_2023}. Moreover, we provide new algebraic lower bounds for the out-sample error.

\setcounter{lem}{0}
\begin{lem}[Algebraic upper bounds of out-sample error] Given invertible $\Sigma_{1:d}$, we have
\label{lem1}
\begin{align}
    \B_{\out}&\leq  2\Vert\theta^*_{1:d}\Vert^2_{\Sigma_{1:d}^{-1}} (\frac{1}{\lambda_1}+\frac{\mu_n(H_dH_d^T)}{\mu_1(A_d)})^{-2}+ \frac{2\mu_1(H_dH_d^T)}{\mu_n^2(H_dH_d^T)} \frac{\mu_1^2(A_d)}{\mu_n^2(A_d)} \Vert X_{(d+1):p}\theta^*_{(d+1):p}\Vert^2\notag \\
    &+3(\Vert M_d\Vert \frac{\mu_1(H_dH_d^T)}{\mu^2_n(A_d)} \Vert\theta^*_{1:d}\Vert_{\Sigma^{-1}}^2 (\frac{1}{\lambda_1}+\frac{\mu_n(H_dH_d^T)}{\mu_1(A_d)})^{-2}  +\Vert\theta^*_{(d+1):p}\Vert^2_{\Sigma_{(d+1):p}} \notag\\
    &+ \frac{\Vert M_d\Vert}{\mu_n(A)^2}\Vert X_{(d+1):p}\theta^*_{(d+1):p}\Vert^2 ), \notag   \\
    \V_{\out} &\leq   \sigma^2 \frac{\mu_1^2(A_d)}{\mu_n^2(A_d)}\frac{\Tr(H_d^\T H_d)}{\mu_d(H_dH_d^\T )^2} +  \sigma^2 \frac{\Tr(M_d)}{\mu_n(A_d)^2} . \notag
\end{align}
\end{lem}

\begin{prf}
\newline
\textbf{Algebraic upper bound for the out-sample bias.} From Section  H.2 in Supplement of \cite{Tsigler_ridge_2023}, we have
\begin{align}
\B_{\out,2}&\leq 3(\Vert X_{(d+1):p}^TA^{-1}X_{1:d}\theta_{1:d}^*\Vert^2_{\Sigma_{(d+1):p}}
+\Vert X_{(d+1):p} A^{-1}X_{(d+1):p}\theta^*_{(d+1):p}\Vert_{\Sigma_{(d+1):p}}^2 +\Vert \theta^*_{(d+1):p}\Vert_{\Sigma_{(d+1):p} }^2 )  \notag \\
&\leq 3(\Vert X_{(d+1):p}^TA^{-1}X_{1:d}\theta_{1:d}^*\Vert^2_{\Sigma_{(d+1):p}}
+\frac{\Vert M_d\Vert}{\mu_n(A)^2}\Vert X_{(d+1):p}\theta^*_{(d+1):p}\Vert^2 +\Vert \theta^*_{(d+1):p}\Vert_{\Sigma_{(d+1):p} }^2 ). \notag
\end{align}
From Lemma \ref{lem_01} (ii), we have
\begin{align*}
    &\Vert X_{(d+1):p}^TA^{-1}X_{1:d}\theta_{1:d}^*\Vert^2_{\Sigma_{(d+1):p}}  \\
    &=  \Vert X_{(d+1):p}^TA_d^{-1}X_{1:d}(I_d+X_{1:d}^TA_d^{-1}X_{1:d})^{-1}\theta_{1:d}^*\Vert^2_{\Sigma_{(d+1):p}} \\
    &=  \Vert X_{(d+1):p}^TA_d^{-1}X_{1:d}\Sigma_{1:d}^{-1/2}(\Sigma_{1:d}^{-1}+\Sigma_{1:d}^{-1/2}X_{1:d}^TA_d^{-1}X_{1:d}\Sigma_{1:d}^{-1/2})^{-1}\Sigma_{1:d}^{-1/2}\theta_{1:d}^*\Vert^2_{\Sigma_{(d+1):p}} \\
    &\leq  \Vert M_d\Vert \frac{\mu_1(H_dH_d^T)}{\mu^2_n(A_d)} \Vert\theta^*_{1:d}\Vert_{\Sigma^{-1}}^2 (\frac{1}{\lambda_1}+\frac{\mu_n(H_dH_d^T)}{\mu_1(A_d)})^{-2}.
\end{align*}
Hence we have
\begin{align}
  \B_{\out,2} &\leq  3(\Vert M_d\Vert \frac{\mu_1(H_dH_d^T)}{\mu^2_n(A_d)} \Vert\theta^*_{1:d}\Vert_{\Sigma^{-1}}^2 (\frac{1}{\lambda_1}+\frac{\mu_n(H_dH_d^T)}{\mu_1(A_d)})^{-2} \notag\\
&+\frac{\Vert M_d\Vert}{\mu_n(A)^2}\Vert X_{(d+1):p}\theta^*_{(d+1):p}\Vert^2 +\Vert \theta^*_{(d+1):p}\Vert_{(d+1):p}^2 ) .  \label{lem_prf_04}
\end{align}

It remains to give an upper bound of $\B_{\out,1}$. From Lemma \ref{lem_01} (i), we have
\begin{align}
    &\hat{\theta}(\tau, X\theta^*)_{1:d} + X_{1:d}^\T A_d^{-1}X_{1:d}\hat{\theta}(\tau, X\theta^*)_{1:d} = X_{1:d}^\T A_d^{-1}X\theta^* .\notag
\end{align}
Denote $\zeta_{1:d}=\hat{\theta}(\tau, X\theta^*)_{1:d}-\theta_{1:d}^*$. Then
\begin{align}
    H_d A_d^{-1}X_{(d+1):p}\theta^*_{(d+1):p}-\Sigma_{1:d}^{-1/2}\theta_{1:d}^*&=\Sigma_{1:d}^{-1/2}\zeta_{1:d}+\Sigma_{1:d}^{-1/2}X_{1:d}^\T A_d^{-1}X_{1:d}\zeta_{1:d} \notag \\
    &=(\Sigma_{1:d}^{-1}+ H_dA_d^{-1} H_d^\T )\Sigma_{1:d}^{1/2}\zeta_{1:d} .\label{lem_prf_05}
\end{align}
By standard manipulations,
\begin{align}
\Vert(\Sigma_{1:d}^{-1}+ H_dA_d^{-1} H_d^\T )\Sigma_{1:d}^{1/2}\zeta_{1:d} \Vert^2&\geq  \mu_n(\Sigma_{1:d}^{-1}+ H_dA_d^{-1} H_d^\T )^2\Vert\Sigma_{1:d}^{1/2}\zeta_{1:d}\Vert^2 \notag\\
&\geq  \B_{\out,1}(\frac{1}{\lambda_1}+\frac{\mu_n(H_dH_d^T)}{\mu_1(A_d)})^2,\label{lem1_1_ea}\\
\Vert H_dA_d^{-1}X_{(d+1):p}\theta^*_{(d+1):p}- \Sigma_{1:d}^{-1/2}\theta_{1:d}^*\Vert^2&\leq (\Vert\theta^*_{1:d}\Vert_{\Sigma_{1:d}^{-1}} + \Vert H_dA_d^{-1}X_{(d+1):p}\theta^*_{(d+1):p}\Vert)^2 .\label{lem1_2_ea}
\end{align}
From (\ref{lem_prf_05}), (\ref{lem1_1_ea}) and (\ref{lem1_2_ea}), we have
\begin{align}
&\B_{\out,1}(\frac{1}{\lambda_1}+\frac{\mu_n(H_dH_d^T)}{\mu_1(A_d)})^2  \leq  (\Vert\theta^*_{1:d}\Vert_{\Sigma_{1:d}^{-1}} + \Vert H_dA_d^{-1}X_{(d+1):p}\theta^*_{(d+1):p}\Vert)^2  .\notag
\end{align}
That is, we have
\begin{align}
   \B_{\out,1}&\leq  (\Vert\theta^*_{1:d}\Vert_{\Sigma_{1:d}^{-1}} + \Vert H_dA_d^{-1}X_{(d+1):p}\theta^*_{(d+1):p}\Vert)^2(\frac{1}{\lambda_1}+\frac{\mu_n^2(H_dH_d^T)}{\mu_1^2(A_d)})^{-2} \notag  \\
 &\leq 2\Vert\theta^*_{1:d}\Vert^2_{\Sigma_{1:d}^{-1}} (\frac{1}{\lambda_1}+\frac{\mu_n^2(H_dH_d^T)}{\mu_1^2(A_d)})^{-2}+ 2\Vert H_dA_d^{-1}X_{(d+1):p}\theta^*_{(d+1):p}\Vert^2\frac{\mu^2_1(A_d)}{\mu^2_n(H_dH_d^\T)}  \notag \\
 &\leq 2\Vert\theta^*_{1:d}\Vert^2_{\Sigma_{1:d}^{-1}} (\frac{1}{\lambda_1}+\frac{\mu_n^2(H_dH_d^T)}{\mu_1^2(A_d)})^{-2}+\frac{2\mu_1(H_dH_d^T)}{\mu_n^2(H_dH_d^T)} \frac{\mu_1^2(A_d)}{\mu_n^2(A_d)} \Vert X_{(d+1):p}\theta^*_{(d+1):p}\Vert^2 .\label{lem_prf_06b}
\end{align}

Combining (\ref{lem_prf_04}) and (\ref{lem_prf_06b}) gives the upper bound of out-sample bias.

\textbf{Algebraic upper bound for the out-sample variance.} From Lemma 27 in \cite{Tsigler_ridge_2023}, we have
\begin{align}
    \V_{\out,1} \leq \sigma^2 \frac{\mu_1^2(A_d)}{\mu_n^2(A_d)}\frac{\Tr(H_d^\T H_d)}{\mu_d(H_dH_d^\T )^2} . \notag
\end{align}
To upper bound $\V_{\out,2}$, we have
\begin{align}
    \V_{\out,2}&= \sigma^2\Tr( A^{-1}X_{(d+1):p}\Sigma_{(d+1):p}X_{(d+1):p}^\T A^{-1})   \notag \\
    &\leq \sigma^2 \frac{\Tr(M_d)}{\mu_n(A)^2}  \notag \\
    &\leq  \sigma^2 \frac{\Tr(M_d)}{\mu_n(A_d)^2}.  \notag
\end{align}
Combining the preceding two displays gives the upper bound of out-sample variance.
\end{prf}

\begin{lem}[Algebraic lower bounds of out-sample error] Given invertible $\Sigma_{1:d}$, we have
\label{lem2}
\begin{align*}
    \V_{\out} &\geq   \sigma^2\frac{d\mu_d(H_dH_d^\T )}{\mu_1(A_d)^2}(\frac{1}{\lambda_d}+\frac{\mu_1(H_dH_d^\T )}{\mu_n(A_d)})^{-2} +  \max\{0,\sigma^2\frac{\Tr(M_d)-d\mu_1(M_d)}{\mu_1(A_d)^2} \}.
\end{align*}
Further if $\Vert\theta^*_{1:d}\Vert_{\Sigma_{1:d}^{-1}}\geq \frac{\mu_1(H_dH_d)^{1/2}\Vert X_{(d+1):p}\theta^*_{(d+1):p}\Vert}{\mu_n(A_d)}$, then
\begin{align}
     \B_{\out}&\geq  (1 - \frac{\mu_1(H_dH_d^T)^{1/2}\Vert X_{(d+1):p}\theta^*_{(d+1):p}\Vert}{\mu_n(A_d)\Vert\theta^*_{1:d}\Vert_{\Sigma_{1:d}^{-1}}})^2\Vert\theta^*_{1:d}\Vert^2_{\Sigma_{1:d}^{-1}}(\frac{1}{\lambda_d}+\frac{\mu_1(H_dH_d^\T )}{\mu_n(A_d)})^{-2} .\label{lem_prf_08a}
\end{align}
\end{lem}
\begin{prf}
$\newline$
\textbf{Algebraic lower bound for the out-sample variance:} From Lemma \ref{lem_01} (i), we have
\begin{align}
    X_{1:d}^\T  A_{d}^{-1}\epsilon &= \hat{\theta}_{1:d}(\tau, \epsilon) + X_{1:d}^\T  A_{d}^{-1}X_{1:d}\hat{\theta}_{1:d}(\tau, \epsilon) . \notag
\end{align}
Multiplying the two sides with $\Sigma^{-1/2}$, we have
\begin{align}
    H_dA_d^{-1}\epsilon = (\Sigma_{1:d}^{-1} + H_dA_d^{-1}\Sigma_{1:d}^{-1/2})\Sigma_{1:d}^{1/2}\hat{\theta}(\tau,\epsilon) ,\notag
\end{align}
and hence
\begin{align*}
    \Vert H_dA_d^{-1}\epsilon\Vert^2 &= \Vert(\Sigma_{1:d}^{-1} +  H_dA_d^{-1}\Sigma_{1:d}^{-1/2})\Sigma_{1:d}^{1/2}\hat{\theta}(\tau,\epsilon)\Vert^2\\
    &\leq  \mu_1(\Sigma_{1:d}^{-1} +  H_dA_d^{-1}\Sigma_{1:d}^{-1/2})^2\Vert\Sigma_{1:d}^{1/2}\hat{\theta}(\tau,\epsilon)\Vert^2.
\end{align*}
Taking the expectations of the two sides with respect to $\epsilon$ yields
\begin{align*}
    \mu_1(\Sigma_{1:d}^{-1} +  H_dA_d^{-1}\Sigma_{1:d}^{-1/2})^2\V_{\out,1}&\geq E_\epsilon[\Vert H_dA_d^{-1}\epsilon\Vert^2] \\
    &= \sigma^2\Tr(A_d^{-1}H_dH_d^\T A_d^{-1}).
\end{align*}
By simple manipulations, we have
 \begin{align}
 (\frac{1}{\lambda_d}+\frac{\mu_1( H_d H_d^\T )}{\mu_n(A_d)})^2\V_{\out,1}&\geq \sigma^2\Tr(H_d A_d^{-2}H_d^\T)  \notag \\
  &\geq \sigma^2\frac{\Tr(H_dH_d^\T)}{\mu_1(A_d)^2}  \notag\\
  &\geq \sigma^2\frac{d\mu_d( H_d H_d^\T )}{\mu_1(A_d)^2} ,  \notag
\end{align}
and hence
\begin{align}
    \V_{\out,1} \geq \sigma^2\frac{d\mu_d( H_d H_d^\T )}{\mu_1(A_d)^2}(\frac{1}{\lambda_d}+\frac{\mu_1( H_d H_d^\T )}{\mu_n(A_d)})^{-2}  .\label{lem_prf_08}
\end{align}

To lower bound $\V_{\out,2}$, we have
\begin{align}
\V_{\out,2}&=\sigma^2\Tr(A^{-1}M_d A^{-1})\notag \\
    &\geq   \sum_{i=1}^{n-d}\mu_i(A^{-2})\mu_{n-i+1}(M_d )  \notag \\
    &\geq  \frac{\sum_{i=1}^{n-d}\mu_{n-i+1}(M_d )}{\mu_1(A_d)^2} \notag\\
    &\geq \frac{\Tr(M_d )-d\mu_1(M_d)}{\mu_1(A_d)^2}  .\notag
\end{align}
The first inequality above is from Lemma \ref{lem:ruhe_inq}. For the second inequality,
because the rank of $A - A_d$ is at most $d$,
\begin{align*}
    \mu_i(A-A_d) = 0 ~~~\forall i\geq d+1.
\end{align*}
From Weyl's inequality (Lemma \ref{lem_S12b}), we have
\begin{align*}
   & \mu_i(A)-\mu_1(A_d) \leq    \mu_i(A-A_d)  = 0 \quad \forall i\geq d+1\\
   \Longrightarrow &\mu_i(A) \leq \mu_1(A_d)\quad \forall i\geq d+1\\
   \Longrightarrow &\frac{1}{\mu_1(A_d)} \leq \frac{1}{\mu_i(A)} \quad \forall i\geq d+1\\
   \Longrightarrow &\frac{1}{\mu_1(A_d)}\leq \mu_i(A^{-1})   \quad\forall i\leq n-d.
\end{align*}
By requiring $\V_{\out,2}\geq 0$, we have
\begin{align}
    \V_{\out,2}&\geq \max\{0,\frac{\Tr(M_d )-d\mu_1(M_d)}{\mu_1(A_d)^2}\}.\label{lem_prf_09}
\end{align}

Combining (\ref{lem_prf_08}) and (\ref{lem_prf_09}) gives the lower bound of out-sample variance.

\noindent \textbf{Algebraic lower bound for the out-sample bias.} The norms of the two sides of (\ref{lem_prf_05}) can be bounded as follows:
\begin{align}
    \Vert(\Sigma_{1:d}^{-1}+ H_dA_d^{-1} H_d^\T )\Sigma_{1:d}^{1/2}\zeta_{1:d} \Vert^2&\leq  \mu_1(\Sigma_{1:d}^{-1}+ H_dA_d^{-1} H_d^\T )^2\Vert\Sigma_{1:d}^{1/2}\zeta_{1:d}\Vert^2 \notag  \\
&\leq  \B_{\out,1}(\frac{1}{\lambda_d}+\frac{\mu_1( H_d H_d^\T )}{\mu_n(A_d)})^2, \label{lem1_1_e} \\
\Vert H_dA_d^{-1}X_{(d+1):p}\theta^*_{(d+1):p}- \Sigma_{1:d}^{-1/2}\theta_{1:d}^*\Vert^2&\geq (\Vert\theta^*_{1:d}\Vert_{\Sigma_{1:d}^{-1}} - \Vert H_dA_d^{-1}X_{(d+1):p}\theta^*_{(d+1):p}\Vert)^2 .\label{lem1_2_e}
\end{align}
From (\ref{lem_prf_05}), (\ref{lem1_1_e}) and (\ref{lem1_2_e}), we have
\begin{align}
&\B_{\out,1}(\frac{1}{\lambda_d}+\frac{\mu_1( H_d H_d^\T )}{\mu_n(A_d)})^2 \geq  (\Vert\theta^*_{1:d}\Vert_{\Sigma_{1:d}^{-1}} - \Vert H_dA_d^{-1}X_{(d+1):p}\theta^*_{(d+1):p}\Vert)^2, \notag
\end{align}
and hence
\begin{align*}
    \B_\out&\geq \B_{\out,1}  \\
    &\geq  (\Vert\theta^*_{1:d}\Vert_{\Sigma_{1:d}^{-1}} - \Vert H_dA_d^{-1}X_{(d+1):p}\theta^*_{(d+1):p}\Vert)^2(\frac{1}{\lambda_d}+\frac{\mu_1( H_d H_d^\T )}{\mu_n(A_d)})^{-2}.
\end{align*}
With $\Vert\theta^*_{1:d}\Vert_{\Sigma_{1:d}^{-1}}\geq \frac{\mu_1(H_dH_d)^{1/2}\Vert X_{(d+1):p}\theta^*_{(d+1):p}\Vert}{\mu_n(A_d)}$, we have
\begin{align}
    \Vert\theta^*_{1:d}\Vert_{\Sigma_{1:d}^{-1}}\geq \frac{\mu_1(H_dH_d)^{1/2}\Vert X_{(d+1):p}\theta^*_{(d+1):p}\Vert}{\mu_n(A_d)}
    \geq \Vert H_dA_d^{-1}X_{(d+1):p}\theta^*_{(d+1):p}\Vert , \notag
\end{align}
and hence
\begin{align}
    \B_\out&\geq (\Vert\theta^*_{1:d}\Vert_{\Sigma_{1:d}^{-1}} - \Vert H_dA_d^{-1}X_{(d+1):p}\theta^*_{(d+1):p}\Vert)^2(\frac{1}{\lambda_d}+\frac{\mu_1( H_d H_d^\T )}{\mu_n(A_d)})^{-2} \notag \\
   &\geq (1 - \frac{\mu_1(H_dH_d^T)^{1/2}\Vert X_{(d+1):p}\theta^*_{(d+1):p}\Vert}{\mu_n(A_d)\Vert\theta^*_{1:d}\Vert_{\Sigma_{1:d}^{-1}}})^2\Vert\theta^*_{1:d}\Vert^2_{\Sigma_{1:d}^{-1}}(\frac{1}{\lambda_d}+\frac{\mu_1(H_dH_d^\T )}{\mu_n(A_d)})^{-2} . \notag 
\end{align}
\end{prf}

\subsubsection{Intermediate bounds of the out-sample error}
\label{supp_subsec1}
We give the intermediate bounds of out-sample error under the event that some random quantities in the algebraic bounds above are controlled. In the event $\Omega_1(\nu)\cap \Omega_2 \cap \Omega_4$ for $0<\nu <\frac{1}{2}$ defined in Lemma \ref{lem3}, substituting the bounds of $\mu_1(H_dH_d^T)$, $\mu_d(H_dH_d^T)$, $\Vert X_{(d+1):p}\theta^*_{(d+1):p}\Vert$ and $\Vert M_d\Vert$ into the algebraic upper bound of $\B_\out$ yields
\begin{align}
    \B_{\out} &\leq \frac{2}{(1-\nu-\eta_1)^4}
\Vert\theta^*_{1:d}\Vert_{\Sigma_{1:d}^{-1}}^2
(\frac{1}{\lambda_1} + \frac{n}{\mu_1(A_d)})^{-2}
    +\frac{\poly_2(\sigma_x)(1+\nu+\eta_1)^2}{(1-\nu-\eta_1)^4} \frac{\mu_1^2(A_d)}{\mu_n^2(A_d)} \Vert\theta^*_{(d+1):p}\Vert_{\Sigma_{(d+1):p}}^2 \notag \\
    &+\frac{\poly_2(\sigma_x)(1+\nu+\eta_1)^2}{(1-\nu-\eta_1)^4} \frac{(n^2\lambda_{d+1}^2+n\sum_{j>d}\lambda_j^2)}{\mu^2_n(A_d)} \Vert\theta^*_{1:d}\Vert_{\Sigma^{-1}}^2 (\frac{1}{\lambda_1} + \frac{n}{\mu_1(A_d)})^{-2} \notag\\
   &+(3+  \poly_4(\sigma_x) \frac{(n^2\lambda_{d+1}^2+n\sum_{j>d}\lambda_j^2)}{\mu_n(A_d)^2})\Vert\theta^*_{(d+1):p}\Vert_{\Sigma_{(d+1):p}}^2  .\label{lem_prf_10}
\end{align}
In the event $\Omega_1(\nu)\cap \Omega_4$ for $0<\nu <\frac{1}{2}$ defined in Lemma \ref{lem3} and with
\begin{align*}
\Vert\theta^*_{1:d}\Vert_{\Sigma_{1:d}^{-1}}\geq \frac{(1+\nu+\eta_1)(1+\sigma_x^2)^{1/2}n\Vert \theta^*_{(d+1):p}\Vert_{\Sigma_{(d+1):p}}}{\mu_n(A_d)},
\end{align*}
substituting the bounds of $\mu_1(H_dH_d^T)$ and $\Vert X_{(d+1):p}\theta^*_{(d+1):p}\Vert$ into the algebraic lower bound of $\B_\out$ yields
\begin{align}
    \B_{\out} &\geq \frac{1}{(1+\nu+\eta_1)^2}(1 - \frac{(1+\nu+\eta_1)(1+\sigma_x^2)^{1/2}n\Vert\theta^*_{(d+1):p}\Vert_{\Sigma_{(d+1):p}}}{\mu_n(A_d)\Vert\theta^*_{1:d}\Vert_{\Sigma_{1:d}^{-1}}})^2\Vert\theta^*_{1:d}\Vert_{\Sigma_{1:d}^{-1}}^2(\frac{1}{\lambda_d}+\frac{n}{\mu_n(A_d)})^{-2} .\label{lem_prf_12}
\end{align}
In the event $\Omega_1(\nu)\cap \Omega_{31}\cap \Omega_{32}$ for $\nu <\frac{1}{2}$ defined in Lemma \ref{lem3}, substituting the bounds of $\mu_d(H_dH_d^T)$, $\Tr(H_d^TH_d)$ and $\Tr(M_d)$ into the algebraic upper bound of $\V_\out$ yields
\begin{align}
   \V_{\out} &\leq \frac{\poly_2(\sigma_x)}{(1-\nu-\eta_1)^4}\frac{\mu_1^2(A_d)}{\mu_n^2(A_d)}\sigma^2\frac{d}{n} + \poly_2(\sigma_x)\sigma^2 \frac{n\sum_{j>d}\lambda_j^2}{\mu_n(A_d)^2} .\label{lem_prf_11}
\end{align}
In the event $\Omega_1(\nu)\cap \Omega_2\cap \Omega_{33}(\nu)$ for $\nu<\frac{1}{2}\min\{1,\sigma_x^2\}$ defined in Lemma \ref{lem3}, substituting the bounds of $\mu_1(H_dH_d^T))$, $\mu_d(H_dH_d^T)$, $\Tr(M_d)$ and $\mu_1(M_d)$ into the algebraic lower bound of $\V_\out$ yields
\begin{align}
    \V_{\out} &\geq  \frac{(1-\nu-\eta_1)^2}{(1+\nu+\eta_1)^4} \sigma^2\frac{nd}{\mu_1(A_d)^2}(\frac{1}{\lambda_d}+\frac{n}{\mu_n(A_d)})^{-2} + \max\{0,\frac{n\sum_{j>d}\lambda_j^2}{\mu_1(A_d)^2}(1-\nu-C_0\sigma_x^2(\frac{2d}{r_d(\Sigma^2)}+\frac{d}{n})\}.\label{lem_prf_13}
\end{align}

\subsubsection{Final upper bounds of out-sample bias and variance}
We give the final upper bounds of out-sample bias and variance as stated in Proposition \ref{pro1} and \ref{pro7}.
From the intermediate bounds in Section~\ref{supp_subsec1}, we derive the final bounds by further controlling $\mu_1(A_d)$ and $\mu_n(A_d)$ (from Lemma \ref{lem3a}) in the intermediate bounds and incorporating Assumption \ref{ass1} to control the terms related to $\Vert\theta^*_{(d+1):p}\Vert^2_{\Sigma_{(d+1):p}}$. We first discuss the small or moderate TER regime and then the large TER regime.

\noindent \textbf{(i) Small or moderate TER}

From Lemma \ref{lem3a}(i), (\ref{lem_prf_10}) and (\ref{lem_prf_11}), in the event $\Omega_1(\nu)\cap \Omega_2 \cap \Omega_{31}\cap \Omega_{32}\cap\Omega_4\cap \Omega_5$ for $\nu<\frac{1}{2}$ defined in Lemma \ref{lem3} and Assumption \ref{ass:3a}, substituting the bounds of $\mu_1(A_d)$ and $\mu_n(A_d)$ in (\ref{lem:bq_11_00})--(\ref{lem:bq_11_0}) into (\ref{lem_prf_10}) and (\ref{lem_prf_11}), we have for $\tau\geq \lambda_{d+1}$,
\begin{align}
   \B_{\out} &\leq \frac{\poly_4(\sigma_x)(1+C_1)^2}{(1-\nu-\eta_1)^4}\Vert\theta^*_{1:d}\Vert_{\Sigma_{1:d}^{-1}}^2(\frac{1}{\lambda_1}+\frac{1}{\tau})^{-2} +\frac{\poly_6(\sigma_x)(1+\nu+\eta_1)^2}{(1-\nu-\eta_1)^4} (1+C_1)^2 \Vert\theta^*_{(d+1):p}\Vert_{\Sigma_{(d+1):p}}^2 \notag \\
    &+\frac{\poly_6(\sigma_x)(1+C_1)^3(1+\nu+\eta_1)^2}{(1-\nu-\eta_1)^4}  \Vert\theta^*_{1:d}\Vert_{\Sigma_{1:d}^{-1}}^2 (\frac{1}{\lambda_1}+\frac{1}{\tau})^{-2} + \poly_4(\sigma_x) (1+C_1) \Vert\theta^*_{(d+1):p}\Vert_{\Sigma_{(d+1):p}}^2 ,\notag  \\
   \V_{\out} &\leq \frac{\poly_6(\sigma_x)(1+C_1)^2}{(1-\nu-\eta_1)^4}(\sigma^2\frac{d}{n} + \sigma^2 \frac{\sum_{j>d}\lambda_j^2}{n\tau^2}) .\notag
\end{align}
Hence, we obtain the upper bounds
\begin{align}
  \B_{\out} &\leq \frac{(1+C_1)^3(1+\nu+\eta_1)^2\poly_6(\sigma_x)}{(1-\nu-\eta_1)^4}(\Vert\theta^*_{1:d}\Vert^2_{\Sigma_{1:d}^{-1}}(\frac{1}{\lambda_1}+\frac{1}{\tau})^{-2} + \Vert \theta_{(d+1):p}^*\Vert^2_{\Sigma_{(d+1):p}}) , \notag \\
   \V_{\out} &\leq \frac{\poly_6(\sigma_x)(1+C_1)^2}{(1-\nu-\eta_1)^4}(\sigma^2\frac{d}{n} + \sigma^2 \frac{\lambda_{d+1}^2}{\tau^2}\frac{r_d( \Sigma^2)}{n} ) .\notag
\end{align}
With Assumption \ref{ass1}(i), we have for $\tau\geq \lambda_{d+1}$,
\begin{align}
\Vert\theta^*_{(d+1):p}\Vert_{\Sigma_{(d+1):p}}^2&\leq \frac{\delta_1}{4}\Vert\theta^*_{1:d}\Vert_{\Sigma_{1:d}^{-1}}\lambda_{d+1}^2\leq \delta_1\Vert\theta^*_{1:d}\Vert_{\Sigma_{1:d}^{-1}}(\frac{1}{\lambda_1}+\frac{1}{\tau})^{-2} .  \label{lem_prf_14a}
\end{align}
Further with (\ref{lem_prf_14a}), we have for $\tau\geq \lambda_{d+1}$,
\begin{align*}
    \B_{\out} &\leq \frac{(1+C_1)^3(1+\nu+\eta_1)^2\poly_6(\sigma_x)}{(1-\nu-\eta_1)^4}\Vert\theta^*_{1:d}\Vert^2_{\Sigma_{1:d}^{-1}}(\frac{1}{\lambda_1}+\frac{1}{\tau})^{-2}.
\end{align*}
This gives Proposition \ref{pro1}.

\noindent \textbf{(ii) Large TER}

From Lemma \ref{lem3a}(ii), (\ref{lem_prf_10}) and (\ref{lem_prf_11}), in the event $\Omega_1(\nu)\cap \Omega_2 \cap \Omega_{31}\cap \Omega_{32}\cap \Omega_4\cap \Omega_6(\nu)$ for $0<\nu<\frac{1}{2}$ defined in Lemma \ref{lem3} and Assumption \ref{ass:4}, substituting the bounds of $\mu_1(A_d)$ and $\mu_n(A_d)$ in (\ref{lem:bq_11_01aa})--(\ref{lem:bq_11_01a}) into (\ref{lem_prf_10}) and (\ref{lem_prf_11}), for $\tau\geq 0$,
\begin{align}
   \B_{\out} &\leq \frac{(1+\nu+\eta_2)^2}{(1-\nu-\eta_1)^4}\Vert\theta^*_{1:d}\Vert_{\Sigma_{1:d}^{-1}}^2 (\frac{1}{\lambda_1}+\frac{1}{\tau+\lambda_{d+1}\frac{r_d(\Sigma)}{n}})^{-2} \notag\\
   &+\frac{\poly_2(\sigma_x)(1+\nu+\eta_1)^2}{(1-\nu-\eta_1)^4} \frac{(1+\nu+\eta_2)^2}{(1-\nu-\eta_2)^2} \Vert\theta^*_{(d+1):p}\Vert_{\Sigma_{(d+1):p}}^2  \notag \\
    &+\frac{\poly_2(\sigma_x)(1+\nu+\eta_1)^2(1+\nu+\eta_2)^2}{(1-\nu-\eta_1)^4(1-\nu-\eta_2)^2} \Vert\theta^*_{1:d}\Vert_{\Sigma^{-1}}^2 (\frac{1}{\lambda_1}+\frac{1}{\tau+\lambda_{d+1}\frac{r_d(\Sigma)}{n}})^{-2} \notag\\
   &+\frac{\poly_4(\sigma_x)}{(1-\nu-\eta_2)^2} \Vert\theta^*_{(d+1):p}\Vert_{\Sigma_{(d+1):p}}^2, \notag \\
   \V_{\out} &\leq \frac{\poly_2(\sigma_x)(1+\nu+\eta_2)^2}{(1-\nu-\eta_1)^4(1-\nu-\eta_2)^2}\sigma^2(\frac{d}{n} +   \frac{\lambda_{d+1}^2}{(\tau+\lambda_{d+1}\frac{r_d(\Sigma)}{n})^2}\frac{r_d(\Sigma^2)}{n} )  .\notag
\end{align}
Hence, we obtain the upper bounds
\begin{align}
     \B_{\out}&\leq  \frac{(1+\nu +\eta_1)^2(1+\nu + \eta_2)^2\poly_4(\sigma_x)}{(1-\nu-\eta_1)^4(1-\nu-\eta_2)^2}(\Vert\theta^*_{1:d}\Vert^2_{\Sigma_{1:d}^{-1}}(\frac{1}{\lambda_1}+\frac{1}{\tau+\frac{\sum_{j>d}\lambda_j}{n}})^{-2} +\Vert\theta^*_{(d+1):p}\Vert_{\Sigma_{(d+1):p}}^2)\notag, \\
    \V_{\out} &\leq  \frac{(1+\nu+\eta_2)^2\poly_2(\sigma_x)}{(1-\nu-\eta_1)^4(1-\nu-\eta_2)^2}\sigma^2(\frac{d}{n}+\frac{\lambda_{d+1}^2}{(\tau+\lambda_{d+1}\frac{r_d(\Sigma)}{n})^2}\frac{r_d(\Sigma^2)}{n}) .  \notag
\end{align}
With Assumption \ref{ass1}(ii), we have for $\tau\geq 0$,
\begin{align}
\Vert\theta^*_{(d+1):p}\Vert_{\Sigma_{(d+1):p}}^2&\leq \frac{\delta_2}{4}\Vert\theta^*_{1:d}\Vert_{\Sigma_{1:d}^{-1}}(\frac{1}{\lambda_d}+\frac{1}{\lambda_{d+1}\frac{r_d(\Sigma)}{n}})^{-2}  \leq \frac{\delta_2}{4}\Vert\theta^*_{1:d}\Vert_{\Sigma_{1:d}^{-1}}(\frac{1}{\lambda_1}+\frac{1}{\tau+\lambda_{d+1}\frac{r_d(\Sigma)}{n}})^{-2} .\label{lem_prf_19a}
\end{align}
Further with (\ref{lem_prf_19a}), we have for $\tau\geq 0$,
\begin{align}
    \B_{\out}&\leq   \frac{(1+\nu +\eta_1)^2(1+\nu + \eta_2)^2\poly_4(\sigma_x)}{(1-\nu-\eta_1)^4(1-\nu-\eta_2)^2}\Vert\theta^*_{1:d}\Vert^2_{\Sigma_{1:d}^{-1}}(\frac{1}{\lambda_1}+\frac{1}{\tau+\lambda_{d+1}\frac{r_d(\Sigma)}{n}})^{-2} .\notag
\end{align}
This gives Proposition \ref{pro7}.

\subsubsection{Final lower bounds of out-sample bias and variance}
\label{sec_I_1_4}
We give the final lower bounds of out-sample bias and variance. We first discuss the small or moderate TER regime and then the large TER regime.

\noindent \textbf{(i) Small or moderate TER}

\noindent \textbf{Lower bound of out-sample bias.} For $\tau\geq \lambda_{d+1}$,  we have $\mu_n(A_d)\geq n\tau\geq n\lambda_{d+1}$.
Then from Assumption \ref{ass1}(i), we have for $0<\nu<\frac{1}{2}$ and $\tau\geq \lambda_{d+1}$,
\begin{align}
\Vert\theta^*_{1:d}\Vert_{\Sigma_{1:d}^{-1}} &\geq \frac{2(1+\sigma_x^2)^{1/2}n\Vert \theta^*_{(d+1):p}\Vert_{\Sigma_{(d+1):p}}}{n\lambda_{d+1}} \notag \\
&\geq \frac{(1+\nu+\eta_1)(1+\sigma_x^2)^{1/2}n\Vert \theta^*_{(d+1):p}\Vert_{\Sigma_{(d+1):p}}}{\mu_n(A_d)}.  \label{lem_prf_15a}
\end{align}
From (\ref{lem_prf_12}), in the event $\Omega_1(\nu)\cap \Omega_4$ for $\nu<\frac{1}{2}$ defined in Lemma \ref{lem3} and Assumption \ref{ass1}(i),
substituting the lower bound of $\mu_n(A_d)$ in (\ref{lem:bq_11_0}) into (\ref{lem_prf_12}), we have for $\tau\geq \lambda_{d+1}$,
\begin{align}
    \B_{\out} &\geq \frac{1}{(1+\nu+\eta_1)^2}(1 - \frac{(1+\nu+\eta_1)(1+\sigma_x^2)^{1/2}\Vert\theta^*_{(d+1):p}\Vert_{\Sigma_{(d+1):p}}}{\tau\Vert\theta^*_{1:d}\Vert_{\Sigma_{1:d}^{-1}}})^2\Vert\theta^*_{1:d}\Vert_{\Sigma_{1:d}^{-1}}^2(\frac{1}{\lambda_d}+\frac{1}{\tau})^{-2} . \notag
\end{align}
With Assumption \ref{ass1}(i), we also have for $\tau\geq\lambda_{d+1}$,
\begin{align}
     \frac{(1+\nu+\eta_1)(1+\sigma_x^2)^{1/2}\Vert\theta^*_{(d+1):p}\Vert_{\Sigma_{(d+1):p}}}{\tau\Vert\theta^*_{1:d}\Vert_{\Sigma_{1:d}^{-1}}} &\leq  \frac{2(1+\sigma_x^2)^{1/2}\Vert\theta^*_{(d+1):p}\Vert_{\Sigma_{(d+1):p}}}{\lambda_{d+1}\Vert\theta^*_{1:d}\Vert_{\Sigma_{1:d}^{-1}}} \leq \sqrt{\delta_1}.  \label{lem_prf_14}
\end{align}
By applying (\ref{lem_prf_14}), in the event $\Omega_1(\nu)\cap \Omega_4$ for $\nu<\frac{1}{2}$ defined in Lemma \ref{lem3} and Assumption \ref{ass1}(i), we have for $\tau\geq \lambda_{d+1}$,
\begin{align}
    \B_{\out} &\geq \frac{(1 - \sqrt{\delta_1})^2}{(1+\nu+\eta_1)^2}\Vert\theta^*_{1:d}\Vert_{\Sigma_{1:d}^{-1}}^2(\frac{1}{\lambda_d}+\frac{1}{\tau})^{-2}  .\notag
\end{align}
This gives the results of Proposition \ref{pro2}.

\noindent \textbf{Lower bound of out-sample variance.}
To deduce the final lower bound of $\V_{\out}$, our strategy is as follows.
We first derive a lower bound for the first term in the intermediate lower bound (\ref{lem_prf_13}). Then we discuss two complementary cases.
The first case is that $\frac{r_d(\Sigma^2)}{d}$ is upper bounded.
The second case is that $\frac{r_d(\Sigma^2)}{d}$ is large enough such that $1-\nu-C_0\sigma_x^2(\frac{2d}{r_d(\Sigma^2)}+\frac{d}{n})>0$.  Lastly, we show
that in these two cases, $\V_\out$ satisfies lower bounds of the same order, which gives the final lower bound for $\V_\out$.

As preparation, we derive an equivalence relationship, which is useful in the following analysis. Given Assumption \ref{ass:3}, we have
\begin{align}
   \frac{2C_0\sigma_x^2d}{r_d(\Sigma^2)}+\frac{C_0\sigma_x^2d}{n} \geq \frac{1}{2}
   &\Longleftrightarrow \frac{2C_0\sigma_x^2d}{r_d(\Sigma^2)} \geq \frac{1}{2} - \frac{C_0\sigma_x^2d}{n} \notag \\
    &\Longleftrightarrow   \frac{2C_0\sigma_x^2d}{r_d(\Sigma^2)} \geq  \frac{1}{2} - \eta_1 \quad \text{(from~Assumption~\ref{ass:3})}  \notag \\
    &\Longleftrightarrow  \frac{2C_0\sigma_x^2d}{\frac{1}{2} - \eta_1} \geq \frac{\sum_{j>d}\lambda_j^2}{\lambda_{d+1}^2}.  \label{lem_prf_16}
\end{align}

Now, we are ready to give the lower bound of $\V_\out$. We first consider the first term of the right-hand side in (\ref{lem_prf_13}). From Lemma \ref{lem3a}(i) and (\ref{lem_prf_13}), in the event $\Omega_1(\nu)\cap \Omega_2\cap \Omega_{33}(\nu)\cap \Omega_5$ for $0<\nu<\frac{1}{2}\min\{1,\sigma_x^2\}$ defined in Lemma \ref{lem3}, substituting the bounds of $\mu_1(A_d)$ and $\mu_n(A_d)$ in  (\ref{lem:bq_11_00})--(\ref{lem:bq_11_0}) into the first term in (\ref{lem_prf_13}), we have for $\tau\geq \lambda_{d+1}$,
\begin{align}
    \frac{(1-\nu-\eta_1)^2}{(1+\nu+\eta_1)^4} \sigma^2\frac{nd}{\mu_1(A_d)^2}(\frac{1}{\lambda_d}+\frac{n}{\mu_n(A_d)})^{-2} \geq   \frac{(1-\nu-\eta_1)^2}{(1+C_1)^2(1+\nu+\eta_1)^4\poly_4(\sigma_x)} \sigma^2\frac{d}{n\tau^2}(\frac{1}{\lambda_d}+\frac{1}{\tau})^{-2} .\notag
\end{align}
Then in the event $\Omega_1(\nu)\cap \Omega_2\cap \Omega_{33}(\nu)\cap \Omega_5$ for $0<\nu<\frac{1}{2}\min\{1,\sigma_x^2\}$ defined in Lemma \ref{lem3} and Assumption \ref{ass:3a}, we have for $\lambda_{d+1}\leq \tau \leq \lambda_d$,
\begin{align}
    \frac{(1-\nu-\eta_1)^2}{(1+\nu+\eta_1)^4} \sigma^2\frac{nd}{\mu_1(A_d)^2}(\frac{1}{\lambda_d}+\frac{n}{\mu_n(A_d)})^{-2} \geq   \frac{(1-\nu-\eta_1)^2}{(1+C_1)^2(1+\nu+\eta_1)^4\poly_4(\sigma_x)} \sigma^2\frac{d}{n} .\label{lem_prf_17}
\end{align}
Then we discuss two complementary cases. The first cases is that $\frac{2C_0\sigma_x^2d}{r_d(\Sigma^2)}+\frac{C_0\sigma_x^2d}{n} \geq \frac{1}{2}$.
From the equivalence in (\ref{lem_prf_16}) and with Assumption \ref{ass:3}, for $\tau\geq \lambda_{d+1}$,
\begin{align*}
    \frac{2C_0\sigma_x^2d}{\frac{1}{2} - \eta_1} \geq \frac{\sum_{j>d}\lambda_j^2}{\lambda_{d+1}^2}\geq \frac{\sum_{j>d}\lambda_j^2}{\tau^2},
\end{align*}
which implies that
\begin{align*}
&    \frac{d}{n} \geq \frac{\frac{1}{2}-\eta_1}{2C_0\sigma_x^2} \frac{\sum_{j>d}\lambda_j^2}{n\tau^2}  \\
\Longleftrightarrow & \frac{d}{n} \geq \frac{\frac{\frac{1}{2}-\eta_1}{2C_0\sigma_x^2}}{1+\frac{\frac{1}{2}-\eta_1}{2C_0\sigma_x^2}}(\frac{d}{n}+\frac{\sum_{j>d}\lambda_j^2}{n\tau^2}).
\end{align*}
Then from (\ref{lem_prf_13}) and (\ref{lem_prf_17}), in the event $\Omega_1(\nu)\cap \Omega_2 \cap \Omega_{33}(\nu)\cap \Omega_5$ for $\nu<\frac{1}{2}\min\{1,\sigma_x^2\}$ defined in Lemma \ref{lem3},  we have for $\lambda_{d+1}\leq \tau \leq \lambda_d$,
\begin{align}
    \V_\out&\geq \frac{(1-\nu-\eta_1)^2}{(1+C_1)^2(1+\nu+\eta_1)^4\poly_4(\sigma_x)} \sigma^2\frac{d}{n} \notag \\
    &\geq  \frac{(1-\nu-\eta_1)^2}{(1+C_1)^2(1+\nu+\eta_1)^4\poly_4(\sigma_x)} \frac{\frac{\frac{1}{2}-\eta_1}{2C_0\sigma_x^2}}{1+\frac{\frac{1}{2}-\eta_1}{2C_0\sigma_x^2}}\sigma^2(\frac{d}{n}+\frac{\sum_{j>d}\lambda_j^2}{n\tau^2}) .\notag
\end{align}
The second case is that $\frac{2C_0\sigma_x^2d}{r_d(\Sigma^2)}+\frac{C_0\sigma_x^2d}{n} < \frac{1}{2}$. From Lemma \ref{lem3a}(i) and (\ref{lem_prf_13}), in the event $\Omega_1(\nu)\cap \Omega_2 \cap \Omega_{33}(\nu)\cap \Omega_5$ for $0<\nu<\frac{1}{2}\min\{1,\sigma_x^2\}$ defined in Lemma \ref{lem3} and Assumption \ref{ass:3a}, for $\tau\geq \lambda_{d+1}$, the second term in (\ref{lem_prf_13}) satisfies
\begin{align}
    \max\{0,\frac{n\sum_{j>d}\lambda_j^2}{\mu_1(A_d)^2}(1-\nu-C_0\sigma_x^2(\frac{2d}{r_d(\Sigma^2)}+\frac{d}{n})\} \geq  \frac{1}{(1+C_1)^2\poly_4(\sigma_x)}\frac{\sum_{j>d}\lambda_j^2}{n\tau^2} (\frac{1}{2}-\nu) .\label{lem_prf_18}
\end{align}
Then from (\ref{lem_prf_13}), (\ref{lem_prf_17}) and (\ref{lem_prf_18}), in the event $\Omega_1(\nu)\cap \Omega_2 \cap \Omega_{33}(\nu)\cap \Omega_5$ for $0<\nu<\frac{1}{2}\min\{1,\sigma_x^2\}$ defined in Lemma \ref{lem3} and Assumption \ref{ass:3a}, we have for $\tau\geq \lambda_{d+1}$,
\begin{align}
    \V_\out\geq \frac{(1-\nu-\eta_1)^2(\frac{1}{2}-\nu)}{(1+C_1)^2(1+\nu+\eta_1)^4\poly_4(\sigma_x)} \sigma^2(\frac{d}{n}+\frac{\sum_{j>d}\lambda_j^2}{n\tau^2}).  \notag
\end{align}
In conclusion, in the event $\Omega_1(\nu)\cap \Omega_2 \cap \Omega_{33}(\nu)\cap \Omega_5$ for $0<\nu<\frac{1}{2}\min\{1,\sigma_x^2\}$ defined in Lemma \ref{lem3} and Assumption \ref{ass:3a}, we have for $\lambda_{d+1}\leq \tau\leq \lambda_d$,
\begin{align}
   \V_\out\geq \frac{(1-\nu-\eta_1)^2(\frac{1}{2}-\nu)}{(1+C_1)^2(1+\nu+\eta_1)^4\poly_4(\sigma_x)} \frac{\frac{\frac{1}{2}-\eta_1}{2C_0\sigma_x^2}}{1+\frac{\frac{1}{2}-\eta_1}{2C_0\sigma_x^2}}\sigma^2(\frac{d}{n}+\frac{\sum_{j>d}\lambda_j^2}{n\tau^2}) . \notag
\end{align}
or equivalently
\begin{align}
     \V_\out\geq \frac{(1-\nu-\eta_1)^2(\frac{1}{2}-\nu)}{(1+C_1)^2(1+\nu+\eta_1)^4\poly_4(\sigma_x)} \frac{\frac{\frac{1}{2}-\eta_1}{2C_0\sigma_x^2}}{1+\frac{\frac{1}{2}-\eta_1}{2C_0\sigma_x^2}}\sigma^2(\frac{d}{n}+\frac{\lambda_{d+1}^2}{\tau^2}\frac{r_d(\Sigma^2)}{n}) . \notag
\end{align}
The result for $\tau \leq \lambda_{d+1}$ follows from the monotonicity of variance in Lemma \ref{lem_12b}. This gives Proposition \ref{pro3}.
$\newline\newline$
\textbf{(ii) Large TER}

\noindent \textbf{Lower bound of out-sample bias.} In the event $\Omega_6(\nu)$ for $0<\nu <\frac{1}{4}$ defined in Lemma \ref{lem3} and Assumption \ref{ass:4} and \ref{ass1}(ii), incorporating the bound of $\mu_n(A_d)$, we have for $\tau\geq 0$,
\begin{align}
\Vert\theta^*_{1:d}\Vert_{\Sigma_{1:d}^{-1}}&\geq \frac{2(1+\sigma_x^2)^{1/2}\Vert \theta^*_{(d+1):p}\Vert_{\Sigma_{(d+1):p}}}{\frac{1}{4}(\lambda_{d+1}\frac{r_d(\Sigma)}{n}+\tau) } \notag \\
&\geq \frac{(1+\nu+\eta_1)(1+\sigma_x^2)^{1/2}n\Vert \theta^*_{(d+1):p}\Vert_{\Sigma_{(d+1):p}}}{\mu_n(A_d)}.  \label{lem_prf_20a}
\end{align}
Then from (\ref{lem_prf_12}) and Lemma \ref{lem3a}(ii), in the event $\Omega_1(\nu)\cap \Omega_4\cap \Omega_6(\nu)$ for $0<\nu<\frac{1}{4}$ defined in Lemma \ref{lem3} and Assumption \ref{ass:4} and \ref{ass1}(ii), substituting the lower bounds of $\mu_n(A_d)$ in (\ref{lem:bq_11_01a}) into (\ref{lem_prf_12}), we have for $\tau\geq 0$,
\begin{align}
    \B_{\out} &\geq \frac{(1-\nu-\eta_2)^2}{(1+\nu+\eta_1)^2}(1 - \frac{(1+\nu+\eta_1)(1+\sigma_x^2)^{1/2}\Vert\theta^*_{(d+1):p}\Vert_{\Sigma_{(d+1):p}}}{(1-\nu-\eta_2)(\lambda_{d+1}\frac{r_d(\Sigma)}{n}+\tau)\Vert\theta^*_{1:d}\Vert_{\Sigma_{1:d}^{-1}}})^2\Vert\theta^*_{1:d}\Vert_{\Sigma_{1:d}^{-1}}^2(\frac{1}{\lambda_d}+\frac{1}{(\lambda_{d+1}\frac{r_d(\Sigma)}{n}+\tau)})^{-2}.  \notag
\end{align}
With Assumption \ref{ass1}(ii), we have for $\tau\geq 0$,
\begin{align}
    \frac{(1+\nu+\eta_1)(1+\sigma_x^2)^{1/2}\Vert\theta^*_{(d+1):p}\Vert_{\Sigma_{(d+1):p}}}{(1-\nu-\eta_2)(\lambda_{d+1}\frac{r_d(\Sigma)}{n}+\tau)\Vert\theta^*_{1:d}\Vert_{\Sigma_{1:d}^{-1}}} &\leq  \frac{2(1+\sigma_x^2)^{1/2}\Vert\theta^*_{(d+1):p}\Vert_{\Sigma_{(d+1):p}}}{\frac{1}{4}\lambda_{d+1}\frac{r_d(\Sigma)}{n}\Vert\theta^*_{1:d}\Vert_{\Sigma_{1:d}^{-1}}} \leq \sqrt{\delta_2} .\label{lem_prf_19}
\end{align}
Moreover, by applying (\ref{lem_prf_19}), in the event $\Omega_1(\nu)\cap \Omega_4\cap \Omega_6(\nu)$ for $0<\nu<\frac{1}{4}$ defined in Lemma \ref{lem3} and Assumption \ref{ass1}(ii), we have for $\tau\geq 0$,
\begin{align}
    \B_{\out} &\geq \frac{(1-\nu-\eta_2)^2}{(1+\nu+\eta_1)^2}(1 - \sqrt{\delta_2})^2\Vert\theta^*_{1:d}\Vert_{\Sigma_{1:d}^{-1}}^2(\frac{1}{\lambda_d}+\frac{1}{(\tau+\lambda_{d+1}\frac{r_d(\Sigma)}{n})})^{-2}.  \notag
\end{align}
This gives Proposition \ref{pro8}.

\noindent \textbf{Lower bound of out-sample variance.} Our strategy for deriving the lower bound of out-sample variance in large TER regime is similar to
that in small or moderate TER regime. With Assumption \ref{ass:4}, we have
\begin{align}
    \frac{1}{2} &> C_0\sigma_x^2\sqrt{\frac{4n^2}{r_d(\Sigma)^2}+\frac{n}{r_d(\Sigma)}}
    >C_0\sigma_x^2 \frac{2n}{r_d(\Sigma)} , \notag
\end{align}
and hence
\begin{align}
    \frac{1}{4C_0\sigma_x^2} \frac{\sum_{j>d}\lambda_j}{n} > \lambda_{d+1} . \label{eq:s_1_1_1_a}
\end{align}
Moreover, with the(\ref{lem_prf_16}) under Assumption \ref{ass:3} and (\ref{eq:s_1_1_1_a}) under Assumption \ref{ass:4}, we have for $\tau\geq 0$,
\begin{align}
   & \frac{2C_0\sigma_x^2d}{r_d(\Sigma^2)}+\frac{C_0\sigma_x^2d}{n} \geq \frac{1}{2} \notag \\
   \Longrightarrow &  \frac{2C_0\sigma_x^2d}{\frac{1}{2} - \eta_1} \geq \frac{\sum_{j>d}\lambda_j^2}{\lambda_{d+1}^2} \notag
   \\ \Longrightarrow &\frac{d}{8C_0\sigma_x^2(\frac{1}{2} - \eta_1) } \geq \frac{\sum_{j>d}\lambda_j^2}{(\lambda_{d+1}\frac{r_d(\Sigma)}{n}+\tau)^2} .\label{lem_prf_20}
\end{align}

We first consider the first term of the right-hand side in (\ref{lem_prf_13}). From Lemma \ref{lem3a}(ii) and (\ref{lem_prf_13}), in the event $\Omega_1(\nu)\cap \Omega_2\cap \Omega_{33}(\nu)\cap \Omega_6(\nu)$ for $0<\nu<\frac{1}{2}\min\{1,\sigma_x^2\}$ defined in Lemma \ref{lem3}, substituting the bounds of $\mu_1(A_d)$ and $\mu_n(A_d)$ in (\ref{lem:bq_11_01aa})--(\ref{lem:bq_11_01a}) into the first term in (\ref{lem_prf_13}), we have for $\tau\geq 0$,
\begin{align}
   & \frac{(1-\nu-\eta_1)^2}{(1+\nu+\eta_1)^4} \sigma^2\frac{nd}{\mu_1(A_d)^2}(\frac{1}{\lambda_d}+\frac{n}{\mu_n(A_d)})^{-2} \notag \\
    &\geq   \frac{(1-\nu-\eta_2)^2(1-\nu-\eta_1)^2}{(1+\nu+\eta_2)^2(1+\nu+\eta_1)^4} \sigma^2\frac{d}{n(\lambda_{d+1}\frac{r_d(\Sigma)}{n}+\tau)^2}(\frac{1}{\lambda_d}+\frac{1}{(\lambda_{d+1}\frac{r_d(\Sigma)}{n}+\tau)})^{-2}. \notag
\end{align}
Then in the event $\Omega_1(\nu)\cap \Omega_2\cap \Omega_{33}(\nu)\cap \Omega_6(\nu)$ for $0<\nu<\frac{1}{2}\min\{1,\sigma_x^2\}$ defined in Lemma \ref{lem3},  we have for $\tau +\lambda_{d+1}\frac{r_d(\Sigma)}{n}\leq \lambda_d$,
\begin{align}
   & \frac{(1-\nu-\eta_1)^2}{(1+\nu+\eta_1)^4} \sigma^2\frac{nd}{\mu_1(A_d)^2}(\frac{1}{\lambda_d}+\frac{n}{\mu_n(A_d)})^{-2} \geq   \frac{(1-\nu-\eta_2)^2(1-\nu-\eta_1)^2}{4(1+\nu+\eta_2)^2(1+\nu+\eta_1)^4} \sigma^2\frac{d}{n} .\label{lem_prf_21}
\end{align}
Then we discuss two complementary cases. The first cases is that $\frac{2C_0\sigma_x^2d}{r_d(\Sigma^2)}+\frac{C_0\sigma_x^2d}{n} \geq \frac{1}{2}$. From (\ref{lem_prf_20}), with Assumption \ref{ass:3} and \ref{ass:4}, we have for $\tau\geq 0$,
\begin{align*}
 &   \frac{d}{8C_0\sigma_x^2(\frac{1}{2} - \eta_1) } \geq \frac{\sum_{j>d}\lambda_j^2}{(\lambda_{d+1}\frac{r_d(\Sigma)}{n}+\tau)^2}  \\
\Longrightarrow &    \frac{d}{n} \geq \frac{8C_0\sigma_x^2(\frac{1}{2}-\eta_1)}{1+8C_0\sigma_x^2(\frac{1}{2}-\eta_1)}(\frac{d}{n}+\frac{\sum_{j>d}\lambda_j^2}{n(\lambda_{d+1}\frac{r_d(\Sigma)}{n}+\tau)^2}).
\end{align*}
Combining with (\ref{lem_prf_21}), we have for $\lambda_{d+1}\frac{r_d(\Sigma)}{n}+\tau\leq \lambda_d$,
\begin{align*}
    \V_{\out}&\geq \frac{(1-\nu-\eta_1)^2}{(1+\nu+\eta_1)^4} \sigma^2\frac{nd}{\mu_1(A_d)^2}(\frac{1}{\lambda_d}+\frac{n}{\mu_n(A_d)})^{-2} \notag\\
    &\geq \frac{(1-\nu-\eta_1)^2(1-\nu-\eta_1)^2}{4(1+\nu+\eta_1)^4(1+\nu+\eta_2)^2} \sigma^2\frac{d}{n}\notag\\
    &\geq   \frac{(1-\nu-\eta_1)^2(1-\nu-\eta_1)^2}{4(1+\nu+\eta_1)^4(1+\nu+\eta_2)^2} \frac{8C_0\sigma_x^2(\frac{1}{2}-\eta_1)}{1+8C_0\sigma_x^2(\frac{1}{2}-\eta_1)}\sigma^2(\frac{d}{n}+\frac{\sum_{j>d}\lambda_j^2}{n(\tau+\lambda_{d+1}\frac{r_d(\Sigma)}{n})^2}).
\end{align*}
The second case is that $\frac{2C_0\sigma_x^2d}{r_d(\Sigma^2)}+\frac{C_0\sigma_x^2d}{n} < \frac{1}{2}$. From (\ref{lem_prf_13}) and Lemma \ref{lem3a}(ii), in the event $\Omega_1(\nu)\cap \Omega_2\cap \Omega_{33}(\nu)\cap\Omega_6(\nu)$ for $0<\nu<\frac{1}{2}\min\{1,\sigma_x^2\}$ defined in Lemma \ref{lem3}, we have for $\tau\geq 0$,
\begin{align}
    \max\{0,\frac{n\sum_{j>d}\lambda_j^2}{\mu_1(A_d)^2}(1-\nu-C_0\sigma_x^2(\frac{2d}{r_d(\Sigma^2)}+\frac{d}{n})\} \geq \frac{\sum_{j>d}\lambda_j^2}{n(\lambda_{d+1}\frac{r_d(\Sigma)}{n}+\tau)^2}(\frac{1}{2}-\nu).\label{lem_prf_22}
\end{align}
From (\ref{lem_prf_13}), (\ref{lem_prf_21}) and (\ref{lem_prf_22}), in the event $\Omega_1(\nu)\cap \Omega_2\cap \Omega_{33}(\nu)\cap\Omega_6(\nu)$ for $\nu<\frac{1}{2}\min\{1,\sigma_x^2\}$ defined in Lemma \ref{lem3}, we have for $\tau\geq 0$,
\begin{align}
    \V_{\out} &\geq  \frac{(1-\nu-\eta_1)^2(1-\nu-\eta_1)^2(\frac{1}{2}-\nu)}{4(1+\nu+\eta_1)^4(1+\nu+\eta_2)^2} \sigma^2(\frac{d}{n}+\frac{\sum_{j>d}\lambda_j^2}{n(\tau+\lambda_{d+1}\frac{r_d(\Sigma)}{n})^2}) . \notag
\end{align}
In conclusion, in the event $\Omega_1(\nu)\cap \Omega_2\cap \Omega_{33}(\nu)\cap\Omega_6(\nu)$ for $0<\nu<\frac{1}{2}\min\{1,\sigma_x^2\}$ defined in Lemma \ref{lem3} and Assumption \ref{ass:3},\ref{ass:4}, we have for $\tau\geq 0$,
\begin{align}
    \V_{\out} &\geq  \frac{(1-\nu-\eta_1)^2(1-\nu-\eta_1)^2(\frac{1}{2}-\nu)}{4(1+\nu+\eta_1)^4(1+\nu+\eta_2)^2} \frac{8C_0\sigma_x^2(\frac{1}{2}-\eta_1)}{1+8C_0\sigma_x^2(\frac{1}{2}-\eta_1)}\sigma^2(\frac{d}{n}+\frac{\sum_{j>d}\lambda_j^2}{n(\tau+\frac{\sum_{j>d}\lambda_j}{n})^2}) , \notag
\end{align}
or equivalently
\begin{align}
    \V_{\out} &\geq  \frac{(1-\nu-\eta_1)^2(1-\nu-\eta_1)^2(\frac{1}{2}-\nu)}{4(1+\nu+\eta_1)^4(1+\nu+\eta_2)^2} \frac{8C_0\sigma_x^2(\frac{1}{2}-\eta_1)}{1+8C_0\sigma_x^2(\frac{1}{2}-\eta_1)}\sigma^2(\frac{d}{n}+\frac{\lambda_{d+1}^2}{(\tau+\lambda_{d+1}\frac{r_d(\Sigma)}{n})^2}\frac{r_d(\Sigma^2)}{n}) . \notag
\end{align}
This gives Proposition \ref{pro9}.

\subsection{Proof of the bounds for in-sample error}

\label{sec_s2}

\subsubsection{Algebraic bounds of the in-sample error}

The bias and variance of the in-sample error can be decomposed or upper bounded as follows:
\begin{align}
    \B_{\In} &= \underbrace{\Vert\theta^*_{1:d} - X_{1:d}^\T A^{-1}X\theta^*\Vert_{\hat{\Sigma}_{1:d}}^2 }_{\B_{\In,1}} + \underbrace{\Vert\theta^*_{(d+1):p} - X_{d+1:p}^\T A^{-1}X\theta^*\Vert_{\hat{\Sigma}_{d+1:p}}^2 }_{\B_{\In,2}} \notag \\
    &+\underbrace{2(\theta^{*T}_{1:d}-\theta^{*T}X^\T A^{-1}X_{1:d})\hat{\Sigma}_{1:d, (d+1):p}(\theta_{(d+1):p}^* - X_{d+1:p}^\T A^{-1}X\theta^*)}_{\B_{\In,12}} ,\notag  \\
    \V_{\In} &\leq \underbrace{2\sigma^2\Tr(A^{-1}X_{1:d}\hat{\Sigma}_{1:d} X_{1:d}^\T A^{-1})}_{\V_{\In,1}} + \underbrace{2\sigma^2\Tr( A^{-1}X_{(d+1):p}\hat{\Sigma}_{(d+1):p}X_{(d+1):p}^\T A^{-1})}_{\V_{\In,2}}  .\label{lem_prf_02a}
\end{align}

The following algebraic bounds are the foundation of the upper bounds and lower bounds for the in-sample error.

\begin{lem}[Algebraic upper bounds of in-sample error] Given $\hat{\Sigma}_{1:d}$ is invertible, we have
\label{lem1a}
\begin{align}
    \B_{\In}&\leq  2\Vert\theta^*_{1:d}\Vert^2_{\hat{\Sigma}_{1:d}^{-1}} (\frac{1}{\lambda_1}+\frac{n^2}{\mu_1^2(A_d)})^{-2} + \frac{2}{n} \frac{\mu_1^2(A_d)}{\mu_n^2(A_d)} \Vert X_{(d+1):p}\theta^*_{(d+1):p}\Vert^2  \notag \\
    &+3(\Vert \hat{M}_d\Vert \frac{n}{\mu^2_n(A_d)} \Vert\theta^*_{1:d}\Vert_{\hat{\Sigma}^{-1}}^2 (\frac{1}{\lambda_1}+\frac{n}{\mu_1(A_d)})^{-2}  +\Vert\theta^*_{(d+1):p}\Vert^2_{\hat{\Sigma}_{(d+1):p}} + \frac{\Vert \hat{M}_d\Vert}{\mu_n(A)^2}\Vert X_{(d+1):p}\theta^*_{(d+1):p}\Vert^2 )  ,\label{lem_prf_23} \\
    \V_{\In} &\leq   2\sigma^2 \frac{\mu_1^2(A_d)}{\mu_n^2(A_d)}\frac{d}{n} +  2\sigma^2 \frac{\frac{1}{n}\Tr(X_{(d+1):p}X_{(d+1):p}^\T)\mu_1(X_{(d+1):p}X_{(d+1):p}^\T )}{\mu_n(A_d)^2}   . \notag
\end{align}
\end{lem}
\begin{prf}
$\newline$
\textbf{Algebraic upper bound for the in-sample bias.} Given invertible $\hat{\Sigma}_{1:d}$, (\ref{lem_prf_23}) can be derived similarly as (\ref{lem_prf_03a}) .
$\newline$
\textbf{Algebraic upper bound for the in-sample variance.} Similarly as the derivation in Lemma 27 in \cite{Tsigler_ridge_2023}, given invertible $\hat{\Sigma}_{1:d}$, we have
\begin{align}
\V_{\In,1}=2\sigma^2\Tr(A^{-1}X_{1:d}\hat{\Sigma}_{1:d} X_{1:d}^\T A^{-1}) &\leq  2\frac{\sigma^2\mu_1(A_d)^2\Tr(X_{1:d}\hat{\Sigma}_{1:d}^{-1}X_{1:d})}{\mu_n(A_d)^2\mu_d(\hat{H}_d  \hat{H}_d^\T )^2} \notag \\
&=2\sigma^2 \frac{\mu_1^2(A_d)}{\mu_n^2(A_d)}\frac{d}{n} . \notag
\end{align}
Moreover, we have
\begin{align}
    \V_{\In,2}=2\sigma^2\Tr( A^{-1}X_{(d+1):p}\hat{\Sigma}_{(d+1):p}X_{(d+1):p}^\T A^{-1}) &\leq  2\sigma^2\frac{\Tr(X_{(d+1):p}\hat{\Sigma}_{(d+1):p}X_{(d+1):p}^\T )}{\mu_n(A_d)^2} \notag  \\
&=2\sigma^2\frac{\frac{1}{n}\Tr((X_{(d+1):p}X_{(d+1):p}^\T)^2 )}{\mu_n(A_d)^2} \notag \\
    &\leq 2\sigma^2 \frac{\frac{1}{n}\Tr(X_{(d+1):p}X_{(d+1):p}^\T)\mu_1(X_{(d+1):p}X_{(d+1):p}^\T )}{\mu_n(A_d)^2} . \notag
\end{align}
By (\ref{lem_prf_02a}), combining the preceding bounds gives the algebraic upper bound of in-sample variance.
\end{prf}

\begin{lem}[Algebraic lower bounds of in-sample error] Given invertible $\hat{\Sigma}_{1:d}$, we have
\begin{align}
    \V_\In
    &\geq \sigma^2 \frac{1}{2}(\frac{1}{n}\sum_{i=1}^d \frac{\mu^2_i(X_{1:d}X_{1:d}^\T )}{(\mu_i(X_{1:d}X_{1:d}^\T )+n\tau)^2} +\frac{1}{n}\sum_{i=1}^n \frac{\mu_i^2(X_{(d+1):p}X_{(d+1):p}^\T )}{(\mu_i(X_{(d+1):p}X_{(d+1):p}^\T )+n\tau)^2}) ,\label{lem_prf_25}
\end{align}
Further if $\Vert\theta^*_{1:d}\Vert_{\hat{\Sigma}_{1:d}^{-1}}\geq \frac{n^{1/2}\Vert X_{(d+1):p}\theta^*_{(d+1):p}\Vert}{\mu_n(A_d)}$, then
\begin{align}
     \B_{\In}&\geq  \max\{0,1-\frac{|\B_{\In,12}|}{\B_{\In}}\}(1 - \frac{n^{1/2}\Vert X_{(d+1):p}\theta^*_{(d+1):p}\Vert}{\mu_n(A_d)\Vert\theta^*_{1:d}\Vert_{\hat{\Sigma}_{1:d}^{-1}}})^2\Vert\theta^*_{1:d}\Vert^2_{\hat{\Sigma}_{1:d}^{-1}}(\frac{1}{\lambda_d}+\frac{n}{\mu_n(A_d)})^{-2}  .\notag
\end{align}
\end{lem}
\begin{prf}
$\newline$
\textbf{Algebraic lower bound for the in-sample variance.} The in-sample variance is
\begin{align*}
    V_{\In} &= \sigma^2  \Tr(A^{-1}X\hat{\Sigma} X^\T  A^{-1}) = \sigma^2\frac{1}{n}\sum_{i=1}^n  \frac{\mu^2_i(XX^\T )}{(\mu_i(XX^\T )+n\tau)^2}.
\end{align*}
From Weyl's inequality (Lemma \ref{lem_S12b}), we have
\begin{align*}
  &  \mu_i(XX^\T ) \geq \mu_i(X_{1:d}X_{1:d}^\T ), \quad i=1, \ldots, d, \\
  &  \mu_i(XX^\T ) \geq \mu_i(X_{(d+1):p}X_{(d+1):p}^\T ). \quad i=1, \ldots, n
\end{align*}
Then $\V_{\In}$ can be lower bounded by
\begin{align*}
    \V_\In &\geq \sigma^2 \max\{\frac{1}{n}\sum_{i=1}^d \frac{\mu^2_i(X_{1:d}X_{1:d}^\T )}{(\mu_i(X_{1:d}X_{1:d}^\T )+n\tau)^2}   ,\frac{1}{n}\sum_{i=1}^n \frac{\mu_i^2(X_{(d+1):p}X_{(d+1):p}^\T )}{(\mu_i(X_{(d+1):p}X_{(d+1):p}^\T )+n\tau)^2}\} \\
    &\geq \frac{\sigma^2}{2}(\frac{1}{n}\sum_{i=1}^d \frac{\mu^2_i(X_{1:d}X_{1:d}^\T )}{(\mu_i(X_{1:d}X_{1:d}^\T )+n\tau)^2} +\frac{1}{n}\sum_{i=1}^n \frac{\mu_i^2(X_{(d+1):p}X_{(d+1):p}^\T )}{(\mu_i(X_{(d+1):p}X_{(d+1):p}^\T )+n\tau)^2}).
\end{align*}
$\newline$
\textbf{Algebraic lower bound for the in-sample bias.} Note that
\begin{align*}
    \B_\In \geq \max\{0,1-\frac{|\B_{\In,12}|}{\B_{\In,1}}\}\B_{\In,1}.
\end{align*}
The result follows from the algebraic lower bound for $\B_{\In,1}$:
\begin{align}
     \B_{\In,1}&\geq  (1 - \frac{n^{1/2}\Vert X_{(d+1):p}\theta^*_{(d+1):p}\Vert}{\mu_n(A_d)\Vert\theta^*_{1:d}\Vert_{\hat{\Sigma}_{1:d}^{-1}}})^2\Vert\theta^*_{1:d}\Vert^2_{\hat{\Sigma}_{1:d}^{-1}}(\frac{1}{\lambda_d}+\frac{n}{\mu_n(A_d)})^{-2} , \notag
\end{align}
which can be derived similarly as (\ref{lem_prf_08a}).
\end{prf}

\subsubsection{Intermediate bounds of the in-sample error}

We give the intermediate bounds of in-sample error under the event that some random quantities in the algebraic bounds above are controlled. In the event $\Omega_1(\nu)\cap \Omega_4\cap \Omega_5$ for $0<\nu <\frac{1}{2}$ defined in Lemma \ref{lem3}, substituting the bounds of $\Vert\theta^*_{1:d}\Vert_{\hat{\Sigma}_{1:d}^{-1}}$, $\Vert X_{(d+1):p}\theta^*_{(d+1):p}\Vert$ and $\Vert \hat{M}_d\Vert$ into the algebraic upper bound of $\B_\In$
yields
\begin{align}
    \B_{\In}&\leq \frac{1}{(1-\nu-\eta_1)^2}\Vert\theta^*_{1:d}\Vert_{\Sigma_{1:d}^{-1}}^2(\frac{1}{\lambda_1}+\frac{n}{\mu_1(A_d)})^{-2}
    + \poly_2(\sigma_x) \frac{\mu_1^2(A_d)}{\mu_n^2(A_d)} \Vert \theta^*_{(d+1):p}\Vert_{\Sigma_{(d+1):p}}^2 \notag \\
    &+ \frac{1}{(1-\nu-\eta_1)^2}\poly_4(\sigma_x)\frac{(n\lambda_{d+1}+\sum_{j>d}\lambda_j)^2}{\mu^2_n(A_d)}  \Vert\theta^*_{1:d}\Vert_{\Sigma^{-1}}^2  (\frac{1}{\lambda_1}+\frac{n}{\mu_1(A_d)})^{-2}\notag \\
    &+(\poly_2(\sigma_x)+\poly_6(\sigma_x)\frac{(n\lambda_{d+1}+\sum_{j>d}\lambda_j)^2}{\mu_n(A_d)^2}\Vert\theta^*_{(d+1):p}\Vert_{\Sigma_{(d+1):p}}^2  . \label{lem_prf_26}
\end{align}
In the event $\Omega_1(\nu)\cap \Omega_4$ for $0<\nu<\frac{1}{2}$ defined in Lemma \ref{lem3} and with
\begin{align*}
\Vert\theta^*_{1:d}\Vert_{\Sigma_{1:d}^{-1}}\geq \frac{(1+\nu+\eta_1)(1+\sigma_x^2)^{1/2}n\Vert \theta^*_{(d+1):p}\Vert_{\Sigma_{(d+1):p}}}{\mu_n(A_d)},
\end{align*}
substituting the bound of $\Vert\theta^*_{1:d}\Vert_{\hat{\Sigma}_{1:d}^{-1}}$ and $\Vert X_{(d+1):p}\theta^*_{(d+1):p}\Vert$ into the algebraic lower bound of $\B_\In$
yields
\begin{align}
    \B_{\In} \geq \frac{\max\{0,1-\frac{|\B_{\In,2}|}{\B_{\In,1}}\}}{(1+\nu+\eta_1)^2}(1 - \frac{(1+\nu+\eta_1)(1+\sigma_x^2)^{1/2}n\Vert \theta^*_{(d+1):p}\Vert_{\Sigma_{(d+1):p}}}{\mu_n(A_d)\Vert\theta^*_{1:d}\Vert_{\Sigma_{1:d}^{-1}}})^2\Vert\theta^*_{1:d}\Vert_{\Sigma_{1:d}^{-1}}^2(\frac{1}{\lambda_d}+\frac{n}{\mu_n(A_d)})^{-2} .\label{lem_prf_26a}
\end{align}
In the event $\Omega_5\cap \Omega_7(\nu)$ for $0<\nu < \min\{1,\sigma_x^2\}$, substituting the $\Tr(X_{(d+1):p}X_{(d+1):p}^\T)$ and $\mu_1(X_{(d+1):p}\cdot X_{(d+1):p}^\T )$ into the algebraic upper bound of $\V_\In$ yields
\begin{align}
    \V_{\In}\leq 2\sigma^2 \frac{\mu_1^2(A_d)}{\mu_n^2(A_d)} \frac{d}{n} + \poly_4(\sigma_x) \sigma^2 \frac{(\sum_{j>d}\lambda_j)(n\lambda_{d+1}+\sum_{j>d}\lambda_j)}{\mu_n(A_d)^2} . \label{lem_prf_27}
\end{align}

\subsubsection{Final upper bounds of in-sample bias and variance}

We give the final upper bounds of in-sample bias and variance. We first discuss the small or moderate TER regime and then the large TER regime.

\noindent \textbf{(i) Small or moderate TER}

From Lemma \ref{lem3a}(i) and (\ref{lem_prf_26}),(\ref{lem_prf_27}), in the event $\Omega_1(\nu)\cap\Omega_4\cap \Omega_5\cap\Omega_7(\nu)$ for $0<\nu<\frac{1}{2}\min\{1,\sigma_x^2\}$ defined in Lemma \ref{lem3} and Assumption \ref{ass:3a}, substituting the bounds of $\mu_1(A_d)$ and $\mu_n(A_d)$ in (\ref{lem:bq_11_00})--(\ref{lem:bq_11_0}) into (\ref{lem_prf_26}) and (\ref{lem_prf_27}), we have for $\tau\geq \lambda_{d+1}$,
\begin{align}
     \B_{\In}&\leq \frac{(1+C_1)^2\poly_4(\sigma_x)}{(1-\nu-\eta_1)^2}\Vert\theta^*_{1:d}\Vert_{\Sigma_{1:d}^{-1}}^2(\frac{1}{\lambda_1}+\frac{1}{\tau})^{-2}
    +\poly_6(\sigma_x) (1+C_1)^2 \Vert \theta^*_{(d+1):p}\Vert_{\Sigma_{(d+1):p}}^2 \notag \\
    &+ \frac{(1+C_1)^2\poly_8(\sigma_x)}{(1-\nu-\eta_1)^2} \Vert\theta^*_{1:d}\Vert_{\Sigma^{-1}}^2  (\frac{1}{\lambda_1}+\frac{1}{\tau})^{-2}\notag \\
    &+(\poly_2(\sigma_x)+\poly_6(\sigma_x))\Vert\theta^*_{(d+1):p}\Vert_{\Sigma_{(d+1):p}}^2 ,\notag \\
     \V_{\In}&\leq\poly_4(\sigma_x)(1+C_1)^2 \frac{d}{n} + \poly_4(\sigma_x) \sigma^2 \frac{(\sum_{j>d}\lambda_j)(n\lambda_{d+1}+\sum_{j>d}\lambda_j)}{n^2\tau^2}  .\notag
\end{align}
Hence, we obtain the upper bounds
\begin{align}
    \B_{\In}&\leq \frac{(1+C_1)^4\poly_8(\sigma_x)}{(1-\nu-\eta_1)^2}(\Vert\theta^*_{1:d}\Vert_{\Sigma_{1:d}^{-1}}^2(\frac{1}{\lambda_1}+\frac{1}{\tau})^{-2}+ \Vert \theta^*_{(d+1):p}\Vert^2_{\Sigma_{(d+1):p}} ), \notag \\
     \V_{\In}&\leq \poly_4(\sigma_x)(1+C_1)^2 \sigma^2(  \frac{d}{n}  +  \frac{\lambda^2_{d+1}}{\tau^2}\frac{r_d(\Sigma)}{n}).\notag
\end{align}
Further with (\ref{lem_prf_14a}), we have for $\tau\geq \lambda_{d+1}$,
\begin{align}
    \B_{\In}&\leq \frac{(1+C_1)^4\poly_8(\sigma_x)}{(1-\nu-\eta_1)^2}\Vert\theta^*_{1:d}\Vert_{\Sigma_{1:d}^{-1}}^2(\frac{1}{\lambda_1}+\frac{1}{\tau})^{-2} .\notag
\end{align}
This gives Proposition \ref{pro4}.

\noindent \textbf{(ii) Large TER}

From Lemma \ref{lem3a}(ii) and (\ref{lem_prf_10}),(\ref{lem_prf_11}), in the event $\Omega_1(\nu)\cap \Omega_4\cap \Omega_5\cap \Omega_6(\nu)\cap\Omega_7(\nu)$ for $0<\nu <\frac{1}{2}\min\{1,\sigma_x^2\}$ defined in Lemma \ref{lem3} and Assumption \ref{ass:4}, substituting the bounds of $\mu_1(A_d)$ and $\mu_n(A_d)$ in (\ref{lem:bq_11_01aa})--(\ref{lem:bq_11_01a}) into (\ref{lem_prf_26}) and
(\ref{lem_prf_27}),  we have
\begin{align}
    \B_{\In}&\leq \frac{(1+\nu+\eta_2)^2}{(1-\nu-\eta_1)^2}\Vert\theta^*_{1:d}\Vert_{\Sigma_{1:d}^{-1}}^2(\frac{1}{\lambda_1}+\frac{1}{\tau+\lambda_{d+1}\frac{r_d(\Sigma)}{n}})^{-2}+\poly_2(\sigma_x) \frac{(1+\nu+\eta_2)^2}{(1-\nu-\eta_2)^2} \Vert \theta^*_{(d+1):p}\Vert_{\Sigma_{(d+1):p}}^2 \notag \\
    &+ \frac{\poly_4(\sigma_x)(1+\nu+\eta_2)^2}{(1-\nu-\eta_1)^2(1-\nu-\eta_2)^2} \Vert\theta^*_{1:d}\Vert_{\Sigma^{-1}}^2  (\frac{1}{\lambda_1}+\frac{1}{\tau+\lambda_{d+1}\frac{r_d(\Sigma)}{n}})^{-2}\notag \\
    &+(\poly_2(\sigma_x)+\frac{\poly_6(\sigma_x)}{(1-\nu-\eta_2)^2})\Vert\theta^*_{(d+1):p}\Vert_{\Sigma_{(d+1):p}}^2 , \notag \\
    \V_{\In}&\leq \frac{2(1+\nu+\eta_2)^2}{(1-\nu-\eta_2)^2}\sigma^2  \frac{d}{n} + \frac{\poly_4(\sigma_x)}{(1-\nu-\eta_2)^2} \sigma^2 \frac{(\lambda_{d+1}\frac{r_d(\Sigma)}{n})^2}{(\lambda_{d+1}\frac{r_d(\Sigma)}{n}+\tau)^2} .  \notag
\end{align}
Hence, we obtain the upper bound
\begin{align}
    \B_{\In}&\leq \frac{(1+\nu+\eta_2)^2\poly_6(\sigma_x)}{(1-\nu-\eta_1)^2(1-\nu-\eta_2)^2}(\Vert\theta^*_{1:d}\Vert_{\Sigma_{1:d}^{-1}}^2 (\frac{1}{\lambda_1}+\frac{1}{\tau+\lambda_{d+1}\frac{r_d(\Sigma)}{n}})^{-2} + \Vert \theta^*_{(d+1):p}\Vert^2_{\Sigma_{(d+1):p}}) ,\notag \\
    \V_{\In}&\leq \frac{(1+\nu+\eta_2)^2\poly_4(\sigma_x)}{(1-\nu-\eta_2)^2}\sigma^2  (\frac{d}{n} +\frac{\lambda_{d+1}^2}{(\tau+\lambda_{d+1}\frac{r_d(\Sigma)}{n})^2}\frac{r^2_d(\Sigma)}{n^2}). \notag
\end{align}
Further with (\ref{lem_prf_19a}), we have for $\tau\geq 0$,
\begin{align}
    \B_{\In}&\leq \frac{(1+\nu+\eta_2)^2\poly_6(\sigma_x)}{(1-\nu-\eta_1)^2(1-\nu-\eta_2)^2}\Vert\theta^*_{1:d}\Vert_{\Sigma_{1:d}^{-1}}^2 (\frac{1}{\lambda_1}+\frac{1}{\tau+\lambda_{d+1}\frac{r_d(\Sigma)}{n}})^{-2} .
\end{align}
This gives Proposition \ref{pro10}.

\subsubsection{Final lower bounds of in-sample bias and variance}

We give the final lower bounds of in-sample bias and variance. As preparation, we give the following lemma to compare $\B_{\In,1}$ and $|\B_{\In,12}|$.

\begin{lem}[Comparison between $\B_{\In,1}$ and $|\B_{\In,12}|$]
\label{lem6a}
$\newline$
(i) Given Assumption \ref{ass:3a} and \ref{ass1}(i) and in the event $\Omega_1(\nu)\cap \Omega_4\cap \Omega_5$ for $0<\nu<\frac{1}{2}$,
we have for $\tau \geq \lambda_{d+1}$,
\begin{align}
    \max\{1-\frac{|\B_{\In,12}|}{\B_{\In,1}},0\}  \geq \kappa_{1}(\tau) ,\label{eq_I_2_4_1}
\end{align}
where $\kappa_1(\tau)=\max\{1-(\frac{2 C_0\sigma_x^2(2+C_1)\lambda_{d+1}}{\tau}(1+16(2C_0\sigma_x^2+1)(1+C_1)\frac{\sqrt{\delta_1}}{1-\sqrt{\delta_1}}) + 64\frac{\sqrt{\delta_1}}{1-\sqrt{\delta_1}}),0\}$.

$\newline$
(ii) Given Assumption \ref{ass:4} and \ref{ass1}(ii) and in the event $\Omega_1(\nu)\cap \Omega_4\cap \Omega_6(\nu)$ for $0<\nu<\frac{1}{4}$, we have for $\tau\geq 0$,
\begin{align}
    \max\{1-\frac{|\B_{\In,12}|}{\B_{\In,1}},0\} \geq \kappa_{2}(\tau) ,
\end{align}
where $\kappa_2(\tau)=\max\{1-(16\frac{\lambda_{d+1}\frac{r_d(\Sigma)}{n}}{\tau+\lambda_{d+1}\frac{r_d(\Sigma)}{n}}(1+112\frac{\sqrt{\delta_2}}{1-\sqrt{\delta_2}})+64\frac{\sqrt{\delta_2}}{1-\sqrt{\delta_2}}),0\}$.
\end{lem}

The proof of Lemma \ref{lem6a} is left to Section~\ref{prf_lem6a}.
In the following, we first discuss the small or moderate TER regime and then the large TER regime.

\noindent \textbf{(i) Small or moderate TER}
$\newline$
\textbf{Lower bound of in-sample bias.} From (\ref{lem_prf_15a}) and Lemma \ref{lem6a}(i), in the event $\Omega_1(\nu)\cap \Omega_4\cap \Omega_5$ for $0<\nu<\frac{1}{2}$ defined in Lemma \ref{lem3} and Assumption \ref{ass1}(i), substituting the lower bound of $\mu_n(A_d)$ and lower bound of $\max\{1-\frac{|\B_{\In,12}|}{\B_{\In,1}},0\}$ in (\ref{eq_I_2_4_1}) into (\ref{lem_prf_26a}), we have for $\tau\geq \lambda_{d+1}$,
\begin{align}
    \B_{\In} \geq \frac{\kappa_1(\tau)}{(1+\nu+\eta_1)^2}(1 - \frac{2(1+\sigma_x^2)^{1/2}\Vert \theta^*_{(d+1):p}\Vert_{\Sigma_{(d+1):p}}}{\tau\Vert\theta^*_{1:d}\Vert_{\Sigma_{1:d}^{-1}}})^2\Vert\theta^*_{1:d}\Vert_{\Sigma_{1:d}^{-1}}^2(\frac{1}{\lambda_d}+\frac{1}{\tau})^{-2}  .\notag
\end{align}
Moreover, by applying (\ref{lem_prf_14}), we have for $\tau\geq \lambda_{d+1}$,
\begin{align}
    \B_{\In} \geq \frac{\kappa_1(\tau)}{(1+\nu+\eta_1)^2}(1 - \sqrt{\delta_1})^2\Vert\theta^*_{1:d}\Vert_{\Sigma_{1:d}^{-1}}^2(\frac{1}{\lambda_d}+\frac{1}{\tau})^{-2}  .\notag
\end{align}
This gives the Proposition \ref{pro5}.

\noindent \textbf{Lower bound of in-sample variance.} We first study the first term in (\ref{lem_prf_25}). In the event $\Omega_1(\nu)$ for $0<\nu <\frac{1}{2}$, substituting the lower bound of $\mu_d(X_{1:d}X_{1:d}^T)=n\mu_d(\hat{\Sigma}_{1:d})$ in (\ref{lem:bq_2}), we have for $\tau\leq \lambda_d$,
\begin{align*}
    \mu_d(X_{1:d}^\T X_{1:d}) &\geq  \lambda_{d} n (1-\nu-\eta_1)^2,\\
    &\geq  \tau n (1-\nu-\eta_1)^2,
\end{align*}
and hence
\begin{align*}
    \frac{1}{(1-\nu-\eta_1)^2}\mu_d(X_{1:d}^\T X_{1:d})\geq  \tau n.
\end{align*}
Then in the event $\Omega_1(\nu)$ for $0<\nu <\frac{1}{2}$, we have for $\tau \leq \lambda_d$,
\begin{align}
    \frac{1}{n}\sum_{i=1}^d \frac{\mu^2_i(X_{1:d}X_{1:d}^\T )}{(\mu_i(X_{1:d}X_{1:d}^\T )+n\tau)^2}\geq (\frac{1}{1+\frac{1}{(1-\nu-\eta_1)^2}})^2\frac{d}{n}. \label{eq:I_2_4}
\end{align}
Then we study the second term in (\ref{lem_prf_25}).
Given Assumption \ref{ass:3a} and in the event $\Omega_5\cap\Omega_7(\nu)$ for $0<\nu<\sigma_x^2$, substituting the bounds of $\mu_1(A_d)$ and $\Tr(X_{(d+1):p}X_{(d+1):p}^T)$ in (\ref{eq_I_3_1})--(\ref{lem:bq_11_00}) into the second term in (\ref{lem_prf_25}), we have for $\tau\geq \lambda_{d+1}$,
\begin{align}
    \sigma^2\frac{1}{n}\sum_{i=1}^n \frac{\mu_i^2(X_{(d+1):p}X_{(d+1):p}^\T )}{(\mu_i(X_{(d+1):p}X_{(d+1):p}^\T )+n\tau)^2}&\geq  \sigma^2\frac{1}{n}\sum_{i=1}^n \frac{\mu_i^2(X_{(d+1):p}X_{(d+1):p}^\T )}{\mu^2_1(A_d)}  \notag \\
    &\geq  \sigma^2\frac{1}{n}\sum_{i=1}^n \frac{\mu_i^2(X_{(d+1):p}X_{(d+1):p}^\T )}{(2C_0\sigma_x^2+1)^2(1+C_1)^2n^2\tau^2 } \notag \\
    &\geq \sigma^2\frac{1}{n} \frac{\Tr(X_{(d+1):p}X_{(d+1):p}^\T )^2}{(2C_0\sigma_x^2+1)^2(1+C_1)^2n^3\tau^2} \notag \\
    &\geq  \frac{(1-\nu)^2}{(2C_0\sigma_x^2+1)^2(1+C_1)^2} \sigma^2\frac{(\frac{\sum_{j>d}\lambda_j}{n})^2}{\tau^2} .\notag
\end{align}
By combining the two terms,
in the event $\Omega_1(\nu)\cap \Omega_5\cap \Omega_7(\nu)$ for $0<\nu < \frac{1}{2}\min\{1,\sigma_x^2\}$, we have for $\lambda_{d+1}\leq \tau \leq \lambda_d$,
\begin{align}
    V_{\In}
    &\geq \frac{(1-\nu)^2}{2(2C_0\sigma_x^2+1)^2(1+C_1)^2(1+\frac{1}{(1-\nu-\eta_2)^2})^2}\sigma^2(\frac{d}{n} +\frac{(\frac{\sum_{j>d}\lambda_j}{n})^2}{\tau^2})  ,\notag
\end{align}
or equivalently,
\begin{align}
    V_{\In}
    &\geq \frac{(1-\nu)^2}{2(2C_0\sigma_x^2+1)^2(1+C_1)^2(1+\frac{1}{(1-\nu-\eta_2)^2})^2}\sigma^2(\frac{d}{n} +\frac{\lambda_{d+1}^2}{\tau^2}\frac{r_d^2(\Sigma)}{n^2})  .\notag
\end{align}
For $\tau > \lambda_d$, we have
\begin{align}
     V_{\In} &\geq \sigma^2\frac{1}{n}\sum_{i=1}^n \frac{\mu_i^2(X_{(d+1):p}X_{(d+1):p}^\T )}{(\mu_i(X_{(d+1):p}X_{(d+1):p}^\T )+n\tau)^2}\notag \\
     &\geq \frac{(1-\nu)^2}{2(2C_0\sigma_x^2+1)^2(1+C_1)^2} \sigma^2\frac{(\frac{\sum_{j>d}\lambda_j}{n})^2}{\tau^2}, \notag
\end{align}
or equivalently
\begin{align}
     V_{\In} \geq \frac{(1-\nu)^2}{2(2C_0\sigma_x^2+1)^2(1+C_1)^2} \sigma^2\frac{\lambda_{d+1}^2}{\tau^2}\frac{r_d^2(\Sigma)}{n^2}. \notag
\end{align}
The result for $\tau \leq \lambda_{d+1}$ follows from the monotonicity of variance in Lemma \ref{lem_12b}. This gives Proposition \ref{pro6}.

\noindent \textbf{(ii) Large TER}

\noindent \textbf{Lower bound of in-sample bias.} From (\ref{lem_prf_20a}) and Lemma \ref{lem6a}(ii), in the event $\Omega_1(\nu)\cap \Omega_4\cap \Omega_6(\nu)$ for $0<\nu<\frac{1}{4}$ defined in Lemma \ref{lem3} and Assumption \ref{ass1}(ii), substituting the bounds of $\mu_n(A_d)$ in (\ref{lem:bq_11_01a}) and the bound of $\max\{0,1-\frac{|\B_{\In,2}|}{\B_{\In,1}}\}$ in (\ref{eq_I_2_4_2}) into (\ref{lem_prf_26a}), we have for $\tau\geq 0$,
\begin{align}
    \B_{\In}
    \geq &\frac{\kappa_2(\tau)(1-\nu-\eta_2)^2}{(1+\nu+\eta_1)^2}(1 - \frac{(1+\nu+\eta_1)(1+\sigma_x^2)^{1/2}n\Vert \theta^*_{(d+1):p}\Vert_{\Sigma_{(d+1):p}}}{(1-\nu-\eta_2)(\frac{\sum_{j>d}\lambda_j}{n}+\tau)\Vert\theta^*_{1:d}\Vert_{\Sigma_{1:d}^{-1}}})^2\\&\cdot \Vert\theta^*_{1:d}\Vert_{\Sigma_{1:d}^{-1}}^2(\frac{1}{\lambda_d}+\frac{n}{(\frac{\sum_{j>d}\lambda_j}{n}+\tau))})^{-2} . \notag
\end{align}
Moreover, by applying (\ref{lem_prf_19}), we have for $\tau\geq 0$,
\begin{align}
    \B_{\In} \geq \frac{\kappa_2(\tau)(1-\nu-\eta_2)^2}{(1+\nu+\eta_1)^2}(1 - \sqrt{\delta_2})^2\Vert\theta^*_{1:d}\Vert_{\Sigma_{1:d}^{-1}}^2(\frac{1}{\lambda_d}+\frac{n}{(\tau+\lambda_{d+1}\frac{r_d(\Sigma)}{n})})^{-2}  .\notag
\end{align}
This gives Proposition \ref{pro11}.
$\newline$
\textbf{Lower bound of in-sample variance.} The first term in (\ref{lem_prf_25}) can be studied for $\tau \le \lambda_d$  similarly as (\ref{eq:I_2_4}) in the small or moderate TER regime. Then we study the second term in (\ref{lem_prf_25}). In the event $\Omega_6(\nu)$ for $0<\nu < \frac{1}{2}$, we have for $\tau\geq 0$,
\begin{align}
      &  (1-\nu-\eta_2) \sum_{j>d}\lambda_j \leq \mu_n(X_{(d+1):p}X_{(d+1):p}^\T ) \leq  (1+\nu+\eta_2)\sum_{j>d}\lambda_j  , \notag
\end{align}
and hence
\begin{align}
    \sigma^2\frac{1}{n}\sum_{i=1}^n \frac{\mu_i^2(X_{(d+1):p}X_{(d+1):p}^\T )}{(\mu_i(X_{(d+1):p}X_{(d+1):p}^\T )+n\tau)^2}&\geq  \sigma^2\frac{1}{n}\sum_{i=1}^n \frac{\mu_i^2(X_{(d+1):p}X_{(d+1):p}^\T )}{\mu^2_1(A_d)}  \notag \\
    &\geq \sigma^2\frac{1}{n^2} \frac{\Tr(X_{(d+1):p}X_{(d+1):p}^\T )^2}{\mu^2_1(A_d)} \notag \\
    &\geq  \frac{(1-\nu-\eta_2)^2}{(1+\nu+\eta_2)^2} \sigma^2\frac{(\frac{\sum_{j>d}\lambda_j}{n})^2}{(\tau+\frac{\sum_{j>d}\lambda_j}{n})^2} \notag \\
    &=\sigma^2\frac{\lambda_{d+1}^2}{(\tau+\lambda_{d+1}\frac{r_d(\Sigma)}{n})^2}\frac{r^2_d(\Sigma)}{n^2}.
\end{align}
By combining the two terms, in the event $\Omega_1(\nu)\cap \Omega_6(\nu)$ for $0<\nu < \frac{1}{2}$, we have for  $0\leq \tau +\lambda_{d+1} \frac{r_d(\Sigma)}{n} \leq \lambda_d$,
\begin{align*}
    V_{\In} &\geq  \frac{(1-\nu-\eta_2)^2}{2(1+\nu+\eta_2)^2(1+\frac{1}{(1-\nu-\eta_1)^2})}\sigma^2(\frac{d}{n} + \frac{\lambda_{d+1}^2}{(\tau+\lambda_{d+1}\frac{r_d(\Sigma)}{n})^2}\frac{r^2_d(\Sigma)}{n^2}).
\end{align*}
For $\tau+\lambda_{d+1} \frac{r_d(\Sigma)}{n} > \lambda_d$, we have
\begin{align*}
   \V_\In &\geq \sigma^2\frac{1}{n}\sum_{i=1}^n \frac{\mu_i^2(X_{(d+1):p}X_{(d+1):p}^\T )}{(\mu_i(X_{(d+1):p}X_{(d+1):p}^\T )+n\tau)^2} \notag\\
   &\geq \frac{(1-\nu-\eta_2)^2}{2(1+\nu+\eta_2)^2} \sigma^2\frac{\lambda_{d+1}^2}{(\tau+\lambda_{d+1}\frac{r_d(\Sigma)}{n})^2}\frac{r^2_d(\Sigma)}{n^2}.
\end{align*}
This gives Proposition \ref{pro12}.

\subsection{Bounds of random quantities}

We give the probability bounds of some random quantities used in the proofs in Supplement Sections \ref{bounds_pro} and \ref{sec_s2}.
We define the related events and give the probability bounds for the events.

\begin{lem}[Bounds of random quantities]
\label{lem3}
$\newline$
\vspace{-25pt}
\begin{itemize}
 \item[(i)]
[Bounding $\mu_1(H_dH_d^\T )$ and $\mu_d(H_dH_d^\T )$] For $\eta_1$ defined in Assumption \ref{ass:3} and $0<\nu < \frac{1}{2}$, denote by $\Omega_1(\nu)$ the event that
\begin{align}
     &(1-\nu - \eta_1)^2n\leq \mu_d(H_dH_d^\T )
    \leq \mu_1(H_dH_d^\T ) \leq  n(1+\nu+\eta_1)^2.  \label{lem:bq_1}
\end{align}
In the event $\Omega_1(\nu)$, we have
\begin{align}
&\lambda_{d}  (1-\nu-\eta_1)^2\leq \mu_d(\hat{\Sigma}_{1:d}) , \label{lem:bq_2} \\
&\frac{\Vert\theta^*_{1:d}\Vert^2_{\Sigma_{1:d}^{-1}}}{(1+\nu+\eta_1)^2}\leq \Vert\theta^*_{1:d}\Vert_{\hat{\Sigma}^{-1}_{1:d}}^2 \leq \frac{\Vert\theta^*_{1:d}\Vert^2_{\Sigma_{1:d}^{-1}}}{(1-\nu-\eta_1)^2}  . \label{lem:bq_3}
\end{align}
Under the sub-gaussianity of $\Sigma^{-\frac{1}{2}}x_i$ and Assumption \ref{ass:3},  $\P(\Omega_1(\nu)) \geq 1-2 \exp\{-\frac{\nu^2 n}{C_0^2\sigma_x^4}\}$.

\item[(ii)] [Bounding $\Vert M_d\Vert$]
Denote by $\Omega_2$ the event that
\begin{align}
\Vert M_d \Vert&\leq C_0\sigma_x^2(2n\lambda_{d+1}^2+\sum_{j>d}\lambda_j^2). \notag
\end{align}
Under the sub-gaussianity of $\Sigma^{-\frac{1}{2}}x_i$, $\P(\Omega_2)\geq 1-6 \exp\{-\frac{ n}{C_0}\}$.

\item[(iii)] [Bounding the trace of $\Tr(H_d^\T H_d)$ and $\Tr(M_d)$]
Denote by $\Omega_{31}$ the event that
\begin{align}
    \Tr(H_d^\T H_d) \leq (1+\sigma_x^2) nd . \notag
\end{align}
Denote by $\Omega_{32}$ the event that
\begin{align}
\Tr(M_d)&\leq (1+\sigma_x^2) n\sum_{j>d}\lambda_j^2 . \notag
\end{align}
For $0<\nu<\sigma_x^2$, denote by $\Omega_{33}(\nu)$ the event that
\begin{align}
    \Tr(M_d) \geq  (1-\nu)n\sum_{j>d} \lambda_j^2.  \notag
\end{align}
Under the sub-gaussianity of $\Sigma^{-\frac{1}{2}}x_i$, $\P(\Omega_{31})\geq 1-2 \exp\{-\frac{n}{C_0}\}$, $\P(\Omega_{32})\geq 1-2 \exp\{-\frac{n}{C_0}\}$ and $\P(\Omega_{33}(\nu))\geq 1-2\exp\{-\frac{\nu^2 n}
{C_0\sigma_x^4}\}$.

\item[(iv)] [Bounding the $\Vert X_{(d+1):p}\theta^*_{(d+1):p}\Vert^2$]
Denote by $\Omega_4$ the event that
\begin{align}
    \Vert X_{(d+1):p}\theta^*_{(d+1):p}\Vert^2 \leq  (1+\sigma_x^2)n\Vert\theta^*_{(d+1):p}\Vert_{\Sigma_{(d+1):p}}^2  . \label{lem:bq_10a}
\end{align}
Under the sub-gaussianity of $\Sigma^{-\frac{1}{2}}x_i$, $\P(\Omega_4)\geq 1-2 \exp\{-\frac{n}{C_0}\}$.

\item[(v)] [Bounding the $\mu_1(X_{(d+1):p}X_{(d+1):p}^\T )$ with $r_d(\Sigma)\leq C_1 n$]
Denote by $\Omega_5$ the event that
\begin{align}
    \mu_1(X_{(d+1):p}X_{(d+1):p}^\T )&\leq  C_0\sigma_x^2(2n\lambda_{d+1}+\sum_{j>d}\lambda_j) .\label{lem:bq_10}
\end{align}
In the event $\Omega_5$, we have
\begin{align}
\Vert \hat{M}_d \Vert&\leq   \frac{1}{n}\Vert X_{(d+1):p}X_{(d+1):p}^\T \Vert^2 \notag \\
&\leq  \frac{4C_0^2\sigma_x^4}{n}  (n\lambda_{d+1}+\sum_{j>d}\lambda_j)^2 .    \notag
\end{align}
Under the sub-gaussianity of $\Sigma^{-\frac{1}{2}}x_i$, $\P(\Omega_5)\geq 1-6 \exp\{-\frac{n}{C_0}\}$.

\item[(vi)] [Bounding the $\mu_1(X_{(d+1):p}X_{(d+1):p}^\T )$ and $\mu_n(X_{(d+1):p}X_{(d+1):p}^\T )$ with Assumption \ref{ass:4}]
For $\eta_2$ defined in Assumption \ref{ass:4} and $0< \nu < \frac{1}{2}$, denote by $\Omega_6(\nu)$ the event that
\begin{align}
  & (1-\nu-\eta_2) \sum_{j>d}\lambda_j \leq \mu_n(X_{d+1:p}X_{d+1:p}^\T ) \leq \mu_1(X_{d+1:p}X_{d+1:p}^\T )\leq (1+\nu+\eta_2) \sum_{j>d}\lambda_j .\label{lem:bq_12}
\end{align}
Under the sub-gaussianity of $\Sigma^{-\frac{1}{2}}x_i$ and Assumption \ref{ass:4}, let $\frac{\nu^2r_d(\Sigma)}{C_0^2\sigma_x^4}>1$, $\P(\Omega_6(\nu))\geq 1-2n\exp\{-\frac{\nu\sqrt{r_d(\Sigma)}}{C_0\sigma_x^2} \}-4 \exp\{-\frac{n}{C_0}\}$.

\item[(vii)] [Bounding of $\Tr(X_{(d+1):p}X_{(d+1):p}^\T)$]
For $0<\nu < \min\{1, \sigma_x^2\}$, denote by $\Omega_7(\nu)$ the event that
\begin{align}
 (1-\nu)  n\sum_{j>d}\lambda_j  \leq  \Tr(X_{(d+1):p}X_{(d+1):p}^\T ) \leq   (1+\nu) n\sum_{j>d}\lambda_j .\label{eq_I_3_1}
\end{align}
Under the sub-gaussianity of $\Sigma^{-1/2}x_i$, $\P(\Omega_7(\nu))\geq 1-2 exp\{-\frac{\nu^2 n}{C_0\sigma_x^4} \}$.
\end{itemize}
\end{lem}

\begin{prf}
$\newline$
(i) Given the sub-gaussianity of $\Sigma^{-\frac{1}{2}}x_i$, from Lemma \ref{lem_05}, we have with probability $1-2 \exp\{-\frac{t}{C_0^2\sigma_x^4}\}$,
\begin{align*}
(\sqrt{n}-C_0\sigma_x^2\sqrt{d}-\sqrt{t})^2\leq \mu_d( H_d H_d^\T )
\leq \mu_1( H_d H_d^\T ) \leq   (\sqrt{n}+C_0\sigma_x^2\sqrt{d}+\sqrt{t})^2,
\end{align*}
which implies that
\begin{align*}
    n (1-C_0\sigma_x^2\sqrt{\frac{d}{n}}-\sqrt{\frac{t}{n}})^2\leq \mu_d( H_d H_d^\T )
\leq \mu_1( H_d H_d^\T ) \leq   n(1+C_0\sigma_x^2\sqrt{\frac{d}{n}}+\sqrt{\frac{t}{n}})^2.
\end{align*}
Under Assumption \ref{ass:3}, with $0<\nu < \frac{1}{2}$, we have $C_0\sigma_x^2\sqrt{\frac{d}{n}} + \nu \leq \nu + \eta_1 < 1$. Taking $t=\nu^2n$, we have with probability at least $1-2 \exp\{-\frac{\nu^2n}{C_0^2\sigma_x^4}\}$,
\begin{align*}
   n (1-\nu-\eta_1)^2\leq \mu_d( H_d H_d^\T )  \leq \mu_1( H_d H_d^\T ) \leq  n(1+\nu+\eta_1)^2.
\end{align*}
In the event $\Omega_1(\nu)$ for $0<\nu<\frac{1}{2}$, we also have
\begin{align*}
\mu_d(X_{1:d}^\T X_{1:d})&= \mu_d(\Sigma_{1:d}^{1/2}(\Sigma_{1:d}^{-1/2}X_{1:d}^\T X_{1:d} \Sigma_{1:d}^{-1/2}) \Sigma_{1:d}^{1/2})   \notag \\
&\geq  \mu_d(\Sigma_{1:d}) \mu_d(\Sigma_{1:d}^{-1/2}X_{1:d}^\T X_{1:d} \Sigma_{1:d}^{-1/2}) \notag \\
&\geq \lambda_{d}n (1-\nu-\eta_1)^2,
\end{align*}
or equivalently
\begin{align*}
\mu_d(\hat{\Sigma}_{1:d}) &= \mu_d(\frac{X_{1:d}^TX_{1:d}}{n})\geq    \lambda_{d} (1-\nu-\eta_1)^2.
\end{align*}
Moreover, in the event $\Omega_1(\nu)$ for $\nu<\frac{1}{2}$,
\begin{align}
&\Vert\theta^*_{1:d}\Vert_{\hat{\Sigma}^{-1}_{1:d}}^2 = \theta^{*\T}_{1:d}  (\frac{X_{1:d}^\T X_{1:d}}{n})^{-1}\theta^*_{1:d}=\theta^{*\T}_{1:d} \Sigma_{1:d}^{-1/2} (\frac{ H_d H_d^\T }{n})^{-1} \Sigma_{1:d}^{-1/2}\theta^{*\T}_{1:d} , \notag
\end{align}
and hence
\begin{align*}
\frac{\Vert\theta^*_{1:d}\Vert^2_{\Sigma_{1:d}^{-1}}}{(1+\nu+\eta_1)^2}\leq \Vert\theta^*_{1:d}\Vert^2_{\Sigma_{1:d}^{-1}}\mu^{-1}_1(\frac{H_dH_d^T}{n}) \leq \Vert\theta^*_{1:d}\Vert_{\hat{\Sigma}^{-1}_{1:d}}^2 \leq \Vert\theta^*_{1:d}\Vert^2_{\Sigma_{1:d}^{-1}}\mu^{-1}_d(\frac{H_dH_d^T}{n}) \leq \frac{\Vert\theta^*_{1:d}\Vert^2_{\Sigma_{1:d}^{-1}}}{(1-\nu-\eta_1)^2}.
\end{align*}
(ii) Given the sub-gaussianity of $\Sigma^{-\frac{1}{2}}x_i$, from Lemma \ref{lem_04}, we have with probability at least $1-6 \exp\{-\frac{ n}{C_0}\}$,
\begin{align*}
\Vert X_{(d+1):p}\Sigma_{(d+1):p}X_{(d+1):p}^\T \Vert&\leq C_0\sigma_x^2(2n\lambda_{d+1}^2+\sum_{j>d}\lambda_j^2).
\end{align*}
(iii) Given the sub-gaussianity of $\Sigma^{-\frac{1}{2}}x_i$, from Lemma \ref{lem_06}, with probability at least $1-2 \exp\{-\frac{n}{C_0}\}$,
\begin{align*}
    \Tr(X_{1:d}\Sigma_{1:d}^{-1}X_{1:d}^\T ) &< (1+\sigma_x^2) nd , \\
    \Tr(X_{d+1:p}\Sigma_{d+1:p}X_{d+1:p})&<(1+\sigma_x^2) n\sum_{j>d}\lambda_j^2.
\end{align*}
(iv) Given the sub-gaussianity of $\Sigma^{-\frac{1}{2}}x_i$, from Lemma \ref{lem_06}, with probability at least $1-2 \exp\{-\frac{n}{C_0}\}$,
\begin{align*}
\Vert X_{(d+1):p}\theta^*_{(d+1):p}\Vert^2 \leq  (1+\sigma_x^2)n\Vert\theta^*_{(d+1):p}\Vert_{\Sigma_{(d+1):p}}^2.
\end{align*}
(v) Given the sub-gaussianity of $\Sigma^{-\frac{1}{2}}x_i$, from Lemma \ref{lem_04}, with probability at least $1-6 \exp\{-\frac{n}{C_0}\}$,
\begin{align*}
\mu_1(X_{(d+1):p}X_{(d+1):p}^\T )&\leq  C_0\sigma_x^2(2n\lambda_{d+1}+\sum_{j>d}\lambda_j).
\end{align*}
(vi) Given the sub-gaussianity of $\Sigma^{-\frac{1}{2}}x_i$, from Lemma \ref{lem_02}, with probability at least

$1-2\exp\{-\min\{r_d(\Sigma)\frac{\delta^2} {C_0^2\sigma_x^4},\sqrt{r_d(\Sigma)}\frac{\delta}{C_0\sigma_x^2} \} \}$,
\begin{align*}
 (1-\delta) \sum_{j>d}\lambda_j \leq \mu_n(\Diag(X_{d+1:p}X_{d+1:p}^\T )) \leq  \mu_1(\Diag(X_{d+1:p}X_{d+1:p}^\T )) \leq  (1+\delta) \sum_{j>d}\lambda_j.
\end{align*}
From Lemma \ref{lem_03}, we have with probability at least $1-4 \exp\{-n/C_0\}$,
\begin{align*}
   \Vert X_{d+1:p}X_{d+1:p}^\T  - \Diag(X_{d+1:p}X_{d+1:p}^\T )\Vert \leq  C_0\sigma_x^2\sqrt{4n^2\lambda^2_{d+1}+2n\sum_{j>d}\lambda_j^2}.
\end{align*}
Then with probability at least $1-2n\exp\{-\min\{r_d(\Sigma)\frac{t^2} {C_0^2\sigma_x^4},\sqrt{r_d(\Sigma)}\frac{t}{C_0\sigma_x^2} \} \}-4 exp\{-\frac{n}{C_0}\}$,
\begin{align*}
    (1-t)\sum_{j>d}\lambda_j - C_0\sigma_x^2\sqrt{4n^2\lambda_{d+1}^2 + 2n\sum_{i>d}\lambda^2_i} \leq \mu_n(X_{d+1:p}X_{d+1:p}^\T ) ~~~~~~~~~~~~~~~~\\
    \leq \mu_1(X_{d+1:p}X_{d+1:p}^\T ) \leq (1+t)\sum_{j>d}\lambda_j + C_0\sigma_x^2\sqrt{4n^2\lambda_{d+1}^2 + 2n\sum_{j>d}\lambda^2_j},
\end{align*}
which can be equivalently written as
\begin{align*}
    \sum_{j>d}\lambda_j(1-t-C_0\sigma_x^2\sqrt{\frac{4n^2}{r^2_d(\Sigma)}+\frac{2n\sum_{j>d}\lambda_j^2}{(\sum_{j>d}\lambda_j)^2}}) \leq \mu_n(X_{d+1:p}X_{d+1:p}^\T ) ~~~~~~~~~~~~~~~~\\
    \leq  \mu_1(X_{d+1:p}X_{d+1:p}^\T )\leq \sum_{j>d}\lambda_j(1+t+C_0\sigma_x^2\sqrt{\frac{4n^2}{r_d(\Sigma)^2}+\frac{2n\sum_{j>d}\lambda_j^2}{(\sum_{j>d}\lambda_j)^2}}).
\end{align*}
Hence
\begin{align*}
    \sum_{j>d}\lambda_j(1-t-C_0\sigma_x^2\sqrt{\frac{4n^2}{r^2_d(\Sigma)}+\frac{2n}{r_d(\Sigma)}}) \leq \mu_n(X_{d+1:p}X_{d+1:p}^\T ) ~~~~~~~~~~~~~~~~~~~~~~~~\\
    \leq  \mu_1(X_{d+1:p}X_{d+1:p}^\T )\leq \sum_{j>d}\lambda_j(1+t+C_0\sigma_x^2\sqrt{\frac{4n^2}{r_d(\Sigma)^2}+\frac{2n}{r_d(\Sigma)}}).
\end{align*}
Under Assumption \ref{ass:4}, with $0<\nu < \frac{1}{2}$ and $\frac{\nu^2r_d(\Sigma)}{C_0^2\sigma_x^4}>1$, we have $C_0\sigma_x^2\sqrt{\frac{4n^2}{r_d(\Sigma)^2}+\frac{2n}{r_d(\Sigma)}}+\nu\leq\nu + \eta_2<1$.
Taking $t=\nu$, we have with probability at least $1-2\exp\{-\frac{\nu\sqrt{r_d(\Sigma)}}{C_0\sigma_x^2} \}-4 \exp\{-\frac{n}{C_0}\}$,
\begin{align*}
   & (1-\nu-\eta_2) \sum_{j>d}\lambda_j \leq \mu_n(X_{d+1:p}X_{d+1:p}^\T ) \leq \mu_1(X_{d+1:p}X_{d+1:p}^\T ) \leq  (1+\nu+\eta_2) \sum_{j>d}\lambda_j.
\end{align*}
(vii) From Lemma \ref{lem_06}, let $t=\frac{\nu^2 n
}{\sigma_x^4}$ for $0<\nu<\min\{\sigma_x^2,1\}$, with probability  at least $1-exp\{-\frac{\nu^2 n}{C_0\sigma_x^4}\}$,
\begin{align*}
(1-\nu)n\sum_{j>d}\lambda_j \leq \Tr(X_{d+1:p}X_{d+1:p}^\T )\leq\ (1+\nu)n\sum_{j>d}\lambda_j.
\end{align*}
\end{prf}

\subsection{Bounds of $\mu_1(A_d)$ and $\mu_n(A_d)$}

\label{bound_A_d}

We give the following lemma to control $\mu_1(A_d)$ and $\mu_n(A_d)$ under the small or moderate TER regime and the large TER regime.
\begin{lem}[Bounds of $\mu_1(A_d)$ and $\mu_n(A_d)$]
\label{lem3a}

$\newline$
(i)(Small or moderate TER regime) Given Assumption \ref{ass:3a} and in the event $\Omega_5$ defined in Lemma \ref{lem3}, we have for $\tau \geq \lambda_{d+1}$,
\begin{align}
     \mu_1(A_d) &\leq  (2C_0\sigma_x^2+1)(1+C_1)n\tau  ,  \label{lem:bq_11_00}\\
     \mu_n(A_d)&\geq n\tau .\label{lem:bq_11_0}
\end{align}

\noindent (ii)(Large TER regime) Given Assumption \ref{ass:4} and in the event $\Omega_6(\nu)$\ for $0<\nu < \frac{1}{2}$ defined in Lemma \ref{lem3}, we have for $\tau\geq 0$,
\begin{align}
\mu_1(A_d) &\leq (1+\nu+\eta_2)\sum_{j>d}\lambda_j +n\tau ,\label{lem:bq_11_01aa} \\
\mu_n(A_d)&\geq  (1-\nu-\eta_2)\sum_{j>d}\lambda_j +n\tau .\label{lem:bq_11_01a}
\end{align}
\end{lem}
\begin{prf}
$\newline$
(i) By the definition of $A_d$,
\begin{align*}
    \mu_n(A_d) &\geq n\tau.
\end{align*}
Given Assumption \ref{ass:3a} and in the event $\Omega_5$, we have
\begin{align*}
\mu_1(X_{(d+1):p}X_{(d+1):p}^\T )&\leq  C_0\sigma_x^2(2n\lambda_{d+1}+\sum_{j>d}\lambda_j).
\end{align*}
Hence
\begin{align}
    \mu_1(A_d)&\leq  2C_0\sigma_x^2(n\lambda_{d+1}+\sum_{j>d}\lambda_j) + n\tau  \notag \\
&\leq 2C_0\sigma_x^2(1+C_1)n\lambda_{d+1} + n\tau. \notag
\end{align}
Further if $\tau \geq \lambda_{d+1}$, then
\begin{align}
    \mu_1(A_d)&\leq  2C_0\sigma_x^2(1+C_1) n\lambda_{d+1} + n\tau   \label{eq:mu1-upper} \\
    &\leq  (2C_0\sigma_x^2+1) (1+C_1) n\tau.  \nonumber
\end{align}
(ii) Given Assumption \ref{ass:4} and in the event $\Omega_6(\nu)$, we have
\begin{align*}
  & (1-\nu-\eta_2) \sum_{j>d}\lambda_j \leq \mu_1(X_{d+1:p}X_{d+1:p}^\T ) \leq \mu_n(X_{d+1:p}X_{d+1:p}^\T )\leq (1+\nu+\eta_2) \sum_{j>d}\lambda_j.
\end{align*}
Then we have for $\tau\geq 0$,
\begin{align*}
   (1-\nu-\eta_2)\sum_{j>d}\lambda_j +n\tau &\leq \mu_n(A_d)\leq\mu_1(A_d) \leq (1+\nu+\eta_2)\sum_{j>d}\lambda_j +n\tau.
\end{align*}
\end{prf}

In the large TER regime, an upper bound on $\mu_1(A_d)$ similar to (\ref{lem:bq_11_01aa}) can be obtained by using (\ref{eq:mu1-upper})
and then Assumption \ref{ass:4}, $\frac{r_d(\Sigma)}{n} \ge c_x$:
\begin{align*}
 \mu_1(A_d)&\leq  2C_0\sigma_x^2(1+C_1) n\lambda_{d+1} + n\tau \\
 & \le  2C_0\sigma_x^2(1+C_1) c_x \sum_{j>d} \lambda_j + n\tau .
\end{align*}
This inequality can be used to achieve a similar purpose as (\ref{lem:bq_11_01aa}) in our proofs.

\subsection{Comparison between $\B_{\In,1}$ and $|\B_{\In,12}|$}

\label{prf_lem6a}
We give a proof of Lemma \ref{lem6a}, which is re-stated as follows.

\setcounter{lem}{4}
\begin{lem}[Comparison between $\B_{\In,1}$ and $|\B_{\In,12}|$]
$\newline$
(i) Given Assumption \ref{ass:3a} and \ref{ass1}(i) and in the event $\Omega_1(\nu)\cap \Omega_4$ for $0<\nu<\frac{1}{2}$,
 we have for $\tau \geq \lambda_{d+1}$,
\begin{align}
    \max\{1-\frac{|\B_{\In,12}|}{\B_{\In,1}},0\}  \geq \kappa_{1}(\tau),
\end{align}
where $\kappa_1(\tau)=\max\{1-(\frac{2 C_0\sigma_x^2(2+C_1)\lambda_{d+1}}{\tau}(1+16(2C_0\sigma_x^2+1)(1+C_1)\frac{\sqrt{\delta_1}}{1-\sqrt{\delta_1}}) + 64\frac{\sqrt{\delta_1}}{1-\sqrt{\delta_1}}),0\}$.

$\newline$
(ii) Given Assumption \ref{ass:4} and \ref{ass1}(ii) and in the event $\Omega_1(\nu)\cap \Omega_4\cap \Omega_6(\nu)$ for $0<\nu<\frac{1}{4}$, we have for $\tau\geq 0$,
\begin{align}
    \max\{1-\frac{|\B_{\In,12}|}{\B_{\In,1}},0\} \geq \kappa_{2}(\tau),\label{eq_I_2_4_2}
\end{align}
where $\kappa_2(\tau)=\max\{1-(16\frac{\lambda_{d+1}\frac{r_d(\Sigma)}{n}}{\tau+\lambda_{d+1}\frac{r_d(\Sigma)}{n}}(1+112\frac{\sqrt{\delta_2}}{1-\sqrt{\delta_2}})+64\frac{\sqrt{\delta_2}}{1-\sqrt{\delta_2}}),0\}$.
\end{lem}

We first give an algebraic bound of $\frac{|\B_{\In,12}|}{\B_{\In,1}}$.

\setcounter{lem}{7}
\begin{lem}[Algebraic bound of $\frac{|\B_{\In,12}|}{\B_{\In,1}}$]
\label{lem_009}
Given $\frac{\mu_1(X_{(d+1):p}X_{(d+1):p}^\T )}{\mu_n(A_d)}\leq 1$, we have
\begin{align*}
    \frac{|\B_{\In,12}|}{\B_{\In,1}}&\leq \frac{2\mu_1(X_{(d+1):p}X_{(d+1):p}^\T )}{\mu_n(A_d)} \frac{\Vert\theta^*_{1:d} -X_{1:d}^\T A^{-1}X_{1:d}\theta_{1:d}^*\Vert_{\hat{\Sigma}_{1:d}}}{\Vert\theta^{*}_{1:d}-X^\T_{1:d}A^{-1}X\theta^*
    \Vert_{\hat{\Sigma}_{1:d}}} \notag \\
    &~~~~~~+4 \frac{\Vert\theta^*_{(d+1):p}\Vert_{\hat{\Sigma}_{(d+1):p}}}{\Vert\theta^{*}_{1:d}-X^\T_{1:d}A^{-1}X\theta^*\Vert_{\hat{\Sigma}_{1:d}}}.
\end{align*}
\end{lem}
\begin{prf}
Given $\frac{\mu_1(X_{(d+1):p}X_{(d+1):p}^\T )}{\mu_n(A_d)}\leq 1$, we have from the Cauchy--Schwartz inequality and simple manipulations,
\begin{align}
    &\Vert\theta_{(d+1):p}^* - X_{d+1:p}^\T A^{-1}X_{(d+1):p}\theta_{(d+1):p}^*\Vert^2_{\hat{\Sigma}_{(d+1):p}} \notag \\
    &\leq 2\Vert\theta_{(d+1):p}^*\Vert_{\hat{\Sigma}_{(d+1):p}}^2  + 2\Vert X_{d+1:p}^\T A^{-1}X_{(d+1):p}\theta_{(d+1):p}^*\Vert^2_{\hat{\Sigma}_{(d+1):p}} \notag \\
    &= 2\Vert\theta_{(d+1):p}^*\Vert^2_{\hat{\Sigma}_{(d+1):p}}  + 2\Vert \frac{1}{n} X_{(d+1):p}X_{(d+1):p}^\T A^{-1}X_{(d+1):p}\theta_{(d+1):p}^*\Vert^2\notag \\
    &\leq  2\Vert\theta_{(d+1):p}^*\Vert^2_{\hat{\Sigma}_{(d+1):p}}  + 2\mu_1(X_{(d+1):p}X_{d+1:p}^\T A^{-1})\Vert \frac{1}{n} X_{(d+1):p}\theta_{(d+1):p}^*\Vert^2 \notag \\
    &\leq  (2 + 2\mu_1(X_{(d+1):p}X_{d+1:p}^\T A^{-1}))  \Vert \theta_{(d+1):p}^*\Vert^2_{\hat{\Sigma}_{(d+1):p}}  \notag \\
    &\leq (2 + 2\frac{\mu_1(X_{(d+1):p}X_{d+1:p}^\T )}{\mu_n(A_d)})  \Vert \theta_{(d+1):p}^*\Vert^2_{\hat{\Sigma}_{(d+1):p}} \notag \\
    &\leq 4 \Vert \theta_{(d+1):p}^*\Vert^2_{\hat{\Sigma}_{(d+1):p}}  .\label{eq:I_5_1}
\end{align}
Then we apply the triangle inequality to $|\B_{\In,12}|$:
\begin{align}
    |\B_{\In,12}|&=|2(\theta^{*T}_{1:d}-\theta^{*T}X^\T A^{-1}X_{1:d})\hat{\Sigma}_{1:d, (d+1):p}(\theta_{(d+1):p}^* - X_{d+1:p}^\T A^{-1}X\theta^*)|\notag\\
    &\leq  |2(\theta^{*T}_{1:d}-\theta^{*T}X^\T A^{-1}X_{1:d})\hat{\Sigma}_{1:d, (d+1):p} X_{d+1:p}^\T A^{-1}X_{1:d}\theta_{1:d}^*)|\notag \\
    &+ 2|(\theta^{*T}_{1:d}-\theta^{*T}X^\T A^{-1}X_{1:d})\hat{\Sigma}_{1:d, (d+1):p}(\theta_{(d+1):p}^* - X_{d+1:p}^\T A^{-1}X_{(d+1):p}\theta_{(d+1):p}^*)|. \notag
\end{align}
The two terms on the right-hand side of the inequality above can be bounded as follows.
First,
\begin{align}
&|2(\theta^{*T}_{1:d}-\theta^{*T}X^\T A^{-1}X_{1:d})\hat{\Sigma}_{1:d, (d+1):p} X_{d+1:p}^\T A^{-1}X_{1:d}\theta_{1:d}^*)|\notag\\
    &=|2(\theta^{*T}_{1:d}-\theta^*X^\T A^{-1}X_{1:d})\frac{X_{1:d}^\T X_{(d+1):p}}{n}  X_{d+1:p}^\T A^{-1}X_{1:d}\theta^*_{1:d}|\notag\\&=|2(\theta^{*T}_{1:d}-\theta^{*T}X^\T A^{-1}X_{1:d})\frac{X_{1:d}^\T X_{(d+1):p}}{n}  X_{d+1:p}^\T A^{-1}_dX_{1:d}(I_d+X_{1:d}^\T A^{-1}_dX_{1:d})^{-1}\theta^*_{1:d}|\notag\\
    &\quad(\text{from~Lemma~}\ref{lem_01}(ii)) \notag \\
    &=|2(\theta^{*T}_{1:d}-\theta^{*T}X^\T A^{-1}X_{1:d})\frac{X_{1:d}^\T X_{(d+1):p}}{n}  X_{d+1:p}^\T A^{-1}_dX_{1:d}(I_d -X_{1:d}^\T A^{-1}X_{1:d})\theta^*_{1:d}|\notag\\&\quad(\text{from~Lemma~}\ref{lem_01}(iii)) \notag\\
    &=|2(\theta^{*T}_{1:d}-\theta^{*T}X^\T A^{-1}X_{1:d})\frac{X_{1:d}^\T X_{(d+1):p}}{n}  X_{d+1:p}^\T A^{-1}_dX_{1:d}(\theta^*_{1:d} -X_{1:d}^\T A^{-1}X_{1:d}\theta_{1:d}^*)| \notag\\
   &\leq \frac{2\mu_1(X_{(d+1):p}X_{(d+1):p}^\T )}{\mu_n(A_d)}\Vert\theta^{*}_{1:d}- X^\T_{1:d}A^{-1}X\theta^*\Vert_{\hat{\Sigma}_{1:d}} \Vert\theta^*_{1:d} -X_{1:d}^\T A^{-1}X_{1:d}\theta_{1:d}^*\Vert_{\hat{\Sigma}_{1:d}}  ,\notag
\end{align}
Second,
\begin{align}
   &2|(\theta^{*T}_{1:d}-\theta^{*T}X^\T A^{-1}X_{1:d})\hat{\Sigma}_{1:d, (d+1):p}(\theta_{(d+1):p}^* - X_{d+1:p}^\T A^{-1}X_{(d+1):p}\theta_{(d+1):p}^*)|  \notag\\
   &\leq 2 \Vert\theta^{*}_{1:d}- X^\T_{1:d}A^{-1}X\theta^*\Vert_{\hat{\Sigma}_{1:d}} \Vert \theta_{(d+1):p}^* - X_{d+1:p}^\T A^{-1}X_{(d+1):p}\theta_{(d+1):p}^*\Vert_{\hat{\Sigma}_{(d+1):p}}\notag\\
   &\quad(\text{Cauchy--Schwartz inequality}) \notag \\
   &\leq  4 \Vert\theta^{*}_{1:d}- X^\T_{1:d}A^{-1}X\theta^*\Vert_{\hat{\Sigma}_{1:d}} \Vert\theta_{(d+1):p}^* \Vert_{\hat{\Sigma}_{(d+1):p}}.\quad\text{(From (\ref{eq:I_5_1}))} \notag
\end{align}
Combining the preceding three displays yields
\begin{align*}
    &~~~|\B_{\In,12}|\leq \frac{2\mu_1(X_{(d+1):p}X_{(d+1):p}^\T )}{\mu_n(A_d)}\Vert\theta^{*}_{1:d}- X^\T_{1:d}A^{-1}X\theta^*\Vert_{\hat{\Sigma}_{1:d}}\Vert_{\hat{\Sigma}_{1:d}} \Vert\theta^*_{1:d} -X_{1:d}^\T A^{-1}X_{1:d}\theta_{1:d}^*\Vert_{\hat{\Sigma}_{1:d}} \notag\\
    &~~~~~~+ 4 \Vert\theta^{*}_{1:d}- X^\T_{1:d}A^{-1}X\theta^*\Vert_{\hat{\Sigma}_{1:d}} \Vert\theta_{(d+1):p}^* \Vert_{\hat{\Sigma}_{(d+1):p}} ,\notag
\end{align*}
or equivalently
\begin{align*}
    \frac{|\B_{\In,12}|}{\B_{\In,1}}\leq \frac{2\mu_1(X_{(d+1):p}X_{(d+1):p}^\T )}{\mu_n(A_d)} \frac{\Vert\theta^*_{1:d} -X_{1:d}^\T A^{-1}X_{1:d}\theta_{1:d}^*\Vert_{\hat{\Sigma}_{1:d}}}{\Vert\theta^{*}_{1:d}- X^\T_{1:d}A^{-1}X\theta^*\Vert_{\hat{\Sigma}_{1:d}}}+4 \frac{\Vert\theta^*_{(d+1):p}\Vert_{\hat{\Sigma}_{(d+1):p}}}{\Vert\theta^{*}_{1:d}- X^\T_{1:d}A^{-1}X\theta^*}\Vert_{\hat{\Sigma}_{1:d}}.
\end{align*}
\end{prf}

Next, it is desired to control the quantities $\frac{\Vert\theta^*_{1:d} -X_{1:d}^\T A^{-1}X_{1:d}\theta_{1:d}^*\Vert_{\hat{\Sigma}_{1:d}}}{\Vert\theta^{*}_{1:d}- X^\T_{1:d}A^{-1}X\theta^*\Vert_{\hat{\Sigma}_{1:d}}}$ and $\frac{\Vert\theta^*_{(d+1):p}\Vert_{\hat{\Sigma}_{(d+1):p}}}{\Vert\theta^{*}_{1:d}- X^\T_{1:d}A^{-1}X\theta^*}\Vert_{\hat{\Sigma}_{1:d}}$, which is stated in Lemma \ref{lem5} below.
We discuss the small or moderate TER regime and the large TER regime, respectively.

\begin{lem}[Bound of $\frac{\Vert\theta^*_{1:d} -X_{1:d}^\T A^{-1}X_{1:d}\theta_{1:d}^*\Vert_{\hat{\Sigma}_{1:d}}}{\Vert\theta^*_{1:d} -X_{1:d}^\T A^{-1}X\theta^*\Vert_{\hat{\Sigma}_{1:d}}}$ and $\frac{\Vert\theta^*_{(d+1):p}\Vert_{\hat{\Sigma}_{(d+1):p}}}{\Vert\theta^*_{1:d} -X_{1:d}^\T A^{-1}X\theta^*\Vert_{\hat{\Sigma}_{1:d}}}$]
\label{lem5}
$\newline$
(i) Given Assumption \ref{ass:3a} and \ref{ass1}(i), in the event $\Omega_1(\nu)\cap \Omega_4\cap \Omega_5$ for $0<\nu<\frac{1}{2}$, we have for $\tau\geq \lambda_{d+1}$,
\begin{align*}
     \frac{\Vert\theta^*_{1:d} -X_{1:d}^\T A^{-1}X_{1:d}\theta_{1:d}^*\Vert_{\hat{\Sigma}_{1:d}}}{\Vert\theta^*_{1:d} -X_{1:d}^\T A^{-1}X\theta^*\Vert_{\hat{\Sigma}_{1:d}}} &\leq  1+\frac{(2C_0\sigma_x^2+1)(1+C_1)}{(1-\nu-\eta_1)^2} \frac{\sqrt{\delta_1}}{1-\sqrt{\delta_1}},\\
\frac{\Vert\theta^*_{(d+1):p}\Vert_{\hat{\Sigma}_{(d+1):p}}}{\Vert\theta^*_{1:d} -X_{1:d}^\T A^{-1}X\theta^*\Vert_{\hat{\Sigma}_{1:d}}} &\leq   \frac{1}{(1-\nu-\eta_1)^2}\frac{\sqrt{\delta_1}}{1-\sqrt{\delta_1}}.
\end{align*}
$\newline$
(ii) Given Assumption \ref{ass:4} and \ref{ass1}(ii) and in the event $\Omega_1(\nu)\cap \Omega_4\cap \Omega_6(\nu)$ for $0<\nu<\frac{1}{4}$, we have for $\tau\geq 0$,
\begin{align*}
     \frac{\Vert\theta^*_{1:d} -X_{1:d}^\T A^{-1}X_{1:d}\theta_{1:d}^*\Vert_{\hat{\Sigma}_{1:d}}}{\Vert\theta^*_{1:d} -X_{1:d}^\T A^{-1}X\theta^*\Vert_{\hat{\Sigma}_{1:d}}} &\leq  1+\frac{(1+\nu+\eta_2)}{(1-\nu-\eta_2)(1-\nu-\eta_1)^2}\frac{\sqrt{\delta_2}}{1-\sqrt{\delta_2}} ,\\
\frac{\Vert\theta^*_{(d+1):p}\Vert_{\hat{\Sigma}_{(d+1):p}}}{\Vert\theta^*_{1:d} -X_{1:d}^\T A^{-1}X\theta^*\Vert_{\hat{\Sigma}_{1:d}}} &\leq  \frac{1}{(1-\nu-\eta_1)^2}\frac{\sqrt{\delta_2}}{1-\sqrt{\delta_2}}.
\end{align*}
\end{lem}
\begin{prf}
As preparation, we derive a useful identity. We have
\begin{align}
    \theta^*_{1:d} -X_{1:d}^\T A^{-1}X\theta^*
    &=\theta^*_{1:d} - X_{1:d}^\T A^{-1}X_{1:d}\theta_{1:d}^*-X_{1:d}^\T A^{-1}X_{(d+1):p}\theta_{(d+1):p}^*\notag \\
    &=(I_d+X^\T _{1:d}A_d^{-1}X_{1:d})^{-1}\theta^*_{1:d}-(I_d+X_{1:d}^\T A_d^{-1}X_{1:d})^{-1}X_{1:d}^\T A_dX_{(d+1):p}\theta_{(d+1):p}^* \notag\\
    & \quad \text{(from~Lemma~\ref{lem_01}(ii) and (iii)} )\notag \\
    &=-(I_d+X_{1:d}^\T A_d^{-1}X_{1:d})^{-1}(X_{1:d}^\T A_dX_{(d+1):p}\theta_{(d+1):p}^*-\theta^*_{1:d}) .\label{lem:prf_l6_2}
\end{align}

(i) We first discuss the bound of $\frac{\Vert\theta^*_{1:d} -X_{1:d}^\T A^{-1}X_{1:d}\theta_{1:d}^*\Vert_{\hat{\Sigma}_{1:d}}}{\Vert\theta^*_{1:d} -X_{1:d}^\T A^{-1}X\theta^*\Vert_{\hat{\Sigma}_{1:d}}}$. We have,
\begin{align}
&\frac{\Vert\theta^*_{1:d} -X_{1:d}^\T A^{-1}X_{1:d}\theta_{1:d}^*\Vert_{\hat{\Sigma}_{1:d}}}{\Vert\theta^*_{1:d} -X_{1:d}^\T A^{-1}X\theta^*\Vert_{\hat{\Sigma}_{1:d}}}\notag\\
    &\leq 1 + \frac{\Vert X_{1:d}^\T A^{-1}X\theta^*-X_{1:d}^\T A^{-1}X_{1:d}\theta^*_{1:d}\Vert_{\hat{\Sigma}_{1:d}}}{\Vert(I_d+X_{1:d}^\T A_d^{-1}X_{1:d})^{-1}(X_{1:d}^\T A_dX_{(d+1):p}\theta_{(d+1):p}^*-\theta^*_{1:d})\Vert_{\hat{\Sigma}_{1:d}}} \quad(\text{Using~(\ref{lem:prf_l6_2})}) \notag \\
    &=1+\frac{\Vert X_{1:d}^\T A^{-1}X_{(d+1):p}\theta^*_{(d+1):p}\Vert_{\hat{\Sigma}_{1:d}}}{\Vert\hat{\Sigma}_{1:d}^{1/2}(I_{1:d}+X_{1:d}^\T A_d^{-1}X_{1:d})^{-1}(X_{1:d}^\T A_d^{-1}X_{(d+1):p}\theta^*_{(d+1):p}-\theta_{1:d}^*)\Vert}\notag
    \\&=1+\frac{\Vert\hat{\Sigma}_{1:d}^{1/2}(I_{1:d}+X_{1:d}^\T A_d^{-1}X_{1:d})^{-1}X_{1:d}^\T A_d^{-1}X_{(d+1):p}\theta^*_{(d+1):p}\Vert}{\Vert\hat{\Sigma}_{1:d}^{1/2}(I_{1:d}+X_{1:d}^\T A_d^{-1}X_{1:d})^{-1}(X_{1:d}^\T A_d^{-1}X_{(d+1):p}\theta^*_{(d+1):p}-\theta_{1:d}^*)\Vert} \quad(\text{from~Lemma~}\ref{lem_01}(ii)) \notag \\
    &=1+\frac{\Vert(\hat{\Sigma}_{1:d}^{-1}+  \hat{H}_dA_d^{-1} \hat{H}_d^\T )^{-1}  \hat{H}_dA_d^{-1}X_{(d+1):p}\theta^*_{(d+1):p}\Vert}{\Vert(\hat{\Sigma}_{1:d}^{-1}+  \hat{H}_dA_d^{-1} \hat{H}_d^\T )^{-1}( \hat{H}_dTA_d^{-1}X_{(d+1):p}\theta^*_{(d+1):p}-\hat{\Sigma}_{1:d}^{-1/2}\theta_{1:d}^*)\Vert} \notag \\
    &\leq 1+\frac{\mu_1(\hat{\Sigma}_{1:d}^{-1}+  \hat{H}_dA_d^{-1} \hat{H}_d^\T )n^{1/2}\mu_1(A_d^{-1})\Vert X_{(d+1):p}\theta^*_{(d+1):p}\Vert}{\mu_{n}(\hat{\Sigma}_{1:d}^{-1}+  \hat{H}_dA_d^{-1} \hat{H}_d^\T )\Vert  \hat{H}_dA_d^{-1}X_{(d+1):p}\theta^*_{(d+1):p}-\hat{\Sigma}_{1:d}^{-1/2}\theta_{1:d}^*\Vert}\notag \\
    &\leq  1+\frac{\mu_1(\hat{\Sigma}_{1:d}^{-1}+  \hat{H}_dA_d^{-1} \hat{H}_d^\T )n^{1/2}\mu_1(A_d^{-1})\Vert X_{(d+1):p}\theta^*_{(d+1):p}\Vert}{n\mu_{n}(A_d^{-1})\Vert  \hat{H}_dA_d^{-1}X_{(d+1):p}\theta^*_{(d+1):p}-\hat{\Sigma}_{1:d}^{-1/2}\theta_{1:d}^*\Vert} \notag  \\
    &\leq 1+ \frac{(\frac{1}{\mu_d(\hat{\Sigma}_{1:d})}+\frac{n}{\mu_n(A_d)})n^{1/2}\mu_1(A_d^{-1})\Vert X_{(d+1):p}\theta^*_{(d+1):p}\Vert}{n\mu_{n}(A_d^{-1})\Vert  \hat{H}_dA_d^{-1}X_{(d+1):p}\theta^*_{(d+1):p}-\hat{\Sigma}_{1:d}^{-1/2}\theta_{1:d}^*\Vert}  \notag\\
    & \leq 1+ \frac{(\frac{1}{\mu_d(\hat{\Sigma}_{1:d})}+\frac{n}{\mu_n(A_d)})n^{1/2}\mu_1(A_d^{-1})\Vert X_{(d+1):p}\theta^*_{(d+1):p}\Vert}{n\mu_{n}(A_d^{-1})|\Vert\theta^*_{1:d}\Vert_{\hat{\Sigma}_{1:d}^{-1}} - \Vert  \hat{H}_dA_d^{-1}X_{(d+1):p}\theta^*_{(d+1):p}\Vert |}. \label{lem:prf_l6_3}
\end{align}
In the event $\Omega_1(\nu)\cap \Omega_4$ and given $\Vert\theta^*_{1:d}\Vert_{\hat{\Sigma}_{1:d}^{-1}}>\Vert  \hat{H}_dA_d^{-1}X_{(d+1):p}\theta^*_{(d+1):p}\Vert$,
substituting (\ref{lem:bq_1})--(\ref{lem:bq_10a}) into (\ref{lem:prf_l6_3}), we have
\begin{align}
   & \frac{\Vert\theta^*_{1:d} -X_{1:d}^\T A^{-1}X_{1:d}\theta_{1:d}^*\Vert_{\hat{\Sigma}_{1:d}}}{\Vert\theta^*_{1:d} -X_{1:d}^\T A^{-1}X\theta^*\Vert_{\hat{\Sigma}_{1:d}}}\notag\\
    &\leq 1+\frac{1}{(1-\nu-\eta_1)^2}\frac{\mu_1(A_d)}{\mu_n{A_d}}\frac{(\frac{1}{\lambda_d}+\frac{n}{\mu_n(A_d)})(1+\nu+\eta_1)(1+\sigma_x^2)^{1/2}\Vert\theta^*_{(d+1):p}\Vert_{\Sigma_{(d+1):p}}}{\Vert\theta^*_{1:d}\Vert_{\Sigma_{1:d}^{-1}}}\notag\\
&~~~~~~~\times  \frac{\Vert\theta^*_{1:d}\Vert_{\hat{\Sigma}_{1:d}^{-1}}}{|\Vert\theta^*_{1:d}\Vert_{\hat{\Sigma}_{1:d}^{-1}} - \Vert  \hat{H}_dA_d^{-1}X_{(d+1):p}\theta^*_{(d+1):p}\Vert|} \notag \\
&\leq 1+\frac{1}{(1-\nu-\eta_1)^2}\frac{\mu_1(A_d)}{\mu_n{A_d}}\frac{(\frac{1}{\lambda_d}+\frac{n}{\mu_n(A_d)})(1+\nu+\eta_1)(1+\sigma_x^2)^{1/2}\Vert\theta^*_{(d+1):p}\Vert_{\Sigma_{(d+1):p}}}{\Vert\theta^*_{1:d}\Vert_{\Sigma_{1:d}^{-1}}}\notag\\
&~~~~~~~\times \frac{1}{|1-\frac{n(1+\nu+\eta_1)(1+\sigma_x^2)^{1/2}\Vert \theta^*_{(d+1):p}\Vert_{\Sigma_{(d+1):p}}}{\mu_n(A_d)\Vert\theta^*_{1:d}\Vert_{\Sigma_{1:d}^{-1}}}|}. \label{inq:prf_lemma_9_part_I}
\end{align}
Then we discuss $\frac{\Vert\theta^*_{(d+1):p}\Vert_{\hat{\Sigma}_{(d+1):p}}}{\Vert\theta^*_{1:d} -X_{1:d}^\T A^{-1}X\theta^*\Vert_{\hat{\Sigma}_{1:d}}}$. Substituting (\ref{lem:prf_l6_2}) in the denominator, we have
\begin{align}
&\frac{\Vert\theta^*_{(d+1):p}\Vert_{\hat{\Sigma}_{(d+1):p}}}{\Vert\theta^*_{1:d} -X_{1:d}^\T A^{-1}X\theta^*\Vert_{\hat{\Sigma}_{1:d}}} \notag\\
    &=\frac{\Vert\frac{X_{(d+1):p}}{n}\theta^*_{(d+1):p}\Vert}{\Vert(I_d+X_{1:d}^\T A_d^{-1}X_{1:d})^{-1}(X_{1:d}^\T A_dX_{(d+1):p}\theta_{(d+1):p}^*-\theta^*_{1:d})\Vert_{\hat{\Sigma}_{1:d}}} \notag \\
    &\leq \frac{\Vert\frac{X_{(d+1):p}}{n}\theta^*_{(d+1):p}\Vert}{\Vert\hat{\Sigma}^{1/2}(I_{1:d}+X_{1:d}^\T A_d^{-1}X_{1:d})^{-1}(X_{1:d}^\T A_d^{-1}X_{(d+1):p}\theta^*_{(d+1):p}-\theta_{1:d}^*)\Vert}\notag  \\
    &= \frac{\Vert\frac{X_{(d+1):p}}{n}\theta^*_{(d+1):p}\Vert}{\Vert(\hat{\Sigma}_{1:d}^{-1}+ \hat{H}_dTA_d^{-1}X_{1:d}\hat{\Sigma}^{-1/2})^{-1}(  \hat{H}_dA_d^{-1}X_{(d+1):p}\theta^*_{(d+1):p}-\hat{\Sigma}_{1:d}^{-1/2}\theta_{1:d}^*)\Vert} \notag \\
    &\leq   \frac{\mu_1(\hat{\Sigma}_{1:d}^{-1}+  \hat{H}_dA_d^{-1} \hat{H}_d^\T )\Vert\frac{X_{(d+1):p}}{n}\theta^*_{(d+1):p}\Vert}{\Vert  \hat{H}_dA_d^{-1}X_{(d+1):p}\theta^*_{(d+1):p}-\hat{\Sigma}_{1:d}^{-1/2}\theta_{1:d}^*)\Vert}.  \label{lem:prf_l6}
\end{align}
In the event $\Omega_1(\nu)\cap \Omega_4$ and given $\Vert\theta^*_{1:d}\Vert_{\hat{\Sigma}_{1:d}^{-1}}>\Vert  \hat{H}_dA_d^{-1}X_{(d+1):p}\theta^*_{(d+1):p}\Vert$,
substituting (\ref{lem:bq_1})--(\ref{lem:bq_10a}) into (\ref{lem:prf_l6}), we have
\begin{align}
\frac{\Vert\theta^*_{(d+1):p}\Vert_{\hat{\Sigma}_{(d+1):p}}}{\Vert\theta^*_{1:d} -X_{1:d}^\T A^{-1}X\theta^*\Vert_{\hat{\Sigma}_{1:d}}} \leq \frac{1}{(1-\nu-\eta_1)^2}\frac{(\frac{1}{\lambda_d}+\frac{n}{\mu_n(A_d) })(1+\sigma_x^2)^{1/2}\Vert\theta^*_{(d+1):p}\Vert_{\Sigma_{(d+1):p}}}{\Vert\theta^*_{1:d}\Vert_{\hat{\Sigma}^{-1}_{1:d}}}\notag\\~~~~~~~~~~~~~~~~~~~\times\frac{\Vert\theta^*_{1:d}\Vert_{\hat{\Sigma}^{-1}_{1:d}}}{\Vert \hat{H}_d^{\T}A_d^{-1}X_{(d+1):p}\theta^*_{(d+1):p}-\hat{\Sigma}^{-1/2}\theta_{1:d}^*\Vert} \notag \\
\leq \frac{1}{(1-\nu-\eta_1)^2}\frac{(\frac{1}{\lambda_d}+\frac{n}{\mu_n(A_d) })(1+\nu+\eta_1)(1+\sigma_x^2)\Vert\theta^*_{(d+1):p}\Vert_{\Sigma_{(d+1):p}}}{\Vert\theta^*_{1:d}\Vert_{\Sigma^{-1}_{1:d}}}\notag\\~~~~~~~~~~~~~~~~~~~\times \frac{1}{|1-\frac{n(1+\nu+\eta_1)(1+\sigma_x^2)^{1/2}\Vert \theta^*_{(d+1):p}\Vert_{\Sigma_{(d+1):p}}}{\mu_n(A_d)\Vert\theta^*_{1:d}\Vert_{\Sigma_{1:d}^{-1}}}|}. \label{inq:prf_lemma_9_part_II}
\end{align}

We control $\frac{\Vert\theta^*_{1:d} -X_{1:d}^\T A^{-1}X_{1:d}\theta_{1:d}^*\Vert_{\hat{\Sigma}_{1:d}}}{\Vert\theta^*_{1:d} -X_{1:d}^\T A^{-1}X\theta^*\Vert_{\hat{\Sigma}_{1:d}}}$ and $\frac{\Vert\theta^*_{(d+1):p}\Vert_{\hat{\Sigma}_{(d+1):p}}}{\Vert\theta^*_{1:d} -X_{1:d}^\T A^{-1}X\theta^*\Vert_{\hat{\Sigma}_{1:d}}}$ in small or moderate TER and large TER respectively.

In small or moderate TER regime, given Assumption \ref{ass:3a} and \ref{ass1}(i),  and in the event $\Omega_1(\nu)\cap \Omega_4\cap\Omega_5$ for $0<\nu<\frac{1}{4}$, we have for $\tau \geq \lambda_{d+1}$,
\begin{align*}
\Vert\theta^*_{1:d}\Vert_{\hat{\Sigma}_{1:d}^{-1}}>\Vert  \hat{H}_dA_d^{-1}X_{(d+1):p}\theta^*_{(d+1):p}\Vert.
\end{align*}
Substituting bounds of $\mu_1(A_d)$ and $\mu_n(A_d)$ in Lemma \ref{lem3a}(i) into (\ref{inq:prf_lemma_9_part_I}) and (\ref{inq:prf_lemma_9_part_II}), we have
\begin{align}
    &\frac{\Vert\theta^*_{1:d} -X_{1:d}^\T A^{-1}X_{1:d}\theta_{1:d}^*\Vert_{\hat{\Sigma}_{1:d}}}{\Vert\theta^*_{1:d} -X_{1:d}^\T A^{-1}X\theta^*\Vert_{\hat{\Sigma}_{1:d}}} \notag \\
    &\leq 1+\frac{(2C_0\sigma_x^2+1)(1+C_1)}{(1-\nu-\eta_1)^2} \frac{(1+\nu+\eta_1)(\frac{1}{\lambda_d}+\frac{1}{\lambda_{d+1}})(1+\sigma_x^2)^{1/2}\Vert\theta^*_{(d+1):p}\Vert_{\Sigma_{(d+1):p}}}{\Vert\theta^*_{1:d}\Vert_{\Sigma_{1:d}^{-1}}}\frac{1}{1-\sqrt{\delta_1}} \notag \\
    &\leq 1+\frac{(2C_0\sigma_x^2+1)(1+C_1)}{(1-\nu-\eta_1)^2} \frac{\sqrt{\delta_1}}{1-\sqrt{\delta_1}} .\notag
\end{align}
For $\tau\geq \lambda_{d+1}$,
\begin{align}
\frac{\Vert\theta^*_{(d+1):p}\Vert_{\hat{\Sigma}_{(d+1):p}}}{\Vert\theta^*_{1:d} -X_{1:d}^\T A^{-1}X\theta^*\Vert_{\hat{\Sigma}_{1:d}}} &\leq  \frac{1}{(1-\nu-\eta_1)^2}\frac{(\frac{1}{\lambda_d}+\frac{1}{\lambda_{d+1}})(1+\nu+\eta_1)(1+\sigma_x^2)\Vert\theta^*_{(d+1):p}\Vert_{\Sigma_{(d+1):p}}}{\Vert\theta^*_{1:d}\Vert_{\Sigma^{-1}_{1:d}}}\frac{1}{1-\sqrt{\delta_1}} \notag \\
&\leq  \frac{1}{(1-\nu-\eta_1)^2}\frac{\sqrt{\delta_1}}{1-\sqrt{\delta_1}}   .\notag
\end{align}

In  large TER regime, under Assumption \ref{ass:4} and \ref{ass1}(ii) and in the event $\Omega_1(\nu)\cap \Omega_4\cap\Omega_6(\nu)$ for $\nu<\frac{1}{4}$, we have for $ \tau\geq 0$,
\begin{align*}
\Vert\theta^*_{1:d}\Vert_{\hat{\Sigma}_{1:d}^{-1}}>\Vert  \hat{H}_dA_d^{-1}X_{(d+1):p}\theta^*_{(d+1):p}\Vert.
\end{align*}
Substituting bounds of $\mu_1(A_d)$ and $\mu_n(A_d)$ in Lemma \ref{lem3a}(ii) into (\ref{inq:prf_lemma_9_part_I}) and (\ref{inq:prf_lemma_9_part_II}), we have
\begin{align}
    &\frac{\Vert\theta^*_{1:d} -X_{1:d}^\T A^{-1}X_{1:d}\theta_{1:d}^*\Vert_{\hat{\Sigma}_{1:d}}}{\Vert\theta^*_{1:d} -X_{1:d}^\T A^{-1}X\theta^*\Vert_{\hat{\Sigma}_{1:d}}} \notag \\
    &\leq 1+\frac{(1+\nu+\eta_2)}{(1-\nu-\eta_1)^2(1-\nu-\eta_2)} \frac{(1+\nu+\eta_1)(\frac{1}{\lambda_d}+\frac{n}{(1-\nu-\eta_2)\sum_{j>d}\lambda_j})(1+\sigma_x^2)^{1/2}\Vert\theta^*_{(d+1):p}\Vert_{\Sigma_{(d+1):p}}}{\Vert\theta^*_{1:d}\Vert_{\Sigma_{1:d}^{-1}}}\frac{1}{1-\sqrt{\delta_2}} \notag \\
    &\leq \frac{(1+\nu+\eta_2)}{(1-\nu-\eta_1)^2(1-\nu-\eta_2)} \frac{\sqrt{\delta_2}}{1-\sqrt{\delta_2}} ,\notag
\end{align}
and
\begin{align}
\frac{\Vert\theta^*_{(d+1):p}\Vert_{\hat{\Sigma}_{(d+1):p}}}{\Vert\theta^*_{1:d} -X_{1:d}^\T A^{-1}X\theta^*\Vert_{\hat{\Sigma}_{1:d}}}&\leq  \frac{1}{(1-\nu-\eta_1)^2}\frac{(1+\nu+\eta_1)(\frac{1}{\lambda_d}+\frac{n}{(1-\nu-\eta_2)\sum_{j>d}\lambda_j })(1+\sigma_x^2)\Vert\theta^*_{(d+1):p}\Vert_{\Sigma_{(d+1):p}}}{\Vert\theta^*_{1:d}\Vert_{\Sigma^{-1}_{1:d}}}\notag \\&~~~~~~~~~\times \frac{1}{1-\sqrt{\delta_2}} \notag \\
    &\leq  \frac{1}{(1-\nu-\eta_1)^2}\frac{\sqrt{\delta_2}}{1-\sqrt{\delta_2}} .  \notag
\end{align}
\end{prf}

\textbf{Proof of Lemma \ref{lem6a}.} We obtain Lemma \ref{lem6a} by combining Lemma \ref{lem_009} and \ref{lem5}.
We discuss the small or moderate TER regime and the large TER regime, respectively.

\noindent \textbf{Small or moderate TER.} We substitute the bounds from Lemma \ref{lem5}(i) into Lemma \ref{lem_009}. Given Assumption \ref{ass:3a}, \ref{ass1}(i) and in the event $\Omega_1(\nu)\cap \Omega_4\cap \Omega_5$ for $0<\nu<\frac{1}{4}$, we have for $\tau\geq \lambda_{d+1}$ and $\frac{\mu_1(X_{(d+1):p}X_{(d+1):p}^\T )}{\mu_n(A_d)}\leq 1$,
\begin{align}
     \frac{|\B_{\In,12}|}{\B_{\In,1}}&\leq \frac{2\mu_1(X_{(d+1):p}X_{(d+1):p}^\T )}{\mu_n(A_d)} (1+16(2C_0\sigma_x^2+1)(1+C_1)\frac{\sqrt{\delta_1}}{1-\sqrt{\delta_1}}) + 64\frac{\sqrt{\delta_1}}{1-\sqrt{\delta_1}} .\notag
\end{align}
Under Assumption \ref{ass:3a} and in the event $\Omega_5$, from (\ref{lem:bq_10}), we have
\begin{align*}
\mu_1(X_{(d+1):p}X_{(d+1):p}^\T )&\leq  C_0\sigma_x^2(2n\lambda_{d+1}+\sum_{j>d}\lambda_j) \leq C_0\sigma_x^2(2+C_1)n\lambda_{d+1}.
\end{align*}
Note that $\mu_n(A_d)\geq n\tau$. We have for $\tau\geq \lambda_{d+1}$ and $\frac{C_0\sigma_x^2(2+C_1)n\lambda_{d+1}}{n\tau}\leq 1$,
\begin{align}
     \frac{|\B_{\In,12}|}{\B_{\In,1}}&\leq \frac{2C_0\sigma_x^2(2+C_1)\lambda_{d+1}}{\tau} (1+16(2C_0\sigma_x^2+1)(1+C_1) \frac{\sqrt{\delta_1}}{1-\sqrt{\delta_1}}) + 64\frac{\sqrt{\delta_1}}{1-\sqrt{\delta_1}} . \notag
\end{align}
Hence we have for $\tau \geq \lambda_{d+1}$,
\begin{align*}
    \frac{\max\{\B_{\In,1}-|\B_{\In,12}|,0\}}{\B_{\In,1}} \geq \kappa_1(\tau),
\end{align*}
where $\kappa_1=\max\{1-(\frac{2C_0\sigma_x^2(2+C_1)\lambda_{d+1}}{\tau} (1+16(2C_0\sigma_x^2+1)(1+C_1) \frac{\sqrt{\delta_1}}{1-\sqrt{\delta_1}}) + 64\frac{\sqrt{\delta_1}}{1-\sqrt{\delta_1}}),0\}$.

\noindent \textbf{Large TER.} We substitute the bounds from Lemma \ref{lem5}(ii) into Lemma \ref{lem_009}.
Given Assumption \ref{ass:4} and \ref{ass1}(ii) and in the event $\Omega_1(\nu)\cap \Omega_4\cap \Omega_6(\nu)$ for $0<\nu<\frac{1}{4}$, we have for $\tau\geq 0$ and $\frac{\mu_1(X_{(d+1):p}X_{(d+1):p}^\T )}{\mu_n(A_d)}\leq 1$,
\begin{align}
     \frac{|\B_{\In,12}|}{\B_{\In,1}}&\leq \frac{2\mu_1(X_{(d+1):p}X_{(d+1):p}^\T )}{\mu_n(A_d)} (1+112\frac{\sqrt{\delta_2}}{1-\sqrt{\delta_2}} + 64\frac{\sqrt{\delta_2}}{1-\sqrt{\delta_2}}) . \notag
\end{align}
In the event $\Omega_6(\nu)$ for $0<\nu<\frac{1}{4}$, from (\ref{lem:bq_12}) and (\ref{lem:bq_11_01a}), we have
\begin{align*}
    \mu_1(X_{d+1:p}X_{d+1:p}^\T )&\leq (1+\nu+\eta_2) \sum_{j>d}\lambda_j \leq 2 \sum_{j>d}\lambda_j\\
    \mu_n(A_d)&\geq (1-\nu-\eta_2) (\sum_{j>d}\lambda_j +n\tau)\geq \frac{1}{4} (\sum_{j>d}\lambda_j +n\tau).
\end{align*}
Then for $\tau\geq 0$ and $\frac{8(\frac{\sum_{j>d}\lambda_j}{n})}{\tau+\frac{\sum_{j>d}\lambda_j}{n}}\leq 1$, we have
\begin{align}
     \frac{|\B_{\In,12}|}{\B_{\In,1}}&\leq \frac{16(\frac{\sum_{j>d}\lambda_j}{n})}{\tau+\frac{\sum_{j>d}\lambda_j}{n}} (1+112\frac{\sqrt{\delta_2}}{1-\sqrt{\delta_2}} + 64\frac{\sqrt{\delta_2}}{1-\sqrt{\delta_2}}). \notag
\end{align}
Hence we have for $\tau \geq 0$,
\begin{align*}
    \frac{\max\{\B_{\In,1}-|\B_{\In,12}|,0\}}{\B_{\In,1}} \geq \kappa_2(\tau),
\end{align*}
where $\kappa_2(\tau)=\max\{1-(16\frac{\lambda_{d+1}\frac{r_d(\Sigma)}{n}}{\tau+\lambda_{d+1}\frac{r_d(\Sigma)}{n}}(1+112\frac{\sqrt{\delta_2}}{1-\sqrt{\delta_2}})+64\frac{\sqrt{\delta_2}}{1-\sqrt{\delta_2}}),0\}$.

\subsection{Useful identities and inequalities}

\begin{lem}[Identities from \cite{Tsigler_ridge_2023}]
\label{lem_01}
$\newline$
Let $\hat{\theta}(\tau,y)=X^\T (XX^\T +n \tau I_n)^{-1}y$ and $\hat{\theta}(\tau,y)^\T =[\hat{\theta}(\tau, y)_{1:d}^\T , \hat{\theta}(\lambda, y)_{d+1:p}^\T ]$. Then\\
(i)
\begin{align*}
    \hat{\theta}(\tau, y)_{1:d} + X_{1:d}^\T A_d^{-1}X_{1:d}\hat{\theta}(\tau, y)_{1:d} = X_{1:d}^\T A_d^{-1}y,
\end{align*}
(ii)
\begin{align*}
     A^{-1}X_{1:d} &= A_d^{-1} X_{1:d}(I_d+X_{1:d}^\T A_d^{-1}X_{1:d})^{-1},
\end{align*}
(iii)
\begin{align*}
    (I_d+X^\T _{1:d}A_d^{-1}X_{1:d})^{-1}  &= I_d -X_{1:d}^\T A^{-1}X_{1:d}.
\end{align*}
\end{lem}
\begin{prf}
\newline
Lemma \ref{lem_01}(i) is from Section F in \cite{Tsigler_ridge_2023}.
Lemma \ref{lem_01}(ii) is from H.2 in \cite{Tsigler_ridge_2023}.
To show Lemma \ref{lem_01}(iii), from Lemma \ref{lem_01}(ii),
\begin{align}
    A^{-1}X_{1:d} &= A_d^{-1} X_{1:d}(I_d+X_{1:d}^\T A_d^{-1}X_{1:d})^{-1}.  \notag
\end{align}
Then we have
\begin{align*}
     X_{1:d}^\T A^{-1}X_{1:d} &= X_{1:d}^\T A_d^{-1} X_{1:d}(I_d+X^\T _{1:d}A_d^{-1}X_{1:d})^{-1} \notag \\
    &= (X_{1:d}^\T A_d^{-1} X_{1:d} +I_d -I_d)(I_d+X^\T _{1:d}A_d^{-1}X_{1:d})^{-1} \notag \\
    &=I_d - (I_d+X^\T _{1:d}A_d^{-1}X_{1:d})^{-1},
\end{align*}
which gives
\begin{align*}
    (I_d+X^\T_{1:d}A_d^{-1}X_{1:d})^{-1}  &= I_d -X_{1:d}^\T A^{-1}X_{1:d}.
\end{align*}
\end{prf}

\begin{lem}[Monotoncity of variance]
\label{lem_12b}
Denote by $\V_\out(\tau)$ the $\V_\out$ in (\ref{eq:3.3}) with the ridge parameter $\tau$.
If $0\leq \tau_1\leq \tau_2$, then $\V_\out(\tau_2)\leq \V_\out(\tau_1)$.
\end{lem}
\begin{prf}
From the definition, we have
\begin{align*}
    \V_\out(\tau_1) &= \sigma^2 \Tr((n\tau_1 I_n + XX^\T)^{-1}X\Sigma X^\T(n\tau_1 I_n + XX^\T)^{-1}) \\
    &= \sigma^2 \Tr(X\Sigma X^\T(n\tau_1 I_n + XX^\T)^{-2}) \\
    &\geq  \sigma^2 \Tr(X\Sigma X^\T(n\tau_2 I_n + XX^\T)^{-2})\\
    &=\sigma^2 \Tr((n\tau_2 I_n + XX^\T)^{-1}X\Sigma X^\T(n\tau_2 I_n + XX^\T)^{-1}) \\
    &=\V_\out(\tau_2).
\end{align*}
The inequality follows because $(n\tau_1 I_n + XX^\T)^{-2}-(n\tau_2 I_n + XX^\T)^{-2}$  is semi-positive definite.
\end{prf}

\begin{lem}[Weyl's inequality]
\label{lem_S12b}
For $M,N,R\in \mathbb{C}^{n\times n}$, suppose that $M=N+R$, $R$ is Hermitian matrices, and their respective eigenvalues are ordered as follows:
\begin{align*}
    (M)&\quad \mu_1\geq\mu_2\geq ...\geq \mu_n, \\
    (N)&\quad \nu_1\geq\nu_2\geq ...\geq \nu_n, \\
    (R)&\quad \rho_1\geq\rho_2\geq...\geq \rho_n.
\end{align*}
Then for $i=1,2, \ldots,n$,
\begin{align*}
    \nu_i+\rho_n \leq \mu_i \leq \nu_i + \rho_1.
\end{align*}
\end{lem}

\begin{lem}[Ruhe's trace inequality in \cite{marshall11}]
\label{lem:ruhe_inq}
If $U$ and $V$ are $n \times n$ positive semidefinite Hermitian matrices, then
\begin{align*}
    \Tr(UV)\geq \sum_{i=1}^n \lambda_i(U)\lambda_{n-i+1}(V)
\end{align*}

\end{lem}

In Lemma \ref{lem_02}--\ref{lem_06} below, $C_0$ is an absolute constant which may vary from lemma to lemma.
For simplicity, we treat $C_0$ as a common absolute constant, by taking the maximum of such constants from the individual lemmas.

\begin{lem}[Corollary 2.8 in \cite{Zajkowski_2020}]
\label{lem_02}
Suppose that $z\in \mathbb{R}^p$, with $Cov(z, z) = I_p$, is a sub-gaussian vector with norm $\sigma_x$.
Let $x = z \Diag(\sqrt{\lambda_1},....,\sqrt{\lambda_p})$. Then
\begin{align*}
     P(|\sum_{j=1}^p x^2_{j}-\sum_{j=1}^p\lambda_j|\geq \sum_{j=1}^p\lambda_j\delta) &\leq 2\exp\{-min\{r_0(\Sigma)\frac{\delta^2} {C_0^2\sigma_x^4},\sqrt{r_0(\Sigma)}\frac{\delta}{C_0\sigma_x^2} \} \}.
\end{align*}
That is, with probability at least $1-2\exp\{-min\{r_0(\Sigma)\frac{ \delta^2} {C_0^2\sigma_x^4},\sqrt{r_0(\Sigma)}\frac{\delta}{C_0\sigma_x^2} \} \}$,
\begin{align*}
    (1-\delta)\sum_{j=1}^p\lambda_j \leq \sum_{j=1}^p x^2_{j}\leq (1+\delta)\sum_{j=1}^p\lambda_j.
\end{align*}
\end{lem}

\begin{lem}[Lemma 23 in \cite{Tsigler_ridge_2023}]
\label{lem_03}
\label{lem:benign_ridge_1}
Suppose that $z_1,\ldots,z_n$ are independent sub-gaussian vectors in $\mathbb{R}^p$, each with sub-gaussian norm $\sigma_x$. Let $\Sigma=\Diag(\lambda_1,...,\lambda_p)$ for some positive non-increasing sequence $\{\lambda_i\}_{j=1}^p$. Denote $Z$ to be the matrix with rows $\{z_i\Sigma^{1/2}\}_{i=1}^n$ and $A=ZZ^\T$. Denote also $\mathring{A}$ to be the matrix $A$ with zeroed out diagonal elements:  $\mathring{A}_{i,j}=(1-\delta_{i,j})A_{i,j}$, where $\delta_{i,j}=0$ if $i\neq j$ or $\delta_{i,j}=1$ if $i=j$.
Then for any $t>0$ with probability at least $1-4\exp\{-t/C_0\}$,
\begin{align*}
    \Vert\mathring{A}\Vert\leq C_0\sigma_x^2\sqrt{(t+n)(\lambda_1^2(t+n)+\sum_{j=1}^p\lambda_j^2)}.
\end{align*}
\end{lem}
\begin{lem}[Lemma 24 in \cite{Tsigler_ridge_2023}]
\label{lem_04}
In the same setting as Lemma \ref{lem:benign_ridge_1}, we have with probability at least $1-6\exp\{-t/C_0\}$,
\begin{align*}
    \Vert A\Vert\leq C_0 \sigma_x^2 (\lambda_1(t+n) + \sum_{j=1}^p\lambda_j).
\end{align*}
\end{lem}

\begin{lem}[Theorem 5.39 in \cite{Vershynin2010}]
\label{lem_05}
Let $d\leq n$ and $Z\in \mathbb{R}^{n\times d}$ whose rows $z_i$ are independent sub-gaussian isotropic vectors in $\mathbb{R}^d$ with sub-gaussian norm $\sigma_x$. Then for $t\geq 0$, with probability at least $1-2\exp\{-\frac{t^2}{C_0^2\sigma_x^4} \}$,
\begin{align*}
    \sqrt{n} - C_0\sigma_x^2\sqrt{d}-t\leq s_{min}(Z)\leq s_{max}(Z)\leq \sqrt{n} + C_0\sigma_x^2\sqrt{d} + t.
\end{align*}
\end{lem}

\begin{lem}[Lemma 21 in \cite{Tsigler_ridge_2023}]
\label{lem_06}
Let $Z\in \mathbb{R}^{n\times p}$ whose rows $z_i$ are independent isotropic sub-gaussian vectors in $\mathbb{R}^p$ with sub-gaussian norm $\sigma_x$. Let $\Sigma=\Diag(\lambda_1,....,\lambda_p)$ for some positive non-increasing sequence $\{\lambda_i\}_{i=1}^p$. Then for any $t\in(0,n)$ with probability at least $1-2\exp\{-t/C_0\}$,
\begin{align*}
    (n-\sqrt{nt}\sigma_x^2)\sum_{j>k}\lambda_j\leq \sum_{i=1}^n \Vert\Sigma_{d:\infty}^{1/2}Z_{i,d:\infty}\Vert^2\leq (n+\sqrt{nt}\sigma_x^2)\sum_{j>d} \lambda_j.
\end{align*}
\end{lem}

\section{Proofs of additional results in Section \ref{main_results}}

We provide proofs of Corollaries \ref{cor1}, \ref{cor2}, \ref{cor3} and \ref{cor4} which are re-stated below for convenience.

\subsection{Sufficient and necessary conditions for $\MSE_\out=O(\frac{d}{n})$ and $\MSE_\In=O(\frac{d}{n})$}

\setcounter{cor}{0}
\begin{cor}[Conditions for $\MSE_{\out} = O(\frac{d}{n})$ with small or moderate TER]
In the setting of Theorem 1, assume further
that $\sigma^2\asymp 1$ and $\Vert\theta_{1:d}^*\Vert^2_{\Sigma_{1:d}^{-1}}\lambda_d^2\asymp 1$.

(i) A sufficient condition for $\MSE_{\out} = O(\frac{d}{n})$ with a probability approaching 1 as $n\to\infty$ is that $\frac{\lambda_{d+1}}{\lambda_{d}}\lesssim  \sqrt{\frac{d}{n}}\min\{1,\sqrt{\frac{d}{r_d(\Sigma^2)}}\}$ and
the ridge parameter $\tau$ is chosen in the range $\A_0^{-1}\lambda_{d+1} \leq \tau \leq \A_0 \lambda_{d+1}$ if
 $r_d(\Sigma^2) \leq d$ or $\A_0^{-1}\lambda_{d+1}\max\{\frac{1}{c}\sqrt{\frac{r_d( \Sigma^2)}{d}},1\} \leq \tau \leq \A_0\lambda_{d}\min\{c\sqrt{\frac{d}{n}}, 1\}$ if $r_d(\Sigma^2) > d$, where $c$ is a constant satisfying $c\geq 1$ and $\frac{\lambda_{d+1}}{\lambda_{d}}\leq c  \sqrt{\frac{d}{n}}\min\{1,\sqrt{\frac{d}{r_d(\Sigma^2)}}\}$.

(ii) Suppose that $n\gg d$ and $r_d( \Sigma^2) \gg d$. Then a necessary condition for $\MSE_{\out} = O(\frac{d}{n})$ with a probability bounded away from 0
is that $\frac{\lambda_{d+1}}{\lambda_d} \lesssim \sqrt{\frac{d}{n}}\sqrt{\frac{d}{r_d(\Sigma^2)}}$ and
the ridge parameter $\tau$ is chosen in the range $\sqrt{\frac{r_d( \Sigma^2)}{d} } \lambda_{d+1} \lesssim \tau \lesssim \sqrt{\frac{d}{n}}\lambda_d$.

\noindent The sufficient and necessary conditions become matched,  $\frac{\lambda_{d+1}}{\lambda_d} \lesssim \sqrt{\frac{d}{n}}\sqrt{\frac{d}{r_d(\Sigma^2)}}$,  if in the case where $n\gg d$ and $r_d(\Sigma^2)\gg d$ in addition to the assumptions stated.
\end{cor}

\begin{prf}
From Theorem \ref{thm1}, for any $0<\epsilon<1$, the bounds in Theorem \ref{thm1}(i)(ii)(iii) hold with probability at least $1-\epsilon$ for $n\geq N$ if $N$ is large enough. From the bounds in Theorem \ref{thm1}(i)(ii)(iii), $\sigma^2\asymp 1$ and $\Vert\theta_{1:d}^*\Vert^2_{\Sigma_{1:d}^{-1}}\lambda_d^2\asymp 1$, we have
\begin{align}
&\MSE_{\out} \gtrsim \frac{d}{n}+ \frac{r_d( \Sigma^2) }{n }, \quad \text{for } \tau \leq \A_0^{-1} \lambda_{d+1} , \label{eq:s_ii_1_1_1}  \\
 & \frac{\tau^2}{\lambda_d^2} + \frac{d}{n}+\frac{\lambda_{d+1}^2}{\tau^2}\frac{r_d( \Sigma^2)}{n}  \gtrsim \MSE_{\out} \gtrsim \frac{\tau^2}{\lambda_d^2} + \frac{d}{n}+\frac{\lambda_{d+1}^2}{\tau^2}\frac{r_d( \Sigma^2)}{n},
  \quad \text{for } \A_0^{-1}\lambda_{d+1} \leq \tau \leq \A_0\lambda_d,  \label{eq:s_ii_1_1_2}\\
 & \MSE_{\out} \gtrsim 1, \quad \text{for }  \tau \geq \A_0 \lambda_d . \label{eq:s_ii_1_1_3}
\end{align}

\noindent \textcircled{1} Proof of Corollary \ref{cor1}(i):

Suppose $\frac{\lambda_{d+1}}{\lambda_d}\lesssim \sqrt{\frac{d}{n}} \min\{1, \sqrt{\frac{d}{r_d(\Sigma^2)}}\}$. Then
there exists a constant $c\geq 1$, such that
\begin{align}
    \frac{\lambda_{d+1}}{\lambda_d}\leq c \sqrt{\frac{d}{n}} \min\{1, \sqrt{\frac{d}{r_d(\Sigma^2)}}\}. \label{eq:s_ii_1_1_3a}
\end{align}
Then we prove the sufficiency of the condition, $\frac{\lambda_{d+1}}{\lambda_d}\lesssim \sqrt{\frac{d}{n}} \min\{1, \sqrt{\frac{d}{r_d(\Sigma^2)}}\}$, in two cases $d\geq r_d(\Sigma^2)$ and $d<r_d(\Sigma^2)$.
\begin{itemize}
    \item  If $d\geq r_d(\Sigma^2)$,
from (\ref{eq:s_ii_1_1_3a}),  then
\begin{align}
    \frac{\lambda_{d+1}}{\lambda_d}\leq c  \sqrt{\frac{d}{n}} . \label{eq:s_ii_1_1_4}
\end{align}
If we let $\A_0^{-1}\lambda_{d+1}\leq \tau\leq \A_0\lambda_{d+1}$, from (\ref{eq:s_ii_1_1_4}) and $d\geq r_d(\Sigma^2)$, then
\begin{align*}
   & \frac{\tau^2}{\lambda_d^2} \leq  A_0^2 \frac{\lambda^2_{d+1}}{\lambda^2_d} \leq  c^2A_0^2 \frac{d}{n} , \\
   & \frac{\lambda^2_{d+1}}{\tau^2}\frac{r_d(\Sigma^2)}{n}\leq A_0^2 \frac{d}{n}.
\end{align*}
From the upper bound in (\ref{eq:s_ii_1_1_2}), we have $\MSE_{\out}=O(\frac{d}{n})$.

\item  If $d<r_d(\Sigma^2)$, from (\ref{eq:s_ii_1_1_3a}), then
\begin{align*}
    \frac{\lambda_{d+1}}{\lambda_d} \leq c  \frac{d}{ \sqrt{n\, r_d( \Sigma^2)} }.
\end{align*}
Let $\tau$ be in the range $A_0^{-1}\lambda_{d+1}\leq\A_0^{-1}\lambda_{d+1}\max\{\frac{1}{c}\sqrt{\frac{r_d( \Sigma^2)}{d}},1\} \leq \tau \leq \A_0\lambda_{d}\min\{c\sqrt{\frac{d}{n}}, 1\}\leq A_0\lambda_d$, then
\begin{align*}
     &\frac{\tau^2}{\lambda_d^2} \leq c^2 A_0^2  \frac{d}{n}  , \\
    &\frac{r_d( \Sigma^2)}{n} \frac{\lambda_{d+1}^2}{\tau^2} \leq c^2 A_0^2\frac{d}{n}.
\end{align*}
From the upper bound in (\ref{eq:s_ii_1_1_2}), we have $\MSE_{\out}=O(\frac{d}{n})$.
\end{itemize}

In  conclusion, $\frac{\lambda_{d+1}}{\lambda_d}\lesssim \sqrt{\frac{d}{n}} \min\{1, \sqrt{\frac{d}{r_d(\Sigma^2)}}\}$ is a sufficient condition for $\MSE_{\out}=O(\frac{d}{n})$ with a probability approaching 1 as $n\rightarrow \infty$. The ridge parameter $\tau$ is chosen in the range $\A_0^{-1}\lambda_{d+1} \leq \tau \leq \A_0 \lambda_{d+1}$ if
 $r_d(\Sigma^2) \leq d$ or $\A_0^{-1}\lambda_{d+1}\max\{\frac{1}{c}\sqrt{\frac{r_d( \Sigma^2)}{d}},1\} \leq \tau \leq \A_0\lambda_{d}\min\{c\sqrt{\frac{d}{n}}, 1\}$ if $r_d(\Sigma^2) > d$, where $c$ is a constant satisfying $c\geq 1$ and $\frac{\lambda_{d+1}}{\lambda_{d}}\leq c  \sqrt{\frac{d}{n}}\min\{1,\sqrt{\frac{d}{r_d(\Sigma^2)}}\}$.

\noindent \textcircled{2} Proof of Corollary \ref{cor1}(ii):

We first show that $\MSE_{\out}=O(\frac{d}{n})$ with a probability bounded away from 0 only when $\A_0^{-1}\lambda_{d+1}\leq \tau \leq \A_0\lambda_d$ by method of exclusion.
\begin{itemize}
    \item  If $\tau \geq \A_0 \lambda_d$, then from lower bound in (\ref{eq:s_ii_1_1_3}), we have $ \MSE_\out\gtrsim 1$, which is contradictory to $\MSE_{\out}=O(\frac{d}{n})$ and $n\gg d$.
    \item  If $\tau \leq A_0^{-1}\lambda_{d+1}$, then from lower bound in (\ref{eq:s_ii_1_1_1}), we have $ \MSE_{\out} \gtrsim \frac{r_d(\Sigma^2)}{n} \gg \frac{d}{n}$ (because $r(\Sigma^2)\gg d$), which is contradictory to $\MSE_\out = O(\frac{d}{n})$.
\end{itemize}

By excluding the above two possibilities, $\MSE_{\out}=O(\frac{d}{n})$ with a probability bounded away from 0 only when $\A_0^{-1}\lambda_{d+1}\leq \tau \leq \A_0\lambda_d$. From the lower bound in (\ref{eq:s_ii_1_1_2}), we have
    \begin{align*}
        \frac{\tau^2}{\lambda_d^2}&=O(\frac{d}{n}),\\
        \frac{r_d(\Sigma^2)}{n}\frac{\lambda^2_{d+1}}{\tau^2}&=O(\frac{d}{n}).
    \end{align*}
    That is,
    \begin{align*}
        &\sqrt{\frac{r_d( \Sigma^2)}{d} } \lambda_{d+1} \lesssim  \tau \lesssim \sqrt{\frac{d}{n}}\lambda_d ,\\
        &\frac{\lambda_{d+1}}{\lambda_d} \lesssim \sqrt{\frac{d}{n}}\sqrt{\frac{d}{r_d(\Sigma^2)}} .
    \end{align*}
    Hence a necessary condition for $\MSE_{out} = O(\frac{d}{n})$ with a probability bounded away from 0
is that $\frac{\lambda_{d+1}}{\lambda_d} \lesssim \sqrt{\frac{d}{n}}\sqrt{\frac{d}{r_d(\Sigma^2)}}$ and $\sqrt{\frac{r_d( \Sigma^2)}{d} } \lambda_{d+1} \lesssim \tau \lesssim \sqrt{\frac{d}{n}}\lambda_d$.

\end{prf}

\begin{cor}[Conditions for $\MSE_{\In} = O(\frac{d}{n})$ with small or moderate TER]
In the setting of Theorem \ref{thm2}, assume further
that $\sigma^2\asymp 1$ and $\Vert\theta_{1:d}^*\Vert^2_{\Sigma_{1:d}^{-1}}\lambda_d^2\asymp 1$.

(i) A sufficient condition for $\MSE_{\In}=O(\frac{d}{n})$ with a probability approaching 1 as $n\to\infty$ is that $\frac{\lambda_{d+1}}{\lambda_{d}} \lesssim  \sqrt{\frac{d}{n}}\min\{1, \sqrt{\frac{d}{r_d(\Sigma)}}\}$ and
the ridge parameter $\tau$ is chosen in the range $\A_0^{-1}\lambda_{d+1} \leq \tau \leq \A_0 \lambda_{d+1}$ if
 $r_d(\Sigma) \leq d$ or $\A_0^{-1}\lambda_{d+1}\max\{\frac{1}{c}\sqrt{\frac{r_d( \Sigma)}{d}},1\} \leq \tau \leq \A_0\lambda_{d}\min\{c\sqrt{\frac{d}{n}}, 1\}$ if $r_d(\Sigma) > d$, where $c$ is a constant satisfying $c\geq 1$ and $\frac{\lambda_{d+1}}{\lambda_{d}}\leq c   \sqrt{\frac{d}{n}}\min\{1,\sqrt{\frac{d}{r_d(\Sigma)}}\}$.

(ii) Suppose that $n\gg d$, $ \frac{r_d(\Sigma)}{n} \sqrt{\frac{n}{d}}\gg 1$  and
$64\frac{\sqrt{\delta_1}}{1-\sqrt{\delta_1}}<1$. Then a necessary condition for $\MSE_{\In} = O(\frac{d}{n})$ with a probability bounded away from 0 is that $\frac{\lambda_{d+1}}{\lambda_{d}}\lesssim  \frac{d}{r_d(\Sigma)}$ and the ridge parameter $\tau$ is chosen in the range $\lambda_{d+1}\frac{r_d(\Sigma)}{n}\sqrt{\frac{n}{d}} \lesssim \tau \lesssim \lambda_d\sqrt{\frac{d}{n}}$.

\noindent The sufficient and necessary conditions become matched, $\frac{\lambda_{d+1}}{\lambda_{d}}\lesssim  \frac{d}{r_d(\Sigma)}$,  if in the case where $n\gg d$ and $r_d(\Sigma)\asymp n$ in addition to the assumptions stated.
\end{cor}
\begin{prf}

From Theorem \ref{thm2}, for $0<\epsilon<1$, Theorem \ref{thm2}(i)(ii)(iii) hold with probability at least $1-\epsilon$ for $n\geq N$ if $N$ is large enough. From the bounds in Theorem \ref{thm2}(i)(ii)(iii), $\sigma^2\asymp 1$ and $\Vert\theta_{1:d}^*\Vert^2_{\Sigma_{1:d}^{-1}}\lambda_d^2\asymp 1$, we have
\begin{align}
&\MSE_{\In} \gtrsim \frac{d}{n}+ \frac{r^2_d( \Sigma) }{n^2 }, \quad \text{for } \tau \leq \A_0^{-1} \lambda_{d+1} , \label{eq:s_ii_1_1_5}  \\
 & \frac{\tau^2}{\lambda_d^2} + \frac{d}{n}+\frac{\lambda_{d+1}^2}{\tau^2}\frac{r_d( \Sigma)}{n}  \gtrsim \MSE_{\In} \gtrsim \kappa_1(\tau)\frac{\tau^2}{\lambda_d^2} + \frac{d}{n}+\frac{\lambda_{d+1}^2}{\tau^2}\frac{r^2_d( \Sigma)}{n^2 }  ,
  \quad \text{for } \A_0^{-1}\lambda_{d+1} \leq \tau \leq \A_0\lambda_d,  \label{eq:s_ii_1_1_6}\\
 & \MSE_{\In} \gtrsim \kappa_1(\tau)  + \frac{\lambda_{d+1}^2}{\tau^2}\frac{r_d^2(\Sigma)}{n^2} , \quad \text{for }  \tau \geq \A_0 \lambda_d .  \label{eq:s_ii_1_1_7}
\end{align}

\noindent \textcircled{1} Proof of Corollary \ref{cor2}(i):

Suppose $\frac{\lambda_{d+1}}{\lambda_d}\lesssim  \sqrt{\frac{d}{n}} \min\{1, \sqrt{\frac{d}{r_d(\Sigma)}}\}$, then there exists a constant $c$ such that $c\geq 1$ and
\begin{align}
    \frac{\lambda_{d+1}}{\lambda_d}\leq c  \sqrt{\frac{d}{n}} \min\{1, \sqrt{\frac{d}{r_d(\Sigma)}}\}. \label{eq:s_ii_1_1_7a}
\end{align}
We prove the sufficiency of the condition, $\frac{\lambda_{d+1}}{\lambda_d}\lesssim  \sqrt{\frac{d}{n}} \min\{1, \sqrt{\frac{d}{r_d(\Sigma)}}\}$, in two cases, $d\geq r_d(\Sigma)$ and $d<r_d(\Sigma)$.
\begin{itemize}
    \item  If $d\geq r_d(\Sigma)$,
from (\ref{eq:s_ii_1_1_7a}),  we have
\begin{align}
    \frac{\lambda_{d+1}}{\lambda_d}\leq c \sqrt{\frac{d}{n}} .  \label{eq:s_ii_1_1_4a}
\end{align}
If we let $\A_0^{-1}\lambda_{d+1}\leq \tau\leq \A_0\lambda_{d+1}$, from (\ref{eq:s_ii_1_1_4a}) and $d\geq r_d(\Sigma)$, we have
\begin{align*}
   & \frac{\tau^2}{\lambda_d^2} \leq  A_0^2 \frac{\lambda^2_{d+1}}{\lambda^2_d} \leq  c^2A_0^2 \frac{d}{n} , \\
   & \frac{\lambda^2_{d+1}}{\tau^2}\frac{r_d(\Sigma)}{n}\leq  A_0^2 \frac{r_d(\Sigma)}{n} \leq A_0^2 \frac{d}{n}.
\end{align*}
Then from the upper bound in (\ref{eq:s_ii_1_1_6}), we have $\MSE_{\In}=O(\frac{d}{n})$.

\item  If $d<r_d(\Sigma)$, from (\ref{eq:s_ii_1_1_7a}), we have
\begin{align*}
    \frac{\lambda_{d+1}}{\lambda_d} \leq c \frac{d}{ \sqrt{n\, r_d( \Sigma)} }.
\end{align*}
Let $\tau$ be in the range $A_0^{-1}\lambda_{d+1}\leq\A_0^{-1}\lambda_{d+1}\max\{\frac{1}{c}\sqrt{\frac{r_d( \Sigma)}{d}},1\} \leq \tau \leq \A_0\lambda_{d}\min\{c\sqrt{\frac{d}{n}}, 1\}\leq A_0\lambda_d$, then we have
\begin{align*}
     &\frac{\tau^2}{\lambda_d^2} \leq c^2 A_0^2 \frac{d}{n} , \\
    &\frac{\lambda_{d+1}^2}{\tau^2}\frac{r_d(\Sigma)}{n}  \leq c^2 A_0^2\frac{d}{n}.
\end{align*}
From the upper bound in (\ref{eq:s_ii_1_1_6}), we have $\MSE_{\In}=O(\frac{d}{n})$.
\end{itemize}

In  conclusion, $\frac{\lambda_{d+1}}{\lambda_d}\lesssim  \sqrt{\frac{d}{n}} \min\{1, \sqrt{\frac{d}{r_d(\Sigma)}}\}$ is a sufficient condition for $\MSE_{\In}=O(\frac{d}{n})$ with a probability approaching 1 as $n\rightarrow \infty$. The ridge parameter $\tau$ is chosen in the range $\A_0^{-1}\lambda_{d+1} \leq \tau \leq \A_0 \lambda_{d+1}$ if
 $r_d(\Sigma) \leq d$ or $\A_0^{-1}\lambda_{d+1}\max\{\frac{1}{c}\sqrt{\frac{r_d( \Sigma)}{d}},1\} \leq \tau \leq \A_0\lambda_{d}\min\{c\sqrt{\frac{d}{n}}, 1\}$ if $r_d(\Sigma) > d$, where $c$ is a constant satisfying $c\geq 1$ and $\frac{\lambda_{d+1}}{\lambda_{d}}\leq c  \sqrt{\frac{d}{n}}\min\{1,\sqrt{\frac{d}{r_d(\Sigma)}}\}$.

\noindent \textcircled{2} Proof of Corollary \ref{cor2}(ii):

First, we point out that with $64\frac{\sqrt{\delta_1}}{1-\sqrt{\delta_1}}<1$, we have $\kappa_1(\tau)\gtrsim 1$ if $\frac{\lambda_{d+1}}{\tau}\ll 1$.
Then we show that $\MSE_{\In}=O(\frac{d}{n})$ with a probability bounded away from 0 only when $\A_0^{-1}\lambda_{d+1}\leq \tau \leq \A_0\lambda_d$ by method of exclusion.
\begin{itemize}
    \item  If $\tau \leq \A_0^{-1} \lambda_{d+1}$, from lower bound in (\ref{eq:s_ii_1_1_5}) and $\frac{r_d(\Sigma)}{n}\gg \sqrt{\frac{d}{n}}$, we have $ \MSE_\In \gg \frac{d}{n}$, which is contradictory to $\MSE_\In = O(\frac{d}{n}).$
    \item If $\tau \geq \A_0\lambda_d$ and $\MSE_{\In}=O(\frac{d}{n})$, from lower bound of (\ref{eq:s_ii_1_1_7}), we have
    \begin{align*}
        \frac{\lambda_{d+1}^2}{\tau^2}\frac{r_d^2(\Sigma)}{n^2}=O(\frac{d}{n}),
    \end{align*}
   hence there exists a constant $c\geq 1$ such that
   \begin{align*}
       \frac{\lambda_{d+1}^2}{\tau^2}\frac{r_d^2(\Sigma)}{n^2}\leq c\frac{d}{n}.
   \end{align*}
   With $\frac{r_d(\Sigma)}{n}\sqrt{\frac{n}{d}}\gg 1$, we have
   \begin{align*}
       \frac{\lambda_{d+1}}{\tau} \leq  \sqrt{c}\sqrt{\frac{d}{n}} \frac{n}{r_d(\Sigma)}\ll 1.
   \end{align*}
   Hence we have
   \begin{align*}
       &\kappa_1(\tau)\gtrsim 1 \\
       \Longrightarrow & \MSE_\In \gtrsim 1,\quad(\text{from~(\ref{eq:s_ii_1_1_7})})
   \end{align*}
   which is contradictory to $\MSE_\In=O(\frac{d}{n})$ if $n\gg d$.
\end{itemize}

By excluding the above two possibilities, we know that $\MSE_{\In}=O(\frac{d}{n})$ with a
probability bounded away from 0 only when $\A_0^{-1}\lambda_{d+1}\leq \tau\leq \A_0 \lambda_d$. From lower bound in (\ref{eq:s_ii_1_1_6}) and $\MSE_\In=O(\frac{d}{n})$, we have
\begin{align*}
    \frac{\lambda_{d+1}^2}{\tau^2}\frac{r^2_d( \Sigma)}{n^2} =O(\frac{d}{n}).
\end{align*}
Similarly to the derivation in the case of $\tau\geq A_0\lambda_d$, we have $\kappa_1(\tau)\gtrsim 1$. Hence $\MSE_{\In}=O(\frac{d}{n})$ only when $\A_0^{-1}\lambda_{d+1}\leq \tau\leq \A_0 \lambda_d$ and $\kappa_1(\tau)\gtrsim 1$. From the lower bound in (\ref{eq:s_ii_1_1_6}), we have  $\MSE_{\In}=O(\frac{d}{n})$ only when $\A_0^{-1}\lambda_{d+1}\leq \tau\leq \A_0 \lambda_d$ and
\begin{align*}
   & \frac{\lambda_{d+1}^2}{\tau^2}\frac{r^2_d( \Sigma)}{n^2} =O(\frac{d}{n}), \\
   & \frac{\tau^2}{\lambda_d^2} = O(\frac{d}{n}).
\end{align*}
That is,
\begin{align*}
  & A_0^{-1}\lambda_{d+1} \ll   \lambda_{d+1}\frac{r_d(\Sigma)}{n}\sqrt{\frac{n}{d}} \lesssim \tau\lesssim \lambda_d\sqrt{\frac{d}{n}} \ll A_0 \lambda_d, \\
  &\frac{\lambda_{d+1}}{\lambda_d}\lesssim \frac{d}{r_d(\Sigma)}.
\end{align*}
Hence a necessary condition for $\MSE_{\In} = O(\frac{d}{n})$ with a probability bounded away from 0 is $\frac{\lambda_{d+1}}{\lambda_{d}}\lesssim  \frac{d}{r_d(\Sigma)}$ and $\lambda_{d+1}\frac{r_d(\Sigma)}{n}\sqrt{\frac{n}{d}} \lesssim \tau \lesssim \lambda_d\sqrt{\frac{d}{n}}$.
\end{prf}

\setcounter{cor}{3}
\begin{cor}[Conditions for $\MSE_{\out} = O(\frac{d}{n})$ with large TER]
In the setting of Theorem \ref{thm3}, assume further that $\sigma^2\asymp 1$ and
$\Vert\theta_{1:d}^*\Vert^2_{\Sigma_{1:d}^{-1}}\lambda_d^2\asymp 1$.

(i) A sufficient condition for $\MSE_{\out}=O(\frac{d}{n})$ with a probability approaching 1 as $n\to\infty$ is
that $\frac{\lambda_{d+1}}{\lambda_{d}} \lesssim  \sqrt{\frac{d}{n}} \min \{ \sqrt{\frac{d}{r_d(\Sigma^2)}}, \frac{n}{r_d(\Sigma)}\}$
and the ridge parameter $\tau$ is chosen satisfying
$\tau + \lambda_{d+1}\frac{r_d(\Sigma)}{n} \lesssim \sqrt{\frac{d}{n}}\lambda_d$ if $\frac{n\sqrt{r_d(\Sigma^2)}}{\sqrt{d} r_d(\Sigma)}\leq 1$ or $\sqrt{\frac{r_d( \Sigma^2)}{d} } \lambda_{d+1} \lesssim \tau +\lambda_{d+1}\frac{r_d(\Sigma)}{n}\lesssim \sqrt{\frac{d}{n}}\lambda_{d}$ if $\frac{n\sqrt{r_d(\Sigma^2)}}{\sqrt{d} r_d(\Sigma)}> 1$.

(ii) Suppose that $n\gg d$. Then a necessary condition for $\MSE_{\out}=O(\frac{d}{n})$ with a probability bounded away from 0 is that
$\frac{\lambda_{d+1}}{\lambda_{d}} \lesssim  \sqrt{\frac{d}{n}} \min \{ \sqrt{\frac{d}{r_d(\Sigma^2)}}, \frac{n}{r_d(\Sigma)}\}$ and $\tau$ is chosen satisfying $\tau+\lambda_{d+1}\frac{r_d(\Sigma)}{n}\lesssim \sqrt{\frac{d}{n}}\lambda_d$ if $\frac{n\sqrt{r_d(\Sigma^2)}}{\sqrt{d}r_d(\Sigma)}\leq 1$ or $\sqrt{\frac{r_d( \Sigma^2)}{d} } \lambda_{d+1} \lesssim \tau +\lambda_{d+1}\frac{r_d(\Sigma)}{n}\lesssim
           \sqrt{\frac{d}{n}}\lambda_{d}$ if $\frac{n\sqrt{r_d(\Sigma^2)}}{\sqrt{d}r_d(\Sigma)}> 1$.

\noindent The sufficient and necessary conditions become matched, $\frac{\lambda_{d+1}}{\lambda_{d}} \lesssim  \sqrt{\frac{d}{n}} \min \{ \sqrt{\frac{d}{r_d(\Sigma^2)}}, \frac{n}{r_d(\Sigma)}\}$,  if in the case where $n\gg d$ in addition to the assumptions stated.
\end{cor}

\begin{prf}

From Theorem \ref{thm3}, for $0<\epsilon<1$, Theorem \ref{thm3}(i)(ii) hold with probability at least $1-\epsilon$ for $n\geq N$ if $N$ is large enough. Because $A_0$ can be any unbounded positive value in Theorem \ref{thm3}, from the bounds in Theorem \ref{thm3}(i)(ii), $\sigma^2\asymp 1$ and $\Vert \theta^*_{1:d}\Vert_{\Sigma^{-1}_{1:d}}^2\lambda_d^2\asymp 1$, we have
\begin{align}
   & \frac{(\tau+ \lambda_{d+1}\frac{r_d(\Sigma)}{n})^2}{\lambda_d^2} + \frac{d}{n}+ \frac{\lambda_{d+1}^2}{(\tau+\lambda_{d+1}\frac{r_d(\Sigma)}{n})^2}\frac{r_d(\Sigma^2)}{n} \gtrsim \MSE_\out \gtrsim \notag \\
   &~~~~~~~~~~~~\frac{(\tau+ \lambda_{d+1}\frac{r_d(\Sigma)}{n})^2}{\lambda_d^2} + \frac{d}{n}+ \frac{\lambda_{d+1}^2}{(\tau+\lambda_{d+1}\frac{r_d(\Sigma)}{n})^2}\frac{r_d(\Sigma^2)}{n}, \quad   \text{for }  \tau + \lambda_{d+1}\frac{r_d(\Sigma)}{n} \lesssim \lambda_d, \label{eq:s_ii_1_1_10}\\
   &\MSE_\out \gtrsim 1 , \quad   \text{for }  \tau + \lambda_{d+1}\frac{r_d(\Sigma)}{n} \gtrsim  \lambda_d.    \label{eq:s_ii_1_1_11}
\end{align}

\noindent \textcircled{1} Proof of Corollary \ref{cor3}(i):

We prove the sufficiency of the condition, $\frac{\lambda_{d+1}}{\lambda_{d}} \lesssim  \sqrt{\frac{d}{n}} \min \{ \sqrt{\frac{d}{r_d(\Sigma^2)}}, \frac{n}{r_d(\Sigma)}\}$, in two cases $\frac{n\sqrt{r_d(\Sigma^2)}}{\sqrt{d}r_d(\Sigma)}\leq 1$ and $\frac{n\sqrt{r_d(\Sigma^2)}}{\sqrt{d}r_d(\Sigma)} >1$.
\begin{itemize}
    \item  If $\frac{n\sqrt{r_d(\Sigma^2)}}{\sqrt{d}r_d(\Sigma)}\leq 1$, from
$\frac{\lambda_{d+1}}{\lambda_{d}} \lesssim  \sqrt{\frac{d}{n}} \min \{ \sqrt{\frac{d}{r_d(\Sigma^2)}}, \frac{n}{r_d(\Sigma)}\}$, we have
\begin{align*}
    \frac{\lambda_{d+1}}{\lambda_d}\lesssim \frac{\sqrt{nd}}{r_d(\Sigma)} .
\end{align*}
Because $\frac{\lambda_{d+1}}{\lambda_d}\lesssim \frac{\sqrt{nd}}{r_d(\Sigma)}$, we have $\lambda_{d+1}\frac{r_d(\Sigma)}{n}\lesssim \sqrt{\frac{d}{n}}\lambda_d$. Then we can choose non-negative $\tau$ such that $\tau+\lambda_{d+1}\frac{r_d(\Sigma)}{n}\lesssim \sqrt{\frac{d}{n}}\lambda_d$ and
\begin{align}
&\frac{(\tau+\lambda_{d+1}\frac{r_d(\Sigma)}{n})^2}{\lambda_d^2}\lesssim \frac{d}{n} .\label{eq:s_ii_1_1_8}
\end{align}
Moreover, with $\frac{n\sqrt{r_d(\Sigma^2)}}{\sqrt{d}r_d(\Sigma)}\leq 1$, we have $\frac{nr_d(\Sigma^2)}{r^2_d(\Sigma)}\leq \frac{d}{n}$. Hence we have for $\tau\geq 0$,
\begin{align}
    \frac{\lambda_{d+1}^2}{(\tau+\lambda_{d+1}\frac{r_d(\Sigma)}{n})^2}\frac{r_d(\Sigma^2)}{n} &\leq  \frac{n^2}{r^2_d(\Sigma)}\frac{r_d(\Sigma^2)}{n} \notag \\
    &\leq  \frac{d}{n} . \label{eq:s_ii_1_1_9}
\end{align}
Then from (\ref{eq:s_ii_1_1_8})--(\ref{eq:s_ii_1_1_9}) and the upper bound in (\ref{eq:s_ii_1_1_10}), we have $\MSE_\out = O(\frac{d}{n})$.

\item  If $\frac{n\sqrt{r_d(\Sigma^2)}}{\sqrt{d}r_d(\Sigma)} >1$, from $\frac{\lambda_{d+1}}{\lambda_{d}} \lesssim  \sqrt{\frac{d}{n}} \min \{ \sqrt{\frac{d}{r_d(\Sigma^2)}}, \frac{n}{r_d(\Sigma)}\}$, we have
\begin{align}
    \frac{\lambda_{d+1}}{\lambda_d} &\lesssim \frac{d}{ \sqrt{n\,r_d( \Sigma^2)}} . \label{eq:s_ii_1_1_12}
\end{align}
With $\frac{n\sqrt{r_d(\Sigma^2)}}{\sqrt{d}r_d(\Sigma)} >1$, we have
\begin{align}
\sqrt{r_d(\Sigma^2)}>\frac{\sqrt{d}r_d(\Sigma)}{n} .  \label{eq:s_ii_1_1_13}
\end{align}
Now we prove $\sqrt{\frac{d}{n}}\lambda_d\gtrsim \lambda_{d+1}\frac{r_d(\Sigma)}{n}$. From (\ref{eq:s_ii_1_1_12}), we have
\begin{align*}
     \sqrt{\frac{d}{n}}\lambda_d &\gtrsim \lambda_{d+1}\sqrt{\frac{r_d(\Sigma^2)}{d}}.
\end{align*}
Further with (\ref{eq:s_ii_1_1_13}), we have
\begin{align}
    \sqrt{\frac{d}{n}}\lambda_d &\gtrsim \lambda_{d+1}\sqrt{\frac{r_d(\Sigma^2)}{d}} ,\notag \\
    &> \lambda_{d+1}\frac{r_d(\Sigma)}{n}  .\label{eq:s_ii_1_1_14}
\end{align}
From (\ref{eq:s_ii_1_1_14}), we let $\tau\geq 0$ such that
\begin{align}
    \sqrt{\frac{r_d( \Sigma^2)}{d} } \lambda_{d+1} \lesssim \tau +\lambda_{d+1}\frac{r_d(\Sigma)}{n}\lesssim \sqrt{\frac{d}{n}}\lambda_{d}\leq \lambda_d  .\label{eq_I6_1}
\end{align}
Then we have
\begin{align*}
   & \frac{r_d(\Sigma^2)}{n} \frac{\lambda_{d+1}^2}{(\tau+\frac{\sum_{j>d}\lambda_j}{n})^2}\lesssim \frac{d}{n},\\
    &\frac{(\tau+ \frac{\sum_{j>d}\lambda_j}{n})^2}{\lambda_d^2}\lesssim \frac{d}{n}.
\end{align*}
From upper bound in (\ref{eq:s_ii_1_1_10}), we have $\MSE_{\out}=O(\frac{d}{n})$.

\end{itemize}

In conclusion, we have
$\MSE_\out=O(\frac{d}{n})$ with a high probability approaching 1 if $n\rightarrow \infty$ given $\frac{\lambda_{d+1}}{\lambda_{d}} \lesssim  \sqrt{\frac{d}{n}} \min \{ \sqrt{\frac{d}{r_d(\Sigma^2)}}, \frac{n}{r_d(\Sigma)}\}$. The ridge parameter $\tau$ is chosen such that
$\tau + \lambda_{d+1}\frac{r_d(\Sigma)}{n} \lesssim \sqrt{\frac{d}{n}}\lambda_d$ if $\frac{n\sqrt{r_d(\Sigma^2)}}{\sqrt{d} r_d(\Sigma)}\leq 1$ or $\sqrt{\frac{r_d( \Sigma^2)}{d} } \lambda_{d+1} \lesssim \tau +\lambda_{d+1}\frac{r_d(\Sigma)}{n}\lesssim \sqrt{\frac{d}{n}}\lambda_{d}$ if $\frac{n\sqrt{r_d(\Sigma^2)}}{\sqrt{d} r_d(\Sigma)}> 1$.

$\newline$
\textcircled{2} Proof of Corollary \ref{cor3}(ii):

We first show that $\MSE_{\out}=O(\frac{d}{n})$ with a probability bounded away from 0 only when $\tau + \lambda_{d+1}\frac{r_d(\Sigma)}{n}\lesssim \lambda_d$ by method of exclusion.

If $\tau \gtrsim  \lambda_d$, from lower bound (\ref{eq:s_ii_1_1_11}), we have $\MSE_{\out}\gtrsim 1$ , which is contradictory to $\MSE_{\out}=O(\frac{d}{n})$ and $n\gg d$.

By excluding the above possibility, $\MSE_{\out}=O(\frac{d}{n})$ with a probability bounded away from 0 only when $\tau + \lambda_{d+1}\frac{r_d(\Sigma)}{n}\lesssim \lambda_d$. Then we prove the necessity of the condition in two cases, $\frac{n\sqrt{r_d(\Sigma^2)}}{\sqrt{d}r_d(\Sigma)}\leq 1$ and $\frac{n\sqrt{r_d(\Sigma^2)}}{\sqrt{d}r_d(\Sigma)} >1$.

\begin{itemize}
    \item If $\frac{n\sqrt{r_d(\Sigma^2)}}{\sqrt{d}r_d(\Sigma)}\leq 1$, from $\tau + \lambda_{d+1}\frac{r_d(\Sigma)}{n}\lesssim \lambda_d$ and lower bound in (\ref{eq:s_ii_1_1_10}), we have
    \begin{align*}
\frac{(\tau+\lambda_{d+1}\frac{r_d(\Sigma)}{n})^2}{\lambda_d^2}&=O(\frac{d}{n}).
    \end{align*}
    Then we have
    \begin{align*}
        &\tau+\lambda_{d+1}\frac{r_d(\Sigma)}{n} \lesssim  \sqrt{\frac{d}{n}}\lambda_d ,
    \end{align*}
    and
    \begin{align}
    \frac{\lambda_{d+1}}{\lambda_d}\lesssim \sqrt{\frac{d}{n}} \frac{\lambda_{d+1}}{\tau+\lambda_{d+1}\frac{r_d(\Sigma)}{n}}  .\label{eq:cor4_1}
    \end{align}
    Obviously, we have
    \begin{align}
        \sqrt{\frac{d}{n}} \frac{\lambda_{d+1}}{\tau+\lambda_{d+1}\frac{r_d(\Sigma)}{n}} \leq  \frac{\sqrt{nd}}{r_d(\Sigma)}   .\label{eq:cor4_2}
    \end{align}
    Combining (\ref{eq:cor4_1}) and (\ref{eq:cor4_2}), we have
    \begin{align*}
      \frac{\lambda_{d+1}}{\lambda_d}\lesssim \frac{\sqrt{nd}}{r_d(\Sigma)}  .
    \end{align*}
    \item  If $\frac{n\sqrt{r_d(\Sigma^2)}}{\sqrt{d}r_d(\Sigma)}> 1$, from  $\tau + \lambda_{d+1}\frac{r_d(\Sigma)}{n}\lesssim \lambda_d$ and the lower bound in (\ref{eq:s_ii_1_1_10}), we have
\begin{align*}
    \frac{r_d(\Sigma^2)}{n} \frac{\lambda_{d+1}^2}{(\tau+\frac{\sum_{j>d}\lambda_j}{n})^2}&=O(\frac{d}{n}), \\
    \frac{(\tau+ \frac{\sum_{j>d}\lambda_j}{n})^2}{\lambda_d^2}&=O(\frac{d}{n}).
\end{align*}
Hence
\begin{align*}
    &\sqrt{\frac{r_d( \Sigma^2)}{d} } \lambda_{d+1} \lesssim \tau +\frac{\sum_{j>d}\lambda_j}{n}\lesssim
           \sqrt{\frac{d}{n}}\lambda_{d},\\
    &\frac{\lambda_{d+1}}{\lambda_d} \lesssim \frac{d}{ \sqrt{n\, r_d( \Sigma^2)}}     .
\end{align*}
\end{itemize}

In conclusion, a necessary condition for $\MSE_{\out}=O(\frac{d}{n})$ with a probability bounded away from 0 is $\frac{\lambda_{d+1}}{\lambda_{d}} \lesssim  \sqrt{\frac{d}{n}} \min \{ \sqrt{\frac{d}{r_d(\Sigma^2)}}, \frac{n}{r_d(\Sigma)}\}$ and the ridge parameter $\tau$ is chosen such that $\tau+\lambda_{d+1}\frac{r_d(\Sigma)}{n}\lesssim \sqrt{\frac{d}{n}}\lambda_d$ if $\frac{n\sqrt{r_d(\Sigma^2)}}{\sqrt{d}r_d(\Sigma)}\leq 1$ or $\sqrt{\frac{r_d( \Sigma^2)}{d} } \lambda_{d+1} \lesssim \tau +\lambda_{d+1}\frac{r_d(\Sigma)}{n}\lesssim
           \sqrt{\frac{d}{n}}\lambda_{d}$ if $\frac{n\sqrt{r_d(\Sigma^2)}}{\sqrt{d}r_d(\Sigma)}> 1$.
\end{prf}

\begin{cor}[Conditions for $\MSE_{\In} = O(\frac{d}{n})$ with large TER]
In the setting of Theorem \ref{thm4}, assume further
that $\sigma^2\asymp 1$ and $\Vert\theta_{1:d}^*\Vert^2_{\Sigma_{1:d}^{-1}}\lambda_d^2\asymp 1$.

(i) A sufficient condition for $\MSE_{\In}=O(\frac{d}{n})$ with a probability approaching to 1 as $n\to\infty$ is $\frac{\lambda_{d+1}}{\lambda_{d}}\lesssim \frac{d}{r_d(\Sigma)}$ and the ridge parameter $\tau$ is chosen such that $\lambda_{d+1}\frac{r_d(\Sigma)}{n}\sqrt{\frac{n}{d}} \lesssim \tau +\lambda_{d+1}\frac{r_d(\Sigma)}{n} \lesssim  \lambda_d\sqrt{\frac{d}{n}}$.

(ii) Suppose that $n\gg d$ and $64\frac{\sqrt{\delta_2}}{1-\sqrt{\delta_2}}<1$. Then a necessary condition for $\MSE_{\In} = O(\frac{d}{n})$ with a probability bounded away from 0 is
$\frac{\lambda_{d+1}}{\lambda_{d}}\lesssim \frac{d}{r_d(\Sigma)}$ and the ridge parameter $\tau$ is chosen in the range $\lambda_{d+1}\frac{r_d(\Sigma)}{n}\sqrt{\frac{n}{d}} \lesssim \tau+\lambda_{d+1}\frac{r_d(\Sigma)}{n}\lesssim \lambda_d\sqrt{\frac{d}{n}}$.

\noindent The sufficient and necessary conditions become matched, $\frac{\lambda_{d+1}}{\lambda_{d}}\lesssim \frac{d}{r_d(\Sigma)}$, in the case where $n\gg d$ in addition to the assumptions stated.
\end{cor}

\begin{prf}

From Theorem \ref{thm4}, for $0<\epsilon<1$, Theorem \ref{thm4}(i)(ii) hold with probability at least $1-\epsilon$ for $n\geq N$ if $N$ is large enough. Because $A_0$ can be any unbounded positive value in Theorem \ref{thm4}, from the bounds in Theorem \ref{thm4}(i)(ii), $\sigma^2\asymp 1$ and $\Vert \theta^*_{1:d}\Vert_{\Sigma^{-1}_{1:d}}^2\lambda_d^2\asymp 1$, we have
\begin{align}
   & \frac{(\tau+ \lambda_{d+1}\frac{r_d(\Sigma)}{n})^2}{\lambda_d^2} + \frac{d}{n}+ \frac{\lambda_{d+1}^2}{(\tau+\lambda_{d+1}\frac{r_d(\Sigma)}{n})^2}\frac{r^2_d(\Sigma)}{n^2}  \gtrsim \MSE_\In \gtrsim \notag \\
   &~~~~~~~~~~~~\kappa_2(\tau)\frac{(\tau+ \lambda_{d+1}\frac{r_d(\Sigma)}{n})^2}{\lambda_d^2} + \frac{d}{n}+ \frac{\lambda_{d+1}^2}{(\tau+\lambda_{d+1}\frac{r_d(\Sigma)}{n})^2}\frac{r^2_d(\Sigma)}{n^2}, \quad   \text{for }  \tau + \lambda_{d+1}\frac{r_d(\Sigma)}{n} \lesssim \lambda_d, \label{eq:s_ii_1_1_18} \\
   &\MSE_\In \gtrsim \kappa_2(\tau)+ \frac{\lambda_{d+1}^2}{(\tau+\lambda_{d+1}\frac{r_d(\Sigma)}{n})^2}\frac{r^2_d(\Sigma)}{n}, \quad   \text{for }  \tau + \lambda_{d+1}\frac{r_d(\Sigma)}{n} \gtrsim  \lambda_d.    \label{eq:s_ii_1_1_19}
\end{align}

\noindent \textcircled{1} Proof of Corollary \ref{cor4}(i):

If $\frac{\lambda_{d+1}}{\lambda_d}\lesssim \frac{d}{r_d(\Sigma)}$, we have $\lambda_d\sqrt{\frac{d}{n}}\gtrsim \lambda_{d+1}\frac{r_d(\Sigma)}{n}\sqrt{\frac{n}{d}}$. Hence we can choose $\tau\geq 0$ such that
\begin{align}
   \lambda_{d+1} \frac{r_d(\Sigma)}{n}\sqrt{\frac{n}{d}}\lesssim \tau+\lambda_{d+1}\frac{r_d(\Sigma)}{n}\lesssim \lambda_d\sqrt{\frac{d}{n}} . \label{eq_I6_2}
\end{align}
Then
\begin{align*}
\frac{(\tau+\lambda_{d+1}\frac{r_d(\Sigma)}{n})^2}{\lambda_d^2}&=O(\frac{d}{n}), \\
\frac{(\lambda_{d+1}\frac{r_d(\Sigma)}{n})^2}{(\tau+\lambda_{d+1}\frac{r_d(\Sigma)}{n})^2}&=O(\frac{d}{n}).
\end{align*}
From the upper bound in (\ref{eq:s_ii_1_1_18}), we have $\MSE_{\In}=O(\frac{d}{n})$.

\noindent \textcircled{2} Proof of Corollary \ref{cor4}(ii):

We point out that with $64\frac{\sqrt{\delta_2}}{1-\sqrt{\delta_2}}<1$, we have $\kappa_2(\tau)\gtrsim 1$ if $\frac{\lambda_{d+1}\frac{r_d(\Sigma)}{n}}{\tau+\lambda_{d+1}\frac{r_d(\Sigma)}{n}}\ll 1$. If $\MSE_{\In}=O(\frac{d}{n})$, from lower bound in (\ref{eq:s_ii_1_1_18})--(\ref{eq:s_ii_1_1_19}), we have
\begin{align*}
    \frac{\lambda_{d+1}^2}{(\tau+\lambda_{d+1}\frac{r_d(\Sigma)}{n})^2}\frac{r_d(\Sigma^2)}{n}=O(\frac{d}{n}).
\end{align*}
With $n\gg d$, we have
\begin{align*}
\tau+\lambda_{d+1}\frac{r_d(\Sigma)}{n}\gtrsim \lambda_{d+1}\frac{r_d(\Sigma)}{n}\sqrt{\frac{n}{d}}\gg \lambda_{d+1}\frac{r_d(\Sigma)}{n}.
\end{align*}
Equivalently, we have
\begin{align*}
    \frac{\lambda_{d+1}\frac{r_d(\Sigma)}{n}}{\tau+\lambda_{d+1}\frac{r_d(\Sigma)}{n}}\ll 1,
\end{align*}
then
\begin{align*}
    \kappa_2(\tau)\gtrsim 1.
\end{align*}
Hence $\MSE_\In=O(\frac{d}{n})$ only when $\kappa_2(\tau)\gtrsim 1$. Then we show that $\MSE_{\In}=O(\frac{d}{n})$ with a probability bounded away from 0 only when $\tau + \lambda_{d+1}\frac{r_d(\Sigma)}{n}\lesssim \lambda_d$ by method of exclusion.

If $\tau +\lambda_{d+1}\frac{r_d(\Sigma)}{n}\gtrsim \lambda_d$, from the lower bound in (\ref{eq:s_ii_1_1_19}), we have
\begin{align*}
    \MSE_{\In}\gtrsim \kappa_2(\tau) \gtrsim 1,
\end{align*}
which is contradictory to $\MSE_\In=O(\frac{d}{n})$ and $n\gg d$.

By excluding the above possibility, $\MSE_{\In}=O(\frac{d}{n})$ with a probability bounded away from 0 only when $\tau + \lambda_{d+1}\frac{r_d(\Sigma)}{n}\lesssim \lambda_d$. From the lower bound in (\ref{eq:s_ii_1_1_18}), $\MSE_{\In}=O(\frac{d}{n})$ only when
\begin{align*}
\frac{(\tau+\lambda_{d+1}\frac{r_d(\Sigma)}{n})^2}{\lambda_d^2}&=O(\frac{d}{n}), \\
\frac{(\lambda_{d+1}\frac{r_d(\Sigma)}{n})^2}{(\tau+\lambda_{d+1}\frac{r_d(\Sigma)}{n})^2}&=O(\frac{d}{n}).
\end{align*}
Then
\begin{align*}
&\lambda_{d+1}\frac{r_d(\Sigma)}{n}\sqrt{\frac{n}{d}} \lesssim \tau+\lambda_{d+1}\frac{r_d(\Sigma)}{n}\lesssim \lambda_d\sqrt{\frac{d}{n}}, \\
    &\frac{\lambda_{d+1}}{\lambda_d}\lesssim \frac{d}{r_d(\Sigma)}.
\end{align*}

\end{prf}

\subsection{Out-sample and in-sample errors with optimal ridge parameters}

\label{sec_II_2}
\setcounter{cor}{2}
\begin{cor}[Optimal error orders with small or moderate TER]
Suppose that Assumption \ref{ass:3}, \ref{ass:3a} and \ref{ass1}(i) are satisfied and further $\sigma^2\asymp 1$, $\Vert\theta_{1:d}^*\Vert^2_{\Sigma_{1:d}^{-1}}\lambda_d^2\asymp 1$, $r_d(\Sigma)\asymp n$,
$\lambda_d \gtrsim  \lambda_{d+1}\sqrt{\frac{n}{r_d(\Sigma^2)}}$, $\lambda_d \gg \lambda_{d+1}$, and $64\frac{\sqrt{\delta_1}}{1-\sqrt{\delta_1}}<1$. Then
\newline
(i) $\MSE_{\out}^*\asymp\max\{\frac{\lambda_{d+1}}{\lambda_d}\sqrt{\frac{r_d(\Sigma^2)}{n}}, \frac{d}{n}\}$ with a probability approaching to 1 and the optimal $\tau$ is chosen as $\tau =  \sqrt{\lambda_d\lambda_{d+1}\sqrt{\frac{r_d(\Sigma^2)}{n}}}\min\{\sqrt{cA_0^{-2}}, \frac{A_0\lambda_d}{\sqrt{\lambda_d\lambda_{d+1}\sqrt{\frac{r_d(\Sigma^2)}{n}}}}\}$ where $c$ is a constant satisfying $\lambda_{d+1}\sqrt{\frac{n}{r_d(\Sigma^2)}}\leq c \lambda_{d}$.
\newline
(ii) $\MSE^*_{\In}\asymp \max\{\frac{\lambda_{d+1}}{\lambda_d},\frac{d}{n}\}$ with a probability approaching to 1 and the optimal $\tau$ is chosen satisfying $\tau\asymp\sqrt{\lambda_{d+1}\lambda_d}$.
\newline
Therefore $\MSE_{\out}^* \lesssim \MSE_{\In}^*$ with a probability approaching to 1, by noting $r_d(\Sigma^2) \le r_d(\Sigma) \asymp n$.
\end{cor}

\begin{prf}
\newline
\textbf{Proof of $\MSE_{\out}^*$: } For $0<\epsilon < 1$, Theorem \ref{thm1}(i)(ii)(iii) hold with probability at least $1-\epsilon$ for $n\geq N$ if $N$ is large enough under Assumption \ref{ass:3}, \ref{ass:3a} and \ref{ass1}(i). From the bounds in Theorem \ref{thm1}(i)(ii)(iii), $\sigma^2\asymp 1$ and $\Vert\theta_{1:d}^*\Vert^2_{\Sigma_{1:d}^{-1}}\lambda_d^2\asymp 1$, (\ref{eq:s_ii_1_1_1})--(\ref{eq:s_ii_1_1_3}) hold. Denote $\MSE_\out$ with optimal $\tau$ chosen from $\tau \leq \A_0^{-1} \lambda_{d+1}$ as $\MSE_{1,\out}^*$, $\MSE_\out$ with optimal $\tau$ chosen from $\A_0^{-1}\lambda_{d+1} \leq \tau \leq \A_0\lambda_d$ as $\MSE_{2,\out}^*$ and $\MSE_\out$ with optimal $\tau$ chosen from $\tau\geq A_0\lambda_d$ as $\MSE_{3,\out}^*$. Then we give orders of bounds of $\MSE_{1,\out}^*$, $\MSE_{2,\out}^*$ and $\MSE_{3,\out}^*$ and give the order of $\MSE_\out^*$ by the integrating the orders of bounds of $\MSE_{1,\out}^*$, $\MSE_{2,\out}^*$ and $\MSE_{3,\out}^*$.

\begin{itemize}
    \item If $\tau \leq \A_0^{-1} \lambda_{d+1}$, from lower bound in (\ref{eq:s_ii_1_1_1}), we have $\MSE_{1,\out}^* \gtrsim (\frac{d}{n}+ \frac{r_d( \Sigma^2) }{n })\asymp \max\{\frac{d}{n},\frac{r_d(\Sigma^2)}{n}\}$. Moreover, we have $\max\{\frac{d}{n},\frac{r_d(\Sigma^2)}{n}\}\gtrsim \max\{\frac{d}{n}, \frac{\lambda_{d+1}}{\lambda_d}\sqrt{\frac{r_d(\Sigma^2)}{n}}\}$ because $\lambda_d\gtrsim \lambda_{d+1}\sqrt{\frac{n}{r_d(\Sigma^2)}}$.
    \item If $\tau \geq \A_0\lambda_d$, from lower bound in (\ref{eq:s_ii_1_1_3}), we have $\MSE_{3,\out}^* \gtrsim 1$.
    \item  If $\A_0^{-1}\lambda_{d+1} \leq \tau \leq \A_0\lambda_d$, from (\ref{eq:s_ii_1_1_2}), we have $\MSE_\out \asymp \frac{\tau^2}{\lambda_d^2}+\frac{r_d(\Sigma^2)}{n}\frac{\lambda_{d+1}^2}{\tau^2}+\frac{d}{n}$ and $\frac{\tau^2}{\lambda_d^2}+\frac{r_d(\Sigma^2)}{n}\frac{\lambda_{d+1}^2}{\tau^2}\geq \frac{\lambda_{d+1}}{\lambda_d}\sqrt{\frac{r_d(\Sigma^2)}{n}}$. Let
        \begin{align*}
        \tau =  \sqrt{\lambda_d\lambda_{d+1}\sqrt{\frac{r_d(\Sigma^2)}{n}}}\min\{\sqrt{cA_0^{-2}}, \frac{A_0\lambda_d}{\sqrt{\lambda_d\lambda_{d+1}\sqrt{\frac{r_d(\Sigma^2)}{n}}}}\},
        \end{align*}
        where $c$ is a constant satisfying $\lambda_{d+1}\sqrt{\frac{n}{r_d(\Sigma^2)}}\leq c \lambda_{d}$, then $\A_0^{-1}\lambda_{d+1} \leq \tau \leq \A_0\lambda_d$ because $\lambda_{d+1}\sqrt{\frac{n}{r_d(\Sigma^2)}}\leq c \lambda_{d}$. With the above choice of $\tau$, we have $\MSE_\out\asymp \frac{\lambda_{d+1}}{\lambda_d}\sqrt{\frac{r_d(\Sigma^2)}{n}}$. Therefore, if $\A_0^{-1}\lambda_{d+1} \leq \tau \leq \A_0\lambda_d$, we have $\MSE_{2,\out}^*\asymp \frac{d}{n} + \frac{\lambda_{d+1}}{\lambda_d}\sqrt{\frac{r_d(\Sigma^2)}{n}}\asymp \max\{\frac{d}{n}, \frac{\lambda_{d+1}}{\lambda_d}\sqrt{\frac{r_d(\Sigma^2)}{n}}\}$.
\end{itemize}
By integrating the orders of bounds of $\MSE_{1,\out}^*$, $\MSE_{2,\out}^*$ and $\MSE_{3,\out}^*$, we have $\MSE_{\out}^*\asymp \max\{\frac{d}{n}, \frac{\lambda_{d+1}}{\lambda_d}\sqrt{\frac{r_d(\Sigma^2)}{n}}\}$.
\newline\newline
\textbf{Proof of $\MSE_{\In}^*$: } For $0<\epsilon < 1$, Theorem \ref{thm2}(i)(ii)(iii) hold with probability at least $1-\epsilon$ for $n\geq N$ if $N$ is large enough under Assumption \ref{ass:3}, \ref{ass:3a} and \ref{ass1}(i). From the bounds in Theorem \ref{thm2}(i)(ii)(iii), $\sigma^2\asymp 1$ and $\Vert\theta_{1:d}^*\Vert^2_{\Sigma_{1:d}^{-1}}\lambda_d^2\asymp 1$, (\ref{eq:s_ii_1_1_5})--(\ref{eq:s_ii_1_1_7}) hold. Denote $\MSE_\In$ with optimal $\tau$ chosen from $\tau \leq \A_0^{-1} \lambda_{d+1}$ as $\MSE_{1,\In}^*$ and $\MSE_\In$ with optimal $\tau$ chosen from $\A_0^{-1}\lambda_{d+1} \leq \tau \leq \A_0\lambda_d$ as $\MSE_{2,\In}^*$ and $\MSE_\In$ with optimal $\tau$ chosen from $\tau\geq A_0\lambda_d$ as $\MSE_{3,\In}^*$. Then we give orders of bounds of $\MSE_{1,\In}^*$, $\MSE_{2,\In}^*$ and $\MSE_{3,\In}^*$ and give the order of $\MSE_\In^*$ by the integrating the orders of bounds of $\MSE_{1,\In}^*$, $\MSE_{2,\In}^*$ and $\MSE_{3,\In}^*$.
\begin{itemize}
    \item If $\tau \leq \A_0^{-1} \lambda_{d+1}$, from lower bound in (\ref{eq:s_ii_1_1_5}) and $r_d(\Sigma)\asymp n$, we have
    \begin{align*}
        \MSE_{1,\In}^* \gtrsim \frac{d}{n}+ \frac{r^2_d( \Sigma) }{n^2}\asymp 1.
    \end{align*}
    \item If $\tau\geq A_0\lambda_d$, then $\tau\gg \lambda_{d+1}$ because $\lambda_d\gg \lambda_{d+1}$. Further with $64\frac{\sqrt{\delta_1}}{1-\sqrt{\delta_1}}<1$, we have $\kappa_1(\tau)\gtrsim 1$. Then from the lower bound in (\ref{eq:s_ii_1_1_7}), we have for $\tau\geq A_0\lambda_d$,
    \begin{align*}
        \MSE_{3,\In}^* \gtrsim \kappa_1(\tau)\gtrsim 1.
    \end{align*}
    \item  If $\A_0^{-1}\lambda_{d+1} \leq \tau \leq \A_0\lambda_d$, we first show that $\MSE_{2,\In}^*\lesssim  \max\{\frac{d}{n}, \frac{\lambda_{d+1}}{\lambda_d}\}$ and then show that $\MSE_{2,\In}^*\gtrsim  \max\{\frac{d}{n}, \frac{\lambda_{d+1}}{\lambda_d}\}$. These give $\MSE_{2,\In}^*\asymp \max\{\frac{d}{n}, \frac{\lambda_{d+1}}{\lambda_d}\}$.

    We first show that $\MSE_{2,\In}^*\lesssim  \max\{\frac{d}{n}, \frac{\lambda_{d+1}}{\lambda_d}\}$. From the upper bound in (\ref{eq:s_ii_1_1_6}),
    we have $\MSE_\In \lesssim \frac{\tau^2}{\lambda_d^2}+\frac{\lambda_{d+1}^2}{\tau^2}+\frac{d}{n}$. Let $\tau\asymp\sqrt{\lambda_{d+1}\lambda_d}$, then $A_0^{-1}\lambda_{d+1}\ll \tau\ll A_0\lambda_d$ and $\frac{\tau^2}{\lambda_d^2}+\frac{\lambda_{d+1}^2}{\tau^2}=\frac{\lambda_{d+1}}{\lambda_d}$. Hence we have $\MSE_\In \lesssim \frac{\lambda_{d+1}}{\lambda_d}+\frac{d}{n}\asymp \max\{\frac{d}{n},\frac{\lambda_{d+1}}{\lambda_d}\}$.

     Next we show that $\MSE_{2,\In}^*\gtrsim  \max\{\frac{d}{n}, \frac{\lambda_{d+1}}{\lambda_d}\}$. From lower bound in (\ref{eq:s_ii_1_1_6}), we have $\MSE_{\In}\gtrsim \kappa_1(\tau)\frac{\tau^2}{\lambda_d^2}+\frac{\lambda_{d+1}^2}{\tau^2}+\frac{d}{n}$. Then we show that $\kappa_1(\tau)\frac{\tau^2}{\lambda_d^2}+\frac{\lambda_{d+1}^2}{\tau^2}\gtrsim \frac{\lambda_{d+1}}{\lambda_d}$ for $A_0^{-1}\lambda_{d+1}\leq \tau \leq A_0\lambda_d$. For $A_0^{-1}\lambda_{d+1}\leq \tau \lesssim \sqrt{\lambda_{d+1}\lambda_d}$, we have
 \begin{align*}
     \kappa_1(\tau)\frac{\tau^2}{\lambda_d^2}+\frac{\lambda_{d+1}^2}{\tau^2}\geq \frac{\lambda_{d+1}^2}{\tau^2}\geq \frac{\lambda_{d+1}}{\lambda_d}.
 \end{align*}
 For $\sqrt{\lambda_{d+1}\lambda_d}\lesssim \tau\leq A_0\lambda_d$, from $\lambda_d\gg \lambda_{d+1}$, we have $\tau\gg \lambda_{d+1}$. Further with $64\frac{\sqrt{\delta_1}}{1-\sqrt{\delta_1}}<1$, we have $\kappa_1(\tau)\gtrsim 1$, and then
 \begin{align*}
     \kappa_1(\tau)\frac{\tau^2}{\lambda_d^2}+\frac{\lambda_{d+1}^2}{\tau^2}\geq \kappa_1(\tau)\frac{\tau^2}{\lambda_d^2}\gtrsim \frac{\lambda_{d+1}}{\lambda_d}.
 \end{align*}
Hence we have $\kappa_1(\tau)\frac{\tau^2}{\lambda_d^2}+\frac{\lambda_{d+1}^2}{\tau^2}\gtrsim \frac{\lambda_{d+1}}{\lambda_d}$ for $A_0^{-1}\lambda_{d+1}\leq \tau \leq A_0\lambda_d$. That is, $\MSE_{2,\In}^*\gtrsim \frac{d}{n}+\frac{\lambda_{d+1}}{\lambda_d}\asymp \max\{\frac{d}{n}, \frac{\lambda_{d+1}}{\lambda_d}\}$.

From $\MSE_{2,\In}^*\lesssim  \max\{\frac{d}{n}, \frac{\lambda_{d+1}}{\lambda_d}\}$ and $\MSE_{2,\In}^*\gtrsim \max\{\frac{d}{n}, \frac{\lambda_{d+1}}{\lambda_d}\}$, we have $\MSE_{2,\In}^*\asymp \max\{\frac{d}{n}, \frac{\lambda_{d+1}}{\lambda_d}\}$.
\end{itemize}
By integrating the orders of bounds of $\MSE_{1,\In}^*$, $\MSE_{2,\In}^*$ and $\MSE_{3,\In}^*$, we have $\MSE_{\In}^*\asymp \max\{\frac{d}{n}, \frac{\lambda_{d+1}}{\lambda_d}\}$.
\end{prf}

\setcounter{cor}{5}
\begin{cor}[Optimal error orders with large TER]
Suppose that Assumption \ref{ass:3}, \ref{ass:4} and \ref{ass1}(ii) are satisfied, and further
$\sigma^2\asymp 1$, $\Vert\theta_{1:d}^*\Vert^2_{\Sigma_{1:d}^{-1}}\lambda_d^2\asymp 1$, $\lambda_d \gg \lambda_{d+1}\frac{r_d(\Sigma)}{n}$ and $64\frac{\sqrt{\delta_2}}{1-\sqrt{\delta_2}}<1$. Then
\newline
(i) The order of $\MSE_{\out}^*$ is $\max\{\frac{\lambda_{d+1}}{\lambda_d}\sqrt{\frac{r_d(\Sigma^2)}{n}},\frac{\lambda_{d+1}^2}{\lambda_d^2}\frac{ r_d(\Sigma)^2}{n^2},\frac{d}{n}\}$ with a probability approaching to 1 and the optimal $\tau$ is chosen as $\tau=0$ if $\frac{n\sqrt{r_d(\Sigma^2)}}{\sqrt{d}r_d(\Sigma)}\leq \frac{\lambda_{d+1}}{\lambda_d}\frac{r_d(\Sigma)}{n}$ or satisfying $\tau+\lambda_{d+1}\frac{r_d(\Sigma)}{n}\asymp\sqrt{\lambda_d\lambda_{d+1}\sqrt{\frac{r_d(\Sigma^2)}{n}}}$ if $\frac{n\sqrt{r_d(\Sigma^2)}}{\sqrt{d}r_d(\Sigma)}> \frac{\lambda_{d+1}}{\lambda_d}\frac{r_d(\Sigma)}{n}$.
\newline
(ii) The order of $\MSE^*_{\In}$ is $\max\{\frac{\lambda_{d+1}}{\lambda_d}\frac{r_d(\Sigma)}{n},\frac{d}{n}\}$ with a probability approaching to 1 and the optimal $\tau$  is chosen satisfying $\tau+\lambda_{d+1}\frac{r_d(\Sigma)}{n} \asymp \sqrt{\lambda_d\lambda_{d+1}\frac{r_d(\Sigma)}{n}}$.

\noindent Therefore $\MSE_{\out}^* \lesssim \MSE_{\In}^*$ with a probability approaching to 1 because $\frac{\lambda_{d+1}}{\lambda_d}\sqrt{\frac{r_d(\Sigma^2)}{n}}\lesssim \frac{\lambda_{d+1}}{\lambda_d}\frac{r_d(\Sigma)}{n}$ by noting $r_d(\Sigma^2)\leq r_d(\Sigma)$ and $r_d(\Sigma)\gtrsim  n$ (by Assumption \ref{ass:4}) and $\frac{\lambda^2_{d+1}}{\lambda_d^2}\frac{r_d(\Sigma)^2}{n^2}\lesssim \frac{\lambda_{d+1}}{\lambda_d}\frac{r_d(\Sigma)}{n}$ by noting $\lambda_d\gg \lambda_{d+1}\frac{r_d(\Sigma)}{n}$.
\end{cor}
\begin{prf}
\newline
\textbf{Proof of $\MSE_{\out}^*$: } For $0<\epsilon < 1$, Theorem \ref{thm3}(i)(ii) hold with probability at least $1-\epsilon$ for $n\geq N$ if $N$ is large enough under Assumption \ref{ass:3}, \ref{ass:4} and \ref{ass1}(ii). Because $A_0$ can be any unbounded positive value in Theorem \ref{thm3}, from the bounds in Theorem \ref{thm3}(i)(ii), $\sigma^2\asymp 1$ and $\Vert \theta^*_{1:d}\Vert_{\Sigma^{-1}_{1:d}}^2\lambda_d^2\asymp 1$, we have (\ref{eq:s_ii_1_1_10})--(\ref{eq:s_ii_1_1_11}) hold. Denote $\MSE_\out$ with optimal $\tau$ chosen from $\tau+\lambda_{d+1}\frac{r_d(\Sigma)}{n}\lesssim \lambda_d$ as $\MSE_{1,\out}^*$ and $\MSE_\out$ with optimal $\tau$ chosen from $\tau+\lambda_{d+1}\frac{r_d(\Sigma)}{n}\gtrsim \lambda_d$ as $\MSE_{2,\out}^*$. Then we give orders of bounds of $\MSE_{1,\out}^*$ and $\MSE_{2,\out}^*$ and give the order of $\MSE_\out^*$ by the integrating the orders of bounds of $\MSE_{1,\out}^*$ and $\MSE_{2,\out}^*$.

\begin{itemize}
    \item  If $\tau+\lambda_{d+1}\frac{r_d(\Sigma)}{n} \gtrsim \lambda_d$, from lower bound in (\ref{eq:s_ii_1_1_11}), we have $\MSE_{2,\out}^* \gtrsim 1$.
    \item If $\tau +\lambda_{d+1}\frac{r_d(\Sigma)}{n} \lesssim \lambda_d$, from (\ref{eq:s_ii_1_1_10}), we
 have $\MSE_\out\asymp \frac{(\tau+ \lambda_{d+1}\frac{r_d(\Sigma)}{n})^2}{\lambda_d^2}+ \frac{r_d(\Sigma^2)}{n} \frac{\lambda_{d+1}^2}{(\tau+\lambda_{d+1}\frac{r_d(\Sigma)}{n})^2} +\frac{d}{n}$. Then we discuss the order of $\MSE_{1,\out}^*$ in two cases, $\frac{n\sqrt{r_d(\Sigma^2)}}{\sqrt{d}r_d(\Sigma)}\leq \frac{\lambda_{d+1}}{\lambda_d}\frac{r_d(\Sigma)}{n}$ and $\frac{n\sqrt{r_d(\Sigma^2)}}{\sqrt{d}r_d(\Sigma)}> \frac{\lambda_{d+1}}{\lambda_d}\frac{r_d(\Sigma)}{n}$.

If $\frac{n\sqrt{r_d(\Sigma^2)}}{\sqrt{d}r_d(\Sigma)}\leq \frac{\lambda_{d+1}}{\lambda_d}\frac{r_d(\Sigma)}{n}$, we have
    \begin{align*}
               \frac{r_d(\Sigma^2)}{n} \frac{\lambda_{d+1}^2}{(\tau+\lambda_{d+1}\frac{r_d(\Sigma)}{n})^2} \leq  \frac{nr_d(\Sigma^2)}{r_d(\Sigma)^2} \leq \frac{\lambda_{d+1}^2}{\lambda_d^2}\frac{r_d(\Sigma)^2}{n^2} \leq  \frac{(\tau+ \lambda_{d+1}\frac{r_d(\Sigma)}{n})^2}{\lambda_d^2}.
    \end{align*}
    Hence
    \begin{align*}
        \frac{(\tau+ \lambda_{d+1}\frac{r_d(\Sigma)}{n})^2}{\lambda_d^2}+ \frac{r_d(\Sigma^2)}{n} \frac{\lambda_{d+1}^2}{(\tau+\lambda_{d+1}\frac{r_d(\Sigma)}{n})^2} &\asymp  \frac{(\tau+ \lambda_{d+1}\frac{r_d(\Sigma)}{n})^2} {\lambda_d^2} \\
        &\geq  \frac{\lambda_{d+1}^2}{\lambda_d^2}\frac{r_d(\Sigma)^2}{n^2}.
    \end{align*}
    If we let $\tau=0$, we have
    \begin{align*}
        \frac{(\tau+ \lambda_{d+1}\frac{r_d(\Sigma)}{n})^2}{\lambda_d^2}+ \frac{r_d(\Sigma^2)}{n} \frac{\lambda_{d+1}^2}{(\tau+\lambda_{d+1}\frac{r_d(\Sigma)}{n})^2} \asymp  \frac{\lambda_{d+1}^2}{\lambda_d^2}\frac{r_d(\Sigma)^2}{n^2}.
    \end{align*}

If $\frac{n\sqrt{r_d(\Sigma^2)}}{\sqrt{d}r_d(\Sigma)}\geq \frac{\lambda_{d+1}}{\lambda_d}\frac{r_d(\Sigma)}{n}$, we have
\begin{align*}
        \frac{(\tau+ \lambda_{d+1}\frac{r_d(\Sigma)}{n})^2}{\lambda_d^2}+ \frac{r_d(\Sigma^2)}{n} \frac{\lambda_{d+1}^2}{(\tau+\lambda_{d+1}\frac{r_d(\Sigma)}{n})^2}\geq \frac{\lambda_{d+1}}{\lambda_d}\sqrt{\frac{r_d(\Sigma^2)}{n}}.
        \end{align*}
        We also have
        \begin{align*}
            \frac{(\tau+ \lambda_{d+1}\frac{r_d(\Sigma)}{n})^2}{\lambda_d^2}+ \frac{r_d(\Sigma^2)}{n} \frac{\lambda_{d+1}^2}{(\tau+\lambda_{d+1}\frac{r_d(\Sigma)}{n})^2}\asymp \frac{\lambda_{d+1}}{\lambda_d}\sqrt{\frac{r_d(\Sigma^2)}{n}},
        \end{align*}
        when
        $\tau+\lambda_{d+1}\frac{r_d(\Sigma)}{n}\asymp\sqrt{\lambda_d\lambda_{d+1}\sqrt{\frac{r_d(\Sigma^2)}{n}}}\lesssim \lambda_d$.

   By combining two cases, $\frac{n\sqrt{r_d(\Sigma^2)}}{\sqrt{d}r_d(\Sigma)}\leq \frac{\lambda_{d+1}}{\lambda_d}\frac{r_d(\Sigma)}{n}$ and $\frac{n\sqrt{r_d(\Sigma^2)}}{\sqrt{d}r_d(\Sigma)}\geq \frac{\lambda_{d+1}}{\lambda_d}\frac{r_d(\Sigma)}{n}$, we have $\MSE_{1,\out}^* \asymp \max\{\frac{\lambda_{d+1}}{\lambda_d}\sqrt{\frac{r_d(\Sigma^2)}{n}},\frac{\lambda_{d+1}^2}{\lambda_d^2}\frac{ r_d(\Sigma)^2}{n^2},\frac{d}{n}\}$.
\end{itemize}
By integrating the order of bounds of $\MSE_{1,\out}^*$ and $\MSE_{2,\out}^*$, we have $\MSE_{\out}^*\asymp \max\{\frac{\lambda_{d+1}}{\lambda_d}\sqrt{\frac{r_d(\Sigma^2)}{n}},$ $
\frac{\lambda_{d+1}^2}{\lambda_d^2}\frac{r_d(\Sigma)^2}{n^2},\frac{d}{n}\}$.
\newline\newline
\textbf{Proof of $\MSE_{\In}^*$: } For $0<\epsilon < 1$, we have Theorem \ref{thm4}(i)(ii) hold with probability at least $1-\epsilon$ for $n\geq N$ if $N$ is large enough under Assumption \ref{ass:3}, \ref{ass:4} and \ref{ass1}(ii). Because $A_0$ can be any unbounded positive value in Theorem \ref{thm4}, from the bounds in Theorem \ref{thm4}(i)(ii), $\sigma^2\asymp 1$ and $\Vert \theta^*_{1:d}\Vert_{\Sigma^{-1}_{1:d}}^2\lambda_d^2\asymp 1$, we have (\ref{eq:s_ii_1_1_18})--(\ref{eq:s_ii_1_1_19}) hold. Denote $\MSE_\In$ with optimal $\tau$ chosen from $\tau+\lambda_{d+1}\frac{r_d(\Sigma)}{n}\lesssim \lambda_d$ as $\MSE_{1,\In}^*$ and $\MSE_\In$ with optimal $\tau$ chosen from $\tau+\lambda_{d+1}\frac{r_d(\Sigma)}{n}\gtrsim \lambda_d$ as $\MSE_{2,\In}^*$. Then we give orders of bounds of $\MSE_{1,\In}^*$ and $\MSE_{2,\In}^*$ and give the order of $\MSE_\In^*$ by the integrating the orders of bounds of $\MSE_{1,\In}^*$ and $\MSE_{2,\In}^*$.

\begin{itemize}
    \item If $\tau+\lambda_{d+1}\frac{r_d(\Sigma)}{n} \gtrsim \lambda_d$, then  $\tau+\lambda_{d+1}\frac{r_d(\Sigma)}{n}\gg \lambda_{d+1}\frac{r_d(\Sigma)}{n}$ because $\lambda_d\gg \lambda_{d+1}\frac{r_d(\Sigma)}{n}$. Further with $64\frac{\sqrt{\delta_2}}{1-\sqrt{\delta_2}}<1$, we have $\kappa_2(\tau)\gtrsim 1$. From the lower bound in (\ref{eq:s_ii_1_1_19}), we have $\MSE^*_{2,\In}\gtrsim \kappa_2(\tau) \gtrsim 1$.
    \item If $\tau +\lambda_{d+1}\frac{r_d(\Sigma)}{n} \lesssim \lambda_d$, we first show that $\MSE_{1,\In}^* \lesssim \max\{\frac{d}{n}+\frac{\lambda_{d+1}}{\lambda_d}\frac{r_d(\Sigma)}{n}\}$ and then show that $\MSE_{1,\In}^* \gtrsim \max\{\frac{d}{n}+\frac{\lambda_{d+1}}{\lambda_d}\frac{r_d(\Sigma)}{n}\}$, which gives the order of $\MSE_{1,\In}^*$.

    We first show that $\MSE_{1,\In}^* \lesssim \max\{\frac{d}{n}+\frac{\lambda_{d+1}}{\lambda_d}\frac{r_d(\Sigma)}{n}\}$. From the upper bound in (\ref{eq:s_ii_1_1_18}), we have
    $\MSE_\In\lesssim \frac{(\tau+\lambda_{d+1}\frac{r_d(\Sigma)}{n})^2}{\lambda_d^2}+\frac{(\lambda_{d+1}\frac{r_d(\Sigma)}{n})^2}{(\tau+\lambda_{d+1}\frac{r_d(\Sigma)}{n})^2}+\frac{d}{n}$. Note that
    \begin{align*}
    \frac{(\tau+\lambda_{d+1}\frac{r_d(\Sigma)}{n})^2}{\lambda_d^2}+\frac{(\lambda_{d+1}\frac{r_d(\Sigma)}{n})^2}{(\tau+\frac{\sum_{j>d}\lambda_j}{n})^2}\geq \frac{\lambda_{d+1}}{\lambda_d}\frac{r_d(\Sigma)}{n}
    \end{align*}
    and the equality can be achieved when $\tau+\lambda_{d+1}\frac{r_d(\Sigma)}{n} \asymp \sqrt{\lambda_d\lambda_{d+1}\frac{r_d(\Sigma)}{n}}\lesssim \lambda_d$. Hence we have $\MSE_{1,\In}^* \lesssim \frac{d}{n}+\frac{\lambda_{d+1}}{\lambda_d}\frac{r_d(\Sigma)}{n}\asymp \max\{\frac{d}{n},\frac{\lambda_{d+1}}{\lambda_d}\frac{r_d(\Sigma)}{n}\}$.

    Next we show that $\MSE_{1,\In}^* \gtrsim \max\{\frac{d}{n}+\frac{\lambda_{d+1}}{\lambda_d}\frac{r_d(\Sigma)}{n}\}$. From the lower bound in (\ref{eq:s_ii_1_1_18}), we have
    $\MSE_\In\gtrsim \kappa_2(\tau)\frac{(\tau+\lambda_{d+1}\frac{r_d(\Sigma)}{n})^2}{\lambda_d^2}+\frac{(\lambda_{d+1}\frac{r_d(\Sigma)}{n})^2}{(\tau+\lambda_{d+1}\frac{r_d(\Sigma)}{n})^2}+\frac{d}{n}$. Then we show $\kappa_2(\tau)\frac{(\tau+\lambda_{d+1}\frac{r_d(\Sigma)}{n})^2}{\lambda_d^2}+\frac{(\lambda_{d+1}\frac{r_d(\Sigma)}{n})^2}{(\tau+\lambda_{d+1}\frac{r_d(\Sigma)}{n})^2}\gtrsim \frac{\lambda_{d+1}}{\lambda_d}\frac{r_d(\Sigma)}{n}$ if $\tau+\lambda_{d+1}\frac{r_d(\Sigma)}{n}\lesssim \lambda_d$. If $\tau+\lambda_{d+1}\frac{r_d(\Sigma)}{n}\lesssim \sqrt{\lambda_{d+1}\lambda_d \frac{r_d(\Sigma)}{n}}$, we have
    \begin{align*}
    \kappa_2(\tau)\frac{(\tau+\lambda_{d+1}\frac{r_d(\Sigma)}{n})^2}{\lambda_d^2}+\frac{(\lambda_{d+1}\frac{r_d(\Sigma)}{n})^2}{(\tau+\lambda_{d+1}\frac{r_d(\Sigma)}{n})^2}\geq \frac{(\lambda_{d+1}\frac{r_d(\Sigma)}{n})^2}{(\tau+\lambda_{d+1}\frac{r_d(\Sigma)}{n})^2}\gtrsim \frac{\lambda_{d+1}}{\lambda_d}\frac{r_d(\Sigma)}{n}.
    \end{align*}
    If $\sqrt{\lambda_{d+1}\lambda_d \frac{r_d(\Sigma)}{n}} \lesssim \tau+\lambda_{d+1}\frac{r_d(\Sigma)}{n}\lesssim \lambda_d$, we have $\tau+\lambda_{d+1}\frac{r_d(\Sigma)}{n}\gg \lambda_{d+1}\frac{r_d(\Sigma)}{n}$ because $\lambda_d\gg \lambda_{d+1}\frac{r_d(\Sigma)}{n}$. Further with $64\frac{\sqrt{\delta_2}}{1-\sqrt{\delta_2}}<1$, we have $\kappa_2(\tau)\gtrsim 1$. Then
    \begin{align*}
    \kappa_2(\tau)\frac{(\tau+\lambda_{d+1}\frac{r_d(\Sigma)}{n})^2}{\lambda_d^2}+\frac{(\lambda_{d+1}\frac{r_d(\Sigma)}{n})^2}{(\tau+\lambda_{d+1}\frac{r_d(\Sigma)}{n})^2}\geq \kappa_2(\tau)\frac{(\tau+\lambda_{d+1}\frac{r_d(\Sigma)}{n})^2}{\lambda_d^2}\gtrsim \frac{\lambda_{d+1}}{\lambda_d}\frac{r_d(\Sigma)}{n}.
    \end{align*}
     Hence we have $\kappa_2(\tau)\frac{(\tau+\lambda_{d+1}\frac{r_d(\Sigma)}{n})^2}{\lambda_d^2}+\frac{(\lambda_{d+1}\frac{r_d(\Sigma)}{n})^2}{(\tau+\lambda_{d+1}\frac{r_d(\Sigma)}{n})^2}\gtrsim \frac{\lambda_{d+1}}{\lambda_d}\frac{r_d(\Sigma)}{n}$ if $\tau+\lambda_{d+1}\frac{r_d(\Sigma)}{n}\lesssim \lambda_d$ and $\MSE_{1,\In}^*\gtrsim \max\{\frac{d}{n}, \frac{\lambda_{d+1}}{\lambda_d}\frac{r_d(\Sigma)}{n}\}$.

    From $\MSE_{1,\In}^* \lesssim \max\{\frac{d}{n},\frac{\lambda_{d+1}}{\lambda_d}\frac{r_d(\Sigma)}{n}\}$ and $\MSE_{1,\In}^*\gtrsim \max\{\frac{d}{n}, \frac{\lambda_{d+1}}{\lambda_d}\frac{r_d(\Sigma)}{n}\}$, we have $\MSE_{1,\out}^*\asymp \max\{\frac{d}{n}, \frac{\lambda_{d+1}}{\lambda_d}\frac{r_d(\Sigma)}{n}\}$.
\end{itemize}

By integrating the orders of bounds of $\MSE_{1,\In}^*$ and $\MSE_{2,\In}^*$, we have $\MSE_{\In}^*\asymp \max\{\frac{d}{n}, \frac{\lambda_{d+1}}{\lambda_d}\frac{r_d(\Sigma)}{n}\}$.

\end{prf}

\section{Error approximation formulas}

\subsection{Convergence of in-sample error approximation formulas}

We provide a proof of Theorem \ref{thm5} in Section \ref{sec:4_1}, which is re-stated below for convenience.

\setcounter{thm}{4}
\begin{thm}[Convergence of in-sample error approximation formulas]
Under Assumption \ref{ass:5}, further assume that $\tau >\frac{1}{M}$ and $n^{-2/3+1/M}<\tau < \frac{M}{2}$. Then for any $D>0$, $\delta>0$, with probability at least $1-C(M,D, \delta)n^{-D}$,
\begin{align*}
    |\mathcal{B}_{\In}(\tau;\hat{H}_n,\hat{G}_n,\gamma) - \B_{\In}| &\leq  C(M)\max\{\frac{1}{\tau^{2/3}n^{(1-\delta)/3}},\frac{8M}{\tau n^{(1-\delta)/2}}\} , \\
    |\mathcal{V}_{\In}(\tau;\hat{H}_n,\gamma) - \V_{\In}| &\leq \sigma^2 C(M)(\max\{\frac{1}{\tau^{2/3}n^{(1-\delta)/3}},\frac{8M}{\tau n^{(1-\delta)/2}}\}  + \frac{1}{n^{(1-\delta)/2}} ).
\end{align*}
where $C(M,D,\delta)$ is a constant depending only on $(M, D, \delta)$ and $C(M)$ is a constant depending only on $M$.
\end{thm}
\begin{prf}
Our proof is inspired by the proof of Theorem 5 in \cite{HasteTrevorAndrea}. We first give the proof for bias and then the proof for variance. In the following proof, $C(M)$ is a constant depending on $M$ and may differ from line to line.
\newline
\textbf{Bias.} We first define two functions, $\bar{F}^{\tau}_n(\eta, \nu)$ and $F^\tau_n(\eta,\nu)$, and control the quantity $|\bar{F}^{\tau}_n(\eta, \nu)-F^\tau_n(\eta,\nu)|$. Moreover, we show that $-\frac{\partial \bar{F}^\tau_n}{\partial \eta}(0,\tau)=\B_\In$ and $-\frac{\partial F^\tau_n}{\partial \eta}(0,\tau)=\mathcal{B}_{\In}(\tau;\hat{H}_n,\hat{G}_n,\gamma)$. Then our objective is to control the quantity $|\frac{\partial \bar{F}^\tau_n}{\partial \eta}(0,\tau)-\frac{\partial F^\tau_n}{\partial \eta}(0,\tau)|$ so that we can control $|\B_\In-\mathcal{B}_{\In}(\tau;\hat{H}_n,\hat{G}_n,\gamma)|$.

Without loss of generality, we let $\Vert \theta^* \Vert_2^2 = 1$. For $\tau\in\mathbb{R}$, $\frac{1}{M}<\tau < \frac{M}{2}$ and $\tau>n^{-2/3+1/M}$, we define
\begin{align*}
    \bar{F}^{\tau}_n(\eta, \nu) &= \nu\langle \theta^*, (\hat{\Sigma}+\nu I + \tau \eta \hat{\Sigma})^{-1}
    \theta^*\rangle     \\
    &=\nu\langle \theta^*, ((1+\tau \eta)\hat{\Sigma}+\nu I)^{-1}\theta^*\rangle .
\end{align*}
Define $\mathbb{D}=\{(\eta,\nu)\in \mathbb{R}\times \mathbb{C}:\eta>-\frac{1}{2M}, \Re(\nu)>0\text{ and }\Im(-\nu)\geq 0\}$. Because $0<\tau < M$, $\bar{F}^{\tau}_n(\eta, \nu)$ is analytical in $\mathbb{D}$ and it can be easily verified that
\begin{align*}
    -\frac{\partial \bar{F}^\tau_n}{\partial \eta}(0,\tau)&=\tau^2\langle \theta^*, (\hat{\Sigma}+\tau I)^{-2}\hat{\Sigma}\theta^*\rangle\\
    &=\langle \theta^*, (\tau(\hat{\Sigma}+\tau I)^{-1})\hat{\Sigma}(\tau(\hat{\Sigma}+\tau I)^{-1})\theta^*\rangle   \\
    &= \langle \theta^*, (I - \hat{\Sigma}(\hat{\Sigma}+\tau I)^{-1})\hat{\Sigma}(I - \hat{\Sigma}(\hat{\Sigma}+\tau I)^{-1})\theta^*\rangle   \\
    &=\B_{\In} .
\end{align*}
Define $\mathbb{D}_0=\{(\eta,\nu)\in \mathbb{R}\times \mathbb{C}:\eta>-\frac{1}{2M},0<\Re(\nu)<M,0<\Im(-\nu) <M\}$.
By using the anisotropic local law for covariance matrices
in Theorem 3.16(i) of \cite{knowles2016anisotropic},
we obtain that for any $\delta>0, \epsilon_0>0, D>0$, with probability at least $1-C(\epsilon_0, D, \delta)n^{-D}$, we have for $(\eta,\nu)\in \mathbb{D}_0$ and $\Re(\nu)>n^{-2/3+\epsilon_0}$,
\begin{align}
|\bar{F}^{\tau}_n(\eta, \nu)-F^{\tau}_n(\eta,\nu)|&\leq \sqrt{\frac{\Im(\tilde{r}_n(\eta,-\nu))}{\Im(-\nu)}n^{-1+\delta}} ,\label{siii_1} \\
    F^{\tau}_n(\eta,\nu)&=\langle\theta^*, (I+\tilde{r}_n(\eta,-\nu)(1+\tau\eta)\Sigma)^{-1}\theta^*\rangle ,  \notag
\end{align}
where $\tilde{r}_n(\eta,z)$ is defined in the demoain $\mathbb{D}_1=\{(\eta,z)\in \mathbb{R}\times \mathbb{C}:\eta>-\frac{1}{2M}, \Re(z)<0\text{ and }\Im(z)\geq 0\}$ and it is defined as the unique solution satisfying $\Im(\tilde{r}_n(\eta,z))>0$ if $\Im(z)>0$ or $\tilde{r}_n(\eta,z)>0$ if $\Im(z)=0$ of
\begin{align*}
    \frac{1}{\tilde{r}_n}&=-z+\gamma \frac{1}{p}\sum_{j=1}^p \frac{(1+\tau\eta)\lambda_j}{1+(1+\tau\eta)\lambda_j \tilde{r}_n}.
\end{align*}

Following a similar process as in Section~A.1.2 of \cite{HasteTrevorAndrea}, we have
\begin{align*}
    |\Im(\tilde{r}_n(\eta,-\nu))|\leq \frac{|\Im(\nu)|}{\Re(\nu)^2}.
\end{align*}
Then taking the limit $\Im(\nu)\rightarrow 0$ in (\ref{siii_1}) and let $\epsilon_0=\frac{1}{M}$, we obtain, with probability at least $1-C( M,D,\delta)n^{-D}$, for $n^{-2/3+1/M}<\nu<M $ and $\eta > -\frac{1}{2M}$,
\begin{align}
    |\bar{F}^\tau_n(\eta,\nu)-F^\tau_n(\eta,\nu)|\leq \frac{1}{n^{(1-\delta)/2}\nu}. \notag
\end{align}
Let $\nu=\tau$, we have
\begin{align}
    |\bar{F}^\tau_n(\eta,\tau)-F^\tau_n(\eta,\tau)|\leq \frac{1}{n^{(1-\delta)/2}\tau}. \label{eq:III_1_1}
\end{align}

Then we reformulate $F^\tau_n(\eta,\tau)$ by substituting $\tilde{r}_n(\eta,-\nu)$ with other quantity and show that $-\frac{\partial F^\tau_n}{\partial \eta}(0,\tau)=\mathcal{B}_{\In}(\tau;\hat{H}_n,\hat{G}_n,\gamma)$. Define $\mathbb{D}_2=\{z\in \mathbb{C}: \Re(z) <0,\Im(z)\geq 0\}$. For $z\in \mathbb{D}_2$, we define $r_n(z) \in \mathbb{R}$ as the unique solution satisfying $\Im(r_n(z))>0$ if $\Im(z)>0$ or $r_n(z)>0$ if $\Im(z)=0$ of
\begin{align}
    \frac{1}{r_n}&=-z+\gamma \frac{1}{p}\sum_{j=1}^p \frac{\lambda_j}{1+\lambda_j r_n}  .\label{eq:s2_1b}
\end{align}
In the following discussion, we consider $(\eta,v)\in \mathbb{D} \cap \mathbb{R}^2$, that is $\eta>-\frac{1}{2M}$ and $\nu > 0$. Then we have,
\begin{align}
     \tilde{r}_n(\eta,-\nu)&=\frac{1}{1+\tau\eta}r_n(-\frac{\nu}{1+\tau\eta}) ,\notag \\
     F^\tau_n(\eta,\nu)&=\langle\theta^*, (I+r_n(-\frac{\nu}{1+\tau\eta})\Sigma)^{-1}\theta^*\rangle .  \label{eq:s2_1c}
\end{align}
Let $m_n(z) = \frac{1-\gamma+z r_n(z)}{\gamma z}$, which by (\ref{eq:s2_1b}) is the unique solution of
\begin{align*}
    m_n = \frac{1}{p}\sum_{j=1}^p \frac{1}{\lambda_j(1-\gamma-\gamma z m_n)-z}.
\end{align*}
Then we have
\begin{align*}
    F^\tau_n(\eta,\nu)&=\langle\theta^*, (I+ (m_n(-\frac{\nu}{1+\tau \eta})\gamma + \frac{(1 - \gamma)(1+\tau \eta) }{\nu})\Sigma)^{-1}\theta^*\rangle,
\end{align*}
and
\begin{align*}
    -\frac{\partial F^\tau_n(0,\tau)}{\partial \eta} &= (\gamma\tau^2 m_n^\prime(-\tau) + 1-\gamma) \langle\theta^*, (I+ (m_n(-\tau)\gamma + \frac{1 - \gamma}{\tau})\Sigma)^{-2}\Sigma\theta^*\rangle  \notag \\
    &= \tau^2(\gamma\tau^2 m_n^\prime(-\tau) + 1-\gamma)||\theta^*||^2\int \frac{s}{[\tau+(1-\gamma+\gamma\lambda m_n(-\tau))s]^2}d\hat{G}_n(s) \\
    &=\mathcal{B}_{\In}(\tau;\hat{H}_n,\hat{G}_n,\gamma).
\end{align*}

In the following discussion, we give upper bounds on $|\frac{\partial \bar{F}^\tau_n}{\partial \eta}(0,\tau)-\frac{\partial F^\tau_n}{\partial \eta}(0,\tau)|$, equivalently $|\B_\In-\mathcal{B}_{\In}(\tau;\hat{H}_n,\hat{G}_n,\gamma)|$. Our strategy is to control $|\frac{\partial^k \bar{F}^\tau_n}{\partial \eta^k}(\eta,\tau)|$, $|\frac{\partial^k F^\tau_n}{\partial \eta^k}(\eta,\tau)|$, then $|\B_\In-\mathcal{B}_{\In}(\tau;\hat{H}_n,\hat{G}_n,\gamma)|$, can be controlled by $|\bar{F}^\tau_n(\eta,\tau)-F^\tau_n(\eta,\tau)|$, $|\frac{\partial^k \bar{F}^\tau_n}{\partial \eta^k}(\eta,\tau)|$ and $|\frac{\partial^k F^\tau_n}{\partial \eta^k}(\eta,\tau)|$ based on the Lemma A.1 in \cite{HasteTrevorAndrea}.

We first give upper bound on $|\frac{\partial^k \bar{F}^\tau_n}{\partial \eta^k}(\eta,\tau)|$.
For $0<\nu < M, \eta > -\frac{1}{2M}$, we have for $k\geq 1$,
\begin{align*}
    \frac{\partial^k \bar{F}^\tau_n}{\partial \eta^k}(\eta,\nu) = k!(-1)^{k+1} \tau^{k}\nu\langle\theta^*, R^k\hat{\Sigma}^kR\theta^*\rangle,
\end{align*}
where $k!=k\times \cdots \times 1$ and $R = ((1+\tau\eta)\hat{\Sigma}+\nu I)^{-1}$. Then we have for $\nu=\tau$, $\eta > -\frac{1}{2M}$ and $k\geq 1$,
\begin{align}
    |\frac{\partial^k \bar{F}^\tau_n}{\partial \eta^k}(\eta,\tau)|\leq k!2^k \tau^k \Vert \theta^*\Vert^2  \leq k! 2^k M^k. \label{eq:III_1_2}
\end{align}

Next we give the upper bound on $|\frac{\partial^k F^\tau_n}{\partial \eta^k}(\eta,\tau)|$.
From (\ref{eq:s2_1c}), it is sufficient to upper bound $|\frac{\tau^{l+1}}{(1+\tau\eta)^{l+1}}|$, $|r^{(l)}_n(-\frac{\tau}{1+\tau\eta})|$ and $|\langle \theta^*,(I+r_n(-\frac{\tau}{1+\tau\eta} )\Sigma)^{-(l+1)}\Sigma^l \theta^* \rangle|$ for $1\leq l\leq k$. We give their upper bounds separately as follows.
\begin{itemize}
    \item \textbf{Upper bound of $|\frac{\tau^{l+1}}{(1+\tau\eta)^{l+1}}|$.} Because $\frac{1}{M}<\tau < \frac{M}{2}$, for $-\frac{1}{2M}<\eta<\frac{1}{2M}$ and $1\leq l \leq k$, we have
\begin{align*}
    |\frac{\tau^{l+1}}{(1+\tau\eta)^{l+1}}|  \leq   M^{l+1}.
\end{align*}
   \item \textbf{Upper bound of $|\langle\theta^*,(I+r_n(-\frac{\tau}{1+\tau\eta} )\Sigma)^{-(l+1)}\Sigma^l \theta^*\rangle|$.}
    We show $r_n(-\frac{\tau}{1+\tau\eta} ) > 0$. Let $u_n(z)=\frac{1}{r_n(-z)}$. Then for $z>0$, $u_n(z)$ is the unique solution of
\begin{align}
    u_n = z + \frac{1}{n}\sum_{j=1}^p\frac{\lambda_ju_n}{u_n+\lambda_j} .  \notag
\end{align}
From Lemma A.2(a)(b) in \cite{HasteTrevorAndrea}, we have for $\frac{1}{2M}<z < M$,
\begin{align}
      u_n(z) &> \frac{1}{4M} \, (>0), \label{s_iii_1_0} \\
     0 < |u_n^\prime(z)| &< C(M). \label{s_iii_1_1}
\end{align}
By (\ref{s_iii_1_0}) and the definition $u_n(z)=\frac{1}{r_n(-z)}$, we have for $\frac{1}{2M}<z < M$,
\begin{align}
    r_n(-z) > 0 .\label{s_iii_1_4}
\end{align}
Because $\frac{1}{M}<\tau < \frac{M}{2}$, then for $-\frac{1}{2M}<\eta <\frac{1}{2M}$, we have $\frac{1}{2M}<\frac{\tau}{1+\tau\eta} < M$ and $r_n(-\frac{\tau}{1+\tau\eta})>0$. Hence for $-\frac{1}{2M}<\eta <\frac{1}{2M}$ and $1\leq l\leq k$,
\begin{align*}
    |\langle\theta^*,(I+r_n(-\frac{\tau}{1+\tau\eta} )\Sigma)^{-(l+1)}\Sigma^l \theta^*\rangle| &\leq  \Vert (I+r_n(-\frac{\tau}{1+\tau\eta})\Sigma)^{-1}\Vert^{l+1}_{\op} \Vert \Sigma\Vert^{l}_{\op} \Vert \theta^*\Vert^2 \\
    &\leq  \Vert \Sigma\Vert^{l}_{\op} \Vert \theta^*\Vert^2 \quad (\text{because $r_n(-\frac{\tau}{1+\tau\eta})>0$}) \\
    &\leq \lambda_1^l \Vert \theta^*\Vert^2  \\
    &\leq  M^l  \Vert \theta^*\Vert^2\quad (\text{from Assumption \ref{ass:5})}\end{align*}

   \item \textbf{Upper bound of $|r^{(l)}_n(-\frac{\tau}{1+\tau\eta})|$.}
   Let $u_n(z)=\frac{1}{r_n(-z)}$. Then $|r_n^{(l)}(-z)|$ can be controlled by the polynomial of $|u_n^{(m)}(z)|$ for $1\leq m\leq l$ and $|\frac{1}{u_n(z)}|$.

   From (\ref{s_iii_1_0}), we have for $\frac{1}{2M}<z < M$,
   \begin{align}
       |u_n^{-1}(z)| \leq 4M.  \label{inq:upper_bound_un}
   \end{align}
   The upper bound of $|u_n^\prime(z)|$ is provided in (\ref{s_iii_1_1}). Then we give the upper bounds of $|u_n^{(m)}(z)|$ for $2\leq m \leq l$. We consider the following function
\begin{align*}
    f(u_n,z) = u_n+z-\frac{1}{n}\sum_{j=1}^p \frac{\lambda_j u_n}{u_n+\lambda_j}.
\end{align*}
From implicit function theorem and (\ref{s_iii_1_1}), we have for $\frac{1}{2M}<z < M$,
\begin{align}
    |\frac{\partial f}{\partial u_n}|=|u^\prime_n(z)|^{-1} \geq \frac{1}{C(M)}.  \label{inq:lower_bound_parital_f_u_n}
\end{align}

To give upper the bounds of $|u_n^{(m)}(z)|$ for $2\leq m \leq l$, from the implicit function theorem, it is sufficient to further give the upper bounds of $|\frac{\partial^{s+t} f}{\partial u_n^s \partial z^t}|$ for all $s+t\leq m$ and $s\geq 1, t \geq 1$ or $s= 0$, $t\geq 2$. Denote
\begin{align*}
    f_1(u_n) = \frac{1}{n}\sum_{j=1}^p \frac{\lambda_j u_n}{u_n+\lambda_j}.
\end{align*}
Then it is sufficient to give upper bounds of $|\frac{\partial^s f_1(u_n)}{\partial u_n^s}| \text{for }1\leq s \leq m$. From (\ref{s_iii_1_0}), we have for $\frac{1}{2M} < z< M$,
\begin{align*}
    u_n(z) > \frac{1}{4M}.
\end{align*}
Then
\begin{align*}
    |\frac{\partial^s f_1(u_n)}{\partial u_n^s} |&=|\frac{s!}{n}\sum_{j=1}^p\frac{\lambda_j^2}{(\lambda_j+u_n)^{s+1}}| \\
    &\leq s!\frac{p}{n} \frac{\lambda_1^2}{(\lambda_p + u_n)^{s+1}} \\
    &\leq s! \frac{p}{n} \frac{\lambda_1^2}{u_n^{s+1}} \\
    &\leq  C(M). \quad\text{(from $\gamma=\frac{p}{n} < M$, $\lambda_1 <M$ and $u_n(z) > \frac{1}{4M}$)}
\end{align*}
Hence from the implicit function theorem, for $\frac{1}{2M} < z < M$, $s+t\leq m$ and $s\geq 0$, $t \geq 1$ or $s= 0$, $t\geq 2$,
\begin{align*}
    |\frac{\partial^{s+t} f}{\partial u_n^s \partial z^t}|\leq C(M) .
\end{align*}
Combining the above with (\ref{inq:lower_bound_parital_f_u_n}), from implicit function theorem, we have for $\frac{1}{2M}<  z < M$ and $2\leq m \leq l$,
\begin{align}
    |u_n^{(m)}(z)| &\leq C(M)  .\label{s_iii_1_2}
\end{align}
With the upper bound of $|u_n^{-1}(z)|$ in (\ref{inq:upper_bound_un}), upper bound of $|u_n^\prime(z)|$ in (\ref{s_iii_1_1}) and upper bound of $|u_n^{(m)}(z)|$ for $2\leq m\leq l$ in (\ref{s_iii_1_2}), we have for $\frac{1}{2M} <  z < M$ and for $1\leq l \leq k$,
\begin{align}
    |r_n^{(l)}(-z)| \leq C(M) .\label{s_iii_1_6}
\end{align}
Because $\frac{1}{M}<\tau < \frac{M}{2}$, for $-\frac{1}{2M}<\eta <\frac{1}{2M}$, we have $\frac{1}{2M}< \frac{\tau}{1+\eta\tau} < M$. From (\ref{s_iii_1_6}), we have for $-\frac{1}{2M}<\eta <\frac{1}{2M}$ and $1\leq l \leq k$,
\begin{align*}
    |r_n^{(l)}(-\frac{\tau}{1+\eta\tau})| \leq  C(M).
\end{align*}
\end{itemize}

From the upper bounds of $|\frac{\tau^{l+1}}{(1+\tau\eta)^{l+1}}|$, $|\langle\theta^*,(I+r_n(-\frac{\tau}{1+\tau\eta} )\Sigma)^{-(l+1)}\Sigma^l \theta^*\rangle|$ and $|r^{(l)}_n(-\frac{\tau}{1+\tau\eta})|$ for $1\leq l\leq k$ above, we have for $-\frac{1}{2M}<\eta < \frac{1}{2M}$ and $k\geq 1$,
\begin{align}
  |\frac{\partial^k F^\tau_n}{\partial \eta^k}(\eta,\tau)|
  &\leq  C(M).\label{eq:III_3}
\end{align}

In the following discussion, we combine the upper bounds of $|\bar{F}^\tau_n(\eta,\tau)-F^\tau_n(\eta,\tau)|$, $|\frac{\partial^k \bar{F}^\tau_n}{\partial \eta^k}(\eta,\tau)|$ and $|\frac{\partial^k F^\tau_n}{\partial \eta^k}(\eta,\tau)|$ from above and apply Lemma A.1 in \cite{HasteTrevorAndrea} to control $|\frac{\partial \bar{F}^\tau_n}{\partial \eta}(0,\tau)-\frac{\partial F^\tau_n}{\partial \eta}(0,\tau)|$. Combining (\ref{eq:III_1_1}),(\ref{eq:III_1_2}) and (\ref{eq:III_3}), from Lemma A.1 in \cite{HasteTrevorAndrea},
and letting $k=3$, we have for $0<\xi <\frac{1}{4M}$, $D>0$ and $\delta > 0$, with probability at least $1-C(\delta, M, D)n^{-D}$,
\begin{align*}
    &|\frac{\partial \bar{F}^\tau_n}{\partial \eta}(0,\tau)-\frac{\partial F^\tau_n}{\partial \eta}(0,\tau)| \leq  C(M)( \frac{1}{\tau n^{(1-\delta)/2}}\frac{1}{\xi}+\xi^2).
\end{align*}
That is, we have
\begin{align*}
    |B_{in}-\mathcal{B}_{\In}(\tau;\hat{H}_n,\hat{G}_n,\gamma)| \leq  C(M)( \frac{1}{\tau n^{(1-\delta)/2}}\frac{1}{\xi}+\xi^2).
\end{align*}
Letting $\xi=\min\{\frac{1}{8M}, \tau^{1/3}n^{(1-\delta)/6}\}$, we have
\begin{align*}
    |\mathcal{B}_{\In}(\tau;\hat{H}_n,\hat{G}_n,\gamma) - \B_{\In}| &\leq  C(M)\max\{\frac{1}{\tau^{2/3}n^{(1-\delta)/3}},\frac{8M}{\tau n^{(1-\delta)/2}}\}.
\end{align*}
\newline
\textbf{Variance.} We have
\begin{align*}
    \V_{\In}&= \frac{\sigma^2}{n} \Tr(\hat{\Sigma}^2(\hat{\Sigma}^2+\tau I)^{-2}) \\
    &=\sigma^2 (\gamma - 2\tau \gamma \frac{1}{p}\Tr((\hat{\Sigma}+\tau I)^{-1}) + \tau^2 \gamma \frac{1}{p}\Tr(\hat{\Sigma}+\tau I)^{-2}).
\end{align*}
From (\ref{eq:s2_2d}) and $m_n(-\tau) = \frac{\gamma-1+\tau r_n(-\tau)}{\gamma \tau}$,
\begin{align*}
     \mathcal{V}_{\In}(\tau;\hat{H}_n,\gamma)=\sigma^2(1 - 2\tau r_n(-\tau)+\tau^2 r_n^\prime(-\tau)).
\end{align*}

We first give the equations and inequality below for the following analysis. From (\ref{eq:s2_1b}), we have
\begin{align}
   \tau r_n(-\tau) + \gamma-1 &=\gamma\frac{1}{p}\sum_{j=1}^n\frac{1}{1+r_n(-\tau)\lambda_j}\notag\\ &= \gamma\frac{1}{p}\Tr((1+r_n(-\tau)\Sigma)^{-1}) .\label{eq:r_n_to_Tr_r_n}
\end{align}
From (\ref{eq:r_n_to_Tr_r_n}), we have
\begin{align*}
  \tau \frac{\partial \gamma\frac{1}{p}\Tr((1+r_n(-\tau)\Sigma)^{-1})}{\partial \tau} = -\tau^2 r_n^\prime(-\tau) + \tau r_n(-\tau).
\end{align*}
Hence we have
\begin{align}
    &\tau\frac{\partial  \gamma \tau \frac{1}{p}\Tr((\hat{\Sigma}+\tau I)^{-1})}{\partial \tau} - \tau\frac{\partial \gamma\frac{1}{p}\Tr((1+r_n(-\tau)\Sigma)^{-1})}{\partial \tau} \notag \\
    &= \tau\gamma\frac{1}{p}\Tr((\hat{\Sigma}+\tau I)^{-1}) - \tau^2\gamma\frac{1}{p}Tr((\hat{\Sigma}+\tau I)^{-2}) +\tau^2 r_n^\prime(-\tau)-\tau r_n(-\tau)  \notag \\
    &= (\tau^2 r_n^\prime(-\tau) +\gamma - 1 - \tau^2\gamma\frac{1}{p}\Tr((\hat{\Sigma}+\tau I)^{-2}) + (\tau \gamma \frac{1}{p}\Tr((\hat{\Sigma}+\tau I)^{-1}) - (\tau r_n(-\tau) + \gamma-1)). \notag
\end{align}
From above, we have
\begin{align}
    &|\tau^2 r_n^\prime(-\tau) + \gamma - 1 - \tau^2\gamma\frac{1}{p}\Tr((\hat{\Sigma}+\tau I)^{-2})| \notag \\
    &\leq |\tau \gamma \frac{1}{p}\Tr((\hat{\Sigma}+\tau I)^{-1}) - (\tau r_n(-\tau) + \gamma-1)| + | \tau\gamma \frac{\partial \tau \frac{1}{p}Tr((\hat{\Sigma}+\tau I)^{-1})}{\partial \tau} - \tau\gamma\frac{\partial \frac{1}{p}\Tr((I+r_n(-\tau)\Sigma)^{-1})}{\partial \tau}|  \notag  \\
    &=  |\tau \gamma \frac{1}{p}\Tr((\hat{\Sigma}+\tau I)^{-1}) - \gamma\frac{1}{p}\Tr((1+r_n(-\tau)\Sigma)^{-1})| + | \tau\gamma \frac{\partial \tau \frac{1}{p}Tr((\hat{\Sigma}+\tau I)^{-1})}{\partial \tau} - \tau\gamma\frac{\partial \frac{1}{p}\Tr((I+r_n(-\tau)\Sigma)^{-1})}{\partial \tau}|. \label{inq:tau_sq_tr}
\end{align}

We give the upper bound of $|\V_{\In}-\mathcal{V}_{\In}(\tau;\hat{H}_n,\gamma)|$ by controlling $|\tau \frac{1}{p}\Tr((\hat{\Sigma}+\tau)^{-1})-\frac{1}{p}\Tr((1+r_n(-\tau)\Sigma)^{-1}) |$ and $|\tau \frac{\partial \tau \frac{1}{p}\Tr((\hat{\Sigma}+\tau I)^{-1})}{\partial \tau} - \tau\frac{\partial \frac{1}{p}\Tr((I+r_n(-\tau)\Sigma)^{-1})}{\partial \tau}|$. In fact, we have
\begin{align}
    &|\V_{\In}-\mathcal{V}_{\In}(\tau;\hat{H}_n,\gamma)| \notag\\
    &= \sigma^2 |(\gamma - 2\tau \gamma \frac{1}{p}\Tr((\hat{\Sigma}+\tau I)^{-1}) + \tau^2 \gamma \frac{1}{p}\Tr((\hat{\Sigma}+\tau I)^{-2}))-(1 - 2\tau r_n(-\tau)+\tau^2 r_n^\prime(-\tau))| \notag \\
    &\leq \sigma^2 |2\tau r_n(-\tau) + 2\gamma-2 - 2\tau\gamma \frac{1}{p}\Tr((\hat{\Sigma}+\tau I)^{-1}) + (\tau^2\gamma \frac{1}{p}\Tr((\hat{\Sigma}+\tau I)^{-2})- (\tau^2r_n^\prime(-\tau) +\gamma-1) | \notag \\
    &\leq  \sigma^2 (2|\gamma\tau \frac{1}{p}\Tr((\hat{\Sigma}+\tau I)^{-1})-(\tau r_n(-\tau) + \gamma-1)| + |\gamma \tau^2 \frac{1}{p}\Tr((\hat{\Sigma}+\tau I)^{-2}) - (\tau^2 r_n^\prime(-\tau)+\gamma-1)|) \notag \\
    &\leq  \sigma^2 \gamma(3|\tau \frac{1}{p}\Tr((\hat{\Sigma}+\tau I)^{-1})-\frac{1}{p}\Tr((1+r_n(-\tau)\Sigma)^{-1}) | \notag \\& \quad\quad\quad\quad\quad\quad+ |\tau \frac{\partial \tau \frac{1}{p}\Tr((\hat{\Sigma}+\tau I)^{-1})}{\partial \tau} - \tau\frac{\partial \frac{1}{p}\Tr((I+r_n(-\tau)\Sigma)^{-1})}{\partial \tau}|) \quad \text{(By (\ref{eq:r_n_to_Tr_r_n}) and (\ref{inq:tau_sq_tr}))}.\label{s_iii_1_11}
\end{align}

We first give the upper bound on $|\tau \frac{1}{p}\Tr((\hat{\Sigma}+\tau)^{-1})-\frac{1}{p}\Tr((1+r_n(-\tau)\Sigma)^{-1}) |$. Using Theorem 3.16(i) in \cite{knowles2016anisotropic}, we have for $D>0$, $\delta>0$, $\Im(-\tau)>0$ and $\Re(\tau)>n^{-2/3+1/M}$, with probablity at least $1-C(M,D,\delta) n^{-D}$,
\begin{align*}
    |\tau \frac{1}{p}\Tr((\hat{\Sigma}+\tau)^{-1}) - \frac{1}{p}\Tr((1+r_n(-\tau)\Sigma)^{-1})| \leq \sqrt{\frac{\Im(r_n(-\tau))}{\Im(-\tau)}}.
\end{align*}
Following a similar process in Section~A.1.2 of \cite{HasteTrevorAndrea}, we have $|\Im(r_n(-\tau))|\leq |\Im(\tau)|/\Re(\tau)^2$.
Letting $\Im(\tau)\rightarrow 0$ shows that
for $D>0$, $\delta > 0$, $\frac{1}{M}<\tau < M$ and $ \tau>n^{-2/3+(1/M)}$, with probability at least $1-C(M,D,\delta)n^{-D}$,
\begin{align}
     &|\tau \frac{1}{p}\Tr((\hat{\Sigma}+\tau I)^{-1}) - \frac{1}{p}\Tr((1+r_n(-\tau)\Sigma)^{-1})| \leq \frac{1}{\tau n^{(1-\delta)/2}} \leq  \frac{M}{n^{(1-\delta)/2}} .\label{s_iii_1_8}
\end{align}

Next, we give the upper bound on $| \tau\frac{\partial \tau \frac{1}{p}\Tr((\hat{\Sigma}+\tau I)^{-1})}{\partial \tau} - \tau\frac{\partial \frac{1}{p}\Tr((I+r_n(-\tau)\Sigma)^{-1})}{\partial \tau}|$. Our strategy is to upper bound $|\frac{\partial^{k} \tau \frac{1}{p}\Tr((\hat{\Sigma}+\tau)^{-1})}{\partial \tau^k}|$ and $|\frac{\partial^k \frac{1}{p}\Tr((1+r_n(-\tau)\Sigma)^{-1})}{\partial \tau^k}|$ for $k\geq 1$,
so that Lemma A.1 in \cite{HasteTrevorAndrea} can be applied again to bound $|\frac{\partial \tau \frac{1}{p}\Tr((\hat{\Sigma}+\tau)^{-1})}{\partial \tau}-\frac{\partial \frac{1}{p}\Tr((1+r_n(-\tau)\Sigma)^{-1})}{\partial \tau}|$ by $|\tau \frac{1}{p}\Tr((\hat{\Sigma}+\tau)^{-1}) - \frac{1}{p}\Tr((1+r_n(-\tau)\Sigma)^{-1})|$, $|\frac{\partial^{k} \tau \frac{1}{p}\Tr((\hat{\Sigma}+\tau)^{-1})}{\partial \tau^k}|$ and $|\frac{\partial^k \frac{1}{p}\Tr((1+r_n(-\tau)\Sigma)^{-1})}{\partial \tau^k}|$ for $k\geq 1$. The upper bound of $|\tau \frac{1}{p}\Tr((\hat{\Sigma}+\tau)^{-1}) - \frac{1}{p}\Tr((1+r_n(-\tau)\Sigma)^{-1})|$ is given in (\ref{s_iii_1_8}). Then we give the upper bound of $|\frac{\partial^{k} \tau \frac{1}{p}\Tr((\hat{\Sigma}+\tau)^{-1})}{\partial \tau^k}|$ and $|\frac{\partial^k \frac{1}{p}\Tr((1+r_n(-\tau)\Sigma)^{-1})}{\partial \tau^k}|$ as below.

We give the upper bound on $|\frac{\partial^k \tau \frac{1}{p}\Tr((\hat{\Sigma}+\tau I)^{-1})}{\partial \tau^k}|$. For $k\geq 1$ and $\frac{1}{M}<\tau< M$,
\begin{align}
    |\frac{\partial^k \tau \frac{1}{p}\Tr((\hat{\Sigma}+\tau I)^{-1})}{\partial \tau^k}| &= |\frac{k!}{p} \Tr(\hat{\Sigma}(\hat{\Sigma}+\tau I)^{-(k+1)}) | \notag \\
    &\leq  \frac{k!}{\tau^k} \leq C(M) . \label{s_iii_1_9}
\end{align}

Then we give the upper bound to $|\frac{\partial^k \frac{1}{p}\Tr((I+r_n(-\tau)\Sigma)^{-1})}{\partial \tau^k}|$. It is sufficient to upper bound $|r_n^{(l)}(-\tau)|$ and $|\frac{1}{p}\Tr(\Sigma^l(I+r_n(-\tau)\Sigma)^{-(l+1)})|$ for $1\leq l\leq k$. From (\ref{s_iii_1_4}), we have for $\frac{1}{2M}<\tau < M$,
\begin{align*}
    r_n(-\tau) > 0,
\end{align*}
and
\begin{align*}
    |\frac{1}{p}\Tr(\Sigma^l(I+r_n(-\tau)\Sigma)^{-(l+1)})| &\leq  \Vert (I+r_n(-\tau)\Sigma)^{-1}\Vert^{l+1}_{\op} \Vert \Sigma\Vert^{l}_{\op} \\
    &\leq  \Vert \Sigma\Vert^{l}_{\op}  \quad (\text{because $r_n(-\tau)>0$}) \\
     &\leq  \lambda_1^l \\
     &\leq  M^l \quad(\text{from Assumption \ref{ass:5}})   \end{align*}
From (\ref{s_iii_1_6}), for $\frac{1}{2M}<\tau < M$ and $1\leq l\leq k$,
\begin{align*}
    |r_n^{(l)}(-\tau)| < C(M).
\end{align*}
Hence for $\frac{1}{2M}<\tau < M$ and $k\geq 1$,
\begin{align}
    |\frac{\partial^k \frac{1}{p}\Tr((I+r_n(-\tau)\Sigma)^{-1})}{\partial \tau^k}| < C(M) .\label{s_iii_1_10}
\end{align}

We combine the upper bounds of $|\tau \frac{1}{p}\Tr((\hat{\Sigma}+\tau)^{-1}) - \frac{1}{p}\Tr((1+r_n(-\tau)\Sigma)^{-1})|$, $|\frac{\partial^{k} \tau \frac{1}{p}\Tr((\hat{\Sigma}+\tau)^{-1})}{\partial \tau^k}|$ and $|\frac{\partial^k \frac{1}{p}\Tr((1+r_n(-\tau)\Sigma)^{-1})}{\partial \tau^k}|$ from above and apply Lemma A.1 in \cite{HasteTrevorAndrea} to control $|\tau\frac{\partial \tau \frac{1}{p}\Tr((\hat{\Sigma}+\tau)^{-1})}{\partial \tau}-\tau\frac{\partial \frac{1}{p}\Tr((1+r_n(-\tau)\Sigma)^{-1})}{\partial \tau}|$. From (\ref{s_iii_1_8})--(\ref{s_iii_1_10}) and Lemma A.1 in \cite{HasteTrevorAndrea}, and letting $k=3$, we have for $D>0$, $\delta > 0$, $\xi < \frac{1}{M}$, for $ \tau>n^{-2/3+(1/M)}$ and $\frac{1}{M}<\tau <M$,  with probability at least $1-C(D,\delta,M)n^{-D}$,
\begin{align}
     |\tau\frac{\partial \tau \frac{1}{p}\Tr((\hat{\Sigma}+\tau)^{-1})}{\partial \tau}-\tau\frac{\partial \frac{1}{p}\Tr((1+r_n(-\tau)\Sigma)^{-1})}{\partial \tau}| \leq  C(M)(\frac{1}{\tau n^{(1-\delta)/2}}\frac{1}{\xi}+\xi^2) .\label{s_iii_1_7}
\end{align}

With the upper bounds of $|\tau \frac{1}{p}\Tr((\hat{\Sigma}+\tau I)^{-1}) - \frac{1}{p}\Tr((1+r_n(-\tau)\Sigma)^{-1})|$ and $|\tau\frac{\partial \tau \frac{1}{p}\Tr((\hat{\Sigma}+\tau)^{-1})}{\partial \tau}-\tau\frac{\partial \frac{1}{p}\Tr((1+r_n(-\tau)\Sigma)^{-1})}{\partial \tau}|$ from above, we give the upper bound of $|\mathcal{V}_{\In}(\tau;\hat{H}_n,\gamma) - \V_{\In}|$. From (\ref{s_iii_1_11})--(\ref{s_iii_1_8}) and (\ref{s_iii_1_7}), we have for any $D>0$, $\delta > 0$, $0<\xi < \frac{1}{M}$, for $ \tau>n^{-2/3+(1/M)}$ and $\frac{1}{M}<\tau <M$, with probability at least $1-C(D,\delta,M)n^{-D}$,
\begin{align*}
    |\mathcal{V}_{\In}(\tau;\hat{H}_n,\gamma) - \V_{\In}|\leq \sigma^2 C(M)(\frac{1}{\tau n^{(1-\epsilon)/2}}\frac{1}{\xi}+\xi^2 + \frac{1}{n^{(1-\delta)/2}}).
\end{align*}
Letting $\xi=\min\{\frac{1}{8M}, \tau^{1/3}n^{(1-\delta)/6}\}$, we have
\begin{align*}
    |\mathcal{V}_{\In}(\tau;\hat{H}_n,\gamma) - \V_{\In}| &\leq \sigma^2 C(M)(\max\{\frac{1}{\tau^{2/3}n^{(1-\epsilon)/3}},\frac{8M}{\tau n^{(1-\epsilon)/2}}\}  + \frac{1}{n^{(1-\epsilon)/2}} ).
\end{align*}

\end{prf}

\subsection{Orders of error approximation formulas}

We provide proofs of Corollaries \ref{cor7_c} and \ref{cor7_d} in Section \ref{sec:4_1}, which are re-stated below for convenience.

\label{derivation_of_equivalence}
\begin{cor}[Matching error approximation formulas with small or moderate TER]
~\newline
\indent(i) Suppose that $\frac{d}{n} < 1$, $r_d(\Sigma)\lesssim n$, and $\Vert\theta^*_{(d+1):p}\Vert^2_{\Sigma_{(d+1):p}}\lesssim  \Vert\theta^*_{1:d}\Vert^2_{\Sigma_{1:d}^{-1}}\lambda_{d+1}^2$.
For $\lambda_{d+1} \lesssim   \tau \lesssim \lambda_{d}$, we have
\begin{align}
    \mathcal{B}_{\out}(\tau;\hat{H}_n,\hat{G}_n,\gamma) +
    \mathcal{V}_{\out}(\tau;\hat{H}_n,\gamma) &\asymp  \Vert\theta^*_{1:d}\Vert_{\Sigma^{-1}_{1:d}}^2 \tau^2+\sigma^2(\frac{d}{n}+\frac{\lambda_{d+1}^2}{\tau^2}\frac{r_d(\Sigma^2)}{n}).  \notag
\end{align}

(ii) Suppose further that $r_d(\Sigma)\asymp n$. For $\lambda_{d+1} \lesssim   \tau \lesssim \lambda_{d}$, we have
\begin{align}
    \mathcal{B}_{\In}(\tau;\hat{H}_n,\hat{G}_n,\gamma) +
    \mathcal{V}_{\In}(\tau;\hat{H}_n,\gamma) &\asymp  \Vert\theta^*_{1:d}\Vert_{\Sigma^{-1}_{1:d}}^2 \tau^2+\sigma^2(\frac{d}{n}+\frac{\lambda_{d+1}^2}{\tau^2}).  \notag
\end{align}
\end{cor}

\begin{prf}
\newline
We first show that $\lambda_d \gtrsim \alpha \tau$. From (\ref{eq:4.1.0}), $\tau\gtrsim \lambda_{d+1}$ and $n> d$, we have
\begin{align}
    &\frac{1}{\alpha}+\frac{1}{n} \frac{\sum_{j>d}^p \lambda_j}{\alpha \tau} \asymp  1 \notag \\
   \Longrightarrow &\alpha \tau \asymp (\tau+\frac{\sum_{j>d}^p \lambda_j}{n}) \notag \\
   \Longrightarrow &\alpha \tau \asymp \tau \quad(\text{from~ $\lambda_{d+1} \lesssim   \tau$~and~$r_d(\Sigma)\lesssim n$}) \notag \\
   \Longrightarrow & \lambda_d \gtrsim \alpha \tau \quad(\text{from~}\tau \lesssim \lambda_d)  .\label{s_iii_2_1}
\end{align}
Then we show $\alpha\asymp 1$. From (\ref{eq:4.1.0}), $\tau\gtrsim \lambda_{d+1}$ and $n> d$, we have
\begin{align}
       &\frac{1}{\alpha}(1+\frac{\sum_{j>d}\lambda_j}{n \tau})   \asymp 1  \notag \\
       & (1+\frac{\sum_{j>d}\lambda_j}{n \tau})   \asymp \alpha \notag \\
       \Longrightarrow & \alpha \asymp 1  \quad (\text{from~}\lambda_{d+1} \lesssim   \tau\text{~and~}r_d(\Sigma)\lesssim n).\label{eq:s2_15e}
\end{align}
Then we give the orders of $1-\frac{1}{n}\sum_{j=1}^p\frac{\lambda_j^2}{(\lambda_j+\alpha\tau)^2}$, $\frac{1}{n}\sum_{j=1}^p \frac{\alpha^2\tau^2 \lambda_j \theta_j^{*2} }{(\lambda_j +\alpha\tau)^2}$ and $\frac{1}{n}\sum_{j=1}^p\frac{\lambda_j^2}{(\lambda_j + \alpha \tau)^2}$, which are important in the formulas (\ref{eq:4.1.1})--(\ref{eq:4.1.4}). We have
\begin{align}
&1-\frac{1}{n}\sum_{j=1}^p\frac{\lambda_j^2}{(\lambda_j+\alpha\tau)^2} \geq 1-\frac{1}{n}\sum_{j=1}^p\frac{\lambda_j}{\lambda_j+\alpha\tau}=\frac{1}{\alpha} \asymp 1 \quad (\text{from~} (\ref{eq:s2_15e}) ), \notag\\
&1-\frac{1}{n}\sum_{j=1}^p\frac{\lambda_j^2}{(\lambda_j+\alpha\tau)^2}\leq 1.   \notag
\end{align}
Hence
\begin{align}
    1-\frac{1}{n}\sum_{j=1}^p\frac{\lambda_j^2}{(\lambda_j+\alpha\tau)^2} \asymp 1 . \label{eq:s2_15f}
\end{align}
From (\ref{s_iii_2_1})--(\ref{eq:s2_15e}) and $\tau \gtrsim \lambda_{d+1}$, we have
\begin{align}
    &\frac{1}{n}\sum_{j=1}^p\frac{\lambda_j^2}{(\lambda_j + \alpha \tau)^2} \asymp \frac{d}{n} + \sum_{j>d}\frac{\lambda_j^2}{n\tau^2}  ,\label{eq:s2_15g} \\
    &\frac{1}{n}\sum_{j=1}^p \frac{\alpha^2\tau^2 \lambda_j \theta_j^{*2} }{(\lambda_j +\alpha\tau)^2} \asymp  \sum_{j=1}^d \frac{\tau^2 \theta_j^{*2}}{\lambda_j} + \sum_{j>d}\lambda_j\theta_j^{*2} .\label{eq:s2_15h}
\end{align}
Substituting (\ref{eq:s2_15f})--(\ref{eq:s2_15h}) into (\ref{eq:4.1.1})--(\ref{eq:4.1.2}), we have
\begin{align}
    \mathcal{B}_{\out}(\tau, \hat{H}_n, \hat{G}_n, \gamma) &\asymp \sum_{j=1}^d \frac{\tau^2 \theta_j^{*2}}{\lambda_j} + \sum_{j>d}\lambda_j\theta_j^{*2}\notag \\
    &\asymp \tau^2 \Vert\theta^*_{1:d}\Vert_{\Sigma^{-1}_{1:d}}^2 +\Vert\theta^*_{(d+1):p}\Vert^2_{\Sigma_{(d+1):p}}\notag \\
    &\asymp \tau^2 \Vert\theta^*_{1:d}\Vert_{\Sigma^{-1}_{1:d}}^2\quad(\text{from~}\Vert\theta^*_{(d+1):p}\Vert^2_{\Sigma_{(d+1):p}}\lesssim  \Vert\theta^*_{1:d}\Vert^2_{\Sigma_{1:d}^{-1}}\lambda_{d+1}^2) , \notag \\
    \mathcal{V}_{\out}(\tau, \gamma, \tilde{\lambda}) &\asymp \sigma^2(\frac{d}{n} + \sum_{j>d}\frac{\lambda_j^2}{n\tau^2}). \notag
\end{align}
Hence we have
\begin{align*}
    \mathcal{B}_{\out}(\tau, \hat{H}_n, \hat{G}_n, \gamma)+\mathcal{V}_{\out}(\tau;\hat{H}_n,\gamma)&\asymp \tau^2 \Vert\theta^*_{1:d}\Vert_{\Sigma^{-1}_{1:d}}^2+\sigma^2(\frac{d}{n} + \sum_{j>d}\frac{\lambda_j^2}{n\tau^2}).
\end{align*}

$\newline$
Note that $\mathcal{V}_{\In}(\tau;\hat{H}_n,\gamma)$ can be also expressed as
\begin{align}
    \mathcal{V}_{\In}(\tau;\hat{H}_n,\gamma) = (1-\frac{1}{\alpha})^2\sigma^2 +\frac{1}{\alpha^2}\frac{\frac{1}{n}\sum_{j=1}^p\frac{\lambda_j^2}{(\lambda_j+\alpha\tau)^2}}{(1-\frac{1}{n}\sum_{j=1}^p\frac{\lambda_j^2}{(\lambda_j+\alpha\tau)^2})} \sigma^2 .\label{eq:s2_15d}
\end{align}
From (\ref{eq:4.1.0}), (\ref{s_iii_2_1}), (\ref{eq:s2_15e}) and $\tau\gtrsim \lambda_{d+1}$, we have
\begin{align}
   1- \frac{1}{\alpha} \asymp  \frac{d}{n} +  \frac{\sum_{j>d}\lambda_j}{n\tau}. \label{eq:s2_15dd}
\end{align}
Substituting (\ref{eq:s2_15e})--(\ref{eq:s2_15h}) into (\ref{eq:4.1.3}) and (\ref{eq:s2_15d}), we have
\begin{align*}
    \mathcal{B}_{\In}(\tau, \hat{H}_n, \hat{G}_n, \gamma) &\asymp \sum_{j=1}^d \frac{\tau^2 \theta_j^{*2}}{\lambda_j} + \sum_{j>d}\lambda_j\theta_j^{*2}\notag \\
    &\asymp \tau^2 \Vert\theta^*_{1:d}\Vert_{\Sigma^{-1}_{1:d}}^2 +\Vert\theta^*_{(d+1):p}\Vert^2_{\Sigma_{(d+1):p}}\notag \\
    &\asymp \tau^2 \Vert\theta^*_{1:d}\Vert_{\Sigma^{-1}_{1:d}}^2 \quad(\text{from~}\Vert\theta^*_{(d+1):p}\Vert^2_{\Sigma_{(d+1):p}}\lesssim  \Vert\theta^*_{1:d}\Vert^2_{\Sigma_{1:d}^{-1}}\lambda_{d+1}^2) , \notag \\
    \mathcal{V}_{\In}(\tau;\hat{H}_n,\gamma) &\asymp \sigma^2 (\frac{d}{n} + \frac{(\frac{\sum_{j>d}\lambda_j}{n})^2}{\tau^2}) \notag \\
    &\asymp  \sigma^2 (\frac{d}{n} + \frac{\lambda_{d+1}^2}{\tau^2})\quad (\text{from~}r_d(\Sigma)\asymp n).  \notag
\end{align*}
Hence we have
\begin{align}
    \mathcal{B}_{\In}(\tau;\hat{H}_n,\hat{G}_n,\gamma) +
    \mathcal{V}_{\In}(\tau;\hat{H}_n,\gamma) &\asymp  \Vert\theta^*_{1:d}\Vert_{\Sigma^{-1}_{1:d}}^2 \tau^2+\sigma^2(\frac{d}{n}+\frac{\lambda_{d+1}^2}{\tau^2})  .
\end{align}

\end{prf}

\begin{cor}[Matching error approximation formulas with large TER]
~$\newline$
\indent (i) Suppose that $\frac{d}{n}<\frac{1}{5}$, $r_d(\Sigma) > c n$ for some $c>10$, $\Vert\theta^*_{(d+1):p}\Vert^2_{\Sigma_{(d+1):p}}\lesssim  \Vert\theta^*_{1:d}\Vert^2_{\Sigma_{1:d}^{-1}}(\frac{1}{\lambda_d}+\frac{n}{\sum_{j>d}\lambda_j})^{-2}$. For
    $\lambda_{d} \gtrsim  \tau + \lambda_{d+1}\frac{r_d(\Sigma)}{n}$ and $\tau>0$, we have
\begin{align}
       \mathcal{B}_{\out}(\tau;\hat{H}_n,\hat{G}_n,\gamma) +
    \mathcal{V}_{\out}(\tau;\hat{H}_n,\gamma) &\asymp \Vert\theta_{1:d}^*\Vert^2_{\Sigma_{1:d}^{-1}}(\tau+\lambda_{d+1}\frac{r_d(\Sigma)}{n})^2 + \sigma^2(\frac{d}{n}+\frac{\lambda_{d+1}^2}{(\tau+\lambda_{d+1}\frac{r_d(\Sigma)}{n})^2}\frac{r_d(\Sigma^2)}{n}) . \notag
\end{align}

(ii) Suppose further that $\tau > \lambda_{d+1}\frac{r_d(\Sigma)}{n}$. For
    $\lambda_{d} \gtrsim  \tau + \lambda_{d+1}\frac{r_d(\Sigma)}{n}$ and $\tau>0$, we have
\begin{align}
    \mathcal{B}_{\In}(\tau;\hat{H}_n,\hat{G}_n,\gamma) +
    \mathcal{V}_{\In}(\tau;\hat{H}_n,\gamma) &\asymp \Vert\theta_{1:d}^*\Vert^2_{\Sigma_{1:d}^{-1}}(\tau+\lambda_{d+1}\frac{r_d(\Sigma)}{n})^2 + \sigma^2(\frac{d}{n}+\frac{\lambda_{d+1}^2}{(\tau+\lambda_{d+1}\frac{r_d(\Sigma)}{n})^2}\frac{r^2_d(\Sigma)}{n^2}) . \notag
\end{align}
\end{cor}

\begin{prf}
$\newline$
We first prove that $\alpha \tau > (c-1)\lambda_{d+1}$. From (\ref{eq:4.1.0}),
\begin{align*}
    n &>  \sum_{j>d}\frac{\lambda_j}{\lambda_j+\alpha\tau} \\
     &> \frac{\sum_{j>d}\lambda_j}{\lambda_{d+1}+\alpha \tau} \\
     &> \frac{\sum_{j>d}\lambda_j}{\lambda_{d+1}} \frac{1}{(1+\frac{\alpha \tau}{\lambda_{d+1}})} \\
     &> cn  \frac{1}{(1+\frac{\alpha \tau}{\lambda_{d+1}})}.
\end{align*}
Hence we have
\begin{align}
   \alpha\tau > (c-1)\lambda_{d+1}.  \label{eq:s2_11a}
\end{align}
Then we prove that $\frac{c}{c-1}>\frac{1}{\alpha}+\frac{1}{n}\frac{\sum_{j>d}\lambda_j}{\alpha\tau}>\frac{4}{5}$. We have,
\begin{align*}
    \frac{1}{\alpha}+\frac{1}{n}\frac{\sum_{j>d}\lambda_j}{\alpha\tau} &>   \frac{1}{\alpha}+\frac{1}{n}\frac{\sum_{j>d}\lambda_j}{\lambda_j+\alpha\tau} \notag \\
    &=1-\frac{1}{n}\frac{\sum_{i=1}^d \lambda_j}{\lambda_j +\alpha\tau}\quad(\text{from }(\ref{eq:4.1.0})) \\
    &\geq 1-\frac{d}{n}\\
    &>\frac{4}{5}\quad(\text{from }\frac{d}{n}<\frac{1}{5}).
\end{align*}
From (\ref{eq:4.1.0}) and (\ref{eq:s2_11a}), we have
\begin{align*}
   \frac{c-1}{c}(\frac{1}{\alpha}+ \frac{1}{n}\frac{\sum_{j>d}\lambda_j}{\alpha\tau})&=\frac{1}{(\frac{1}{c-1}+1)\alpha}+\frac{1}{n}\frac{\sum_{j>d}\lambda_j}{(\frac{1}{c-1}+1)\alpha\tau}  \\
    &<  \frac{1}{\alpha}+\frac{1}{n}\frac{\sum_{j>d}\lambda_j}{\lambda_{d+1}+\alpha\tau} \\
    &\leq \frac{1}{\alpha}+\frac{1}{n}\frac{\sum_{j>d}\lambda_j}{\lambda_{j}+\alpha\tau} \\
    &=1-\frac{1}{n}\frac{\sum_{i=1}^d \lambda_j}{\lambda_j +\alpha\tau}\\
    &\leq 1.
\end{align*}
That is, we have
\begin{align*}
    \frac{1}{\alpha}+\frac{1}{n}\frac{\sum_{j>d}\lambda_j}{\alpha\tau} < \frac{c}{c-1}.
\end{align*}
Hence we have
\begin{align}
    \frac{c}{c-1}>\frac{1}{\alpha}+\frac{1}{n}\frac{\sum_{j>d}\lambda_j}{\alpha\tau}>\frac{4}{5} ,  \label{eq:s2_17a}
\end{align}
and
\begin{align}
   \frac{1}{\alpha}+\frac{1}{n}\frac{\sum_{j>d}\lambda_j}{\alpha\tau}\asymp 1.\label{eq:s2_17ab}
\end{align}
Now we prove that $\lambda_d\gtrsim \alpha\tau$. From (\ref{eq:s2_17ab}), we have
\begin{align}
    &\frac{1}{\alpha}+\frac{1}{n} \frac{\sum_{j>d}^p \lambda_j}{\alpha \tau} \asymp  1   \notag\\
   \Longrightarrow &\alpha \tau \asymp (\tau+\frac{\sum_{j>d}^p \lambda_j}{n}) \notag \\
   \Longrightarrow & \lambda_d \gtrsim \alpha \tau \quad(\text{from~}\lambda_{d} \gtrsim  \tau + \frac{\sum_{j>d}\lambda_j}{n}). \label{eq:s2_11}
\end{align}

We prove the order matching in two cases, $\sum_{j>d}\lambda_j< n\tau$ and $\sum_{j>d}\lambda_j\geq n\tau$. We discuss the two cases separately.

\noindent \textcircled{1} $\sum_{j>d}\lambda_j< n\tau$.

From (\ref{eq:s2_17ab}) and $\sum_{j>d}\lambda_j< n\tau$, we have
\begin{align}
    \alpha \asymp  1 . \label{eq:s2_14a}
\end{align}
From (\ref{eq:s2_11a}), (\ref{eq:s2_11}) and (\ref{eq:s2_14a}), by a similar process as the proof in Corollary \ref{cor7_c}, (\ref{eq:s2_15f}), (\ref{eq:s2_15g}),  (\ref{eq:s2_15h}) and (\ref{eq:s2_15dd}) hold. Substituting (\ref{eq:s2_15f}), (\ref{eq:s2_15g}) and (\ref{eq:s2_15h}) into (\ref{eq:4.1.1}), (\ref{eq:4.1.2}), we have
\begin{align*}
    \mathcal{B}_{\out}(\tau, \hat{H}_n, \hat{G}_n, \gamma)+\mathcal{V}_{\out}(\tau, \hat{H}_n, \hat{G}_n, \gamma)&\asymp \tau^2 \Vert\theta_{1:d}\Vert_{\Sigma^{-1}_{1:d}}^2+\sigma^2(\frac{d}{n} + \sum_{j>d}\frac{\lambda_j^2}{n\tau^2}) \\
    &\asymp  (\tau+\lambda_{d+1}\frac{r_d(\Sigma)}{n})^2 \Vert\theta^*_{1:d}\Vert_{\Sigma^{-1}_{1:d}}^2+\sigma^2(\frac{d}{n} + \sum_{j>d}\frac{\lambda_j^2}{n(\tau+\lambda_{d+1}\frac{r_d(\Sigma)}{n})^2}) \\&(\text{from~}\sum_{j>d}\lambda_j < n\tau).
\end{align*}
Substituting  (\ref{eq:s2_15f}), (\ref{eq:s2_15g}), (\ref{eq:s2_15h}) ,(\ref{eq:s2_15dd}) and (\ref{eq:s2_14a}) into (\ref{eq:4.1.3}) and (\ref{eq:s2_15d}), we have
\begin{align*}
    \mathcal{B}_{\In}(\tau, \hat{H}_n, \hat{G}_n, \gamma) &\asymp \sum_{j=1}^d \frac{\tau^2 \theta_j^{*2}}{\lambda_j} + \sum_{j>d}\lambda_j\theta_j^{*2}\notag \\
    &\asymp (\tau+\lambda_{d+1}\frac{r_d(\Sigma)}{n})^2 \Vert\theta^*_{1:d}\Vert_{\Sigma^{-1}_{1:d}}^2 +\Vert\theta^*_{(d+1):p}\Vert^2_{\Sigma_{(d+1):p}}\quad(\text{from~}\sum_{j>d}\lambda_j < n\tau)  \notag \\
    &\asymp (\tau+\lambda_{d+1}\frac{r_d(\Sigma)}{n})^2 \Vert\theta^*_{1:d}\Vert_{\Sigma^{-1}_{1:d}}^2 \quad(\text{from~}\Vert\theta^*_{(d+1):p}\Vert^2\lesssim  \Vert\theta^*_{1:d}\Vert^2_{\Sigma_{1:d}^{-1}}(\frac{1}{\lambda_d}+\lambda_{d+1}\frac{1}{\frac{\sum_{j>d}\lambda_j}{n}})^{-2}) , \notag \\
    \mathcal{V}_{\In}(\tau;\hat{H}_n,\gamma) &\asymp \sigma^2 (\frac{d}{n} + \frac{(\frac{\sum_{j>d}\lambda_j}{n})^2}{\tau^2}) \notag \\
    &\asymp  \sigma^2 (\frac{d}{n} + \frac{\lambda_{d+1}^2}{(\tau+\lambda_{d+1}\frac{r_d(\Sigma)}{n})^2}\frac{r^2_d(\Sigma)}{n^2}) \quad(\text{from~}\sum_{j>d}\lambda_j < n\tau).
\end{align*}
Hence
\begin{align}
    \mathcal{B}_{\In}(\tau;\hat{H}_n,\hat{G}_n,\gamma) +
    \mathcal{V}_{\In}(\tau;\hat{H}_n,\gamma) &\asymp \Vert\theta_{1:d}^*\Vert^2_{\Sigma_{1:d}^{-1}}(\tau+\lambda_{d+1}\frac{r_d(\Sigma)}{n})^2 + \sigma^2(\frac{d}{n}+\frac{\lambda_{d+1}^2}{(\tau+\lambda_{d+1}\frac{r_d(\Sigma)}{n})^2}\frac{r^2_d(\Sigma)}{n^2}) . \notag
\end{align}

\noindent \textcircled{2} $\sum_{j>d}\lambda_j\geq n\tau$.

We first give the order of the term $1-\frac{1}{n}\sum_{j=1}^p\frac{\lambda_j^2}{(\lambda_j+\alpha\tau)^2}$, which is important in the following analysis.
We have
\begin{align}
    \frac{1}{n}\sum_{j=1}^p \frac{\lambda_j^2}{(\lambda_j+\alpha \tau)^2} &\leq \frac{d}{n} + \frac{\sum_{j>d}\lambda_j^2}{n\alpha^2\tau^2} \notag \\
    &= \frac{d}{n} + \frac{n\sum_{j>d}\lambda_j^2}{n^2\alpha^2\tau^2} \notag \\
    &=  \frac{d}{n} +   \frac{(\sum_{j>d}\lambda_j)^2}{n^2\alpha^2\tau^2}\frac{n\sum_{j>d}\lambda_j^2}{(\sum_{j>d}\lambda_j)^2}               \notag \\
    &\leq  \frac{d}{n} +   \frac{c^2}{(c-1)^2}\frac{n\sum_{j>d}\lambda_j^2}{(\sum_{j>d}\lambda_j)^2}  \quad(\text{from (\ref{eq:s2_17a})})   \notag \\
    &\leq \frac{d}{n} + \frac{c^2}{(c-1)^2} \frac{n}{r_d(\Sigma)} \notag \\
    &\leq \frac{d}{n} + \frac{c}{(c-1)^2}\quad(\text{from $r_d(\Sigma)>cn$})  \\
    &<  2/5\quad(\text{from $\frac{d}{n}<\frac{1}{5}$ and $c>10$}). \label{eq:s2_14b}
\end{align}
From (\ref{eq:s2_14b}), we have
\begin{align}
     1>1-\frac{1}{n}\sum_{j=1}^p \frac{\lambda_j^2}{(\lambda_j+\alpha \tau)^2} > \frac{3}{5},   \label{eq:s2_15b}
\end{align}
and
\begin{align}
    1-\frac{1}{n}\sum_{j=1}^p \frac{\lambda_j^2}{(\lambda_j+\alpha \tau)^2} \asymp 1   .\label{eq:s2_15bb}
\end{align}

Now we give the orders of $\frac{1}{n}\sum_{j=1}^p \frac{\alpha^2\tau^2 \lambda_j \theta_j^{*2} }{(\lambda_j +\alpha\tau)^2}$ and $\frac{1}{n}\sum_{j=1}^p\frac{\lambda_j^2}{(\lambda_j + \alpha \tau)^2}$, which are important in formulas (\ref{eq:4.1.1}) and (\ref{eq:4.1.2}). From (\ref{eq:s2_17ab}), we have
\begin{align}
   & \frac{1}{\alpha}+\frac{1}{n}\frac{\sum_{j>d}\lambda_j}{\alpha \tau} \asymp 1 \notag \\
    \Longrightarrow  &\frac{1}{n}\frac{\sum_{j>d}\lambda_j}{\alpha \tau} \asymp 1 \quad(\text{from~}\frac{\sum_{j>d}\lambda_j}{n}\geq n\tau) \notag \\
    \Longrightarrow & \alpha \tau \asymp \frac{\sum_{j>d}\lambda_j}{n}.  \label{eq:s2_15c}
\end{align}
From (\ref{eq:s2_11a}), (\ref{eq:s2_11}), (\ref{eq:s2_15c}) and $\sum_{j>d}\lambda_j\geq n\tau$, we have
\begin{align}
 &\frac{1}{n}\sum_{j=1}^p\frac{\lambda_j^2}{(\lambda_j + \alpha \tau)^2} \asymp \frac{d}{n} + \sum_{j>d}\frac{\lambda_j^2}{n(\frac{\sum_{j>d}\lambda_j}{n})^2}  , \label{eq:s2_15o}  \\
  &  \frac{1}{n}\sum_{j=1}^p \frac{\alpha^2\tau^2 \lambda_j \theta_j^{*2} }{(\lambda_j +\alpha\tau)^2} \asymp  \sum_{j=1}^d \frac{(\frac{\sum_{j>d}\lambda_j}{n})^2 \theta_j^{*2}}{\lambda_j} + \sum_{j>d}\lambda_j\theta_j^{*2} .  \label{eq:s2_15p}
\end{align}
Substituting (\ref{eq:s2_15b}), (\ref{eq:s2_15c}), (\ref{eq:s2_15o}) and (\ref{eq:s2_15p}) into (\ref{eq:4.1.1})--(\ref{eq:4.1.2}), we have
\begin{align}
    \mathcal{B}_{\out}(\tau, \hat{H}_n, \hat{G}_n, \gamma) &\asymp \sum_{j=1}^d \frac{(\frac{\sum_{j>d}\lambda_j}{n})^2 \theta_j^{*2}}{\lambda_j} + \sum_{j>d}\lambda_j\theta_j^{*2}\notag\\
    &\asymp (\frac{\sum_{j>d}\lambda_j}{n})^2 \Vert\theta^*_{1:d}\Vert_{\Sigma^{-1}_{1:d}}^2 +\Vert\theta^*_{(d+1):p}\Vert^2_{\Sigma_{(d+1):p}}\notag \\
    &\asymp (\frac{\sum_{j>d}\lambda_j}{n})^2 \Vert\theta^*_{1:d}\Vert_{\Sigma^{-1}_{1:d}}^2\quad(\text{from~}\Vert\theta^*_{(d+1):p}\Vert^2_{\Sigma_{(d+1):p}}\ll  \Vert\theta^*_{1:d}\Vert^2_{\Sigma_{1:d}^{-1}}(\frac{1}{\lambda_d}+\frac{1}{\frac{\sum_{j>d}\lambda_j}{n}})^{-2}) ,\notag \\
    \mathcal{V}_{\out}(\tau;\hat{H}_n,\gamma) &\asymp \sigma^2(\frac{d}{n} + \sum_{j>d}\frac{\lambda_j^2}{n(\frac{\sum_{j>d}\lambda_j}{n})^2}) .\notag
\end{align}
That is,
\begin{align*}
    \mathcal{B}_{\out}(\tau, \hat{H}_n, \hat{G}_n, \gamma)+\mathcal{V}_{\out}(\tau;\hat{H}_n,\gamma)&\asymp \Vert\theta^*_{1:d}\Vert_{\Sigma^{-1}_{1:d}}^2(\frac{\sum_{j>d}\lambda_j}{n})^2 + \sigma^2(\frac{d}{n} + \sum_{j>d}\frac{\lambda_j^2}{n(\frac{\sum_{j>d}\lambda_j}{n})^2}) \\
    &\asymp  \Vert\theta^*_{1:d}\Vert_{\Sigma^{-1}_{1:d}}^2(\frac{\sum_{j>d}\lambda_j}{n}+\tau)^2 + \sigma^2(\frac{d}{n} + \sum_{j>d}\frac{\lambda_j^2}{n(\frac{\sum_{j>d}\lambda_j}{n}+\tau)^2}) \notag \\
    &(\text{from~}\sum_{j>d}\lambda_j \geq n\tau).
\end{align*}

\end{prf}

\subsection{Alternative calculation of error approximation formulas}
\label{sec_supp_II4}

The asymptotic out-sample and in-sample errors can also be calculated using a distributional approximation method in \cite{HanShen2022} under the independent components assumption.
By letting $\mu = \sqrt{n} \Sigma^{1/2}\theta$,
the ridge estimator in (\ref{eq:2}) can be equivalently formulated as $\hat{\theta}(\tau)=\frac{1}{\sqrt{n}} \Sigma^{-1/2}\hat{\mu}(\tau)$ with
\begin{align}
    \hat{\mu}(\tau) &= \arg \max_{\mu \in\mathbb{R}^p}  \left\{ \Vert Y-Z\mu\Vert^2 + \tau \Vert \Sigma^{-1/2} \mu\Vert^2 \right\} , \notag
\end{align}
where $Z=\frac{1}{\sqrt{n} }X \Sigma^{-1/2}$.
The rows in $Z$ are covariate vectors with covariance matrix $\frac{1}{n}I_p$.
Let $\mu^* = \sqrt{n} \Sigma^{1/2}\theta^*$.
Following \cite{HanShen2022}, for $z\in\mathbb{R}^p,\tau>0$, we define
\begin{align}
    \psi_{\tilde{\lambda}}(z,\tau) = \argmin_{x\in \mathbb{R}^p}
    \left\{\frac{1}{2}\Vert x - z\Vert^{2} + \frac{\tau}{2}\Vert \Sigma^{-1/2} x\Vert^2\right\}
    =(\frac{z_1}{1+\frac{\tau}{\lambda_1}},..., \frac{z_p}{1+\frac{\tau}{\lambda_p}})^\T .\notag
\end{align}
For an isotropic random vector $z_0\in \mathbb{R}^p$, suppose $(\alpha,\beta)$ is
a unique solution in $ (0,\infty)^2$ to the following equations:
\begin{align}
    \beta^2 -\sigma^2&=  \frac{1}{n}\mathbb{E}\Vert \psi_{\tilde{\lambda}}(\mu^*+\beta z_0,\alpha \tau) -\mu^*\Vert^2 = \frac{1}{n}\sum_{j} \frac{\alpha^2\tau^2\mu_j^{*2} }{(\lambda_j+\alpha\tau)^2} + (\frac{1}{n}\sum_{j}\frac{\lambda_j^2}{(\lambda_j+\alpha \tau)^2})\beta^2, \label{eq:4_2_a}\\
    \frac{1}{\alpha}&= 1-\frac{1}{n}\mathbb{E}\div \psi_{\tilde{\lambda}}(\mu^* + \beta z_0,\alpha \tau)=\gamma\frac{1}{p}\sum_{j=1}^p\frac{1}{1+\frac{\alpha \tau}{\lambda_j}}. \notag
\end{align}
where $\div f(x_1,\ldots,x_p)=\sum_{j=1}^p\frac{\partial f}{\partial x_j}$. Note that from (\ref{eq:4_2_a}), we have
\begin{align}
    \beta^2 &= (1-\frac{1}{n}\sum_{j}\frac{\lambda_j^2}{(\lambda_j+\alpha \tau)^2})^{-1}(\sigma^2 + \frac{1}{n}\sum_{j} \frac{\alpha^2\tau^2\mu_j^{*2} }{(\lambda_j+\alpha\tau)^2} ) \notag \\
    &=(1-\frac{1}{n}\sum_{j}\frac{\lambda_j^2}{(\lambda_j+\alpha \tau)^2})^{-1}(\sigma^2 + \sum_{j} \frac{\alpha^2\tau^2\lambda_j^2\theta^{*2}_j}{(\lambda_j+\alpha\tau)^2} ) .\label{eq:4_2_a_1a}
\end{align}
Then under suitable conditions (see conditions (R1)-(R3) in \cite{HanShen2022}), the distributions of $\hat{\mu}-\mu^*$ and $Z(\hat{\mu}-\mu^*)$ can be approximated as follows:
\begin{equation}
    \begin{aligned}
         \hat{\mu}-\mu^* &\overset{d}{\approx}  \psi_{\tilde{\lambda}}(\mu^*+\beta z_0,\alpha \tau) - \mu^* ,  \\
       Z(\hat{\mu}-\mu^*) &\overset{d}{\approx} (1-\frac{1}{\alpha})\xi +\frac{\sqrt{\gamma^2-\sigma^2}}{\alpha}h,
    \end{aligned} \label{eq:III_4}
\end{equation}
where $\xi\sim N(0,\sigma^2)$ and $h\sim N(0,\sigma^2I_n)$.

Given the approximation results (\ref{eq:III_4}), the asymptotic error formulas in Corollary \ref{cor5b} can also be calculated as follows.
The out-sample error can be approximated by
\begin{align}
    \Vert \hat{\theta}(\tau) - \theta^*\Vert^2_\Sigma &= \frac{1}{n} \Vert \hat{\mu}-\mu^*\Vert^2 \notag \\
    &\approx  \frac{1}{n} \Vert \psi_{\tilde{\lambda}}(\mu^*+\beta z_0,\alpha \tau) - \mu^*\Vert^2 \notag \\
    &\approx \frac{1}{n}\mathbb{E}\Vert \psi_{\tilde{\lambda}}(\mu^*+\beta z_0,\alpha \tau) -\mu^*\Vert^2 \notag \\
    &= \beta^2-\sigma^2 . \label{eq:4_2_a_3}
\end{align}
Substituting (\ref{eq:4_2_a_1a}) into (\ref{eq:4_2_a_3}) yields the sum of (\ref{eq:4.1.1}) and (\ref{eq:4.1.2}).
The in-sample error can be approximated by
\begin{align}
    \Vert \hat{\theta}(\tau) - \theta^*\Vert^2_{\hat{\Sigma} }&= \frac{1}{n} \Vert Z(\hat{\mu}-\mu^*)\Vert^2 \notag \\
    &\approx  (1-\frac{1}{\alpha})^2\sigma^2+\frac{\beta^2-\sigma^2}{\alpha^2} .\label{eq:4_2_a_4}
\end{align}
Substituting (\ref{eq:4_2_a_1a}) into (\ref{eq:4_2_a_4}) yields the sum of (\ref{eq:4.1.3}) and (\ref{eq:4.1.4}).

\section{Comparison with \cite{BuneaFlorentinaStrimas2020}}
\label{sec:approximation_com}

\subsection{Approximations of terms}

\label{sec:terms_approximation}

In this section, we give the approximations of certain terms used in the comparison in Section \ref{sec:4_3_main} between our Theorem \ref{thm3} and Theorem 16 of \cite{BuneaFlorentinaStrimas2020}.
In the setting described in Section \ref{sec:4_3_main}, we show that the following approximations hold:
$\Vert \theta_{1:d}^*\Vert^2_{\Sigma_{1:d}^{-1}}\lambda_d^2\asymp \Vert \theta_{1:d}^*\Vert^2_{\Sigma_{1:d}}$, $\Vert \beta\Vert^2_{\Sigma_Z}\asymp \Vert \theta_{1:d}^*\Vert^2_{\Sigma_{1:d}}$, $\frac{\lambda_{d+1}}{\lambda_d-\lambda_{d+1}}\asymp \frac{\lambda_{d+1}}{\lambda_d}$ and $r_0(\Sigma_E)\asymp r_d(\Sigma)$.

\noindent \textcircled{1} $\Vert \theta_{1:d}^*\Vert^2_{\Sigma_{1:d}^{-1}} \lambda_d^2\asymp \Vert \theta_{1:d}^*\Vert^2_{\Sigma_{1:d}}$

\noindent First, we have
\begin{align*}
 \Vert \theta_{1:d}^*\Vert^2_{\Sigma_{1:d}}  \geq  \Vert \theta_{1:d}^*\Vert^2_{\Sigma_{1:d}^{-1}} \lambda_d^2 \geq    \Vert \theta_{1:d}^*\Vert^2_{\Sigma_{1:d}} \frac{\lambda_d^2}{\lambda_1^2}.
\end{align*}
Then from $\lambda_1\asymp \lambda_d$, we have
$\Vert \theta_{1:d}^*\Vert^2_{\Sigma_{1:d}^{-1}} \lambda_d^2\asymp \Vert \theta_{1:d}^*\Vert^2_{\Sigma_{1:d}}$.\\

\noindent \textcircled{2} $\Vert \theta_{1:d}^* \Vert^2_{\Sigma_{1:d}}\asymp \Vert \beta\Vert_{\Sigma_Z}^2$ and $\frac{\lambda_{d+1}}{\lambda_d-\lambda_{d+1}}\asymp \frac{\lambda_{d+1}}{\lambda_d}$

\noindent From $\lambda_1\geq \ldots \geq \lambda_d\geq c_1 \lambda_{d+1}$ for some $c_1>1$, we have for $i=1,\ldots, d$,
\begin{align*}
    &1\geq  \frac{\lambda_i -\lambda_{d+1}}{\lambda_i} \geq  1-\frac{1}{c_1} , \\
    & \frac{\lambda_i}{\lambda_{d+1}}\geq  \frac{\lambda_i -\lambda_{d+1}}{\lambda_{d+1}} \geq  (1-\frac{1}{c_1}) \frac{\lambda_i}{\lambda_{d+1}}.
\end{align*}
Then for $i=1,\ldots, d$,
\begin{align}
    \frac{\lambda_i -\lambda_{d+1} }{\lambda_i} &\asymp 1, \label{eq:ss_0} \\
     \frac{\lambda_i -\lambda_{d+1} }{\lambda_i} &\asymp \frac{\lambda_i}{\lambda_{d+1}} .\label{eq:ss_0a}
\end{align}
From (\ref{eq:ss_0a}), we have $\frac{\lambda_{d+1}}{\lambda_d-\lambda_{d+1}}\asymp \frac{\lambda_{d+1}}{\lambda_d}$.
From (\ref{eq:ss_0}), we have
\begin{align}
    \Vert \theta_{1:d}^* \Vert^2_{\Sigma_{1:d}} &= \beta^T\Diag(\frac{(\lambda_1-\lambda_d)}{\lambda_1},\ldots, \frac{\lambda_d-\lambda_{d+1}}{\lambda_d})\beta \notag \\
    &\asymp  \beta^T\beta \notag \\
    &=\Vert \beta\Vert_{\Sigma_Z}^2 \quad(\text{with $\Sigma_Z=I_d$}).  \notag
\end{align}

\noindent \textcircled{3} $r_0(\Sigma_E)\asymp r_d(\Sigma)$

\noindent By definition, $r_0(\Sigma_E) = d + r_d(\Sigma)$.
With $r_d(\Sigma) > c_2 d$, we have
\begin{align*}
  (1+c_2)r_d(\Sigma)  \geq   r_0(\Sigma_E) \geq r_d(\Sigma)
\end{align*}
That is, $r_d(\Sigma) \asymp r_0(\Sigma_E)$.

\subsection{Error bounds incorporated with approximations}

In this section, we give the error order of our Theorem 3 and the error upper bound in Theorem 16 of \cite{BuneaFlorentinaStrimas2020} for the min-norm interpolator and incorporate the approximations of terms in Supplement Section \ref{sec:terms_approximation}.

Based on Theorem 3, for $\tau=0$, we have
\begin{align*}
  &  \MSE_\out \asymp \underbrace{\Vert \theta_{1:d}^*\Vert_{\Sigma_{1:d}^{-1}} \lambda_d^2 \frac{\lambda_{d+1}^2}{\lambda_d^2}\frac{r^2_d(\Sigma)}{n^2}}_{\B_\out} + \underbrace{\sigma^2(\frac{d}{n}+ \frac{n r_d(\Sigma^2)}{r^2_d(\Sigma)})}_{\V_\out}, \quad \text{for } \lambda_{d+1}\frac{r_d(\Sigma)}{n}\leq \lambda_{d}, \\
  &  \MSE_\out \gtrsim  \underbrace{\Vert \theta_{1:d}^*\Vert_{\Sigma_{1:d}^{-1}} \lambda_d^2 }_{\B_\out} , \quad
 \text{for } \lambda_{d+1}\frac{r_d(\Sigma)}{n}>\lambda_d.
\end{align*}
Incorporating $\Vert \theta_{1:d}^*\Vert^2_{\Sigma_{1:d}^{-1}}\lambda_d^2\asymp \Vert \theta_{1:d}^*\Vert^2_{\Sigma_{1:d}}$, we have
\begin{align*}
  &  \MSE_\out \asymp \underbrace{\Vert \theta_{1:d}^*\Vert^2_{\Sigma_{1:d}} \frac{\lambda_{d+1}^2}{\lambda_d^2} \frac{r^2_d(\Sigma)}{n^2}}_{\B_\out} + \underbrace{\sigma^2(\frac{d}{n}+ \frac{n r_d(\Sigma^2)}{r^2_d(\Sigma)})}_{\V_\out}, \quad \text{for } \lambda_{d+1}\frac{r_d(\Sigma)}{n}\leq \lambda_{d}, \\
  &  \MSE_\out \gtrsim  \underbrace{\Vert \theta_{1:d}^*\Vert^2_{\Sigma_{1:d}} }_{\B_\out} , \quad
 \text{for } \lambda_{d+1}\frac{r_d(\Sigma)}{n}>\lambda_d.
\end{align*}

Based on Theorem 16 of \cite{BuneaFlorentinaStrimas2020}, we have
\begin{align*}
    \MSE_\out \lesssim \underbrace{ \Vert \beta \Vert^2_{\Sigma_Z} \frac{\lambda_{d+1}}{\lambda_d(A\Sigma_Z A^T)}\frac{r_0(\Sigma_E)}{n}}_{\B_\out} +\underbrace{\sigma^2(\frac{d}{n}+\frac{n}{r_0(\Sigma_E)})}_{\V_\out}.
\end{align*}
Incorporating $\Vert \beta\Vert^2_{\Sigma_Z}\asymp \Vert \theta_{1:d}^*\Vert^2_{\Sigma_{1:d}}$, $\frac{\lambda_{d+1}}{\lambda_d-\lambda_{d+1}}\asymp \frac{\lambda_{d+1}}{\lambda_d}$ and $r_0(\Sigma_E)\asymp r_d(\Sigma)$, we have
\begin{align*}
    \MSE_\out \lesssim \underbrace{\Vert \theta_{1:d}^*\Vert^2_{\Sigma_{1:d}} \frac{\lambda_{d+1}}{\lambda_d}\frac{r_d(\Sigma)}{n}}_{\B_\out} +\underbrace{\sigma^2 (\frac{d}{n}+\frac{n}{r_d(\Sigma)})}_{\V_\out}.
\end{align*}

\end{document}